\documentclass[preprint,12pt]{elsarticle}

\usepackage{amsmath}
\usepackage{amssymb}
\usepackage[margin=2.5cm]{geometry}
\usepackage{subcaption}
\usepackage{algorithm}
\usepackage{algpseudocode}

\newcommand{\dd}{\mathrm{d}}    
\newcommand{\pd}{\partial}  
\renewcommand{\vec}[1]{\boldsymbol{\mathrm{#1}}}     
\newcommand{\vecf}[1]{\boldsymbol{#1}}

\begin{document}

\begin{frontmatter}

\title{Positional Embeddings for Solving PDEs with Evolutional Deep Neural Networks}

\author[inst1]{Mariella Kast \corref{cor1}}
\cortext[cor1]{Corresponding author, email: mariella.kast@epfl.ch}
\affiliation[inst1]{country={Chair of Computational Mathematics and Simulation Science, École Polytechnique Fédérale de Lausanne, 1015 Lausanne, Switzerland}%
            }

\author[inst1]{Jan S. Hesthaven}

\begin{abstract}
 This work extends the paradigm of evolutional deep neural networks (EDNNs) to solving parametric time-dependent partial differential equations (PDEs) on domains with geometric structure. By introducing positional embeddings based on eigenfunctions of the Laplace-Beltrami operator, geometric properties are encoded intrinsically and Dirichlet, Neumann and periodic boundary conditions of the PDE solution are enforced directly through the neural network architecture. The proposed embeddings lead to improved error convergence for static PDEs and extend EDNNs towards computational domains of realistic complexity.
 
 Several steps are taken to improve performance of EDNNs:
Solving the EDNN update equation with a Krylov solver avoids the explicit assembly of Jacobians and enables scaling to larger neural networks. Computational efficiency is further improved by an ad-hoc active sampling scheme that uses the PDE dynamics to effectively sample collocation points.  A modified linearly implicit Rosenbrock method is proposed to  alleviate the time step requirements of stiff PDEs. Lastly,  a completely  training-free approach, which  automatically enforces initial conditions and only requires time integration, is compared against EDNNs that are trained on the initial conditions.
We report results for the Korteweg-de Vries equation,  a nonlinear heat equation and (nonlinear) advection-diffusion problems on domains with and without holes and various boundary conditions, to demonstrate the effectiveness of the method. The numerical results highlight EDNNs as a promising surrogate model for parametrized PDEs with slow decaying Kolmogorov n-width.
\end{abstract}

\begin{keyword}
Scientific machine learning \sep reduced order models \sep many query applications \sep partial differential equations
\end{keyword}

\end{frontmatter}

\section{Introduction}
In many science and engineering applications, partial differential equations (PDEs) are used to model the dynamic behavior of the problem. 
The efficient numerical solution of such PDEs remains an active research topic, and is especially relevant for many-query applications that require repeated evaluation for different parameter settings. Reduced basis (RB) methods \cite{benner_survey_2015,hesthaven_reduced_2022} seek to approximate PDE solutions in a  low-dimensional linear subspace of the solution manifold, leading to computational models of much smaller cost. For  dynamic problems with slowly decaying Kolmogorov n-width, e.g., PDEs with strong advective terms,  RB methods typically require a basis of large size and  they may lose their efficiency. \\
Using deep neural networks (DNNs) for surrogate modelling of PDEs has produced exciting results during the past few years. When experimental data or high fidelity (HF) solution data from a traditional solver is available, neural networks  can be trained to approximate the PDE operator, e.g. the  DeepONet \cite{lu_learning_2021}, neural ODEs \cite{chen_neural_2018}, (Fourier) neural operators \cite{li_fourier_2020},  or to recover the parameter-to-solution map, e.g. POD-NN \cite{hesthaven_non-intrusive_2018}. When HF data is scarce or absent, automatic differentiation tools can be used to train directly on the fit of the PDE residual without requiring access to a classical numerical solver, as utilized in physics-informed neural networks (PINNs) \cite{dissanayake_neuralnetworkbased_1994,raissi_physics-informed_2019} or the Deep Ritz method \cite{yu_deep_2018}. Such neural networks typically learn a mapping from the spatio-temporal input coordinate to the solution output. 

While DNNs offer universal function approximation, it is well known that training PINNs can be challenging and  sensitive to hyperparameter tuning \cite{wang_understanding_2021}, especially when balancing different terms in the residual loss, e.g. the PDE dynamics,boundary conditions, and initial conditions. Recent insights through the lense of neural tangent kernel (NTK) theory \cite{jacot_neural_2018, wang_eigenvector_2021} confirm that the optimization problem in PINNs suffers from stiffness and may not converge to the desired solution.  In \cite{wang_eigenvector_2021, tancik_fourier_2020}, Fourier Features are used to embed the spatial coordinates in a more favorable solution space with a rotation-and translation invariant NTK, which can alleviate the stiffness problem and  enable faster training and better accuracy. 

Building upon these insights, we propose a positional embedding via variational harmonic features. Previously used in Gaussian process regression \cite{solin_hilbert_2020,solin_know_2019}, harmonic features encode geometric information about the computational domain and take into account the  gaps and holes that most engineering computational domains possess. The advantage of the proposed embeddings are two-fold, 1) the solution space of the neural network is informed by the geometry in an intrinsic way, 2) it is possible to  enforce boundary conditions directly on the neural network, hence  eliminating the loss term for the boundary conditions. In \cite{lagaris_artificial_1998, berg_unified_2018}, the authors enforce boundary conditions exactly through the multiplication of the output with specially constructed (polynomial) functions. Similarly, in \cite{du_evolutional_2021} a linear combination of neural network outputs is constructed to enforce homogenous Dirichlet boundary conditions on convex domains. To our knowledge, this work is the first example of  a unified framework to enforce both Dirichlet and Neumann boundary conditions via embedding functions on a neural network. 

When training DNNs with time as an input coordinate, ensuring causality is  another major challenge in solving time-dependent problems with PINNs. 
Recently, exciting progress has been made  for directly evolving the parameters of a neural network in time \cite{du_evolutional_2021, bruna_neural_2022,anderson_evolution_2022}, via an update equation that is derived from the PDE dynamics. The "training" is thus done via time-integration over the PDE dynamics, similar to classical ODE solvers. Building on these ideas, we propose a numerically efficient computation of the update step by employing Krylov solvers. We further describe how linearly-implicit solvers based on Rosenbrock methods \cite{hairer_ii_1991,shampine_matlab_1997} can be used to deal with stiff problems at a small additional computational cost. Constructing neural networks that automatically satisfy the initial condition allows one to completely remove the need for gradient-descent based training. Combining these advances with the feature embeddings, allows for an to efficient solution of  parametrized transport-dominated problems on complicated domains.

In Section 2, we briefly discuss the problem setup for parametrized time-dependent PDEs and the different model paradigms that can be used to approximate the solution. In Section 3, we discuss the construction and properties of the positional embedding and show how boundary conditions are enforced. In Section 4, the time-stepping scheme and "training-free" approach is explained and we introduce the use of linearly-implicit solvers, as well as an active sampling scheme. In Section 5, we demonstrate the methods on numerical examples to highlight their capabilities and limits. We  summarize our results and offer concluding remarks in Section 6.
\section{Problem setting}
\label{sec:PDEprob}
\subsection{PDE formulation}
We consider time dependent PDEs, which may depend on a parameter vector $\vec{\alpha} \in \Omega_\alpha$,  defined  on a spatial domain $\Omega_x \subset \mathbb{R}^d$ with boundary $\pd \Omega_x$:
\begin{equation}
\label{eq:PDE}
\begin{split}
   & \dfrac{\pd u(t,\vec{x}, \vec{\alpha})}{\pd t} = f(u,t,\vec{x}, \vec{\alpha})  \quad \text{for } (t,\vec{x}, \vec{\alpha}) \in [0, \infty ) \times \Omega_x \times \Omega_\alpha, \\
   & \mathcal{B} u(t,\vec{x},\vec{\alpha}) = g(t,\vec{x},\vec{\alpha}) \quad \text{for } (t,\vec{x}, \vec{\alpha}) \in [0, \infty ) \times \partial \Omega_x \times \Omega_\alpha,\\
   & u(0,\vec{x},\vec{\alpha}) =  u_0(\vec{x},\vec{\alpha}) \quad  \text{for } (\vec{x}, \vec{\alpha}) \in \Omega_x \times \Omega_\alpha,
    \end{split}
\end{equation}
where $\vec{x}$ is the spatial coordinate, $u(t,\vec{x}, \vec{\alpha})$ is the time-dependent solution field, $f(u,t,\vec{x}, \vec{\alpha})$ is a nonlinear differential operator and $u_0(\vec{x},\vec{\alpha})$ is a compatible initial condition. We assume that the PDE problem \eqref{eq:PDE} is equipped with an appropriate boundary condition operator $\mathcal{B}$, such that the solution $u$ is well posed for all times $t$ and parameter instances $\vec{\alpha}$. We refer to the collection of all solution instances as the solution manifold $\mathcal{M} = \{ u(t,\vec{x}, \vec{\alpha}): (t,\vec{x}, \vec{\alpha}) \in [0, \infty ) \times \Omega_x \times \Omega_\alpha \}$.

We seek to build a computational model which accurately approximates all solution instances $u$ in $\mathcal{M}$ by their numerical counterparts $\hat{u}$, which are part of the numerical solution manifold $\hat{\mathcal{M}} =\{ \hat{u}(t,\vec{x}, \vec{\alpha}): (t,\vec{x}, \vec{\alpha}) \in [0, T]  \times \Omega_x \times \Omega_\alpha \}$, where $T$ is the final computation time. The choice of numerical method for computing $\hat{u}$ influences the properties of $\hat{\mathcal{M}}$, in particular how well the true solution manifold can be approximated and at what computational cost.

\subsection{FE solution and reduced basis methods}
Classic numerical techniques such as the Finite Element (FE)  method, discretize the physical domain to recover a coefficient vector $\vec{c}$, reflecting a finite number of degrees of freedom, which allows to write the solution as a linear combination of $N$ basis functions $\varphi(\vec{x})$:
\begin{equation}
    \label{eq:FEAnsatz}
    \hat{u}_\text{FE} (t, \vec{x},\vec{\alpha}) = \sum_{i=1}^{N} c_i (t, \vec{\alpha}) \varphi_i (\vec{x})
\end{equation}
This coefficient vector is then evolved in time for each parameter instance separately. The computational cost of approximating the whole solution manifold scales with the dimension of the physical and parameter space and suffers from the curse of dimensionality. Reduced basis methods are an efficient way of decreasing the computational burden, when the problem has a fast-decaying Kolmogorov n-width \cite{pinkus_n-widths_2012, hesthaven_certified_2016, lassila_generalized_2013} -- in such cases it is possible to approximate the true solution manifold by a linear combination of $n$ basis functions $\Psi(\vec{x})$, where $n<< N$:
\begin{equation*}
    \hat{\mathcal{M}}_\text{RB} = \text{span} \{  \Psi_1, \Psi_2, ..., \Psi_n \} \approx \mathcal{\hat{M}}_\text{FE}
\end{equation*}
and the solutions can be expressed  with a reduced coefficient vector $\vec{\beta}$
\begin{equation}
\label{eq:RB}
    \hat{u}_\text{RB} (t, \vec{x},\vec{\alpha}) = \sum_{i=1}^{n} \beta_i (t, \vec{\alpha}) \Psi_i (\vec{x}).
\end{equation}
\subsection{Neural network solvers}
Neural networks have appeared as an appealing alternative of nonlinear models as they promise universal function approximation without suffering from the curse of dimensionality.  In this work, we consider fully connected feed forward neural networks, recursively defined as

\begin{equation}
\label{eq:NN}
    \vec{v_{l+1}}(\vec{v_l})  = \sigma (\vec{W_l} \vec{v_l} +\vec{b_l}),
\end{equation}
where $\sigma(\cdot)$ is the activation function, $\vec{W_l}$ and $\vec{b_l}$ are the weight matrix and bias of the $l$th layer, and inputs to the NN are specified via the 0th layer $\vec{v_0}$. The activation function for the output is taken to be the identity function leading to a linear last layer $v_L = \vec{W_{L-1}} \vec{v_{L-1}} +\vec{b_{L-1}}$.\\
Physics-informed neural networks (PINNs) provide a possible solution approach, which trains the NN to learn the map from inputs $\vec{v_0} = (t,\vec{x},\vec{\alpha})$ to obtain the NN approximation $\hat{u}_\text{PINN}(t,\vec{x}, \vec{\alpha};\vec{\theta})= v_L$. The challenge lies in finding the NN parameter vector $\vec{\theta}$, which is the collection of all weights and biases, such that the PINN  captures the full solution manifold for all time and parameter instances.

\begin{figure}[h]
    \centering
 \subcaptionbox{PINN}
         {\includegraphics[width=0.32\textwidth]{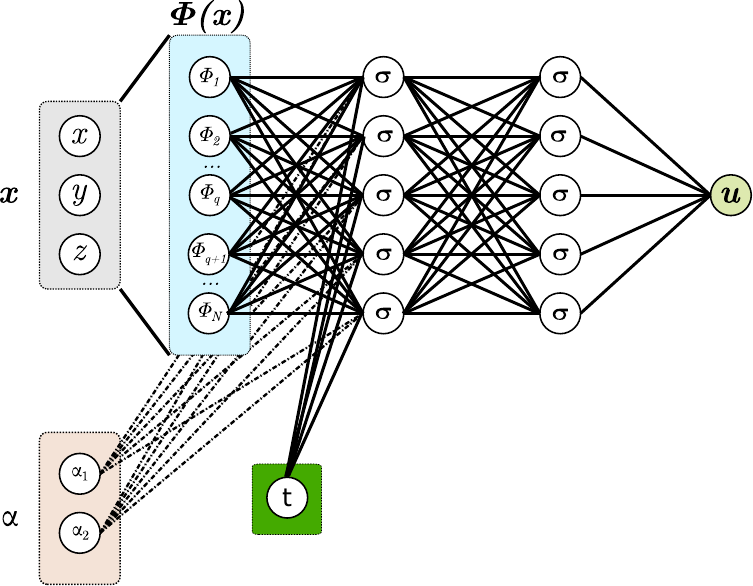}}
                           \hfill
          \subcaptionbox{EDNN}
         {\includegraphics[width=0.32\textwidth]{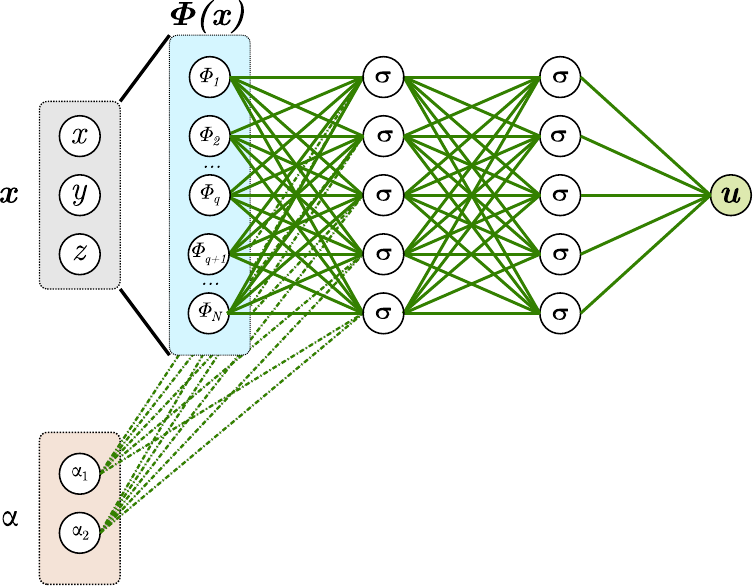}}
                  \hfill
         \subcaptionbox{RB/Linear}
         {\includegraphics[width=0.32\textwidth]{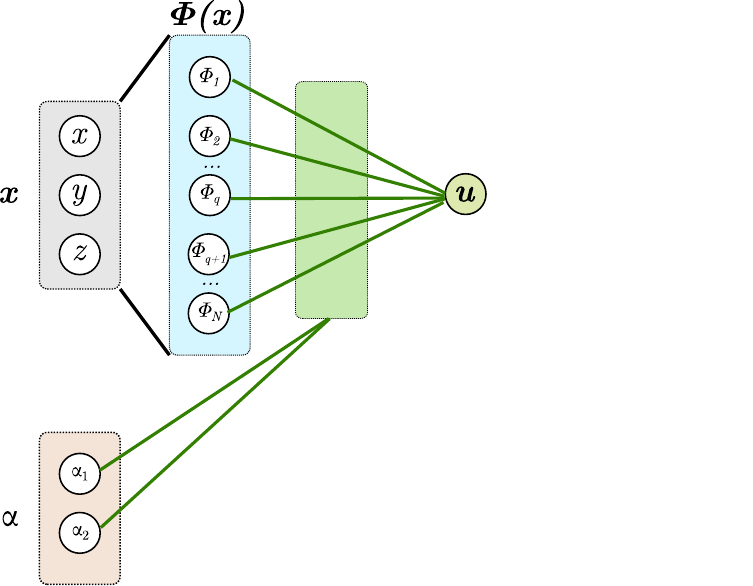}}
    \caption{Overview of different numerical models. Green color indicates that a NN parameter or RB coefficient evolves in time. The green box in the RB model (c) symbolizes the computation of the RB coefficient vector in \eqref{eq:RB}. }
    \label{fig:overview}
\end{figure}

This is usually achieved via the optimization of a  loss function, which incorporates the PDE dynamics, the fit of the initial condition and the boundary conditions as defined in \eqref{eq:PDE}, and is evaluated at (randomly sampled) $n_c = n_{\text{PDE}} +n_{\text{init}} +n_{\text{BC}}$ collocation points:
\begin{equation}
    \label{eq:PINN_loss}
    \begin{split}
    \mathcal{L} & = \sum_{i=1}^{n_{\text{PDE}}} \left(\dfrac{\pd \hat{u}(t_i,\vec{x}_i, \vec{\alpha}_i;\vec{\theta})}{\pd t} - f(\hat{u},t,\vec{x}_i, \vec{\alpha}_i)\right)^2  +  \sum_{i=1}^{n_{\text{init}}}  \left(\hat{u}(0,\vec{x}_i,\vec{\alpha}_i;\vec{\theta} ) -  u_0(\vec{x}_i,\vec{\alpha}_i) \right)^2 \\
    &+ \sum_{i=1}^{n_{\text{BC}}}  \left(\mathcal{B} \hat{u}(t_i,\vec{x}_i,\vec{\alpha}_i;\vec{\theta}) - g(t_i,\vec{x}_i,\vec{\alpha}_i)\right)^2 
    \end{split}
\end{equation}
The loss in \eqref{eq:PINN_loss} is usually minimized via gradient-based methods such as ADAM, where  automatic differentiation (AD) is used to evaluate \eqref{eq:PINN_loss} and the  loss gradient $ \frac{\partial \mathcal{L}}{\partial \vec{\theta}}$. It is well known that this training process requires careful balancing of the loss terms and may fail to converge to acceptable results \cite{wang_understanding_2021, wang_eigenvector_2021} unless all hyperparameters are tuned carefully.  \\
We propose two major improvements to facilitate the learning of parameters $\vec{\theta}$: 
Firstly, we introduce positional embeddings $\Phi(\vec{x}) = \left[ \phi_1(\vec{x}),\phi_2(\vec{x}),...,\phi_{n_{\phi}}(\vec{x}) \right] $ for the spatial coordinate, which adapt the NN to the targeted solution manifolds and automatically enforce the boundary conditions.
Secondly, we use evolutional deep neural  networks (EDNNs), where only the spatial and parametric coordinates form the input vector of the NN,  $\vec{v_0}= (\Phi(\vec{x}),\vec{\alpha})$,  and the neural network parameters are updated in time: $\hat{u}_\text{EDNN}(\Phi(\vec{x}), \vec{\alpha};\vec{\theta}(t))$.  This turns the learning problem into finding the time trajectory of NN parameters $\vec{\theta}(t)$. Figure \ref{fig:overview} provides a visual comparison of the RB, PINN and EDNN approach. \\

\textbf{Remark:} We note an interesting connection between the  EDNN approach and the reduced basis formulation  in \eqref{eq:RB}. Let $ (\vec{v_{L-1}})_i=\Gamma_i(\Phi(\vec{x}), \vec{\alpha}; \vec{\theta} (t))$ be the intermediate values at the second to last layer in \eqref{eq:NN}, and write the NN output as 
\begin{equation}
    \hat{u}_\text{EDNN}(\Phi(\vec{x}), \vec{\alpha};\vec{\theta}(t)) = \sum_{i=1}^{n_\text{hid}} \beta_i (t) \Gamma_i(\Phi(\vec{x}), \vec{\alpha}; \vec{\theta} (t)).
\end{equation}
One can then interpret the positional embedding $\Phi$ of the EDNN approach as a reduced basis, which is evolved in time with the transformation $\Gamma$. This time adaptation of the basis allows to handle problems with slow-decaying Kolmogorov n-width.

\section{Positional embeddings}
In this section, we will show how a positional embedding layer for a NN can be used to directly enforce the boundary conditions.
We extend the work in \cite{du_evolutional_2021} as our framework can be used for Dirichlet and Neumann boundary conditions on spatial domains with arbitrary geometry.
The positional embeddings are  constant in time, consequently their computation can be treated as a pre-processing step and all necessary numerical operations can happen in an "offline" phase.

\subsection{Variational harmonic features}
Variational harmonic features have previously been used for Gaussian process regression \cite{solin_hilbert_2020, solin_know_2019} and are computed via the eigenfunctions of the Laplace (Beltrami) operator:
\begin{equation}
\label{eq:LEigs}
    \begin{split}
         -\nabla^2 \phi(\vec{x}) &=\lambda \phi(\vec{x}),\quad \vec{x} \in \Omega_x, \\
         \mathcal{B} u(\vec{x}) &= 0, \quad \vec{x} \in \pd \Omega_x.
    \end{split}
\end{equation}
When solving PDEs, the geometric shape and boundary conditions often encode key features of the solution, and variational harmonic features allow us to adapt the neural network to these characteristics. For certain geometries,  the solutions to \eqref{eq:LEigs} are known analytically, while for domains with arbitrary shapes,  no closed form solutions to \eqref{eq:LEigs} are known and the eigenfunctions need to be computed numerically.  \\
The embedding with  (Gaussian) random Fourier features \cite{wang_eigenvector_2021,tancik_fourier_2020} takes the form
\begin{equation}
\label{eq:Fourierfeat}
\phi_i(\vec{x}) = [\cos(\vec{b_i}^\top \vec{x}), \sin(\vec{b_i}^\top \vec{x})], \quad \text{with } (\vec{b}_i)_{l} \sim \mathcal{N}(0,\sigma^2),
\end{equation}
where the vectors $\vec{b_i}$ are sampled from a normal distribution according to a lengthscale parameter $\sigma$. This  also fits within this more general framework: Fourier features correspond to the solutions of \eqref{eq:LEigs} for the infinite domain ("no boundary conditions"). Indeed, for computer vision tasks, Fourier features and the induced spatially invariant NTK, are a natural choice as images and videos are usually taken of a continuous real world scene. 
In contrast to this, variational harmonic features impose a spatially-varying NTK, accounting for general features, e.g. holes, in the geometry.\\

\textbf{Remark:} By restricting to a finite  domain, the solution spectrum of \eqref{eq:LEigs} becomes discrete, and the features are not sampled from a continuous probability measure in contrast to \eqref{eq:Fourierfeat}. We chose to select the $n_\Phi$ eigenfunctions associated with the $n_\Phi$ lowest eigenvalues. In certain applications e.g. multi-scale phenomena, other choices may be more appropriate \cite{wang_eigenvector_2021}.

\subsection{Approximation space and boundary conditions}
In this section, we seek to  qualitatively understand which functions we can approximate for a NN with embedding $\Phi(\vec{x})$, i.e. %
we wish to characterize the solution manifold  $\hat{\mathcal{M}}$ of $\hat{u} (\Phi(x);\vec{\theta})$.
It has previously been demonstrated in \cite{du_evolutional_2021, bruna_neural_2022}, that $\hat{u}(\Phi(x))$ inherits the periodicity when the embedding functions $\Phi (x)$ are periodic. We now extend this to embeddings with Neumann and Dirichlet boundary conditions.

\subsubsection{Imposing boundary conditions}
Let $\mathcal{B}_{\text{NM}} u(\vec{x}) = \nabla_{\vec{x}} u \cdot \vec{n}(\vec{x})$, be the  operator acting on the Neumann boundary $\pd \Omega_{\text{NM}} \subset \pd \Omega_x$, where $\vec{n}(\vec{x})$ is the unit vector normal to the boundary. Let  $\mathcal{B}_{\text{DC}} u(\vec{x})= u(\vec{x}) $ be the operator acting on the Dirichlet boundary  $\pd \Omega_{\text{DC}} \subset \pd \Omega_x$. Let $\Phi(x)$ be an embedding, where each $\phi_i(\vec{x})$ is a solution to \eqref{eq:LEigs} such that
\begin{equation}
\begin{split}
\mathcal{B}_{\text{NM}} \phi_i(\vec{x}) = 0 \quad   & \forall \vec{x} \in \pd \Omega_{\text{NM}}, \\
\mathcal{B}_{\text{DC}} \phi_i(\vec{x}) = 0 \quad   &\forall \vec{x} \in \pd \Omega_{\text{DC}}.
\end{split}
\end{equation}
For the NN surrogate $\hat{u} (\Phi(\vec{x}))$, it then holds that
\begin{itemize}
\item $\mathcal{B}_{\text{NM}} \hat{u}(\Phi(\vec{x}))=0$ $\forall \vec{x} \in \pd \Omega_{\text{NM}}$, i.e. the neural network  will inherit the Neumann boundary conditions of the Laplace eigenfunctions. The proof follows directly from applying the chain rule:
\begin{equation}
    \mathcal{B}_{\text{NM}} \hat{u}(\Phi(\vec{x})) =  \nabla_{\vec{x}} \hat{u} \cdot \vec{n}(\vec{x}) = \sum_{i=1}^{n_\Phi} \dfrac{\pd \hat{u}(\Phi(\vec{x}))}{\pd  \phi_i(\vec{x})}  \underbrace{\nabla_{\vec{x}}\phi_i(\vec{x}) \cdot \vec{n}(\vec{x})}_{=0} = 0  \quad \text{for } \vec{x} \in \pd \Omega_{\text{NM}}.
\end{equation}
\item $\mathcal{B}_{\text{DC}} \hat{u}(\Phi(\vec{x}))= \mathcal{B}_{\text{DC}} \hat{u}(\vec{0}) =c$, for $\vec{x} \in \pd \Omega_{\text{DC}}$, i.e. all points on the Dirichlet boundary are mapped to the same constant $c$.
\end{itemize}
We can use this to our advantage and construct an auxiliary network $\tilde{u}$, that automatically enforces homogeneous Dirichlet (and Neumann) boundary conditions :
\begin{equation}
\label{eq:homoBC}
    \tilde{u}(\Phi(\vec{x})) =  \hat{u}(\Phi(\vec{x})) - \hat{u}(\vec{0}),
\end{equation}
so that
\begin{equation}
    \mathcal{B}_{\text{DC}} \tilde{u}(\Phi(\vec{x})) = c-c=  0 \quad \forall \vec{x} \in \pd \Omega_{\text{DC}}.
\end{equation} \\
When inhomogeneous boundary conditions need to be prescribed, one can add a lifting term to \eqref{eq:homoBC}. An appropriate embedding can thus eliminate the boundary condition loss term in \eqref{eq:PINN_loss} and simplify the optimization problem.

It is essential to acknowledge that the proposed embedding restrict the solution space of the neural network. Using an embedding that automatically enforces Neumann boundary conditions would be a poor choice to represent a PDE problem with homogeneous Dirichlet boundary conditions as the gradients along the boundary will be forced to 0.
 Fourier features \eqref{eq:Fourierfeat} form a larger solution space and can be used to approximate any boundary condition. We will further illustrate this generality vs. specificity trade-off in an example in Section \ref{sec:illu}. \\
\textbf{Remark:} One can compute Laplace eigenfunctions with Robin boundary conditions, but this does not eliminate a loss term as the NN does not inherit the Robin boundary condition of the embedding. %

\subsubsection{Selecting basis functions}
We provide some insights into how the size of the embedding affects the approximation space/solution manifold of the NN $\hat{u}\left(\Phi (\vec{x})\right)$.
For a bounded domain $\Omega_x$, it is well established that the eigenfunctions of the Laplace operator form an orthogonal adn complete basis  of the Hilbert space $L_2(\Omega_x)$ \cite{solin_hilbert_2020}. In practice, we will need to truncate this basis to a finite number of embedding functions $n_\Phi$, and $\Phi(\vec{x})$ will not be a complete basis for $L_2(\Omega)$. For reduced basis methods, this truncation error is well understood \cite{hesthaven_reduced_2022}, but the inclusion of nonlinear layers in the neural network enhances the expressivity of its solution manifold and demands a shift in perspective: 

To enforce boundary conditions, it is sufficient to include a single embedding function $\phi_i(\vec{x})$. This may, however, restrict the solution manifold of the NN due to symmetries in $\phi_i(\vec{x})$.  
A trivial example is given by the first Laplace eigenfunctions on the unit line with Dirichlet boundary conditions: $\phi_1(x)= \sin(\pi x) ,\phi_2(x) =\sin(2 \pi x)$. It is impossible to find a mapping $\hat{f}$, such that
\begin{equation}
    \hat{f}\left(\phi_1(x)\right)=\phi_2(x)
\end{equation}
 due to the reflection symmetry around $x=0.5$, e.g. $\phi_1 \left(\frac{1}{4}\right)= \phi_1 \left(\frac{3}{4}\right)$, but $\phi_2(\frac{1}{4}) = -\phi_2(\frac{3}{4})$. This implies that  a NN $\hat{u}(\phi_1(x))$ can only approximate functions with the same reflection symmetry. Conversely, a NN with only two embedding functions $\Phi_2(\vec{x})= [\phi_1(x), \phi_2(x)]$ uniquely identifies each point on the unit line, and consequently  the approximation of 
 $L_2(\Omega)$ via the solution manifold of the NN is only limited by the size of the parameter space of the NN.
 
The Neumann case on the unit line does not suffer from symmetries, in fact we know from the definition of the Chebyshev polynomials that there exists a transformation $T_n$, such that
\begin{equation}
    T_n \left(\cos(\pi x) \right) = \cos(n\pi x),
\end{equation}
i.e. an embedding of size one is sufficient to recover the solution space of interest. 

In higher dimensions, it is more difficult to identify symmetries, but similar considerations apply: as long as the embedding is rich enough to uniquely identify each point in the spatial domain $\Omega_x$, the nonlinear approximation capabilities of the neural network allow for general function approximation. This also implies that  any collection of functions $\phi_i$ can be used as a positional embedding $\Phi(x)$ as long as they encode the boundary conditions and form a rich enough space. It is thus  possible to use functions obtained from a reduced basis approximation, the eigenfunctions of another operator, or directly learn the embeddings from data. We focus on the eigenfunctions of the Laplace operator due to the close connection with Fourier features and their favorable influence on training \cite{wang_eigenvector_2021}.

 \subsection{Numerical approximation}
Using the Ansatz introduced in \eqref{eq:FEAnsatz}  we compute solutions to \eqref{eq:LEigs} using a FE method. Let
\begin{equation}
\label{eq:hatphi}
    \hat{\phi}(\vec{x}) = \sum_{i=1}^{N} c_i  \varphi_i (\vec{x}),
\end{equation}
be the trial functions and $v$ be test functions of the same FE space with $C^1$ continous elements for the weak formulation of \eqref{eq:LEigs}:
\begin{equation}
\label{eq:weakform}
    \int_{\Omega_x } \nabla  \hat{\phi} \cdot  \nabla v \dd \vec{x}=  \lambda \int_{\Omega_x}  \hat{\phi} v \dd \vec{x}, \forall v.
\end{equation}
\\
To utilize \eqref{eq:hatphi} as a positional embedding for a PINN or EDNN, we need to compute its gradients with respect to the spatial coordinate $\vec{x}$.
As the gradient of the FE function $\hat{\phi}(\vec{x})$
\begin{equation}
\label{eq:hatphigrad}
   \nabla_{\vec{x}} \hat{\phi}(\vec{x}) = \sum_{i=1}^{N} c_i    \nabla_{\vec{x}} \varphi_i (\vec{x}),
\end{equation}
is no longer $C^1$ continuous,  we project the gradient \eqref{eq:hatphigrad} back to the original FE space component-wise:
\begin{equation}
    \int_{\Omega_x } g_l   v \dd \vec{x}=   \int_{\Omega_x}  (\nabla_x \hat{\phi})_l v \dd \vec{x},
\end{equation}
so that
\begin{equation}
     \nabla_{\vec{x}} \hat{\phi_i}(\vec{x}) \approx \vec{g_i} (\vec{x}).
\end{equation}
Higher order derivatives can then be evaluated recursively. To control the projection error, we solve \eqref{eq:weakform} on a fine mesh and use P3 Lagrange elements. Assembling the system  leads to the discrete linear eigenvalue problem
\begin{equation}
    K \vec{c} = \lambda M \vec{c},
\end{equation}
which can be solved  with standard linear algebra routines. \\

To compute spatial derivatives of the neural network, we  use the chain rule to combine the spatial FE gradient of the embedding function  and the automatic differentiation  gradient (AD) of the neural network with respect to its inputs $\phi_i(\vec{x})$:
 \begin{equation}
 \label{eq:chainrule}
      \dfrac{\pd \hat{u} \left(\Phi(\vec{x})\right)}{\pd x_l} = \sum_{i=1}^{n_\Phi} \underbrace{\dfrac{\pd \hat{u}(\Phi(\vec{x}))}{\pd  \phi_i(\vec{x})}}_{\text{AD}} \underbrace{\dfrac{\pd \phi_i(\vec{x}) }{ \pd x_l}}_{\text{FE}} \approx \sum_{i=1}^{n_\Phi} \underbrace{\dfrac{\pd \hat{u}(\hat{\Phi}(\vec{x}))}{\pd  \hat{\phi}_i(\vec{x})}}_{\text{AD}} \left(\vec{g_i} (\vec{x})\right)_l .
 \end{equation}
Computing second order spatial derivatives requires slightly more effort -- applying the chain rule twice leads to
\begin{equation}
\label{eq:doublechain}
       \dfrac{\pd^2 \hat{u}}{\pd x_l \pd x_k}= \sum_{i=1}^{n_\Phi} \left(   \underbrace{\dfrac{\pd \hat{u}(\Phi(\vec{x}))}{\pd  \phi_i(\vec{x})}}_{\text{AD}} \underbrace{\frac{\pd^2\phi_i(\vec{x})}{\pd x_l \pd x_k}}_{\text{FE}}  + \sum_{j=1}^{n_\Phi}  \left(\underbrace{\dfrac{\pd^2 \hat{u}(\Phi(\vec{x}))}{\pd  \phi_i(\vec{x}) \pd \phi_j(\vec{x})} \frac{\pd\phi_i(\vec{x})}{\pd x_l}}_{\text{AD of \eqref{eq:chainrule}}} \frac{\pd\phi_j(\vec{x})}{\pd x_k}   \right) \right).
\end{equation}
The first term of \eqref{eq:doublechain} can be computed similarly to \eqref{eq:chainrule}, where  the spatial Hessian $\left((H_i (\vec{x})\right)_{lk} \approx \frac{\pd^2\phi_i(\vec{x})}{\pd x_l \pd x_k}$ is computed with a FE method offline. To avoid an explicit online assembly of the Hessian $\dfrac{\pd^2 \hat{u}(\Phi(\vec{x}))}{\pd  \phi_i(\vec{x}) \phi_j(\vec{x})}$ in the second term and realizing that  AD  treats the FE terms  $ \left(\vec{g_i} (\vec{x})\right)_l$ as constant, we rewrite \eqref{eq:doublechain} as:

\begin{equation}
\label{eq:doublechainre}
       \dfrac{\pd^2 \hat{u}}{\pd x_l \pd x_k}= \sum_{i=1}^{n_\Phi} \left(  \underbrace{ \dfrac{\pd \hat{u}(\Phi(\vec{x}))}{\pd  \phi_i(\vec{x})}}_{\text{AD}} \left((H_i (\vec{x})\right)_{lk} + \sum_{j=1}^{n_\Phi}  \left( \underbrace{\dfrac{\pd  \left(\frac{\pd \hat{u} \left(\Phi(\vec{x})\right)}{\pd x_l}\right) }{\pd \phi_j(\vec{x})}}_{\text{AD}} \left(\vec{g_i} (\vec{x})\right)_k \right) \right),
\end{equation}
which reduces the number of AD calls from $O(n_\Phi^2)$ to $ O(n_\Phi)$ (albeit nested) calls.
Similar nesting strategies can be exploited for higher order derivatives. 
 We stress that using a mesh-based approach does not restrict our function evaluation to the nodes of the mesh: FE methods allow for solution evaluation at any point in the spatial domain. In practice, we compute the eigenfunctions and their derivatives at a pre-selected number of collocation points during the offline phase and then subsample this set at the desired resolution in the online phase. \\
 \subsection{An illustrative example}
 \label{sec:illu}
We shortly illustrate the benefits of the harmonic feature embeddings on a 2D example.  Consider the static PDE
\begin{equation}
\begin{split}
& \nabla \cdot \left( a(\vec{x}) \nabla u \right) =1,\\
&    u(\vec{x}) =0 \quad\forall \vec{x} \in \pd\Omega_{DC}, \\
&    a(\vec{x}) = \exp^{ - (x-0.25)^2 - (y-0.25)^2},
\end{split}
\end{equation}
on the 2D unit square with a hole. We solve the PDE with a NN without an embedding layer, an embedding layer with 10 Fourier features ("FF") and an embedding layer with 10 harmonic features ("HF"). The NN has 4 layers with 10 hidden units each ($n_\theta$=3760), uses a  tanh activation function and is trained on the PDE residual (and boundary condition loss) with ADAM.

\begin{figure}[h]
 \centering
 \subcaptionbox{mesh}
         {\includegraphics[width=0.2\textwidth]{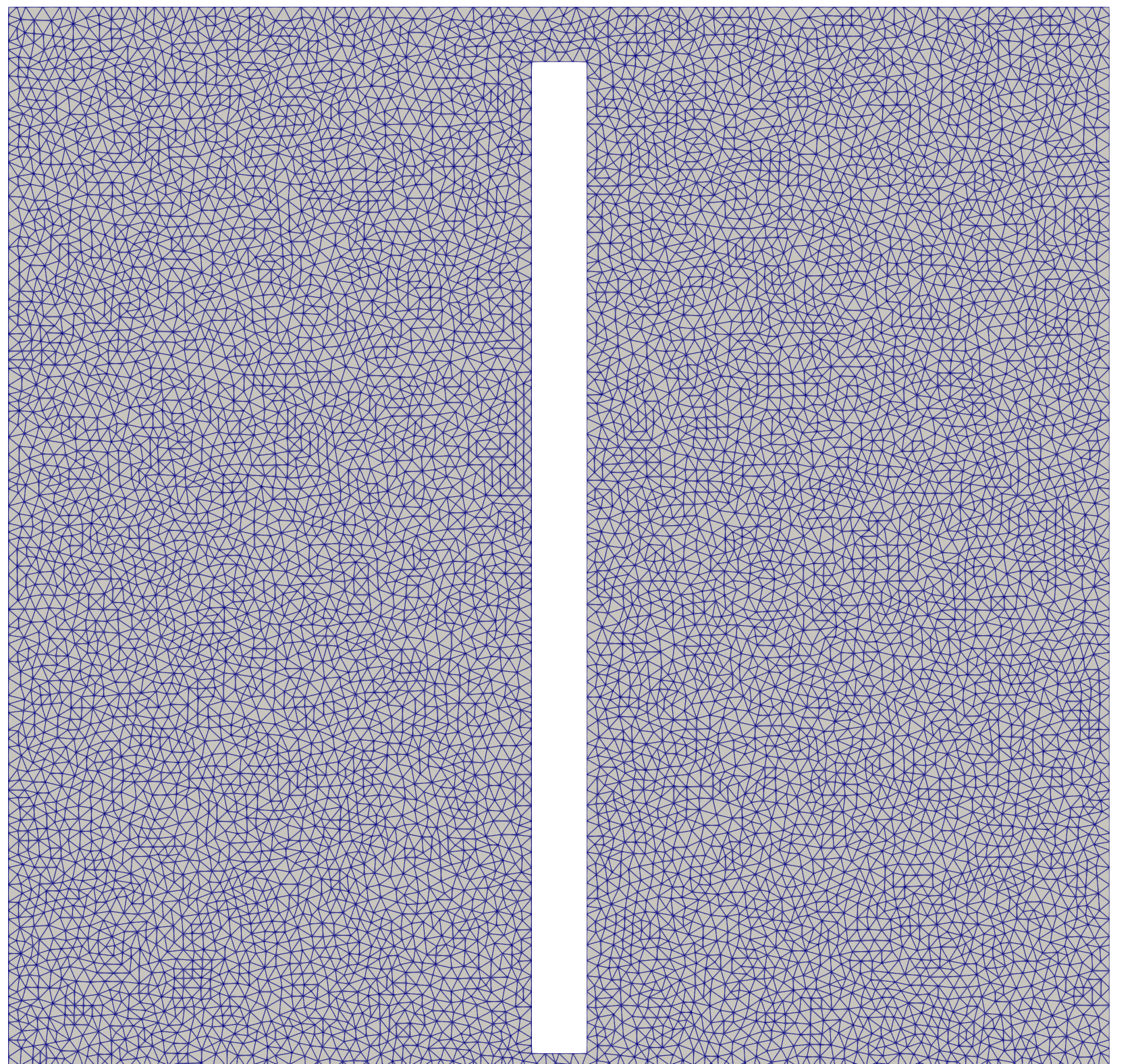}}
  \subcaptionbox{$\hat{\phi}_1$}
         {\includegraphics[width=0.2\textwidth]{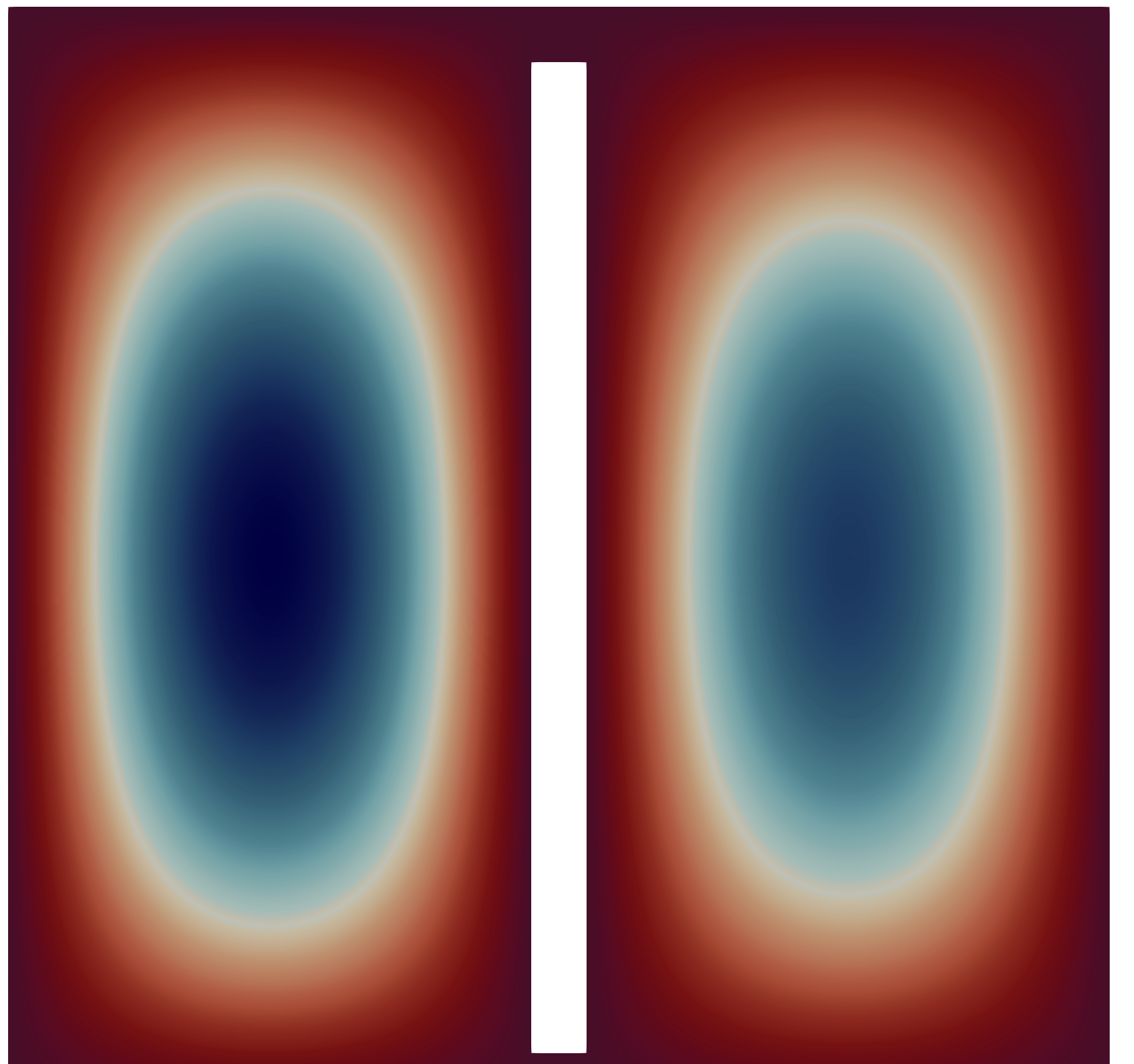}}
           \subcaptionbox{$\hat{\phi}_2$}
         {\includegraphics[width=0.2\textwidth]{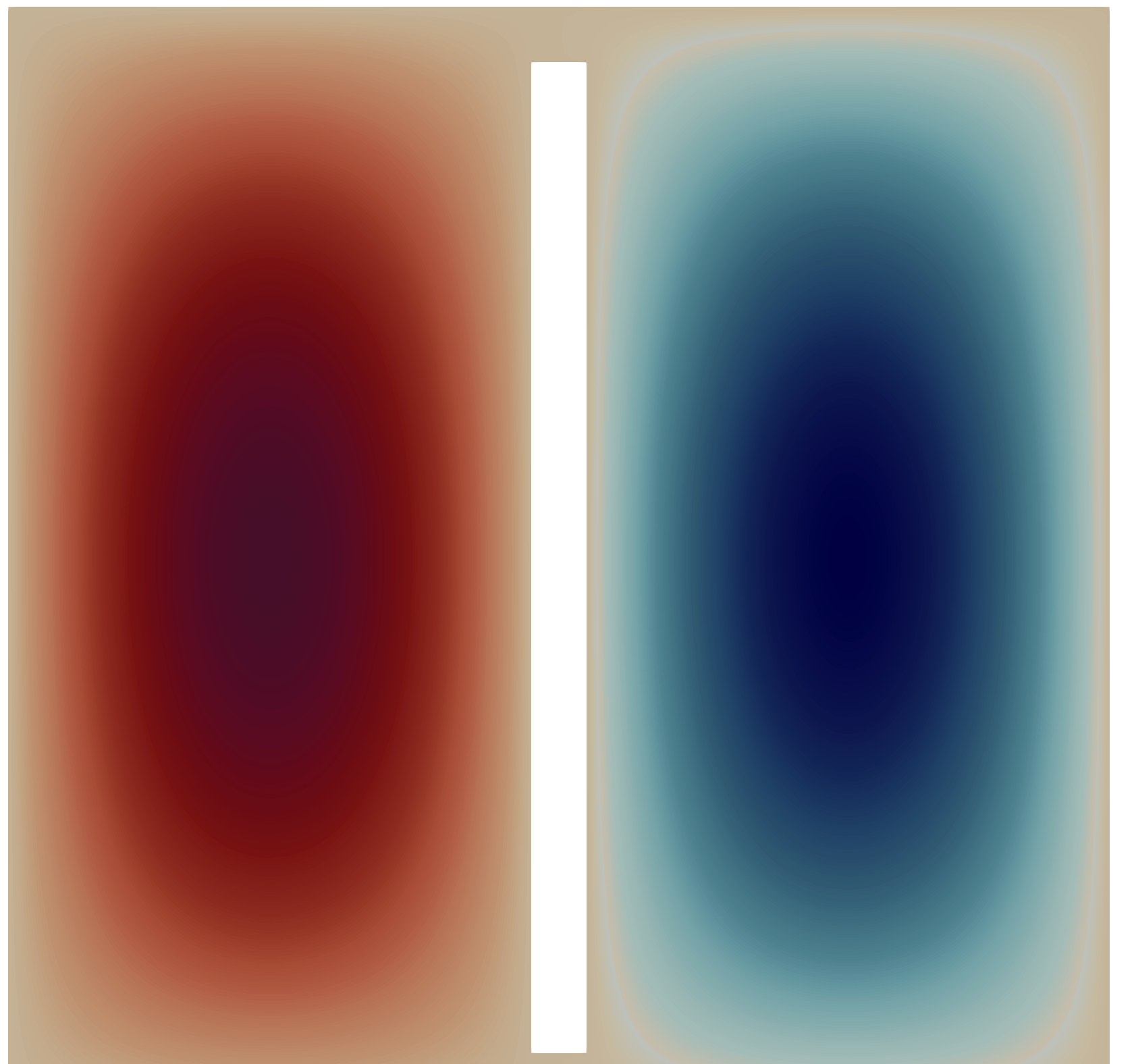}} \\
          \subcaptionbox{$\hat{\phi}_3$}
        {\includegraphics[width=0.2\textwidth]{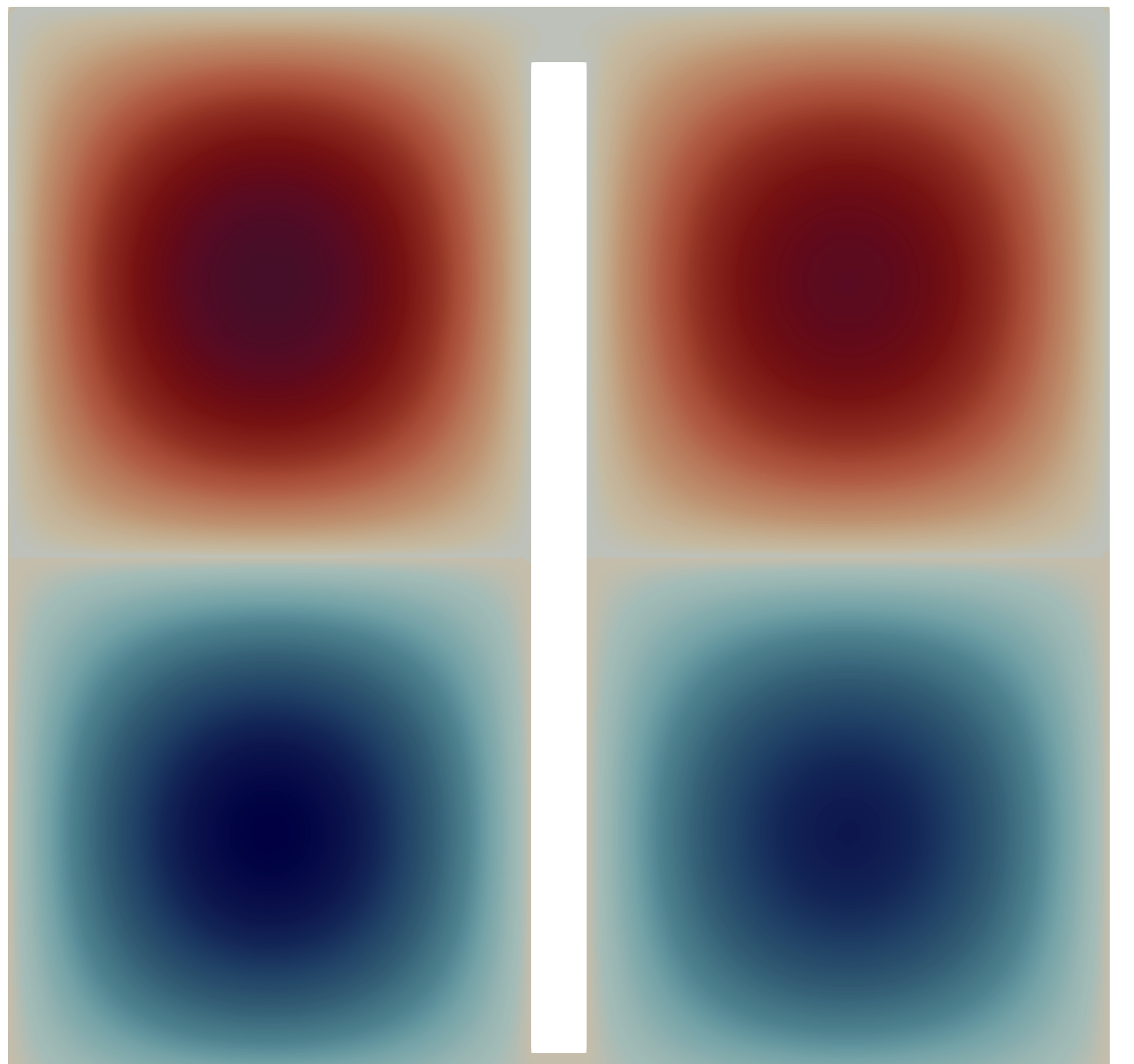}}
                   \subcaptionbox{$\dfrac{ \pd\hat{\phi}_3}{\pd x}$}
         {\includegraphics[width=0.2\textwidth]{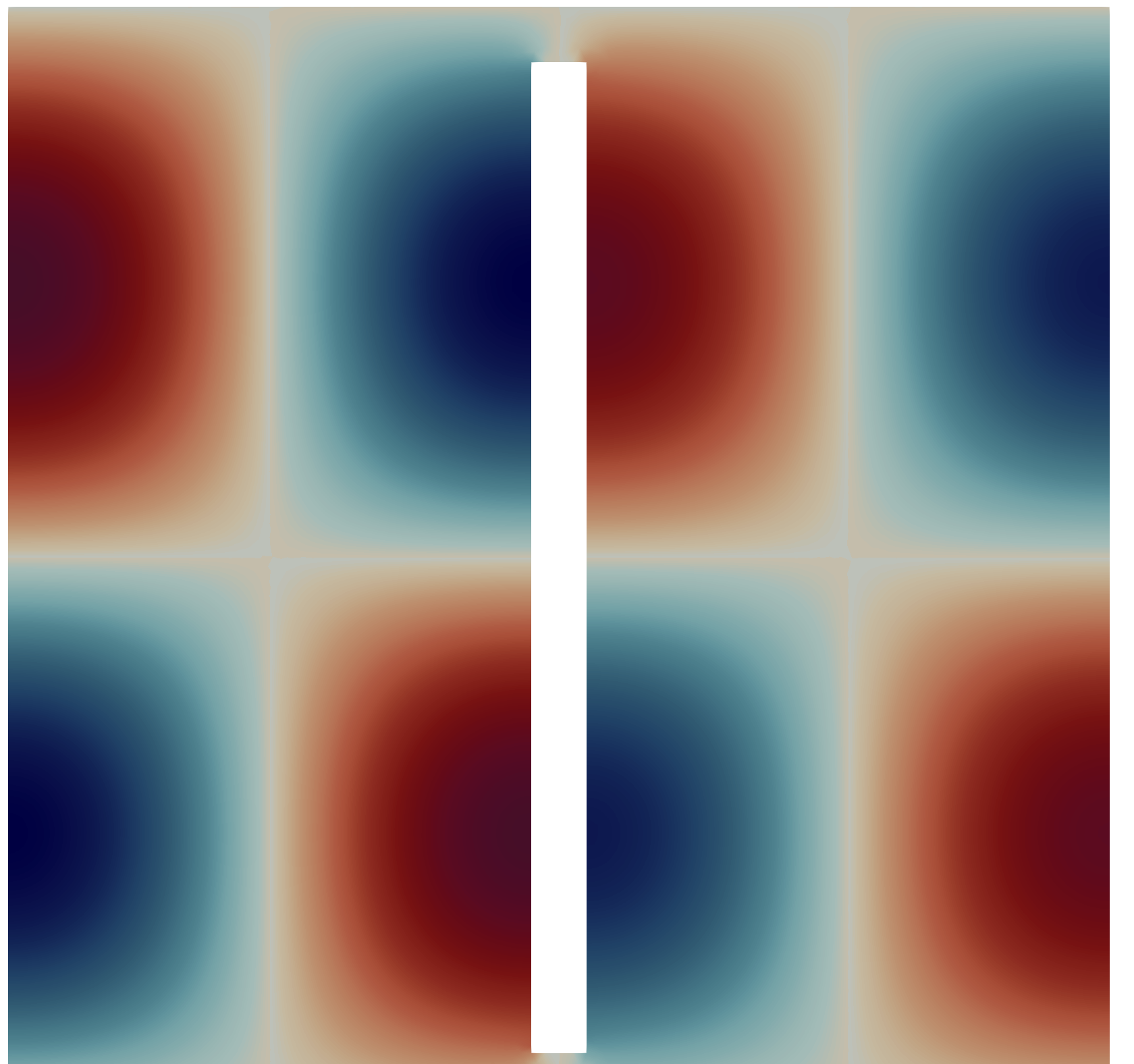}}
                   \subcaptionbox{$\dfrac{ \pd\hat{\phi}_3}{\pd y}$}
         {\includegraphics[width=0.2\textwidth]{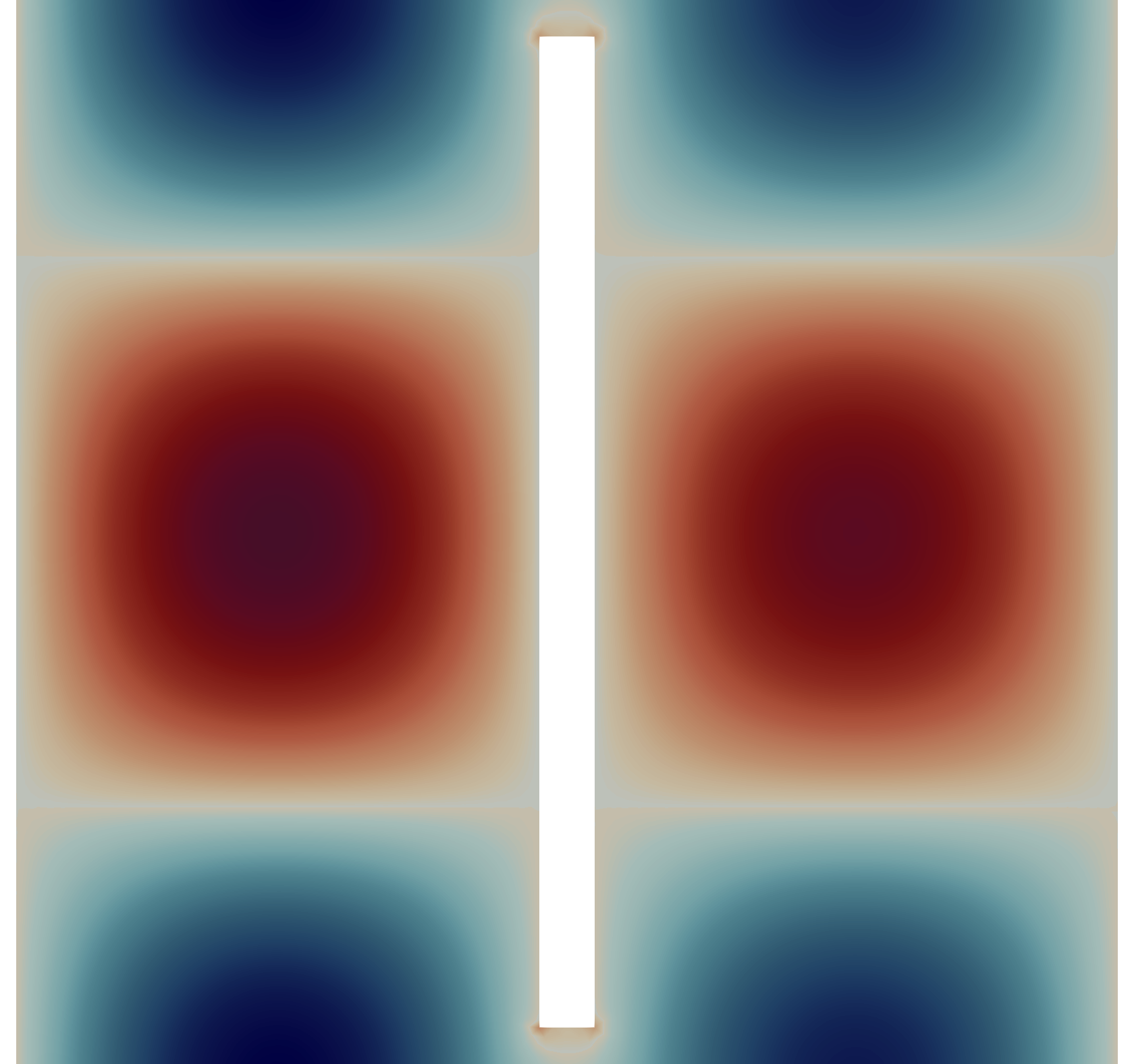}}
\caption{Variational harmonic features computed with a FE method for the static PDE example.}
\label{fig:2Dmodes}
\end{figure}

\begin{figure}[h]
 \centering 
   \begin{minipage}[b]{0.5\linewidth}
  \subcaptionbox{Convergence plot}{\includegraphics[width=1\textwidth]{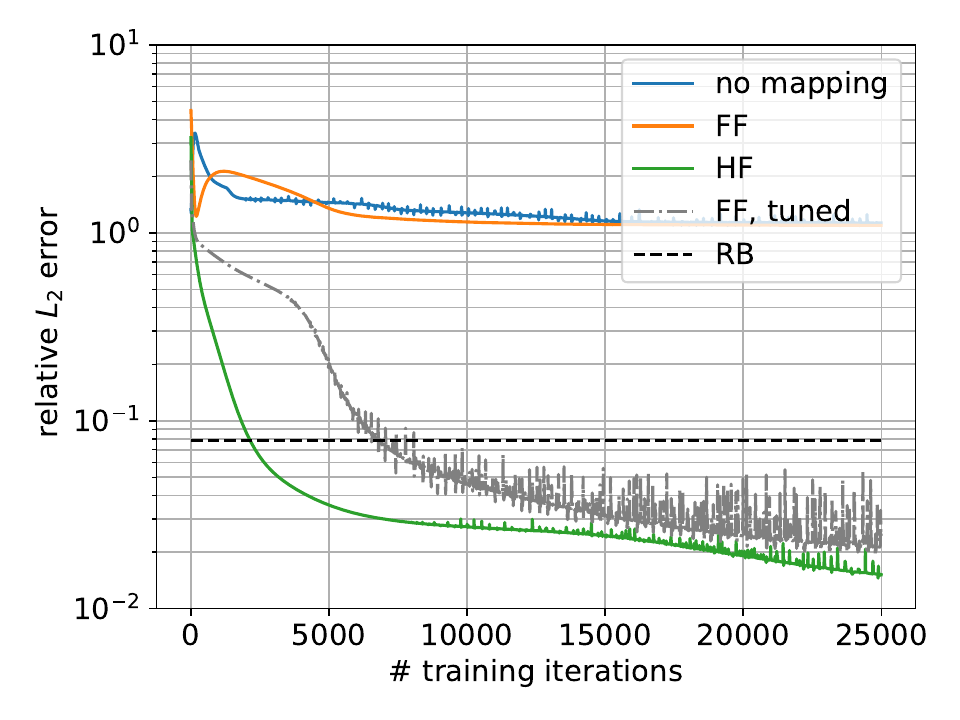}}
  \vfill
\end{minipage}
           \begin{minipage}[b]{0.37\linewidth}
 \subcaptionbox{FE solution}
          {\includegraphics[width=0.49\textwidth]{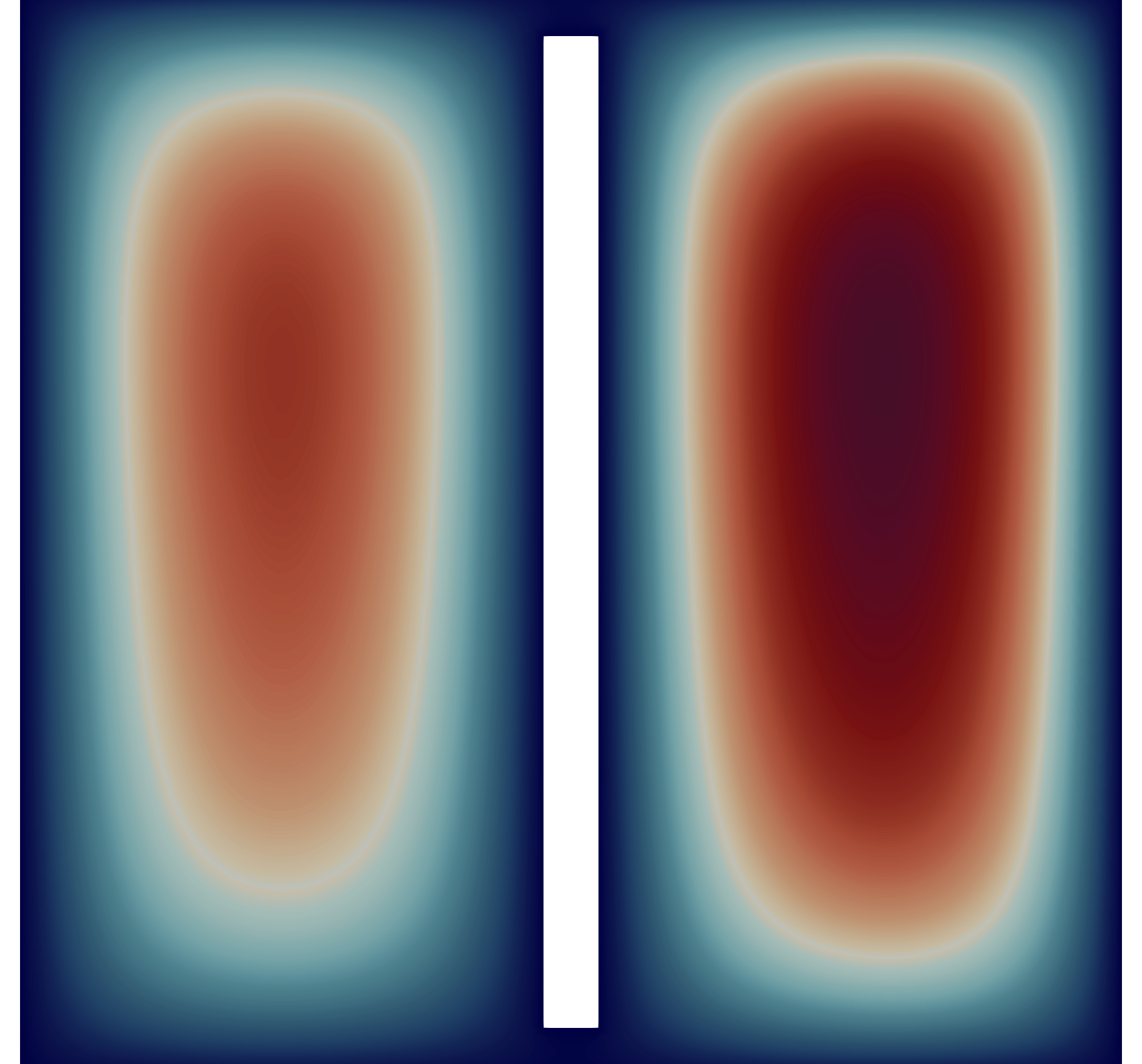}}
  \subcaptionbox{NN with FF}
         {\includegraphics[width=0.49\textwidth]{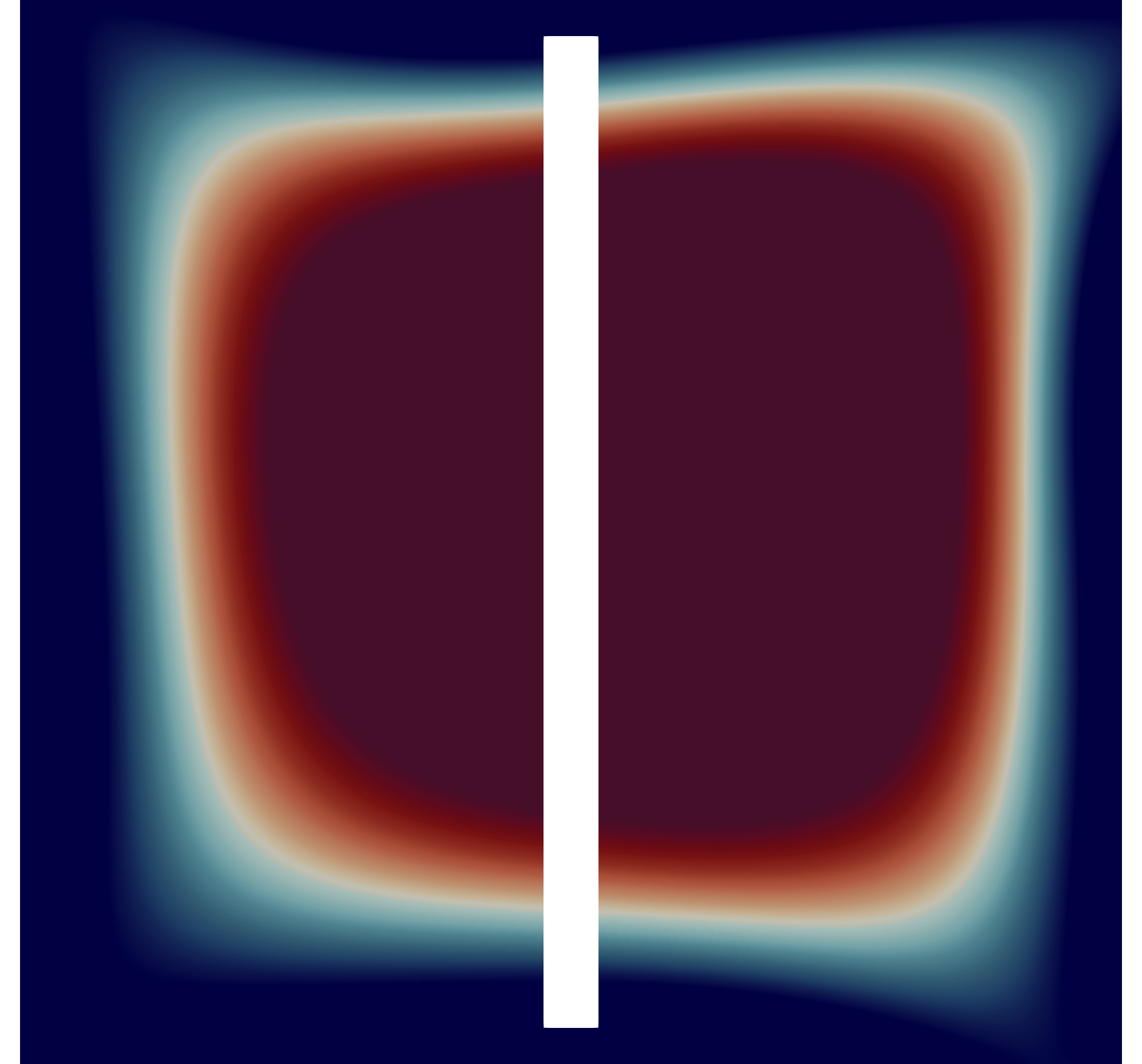}}
           \subcaptionbox{NN with HF}
         {\includegraphics[width=0.49\textwidth]{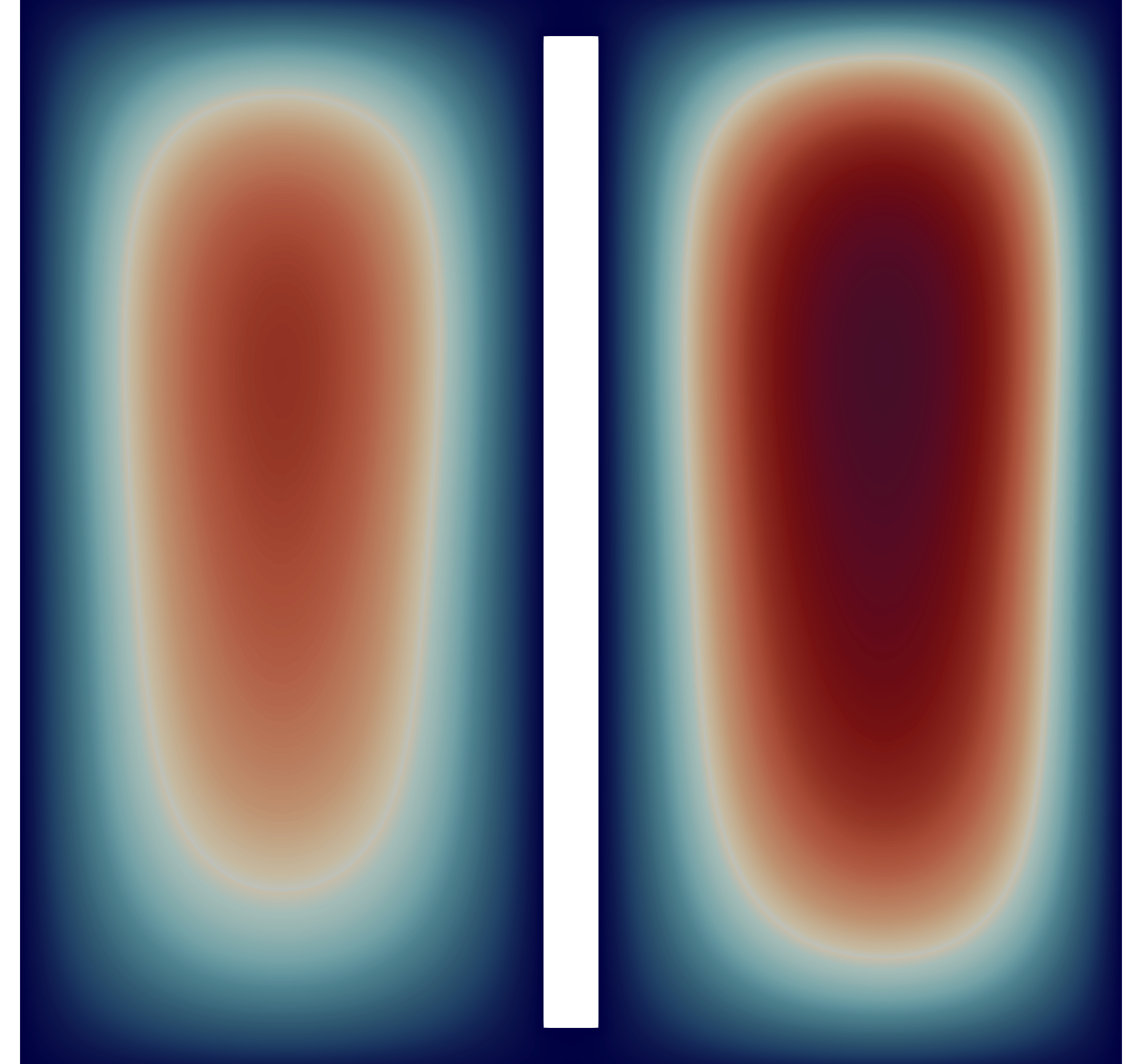}} 
         \centering
         \includegraphics[width=0.15\textwidth]{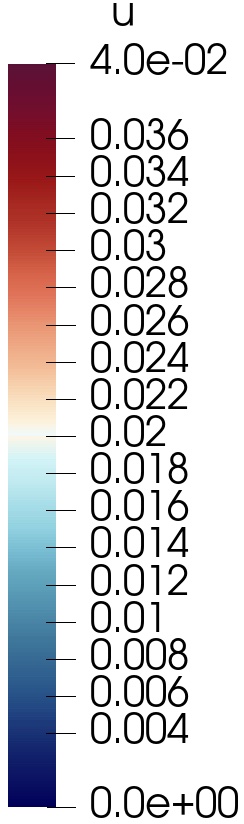}
  \end{minipage}
\caption{Solution plots and error for the static PDE example.}
\label{fig:static}
\end{figure}

In  Figure \ref{fig:2Dmodes}, we visualize the mesh and the computed embeddings.  In Figure \ref{fig:static}, we  visualize the solution and convergence of the relative $L_2$ error 
\begin{equation}
\label{eq:L2error}
    \varepsilon = \dfrac{|| \hat{u}_{\text{NN}}  -\hat{u}_{\text{FE}} ||_2}{||  \hat{u}_{\text{FE}} ||_2}.
\end{equation}

We also provide the projection error of the embedding as a reference, i.e. we compute the error
\begin{equation}
    \varepsilon_\text{proj} = ||\vec{\Phi} {\vec{\Phi}}^\top \hat{u}_\text{FE} - \hat{u}_{\text{FE}} ||_2 
\end{equation}
to confirm that the use of  a neural network enriches the solution manifold compared to the linear subspace.
We observe that variational harmonic features have the fastest convergence and achieve the lowest errors. As expected, the NN without the input mapping fails to converge to the true solution. Notably, the same happens for the NN with Fourier features, unless one carefully balances the weights in the loss term ("tuned" case). 
These results demonstrate that the training process is a main challenge for obtaining good solutions.

\section{Time-stepping with neural networks}
In this section, we derive the time stepping scheme. We drop the explicit parameter dependence of the problem for conciseness: without loss of generality we can write $\tilde{\vec{x}} = [\vec{x}, \vec{\alpha}]$, $\Omega = \Omega_x \times \Omega_\alpha$, and continue the derivation with respect to the combined vector $\tilde{\vec{x}}$. 
\subsection{Deriving an explicit update equation}
To extract the temporal dependence of the neural network on the NN parameter vector $\vec{\theta}$, we apply the chain rule 
\begin{equation}
    \dfrac{\pd \hat{u}(\tilde{\vec{x}}; \vec{\theta} (t)) }{\pd t} = \underbrace{\dfrac{\pd \hat{u}(\tilde{\vec{x}}; \vec{\theta} (t)) }{\pd \vec{\theta} }}_{ J(\tilde{\vec{x}}) } \cdot  \dfrac{\pd \vec{\theta}}{\pd t} = J(\tilde{\vec{x}}) \vec{\dot{\theta}}, \quad \forall \tilde{\vec{x}} \in \Omega.
\end{equation}
Simply speaking, we seek the time derivative $\vec{\dot{\theta}}$, such that the  PDE dynamics \eqref{eq:PDE} are observed at every collocation point $\tilde{\vec{x}}$ of the computational domain. More formally, at each time $t$ we find $\vec{\dot{\theta}}$ as the solution to the minimization problem
\begin{equation}
\label{eq:minprob}
    \begin{split}
      &  \vec{\dot{\theta}} = \text{argmin}_{\vec{\gamma}} \underbrace{\frac{1}{2} \int_\Omega \left(J(\tilde{\vec{x}})\vec{\gamma} -  f\left(\hat{u}(\tilde{\vec{x}}; \vec{\theta}(t)),t, \tilde{\vec{x}}  \right)\right)^2 \dd \tilde{\vec{x}}}_{g(\gamma)}.
        \end{split}
\end{equation}
Applying the first order optimality condition $\nabla_{\vec{\gamma}} g(\vec{\gamma})=0$, yields an equation of  linear operators  to be solved for the parameter update:

\begin{equation}
    M\vec{\gamma}^* = F,
\end{equation}
where the matrix operator $M$ and vector operator $F$ are given by the  inner products:
\begin{equation}
\begin{split}
    M &= \int_\Omega  J(\tilde{\vec{x}})^\top  \cdot J(\tilde{\vec{x}}) \dd \tilde{\vec{x}}, \\
    F &= \int_\Omega  J(\tilde{\vec{x}})^\top  \cdot  f\left(\hat{u}(\tilde{\vec{x}}; \vec{\theta}(t)),t, \tilde{\vec{x}}\right) \dd \tilde{\vec{x}}.
    \end{split}
\end{equation}
We refer the reader to \cite{bruna_neural_2022} for a more thorough discussion of the theoretical background of this optimization problem and its relation to the Dirac-Frenkel variational principle.
In practice, $M$ and $F$ are approximated by randomly sampling $n_x$ collocation points, to obtain a Monte Carlo approximation of the integrals and a finite dimensional linear system
\begin{equation}
    M \approx \frac{1}{n_x}   \vec{J}^\top \vec{J} ,\quad    F \approx \frac{1}{n_x}   \vec{J}^\top \vec{f},
 \end{equation}
\begin{equation}
    \label{eq:linsys}
 \vec{J}^\top \vec{J} \gamma^* = \vec{J}^T \vec{f},
\end{equation}
where $\vec{J}$ is the discrete neural network Jacobian $(\vec{J})_{ik} = \dfrac{\hat{u}(\tilde{\vec{x}}_i; \vec{\theta} (t))}{{\pd \theta_k}}$ and $\vec{f}$ is vector of right hand side evaluations of the differential operator at the collocation points $(\vec{f})_i=f\left(\hat{u}(\tilde{\vec{x}}_i; \vec{\theta}(t)),t, \tilde{\vec{x}}_i\right)  $. 

We remark that there is no a priori guarantee that the linear system \eqref{eq:linsys} has a unique solution as $\vec{J}^\top \vec{J}$ may not be invertible. In such a case, we seek the solution of minimum norm, i.e. the smallest possible update to $\vec{\theta}$. We can thus express the time derivative of the parameters via the Moore-Penrose pseudo-inverse
\begin{equation}
\label{eq:projection}
 \vec{\dot{\theta}} \approx    \vec{\gamma}^* =  (\vec{J}^\top \vec{J})^\dagger \vec{J}^T \vec{f}.
\end{equation}
Equation \eqref{eq:projection} can also be interpreted geometrically - the evolution update equation is  a projection of the PDE right-hand side operator onto the tangent space of the neural network. A succesful EDNN model thus needs to
\begin{itemize}
    \item accurately approximate the true solution trajectory via the output of the neural network,
    \item have a rich enough tangent space along this solution trajectory, such that the PDE dynamics can be accurately represented in the update equation.
\end{itemize}

\subsection{Efficient solution of the parameter update equation}
In practice, direct inversion of the linear system  \eqref{eq:linsys} with the generally dense matrix $\vec{J}^\top \vec{J}$ is costly, especially if  the  Moore-Penrose pseudo-inverse is computed. We  avoid direct assembly of the neural network Jacobian $\vec{J}$, as it requires $n_\theta$ derivative evaluations of size $n_x$. Instead, we propose the use of a Krylov solver, which only needs Jacobian-vector products: \\
To compute the action $\vec{w} =\vec{J}^\top \vec{J} \vec{\gamma}$, we  first compute $\vec{v} =\vec{J} \vec{\gamma}$ and then evaluate $\vec{w}= \vec{J}^\top \vec{v}$ using automatic differentiation. 

Noting that \eqref{eq:linsys} is symmetric, MINRES \cite{paige_solution_1975} would be an appropriate solver choice. It is algebraically equivalent and more numerically stable, especially for ill-conditioned $\vec{J}$, to directly apply LSMR \cite{fong_lsmr_2011} to the following least squares problem:
\begin{equation}
\label{eq:lsqrproblem}
    \text{argmin}_{\vec{\gamma}} ||  \vec{J} \vec{\gamma} -\vec{f} ||_2^2,
\end{equation}
which can be seen as a discrete version of the minimization problem in \eqref{eq:minprob}.

When the least-squares problem in \eqref{eq:lsqrproblem} has multiple solutions, LSMR actually solves the problem
\begin{equation}
\begin{split}
& \vec{\dot{\theta}} =\text{argmin}_{\vec{\gamma}} || \vec{\gamma}||_2^2,     \\
   &\text{subject to } \vec{J} \vec{\gamma} =\vec{f}.
   \end{split}
\end{equation}
Therefore, LSMR is  guaranteed to return the solution of minimum norm. Additionally the norm $||\vec{\gamma}||$ of the current solution iterate $\vec{\gamma}$ has numerically been shown to increase monotonically. 
LSMR is also preferable over the more common LSQR solver as it directly minimizes the residual of \eqref{eq:linsys}, $||  \vec{J}^\top\vec{J} \vec{\gamma} - \vec{J}^\top \vec{f} ||_2^2$, in every step, allowing for faster termination of the procedure. We can then treat underdetermined problems, where the number of NN parameters $n_\theta$ exceeds the number of spatial collocation points $n_x$, identically to the overdetermined case, hence decoupling the spatial resolution from the modeling capabilities of the neural network. 

\textbf{Remark:} Numerically, we terminate the LSMR solve at a finite precision tolerance, such that \eqref{eq:linsys} is not solved exactly even when a compatible solution exists. This can be considered as a form of regularization as LSMR will produce parameter updates with smaller norm for higher error tolerances.

\subsection{Overview time stepping procedure}
In this section, we summarize the time stepping procedure for the neural network parameters.
As a first step, we need to ensure that at time $t=0$, the neural network output $\hat{u}(\vec{\tilde{x}};\vec{\theta}(0))$ reproduces the initial condition $u_0(\vec{\tilde{x}})$. In \cite{du_evolutional_2021, bruna_neural_2022}, this is achieved by training on the residual of the initial condition at a finite number of collocation points $\vec{\tilde{X}}= [\vec{\tilde{x}_1},\vec{\tilde{x}_2}, ..., \vec{\tilde{x}_{n_\text{train}}}]$:
\begin{equation}
\label{eq:inittrain}
    \vec{\theta}_0 = \text{argmin}_{\vec{\gamma}} || \hat{u}(\vec{\tilde{X}};\vec{\gamma}) - u_0(\vec{\tilde{X}}) ||_2^2,
\end{equation}
which may require a large number of epochs until a sufficient accuracy is reached. We propose an additional "trainingfree" approach, where the neural network output is modified to satisfy the initial condition by construction
\begin{equation}
\label{eq:trainfree}
    \tilde{u}(\vec{\tilde{x}};\vec{\theta}(t)) = u_0(\vec{\tilde{x}}) + \hat{u}(\vec{\tilde{x}};\vec{\theta}(t))- \hat{u}(\vec{\tilde{x}};\vec{\theta}(0)),
\end{equation}
which allows skipping the training step, when $u_0(\vec{\tilde{x}})$ is available analytically. This comes at the cost of doubling the number of neural network evaluations,which consequently also increases the cost of evaluating the differential operator and the Jacobian.\\
For the time integration itself, it is possible to rely on established ODE solvers. We consider Runge-Kutta methods which allow flexible adaptation of the time step, specifically the forward Euler scheme and the 5(4) stage method by Tsitouras (Tsit5) \cite{tsitouras_rungekutta_2011}.
The procedure is summarized in  Algorithm \ref{alg:pseudoalgtime}.

\begin{algorithm} [h]
\caption{Time integration of PDEs with neural networks}
\label{alg:pseudoalgtime}
\begin{algorithmic}[1]
\Require Neural network $\hat{u}$, PDE rhs $f$, initial condition $u_0$, collocation points $\vec{\tilde{X}_c}$, RK method with stage weights $a,b,c$, stopping tolerance $\texttt{tol}$ for $\texttt{LSMR}$.
\Ensure PDE solution $\tilde{u}$ at test locations $\vec{\tilde{X}}_\text{test}$ at the final computation time $T$.
\State Randomly initialize $\vec{\theta}_0$.
\State $\vec{\theta}_n \gets \vec{\theta}_0$
\If{trainingfree}
    \State Assign NN: $\tilde{u} \gets u_0(\vec{\tilde{x}}) + \hat{u}(\vec{\tilde{x}};\vec{\theta}_n)- \hat{u}(\vec{\tilde{x}};\vec{\theta}_0)$
\Else
    \State Train $\vec{\theta}_n$ to fit the intial condition according to \eqref{eq:inittrain}.
    \State  Assign NN: $\tilde{u} \gets \hat{u}(\vec{\tilde{x}};\vec{\theta}_n)$
\EndIf
\State $t\gets 0 $
\While{$t <T$}
\For{each RK stage $i$}
\State $t_\text{loc} \gets t+ c_i h$
\State $\vec{\theta}_\text{loc} \gets  \vec{\theta}_n + \sum_{j=1}^{i-1} a_{ij} \vec{\gamma}_j$ 
\State  $\vec{f}\gets f\left(\tilde{u}(\vec{\tilde{X}_c};\vec{\theta}_\text{loc}),t_\text{loc}, \vec{\tilde{X}_c}\right)$
\State $\vec{J} \gets \dfrac{\tilde{u}(\vec{\tilde{X}_c}; \vec{\theta}_\text{loc})}{{\pd \vec{\theta}_\text{loc}}}$  \Comment{Initialize Jacobian action.}
\State $\vec{\gamma}_i \gets \texttt{LSMR}(\vec{J}, \vec{f}, \texttt{tol})$
\EndFor
\State $\Delta\vec{\theta} \gets  h \sum_{i=1}^s b_i \vec{\gamma}_i$ \Comment{Error control and time step adaptation may happen here.}
\State $t \gets t+ h$
\State $\vec{\theta}_n \gets \vec{\theta}_n +  \Delta\vec{\theta}$

\EndWhile
\State Evaluate solution $\tilde{u}(\vec{\tilde{X}}_\text{test},\vec{\theta}_n)$.
\end{algorithmic}
\end{algorithm}

\subsection{Dealing with stiffness}
The parameter update equation \eqref{eq:projection} that we have derived is explicit in nature. This severely limits the time-step choice as some physical systems have natural stiffness, resulting in the time-stepping of the EDNN being expensive. 
In \cite{bruna_neural_2022}, an implicit Euler scheme is proposed. This comes at significant computational cost as  second order derivatives of the neural network with respect to its parameters $\vec{\theta}$ need to be computed repeatedly in an iterative training procedure.  Here, we propose a "linearly implicit" scheme based on Rosenbrock methods. 
Similarly to the previous application of the chain rule, we note that at first order we can write
\begin{equation}
\label{eq:firstorder}
\begin{split}
    u(\theta_{n+1} ) = u(\theta_n +\Delta \theta) \approx u(\theta_n) +  J_u \Delta \theta, \\
    f(\theta_{n+1} ) = f(\theta_n +\Delta \theta) \approx f(\theta_n) +  J_f \Delta \theta
    \end{split}
\end{equation}
where $J_u$ is the previously derived neural network Jacobian $\vec{J}$ and
\begin{equation}
    J_f = \dfrac{\partial f}{\partial \theta}= \dfrac{\partial f}{\partial u} J_u
\end{equation}
is the Jacobian of the right hand side. We consider the classic Rosenbrock triple with time step $h$ and parameters $\gamma = \dfrac{1}{\sqrt{2}+2}$, $e_{32}= 6 + \sqrt{2}$. To integrate from $u_0$ to $u_1$, we compute
\begin{equation}
    \begin{split}
        (I- h \gamma \dfrac{\partial f}{\partial u} ) k_1 &= f(u_0),  \\
        (I- h \gamma \dfrac{\partial f}{\partial u} ) k_2 &= f(u_0 +\frac{1}{2}\gamma k_1) - h\gamma \dfrac{\partial f}{\partial u} k_1, \\
        u_1 &= u_0 + h k_2,
    \end{split}
\end{equation}
and  evaluate the third stage for error control  with error $\varepsilon$ as
\begin{equation}
    \begin{split}
        (I- h \gamma \dfrac{\partial f}{\partial u} ) k_3 &= f(u_1) -e_{32} \left(k_2 - f(u_0 +\frac{1}{2}\gamma k_1)\right) -2 \left(k_1 -f(u_0)\right) , \\
        \varepsilon& = \frac{h}{6} (k_1-2k_2 +k_3).
    \end{split}
\end{equation}

We insert $k_i \approx J^{(i)}_u \kappa_i$,  to obtain 
\begin{equation}
\label{eq:rbparam}
    \begin{split}
        (J_u^{(1)}- h \gamma J_f ) \kappa_1 &= f(u (\theta_0)),  \\
        (J_u^{(2)}- h \gamma J_f ) \kappa_2 &= f(u( \theta_0 +\frac{1}{2}\gamma \kappa_1)) - h\gamma J_f \kappa_1, \\
        \theta_1 &= \theta_0 + h \kappa_2.
    \end{split}
\end{equation}
We solve the first two lines of \eqref{eq:rbparam} in a collocation sense with LSMR, i.e. we first compute 
\begin{equation}
    \text{argmin}_{\kappa_1} || (J_u^{(1)}- h\gamma J_f) \kappa_1 - f(u (\theta_0))||_2,
\end{equation}
and then 
\begin{equation}
\label{eq:2lsmr}
    \text{argmin}_{\kappa_2}||(J_u^{(2)}- h \gamma J_f ) \kappa_2 -f(u( \theta_0 +\frac{1}{2}\gamma \kappa_1) + h\gamma J_f \kappa_1 ||_2
\end{equation}
We then find the final update in \eqref{eq:rbparam}. For error control, we  proceed similarly for the third stage of the Rosenbrock triple and error metric:
\begin{equation}
    \begin{split}
        (J_u^{(3)}- h \gamma J_f ) \kappa_3 &= f(\theta_1) -e_{32} \left(J_u^{(2)}\kappa_2 - f(\theta_0 +\frac{1}{2}\gamma \kappa_1)\right) -2 \left(J_u^{(1)} \kappa_1 -f(u_0)\right) , \\
        \varepsilon& = \frac{h}{6} (\kappa_1-2\kappa_2 +\kappa_3),
    \end{split}
\end{equation}
     which requires one additional LSMR solve as in \eqref{eq:2lsmr}. All remaining computations follow  the usual ODE-based time-stepping described in Algorithm \ref{alg:pseudoalgtime}.

\textbf{Remark:} Rosenbrock methods usually come at severely increased computational cost as they need to form the Jacobian, which is the price to pay for using a larger time step. In this setting, we need to form the costly Jacobian action of the neural network anyway, so the cost increase of forming the Jacobian of the right hand side is limited. 

\subsection{Active sampling of collocation points}
We already indicated that the full resolution  of the FE mesh may not be required for the time-stepping procedure. Subsampling is especially beneficial when only part of the computational domain contains relevant information about the PDE dynamics, and we can concentrate collocation points in relevant areas.
In \cite{bruna_neural_2022}, the authors propose to dynamically  adjust the collocation points based on structural properties of the PDE. %

Inspired by this, we propose an ad-hoc sampling strategy, which can be  applied to any PDE:
 we  choose an indicator criterion, e.g. the magnitude of the right-hand-side
 $$ \omega_i = |f(\tilde{\vec{x}_i};\vec{\theta})|$$
  which informs us whether a collocation point $\tilde{\vec{x}}_i$ contributes to the PDE dynamics residual. 
 As evaluating the neural network is cheap compared to solving the update equation, we can do so on a large set of  $n$ candidate points.
 We then sample the actual collocation points according to a  discrete distribution  on the candidate set with normalized probabilities
 \begin{equation}
 \label{eq:probweight}
     p(\tilde{\vec{x}}_i) = \dfrac{\omega_i}{\sum_{i=1}^{n_s} \omega_i}
 \end{equation}
 $\omega_i$ for each point. The procedure for sampling from a spatial-parametric domain $\Omega_x \times \Omega_\alpha$ is summarized in  Algorithm \ref{alg:pseudoalgactive} and can lead to to a significant speedup for problems with localized features.
 \begin{algorithm} [h]
\caption{Active sampling strategy}
\label{alg:pseudoalgactive}
\begin{algorithmic}[1]
\Require Set of spatial points  $S_{\text{FE}} = \left\{ \vec{x}_k\right\}_{k=1}^{n_s}$ , parameter domain $\Omega_\alpha$, sampling criterion $\tilde{f}(\tilde{\vec{x}})$, size of candidate set $n$.
\Ensure At each time step: Set of collocation points $S_c$ of size $n_c$.
 \State Compute candidate point set $S = \left\{ \tilde{\vec{x}}_i\right\}_{i=1}^{n}$:
\For{ $i=1$ to $n$}
\State Uniformly sample  $ \vec{x}_i \in S_{\text{FE}}$  and independently $\vec{\alpha}_i \in \Omega_\alpha$ to form $\tilde{\vec{x}}_i = [\vec{x}_i,\vec{\alpha}_i]$.
\EndFor

\For{each time step $t_j$}
 \For{ $i=1$ to $n$}
\State Evaluate $\omega_i= \tilde{f} (\tilde{\vec{x}}_i,\vec{\theta}(t_j))$
\EndFor
\State Compute probability weights $p_i = \dfrac{\omega_i}{\sum_{i=1}^{n_s} \omega_i}$ to form discrete probability distribution $\mathcal{P}(\vec{\theta}(t_j)): P(X= \tilde{\vec{x}}_i ) = p_i $.
\State Sample $n_c$ points $X \sim \mathcal{P}(\vec{\theta}(t_j))$ to form $S_c$.
\EndFor

\end{algorithmic}
\end{algorithm}

\newpage
\section{Numerical examples}
We first highlight key components of the developed methods in two 1D cases. On the two-soliton solution of the Korteweg-de-Vries (kdV) equation, we demonstrate that a positional encoding with two basis functions can be sufficient to solve a transport dominated problem with high accuracy. A nonlinear heat equation with a parameterized initial condition then show-cases how time-stepping can be used for many-query applications. We then focus on advection-diffusion problems on a square domain with and without holes, to show the flexibility over a variety of geometries.
\subsection{Implementation and experiment settings}
All FE calculations are performed with FEnics \cite{logg_automated_2012}.
The time stepping scheme is implemented in Julia and builds on several packages: We use \texttt{Flux.jl} \cite{innes_flux_2018} for the NN presentation, \texttt{Zygote.jl} \cite{innes_differentiable_2019}, \texttt{TaylorSeries.jl} \cite{benet_taylorseries_2019} and \texttt{ForwardDiff.jl} \cite{revels_forward-mode_2016} for automatic forward and backward differentiation and \texttt{DifferentialEquations.jl} \cite{rackauckas_differentialequations_2017} for the explicit  Runge Kutte ODE solvers. We implement the Rosenbrock method and its adaptive time stepping scheme, but use the same PI controller settings as  $\texttt{DifferentialEquations.jl}.$ To compute the Jacobian of the right hand-side $J_u$, required for the Rosenbrock method, we use finite differences as the available AD packages do not succeed in this case. We use \texttt{TaylorSeries.jl} to ascertain that the approximation error of this step is below 1e-5. We further use the LSMR implementation of \texttt{IterativeSolvers.jl}, with termination tolerances set to "atol/btol=5e-5" unless specified otherwise. 

We use the same prototype of fully-connected neural network architectures for all problems, where we only vary the input embedding and size and number of hidden layers.
We omit the bias on the last layer as it cancels out when Dirichlet BCs are enforced and we did not observe degraded performance for the other cases.
We report the relative error as defined in \eqref{eq:L2error} at each time instance. For the KdV equation, an analytic solution is used. In the other examples we treat the FE solution as the ground truth. In the parametrized setting, we compute the average $L_2$ -error over a grid of the parametric domain with $N$ points as:
\begin{equation}
\label{eq:L2errorparam}
    \varepsilon = \frac{1}{N} \sum_{i=1} ^N \dfrac{|| \hat{u}_{\text{NN}} (\vec{\alpha}_i) -\hat{u}_{\text{FE}}(\vec{\alpha}_i) ||_2}{||  \hat{u}_{\text{FE}} (\vec{\alpha}_i) ||_2}.
\end{equation}
Another useful metric is the mean solution deviation $\delta$ as it measures how much the solution varies across the parameter space.
It is computed on the ground truth as
\begin{equation}
\label{eq:meandev}
    \begin{split}
        \mu &=  \frac{1}{N} \sum_{i=1} ^N \hat{u}_{\text{FE}}(\vec{\alpha}_i) \\
        \delta & =  \frac{1}{N} \sum_{i=1} ^N \dfrac{|| \hat{u}_{\text{FE}} (\vec{\alpha}_i) -\mu  ||_2} {||  \mu ||_2}
    \end{split}
\end{equation}

\subsection{Korteweg-de Vries equation}
As a first test case, we consider the Korteweg-de Vries equation:
\begin{equation}
\begin{split}
   \dfrac{\partial u}{\partial t} &= - \dfrac{\partial^3 u}{\partial x^3} - 6u \dfrac{\partial u}{\partial x}, \quad  x \in [-20,20] \\
    u(-20,t) &= u(20,t).
    \end{split}
\end{equation}
We choose an initial condition that leads to a travelling two-soliton solution as described in \cite{taha_analytical_1984}, see \ref{app:kdV} for the exact construction of the initial condition and reference solution. We use a NN with a simple periodic embedding 
\begin{equation*}
    \Phi =[ \cos(\dfrac{2\pi}{L}), \sin(\dfrac{2\pi}{L})],
\end{equation*}
where $L=40$ is the size of the domain, and vary the number of layers $n_l$ and number of hidden units $n_n$. The computational domain is discretized by randomly sampling $n_x=1000$ points.
We train each architecture to the initial condition on 5000 randomly sampled points for 100000 iterations. The same trained neural network is then used for time-stepping with a fixed time step Euler scheme, Tsit5 with automatic error control  and the proposed Rosenbrock 2-step method.

\begin{figure}

 \centering
 \subcaptionbox{t=0}
         {\includegraphics[width=0.45\textwidth]{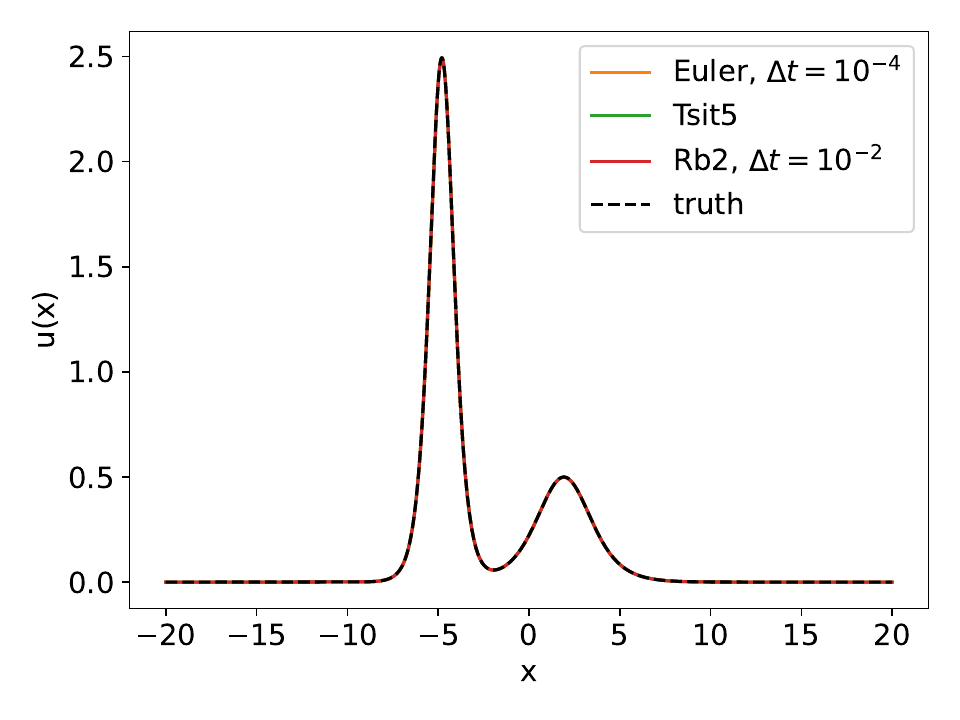}}
  \subcaptionbox{t=1}
         {\includegraphics[width=0.45\textwidth]{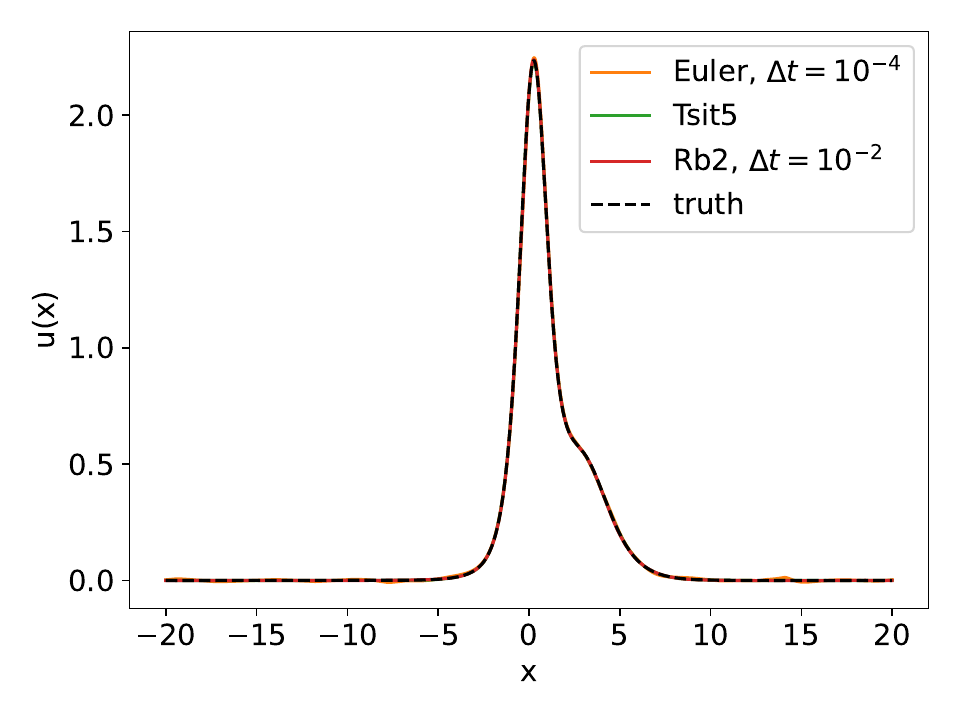}}
  \subcaptionbox{t=2}
       {  \includegraphics[width=0.45\textwidth]{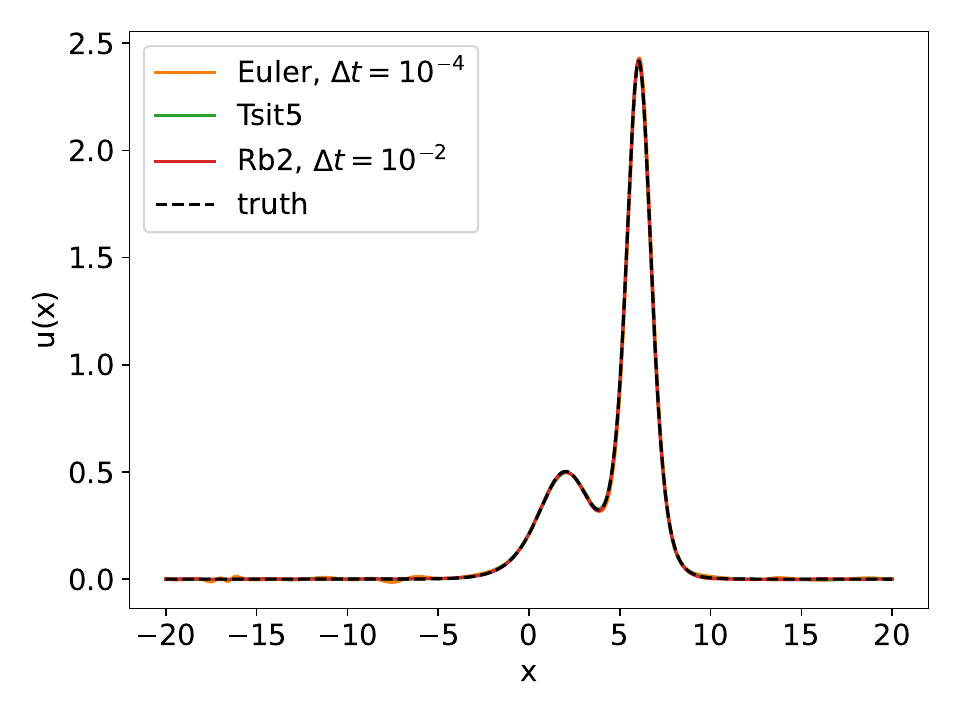}}
    \subcaptionbox{t=3}
       {  \includegraphics[width=0.45\textwidth]{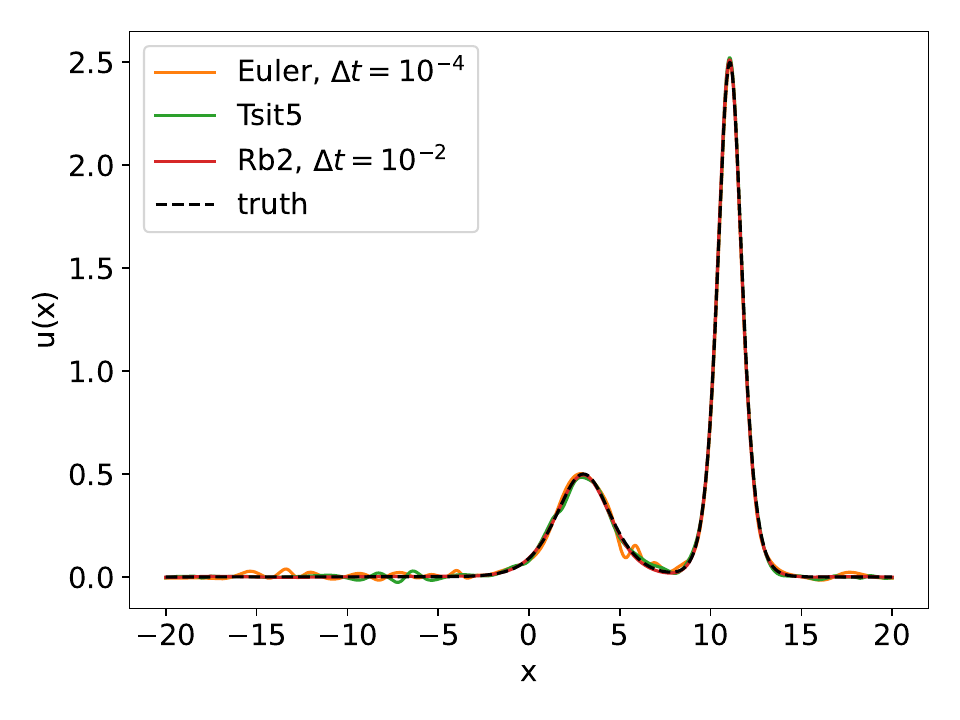}}
\caption{Time evolution of the Korteweg-de Vries equation.}
\label{fig:kdVsol}
\end{figure}

\begin{figure}

 \centering
 \subcaptionbox{Performance of the Rosenbrock (Rb2) method with $\Delta t=0.01$ for different architectures.}
         {\includegraphics[width=0.45\textwidth]{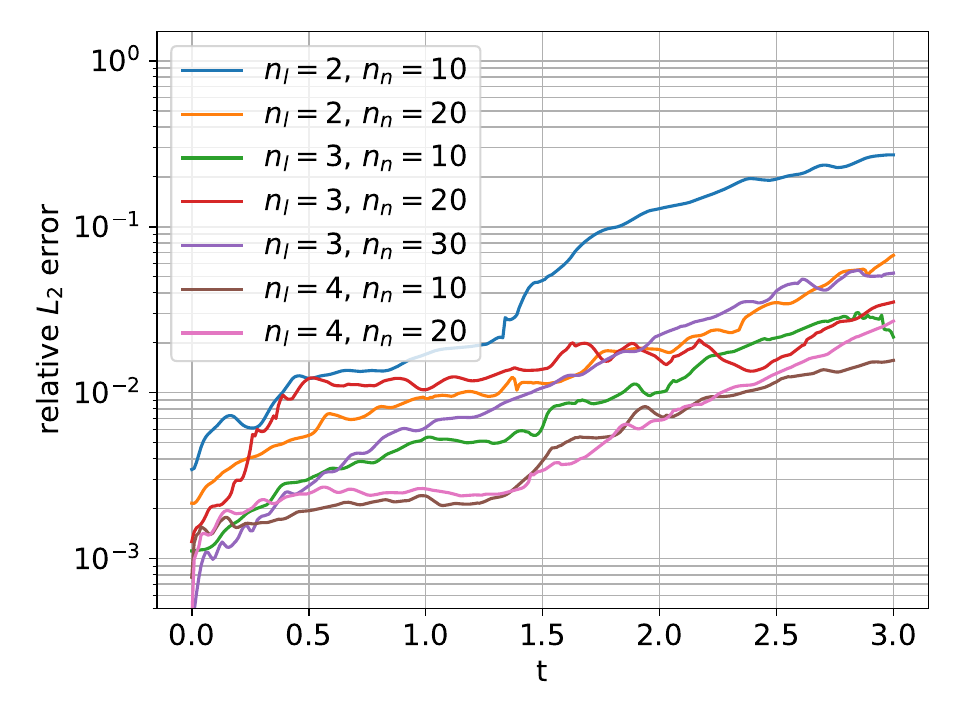}}
  \subcaptionbox{Performance of time steppers for a NN with $n_l=4$, $n_n=10$.}
         {\includegraphics[width=0.45\textwidth]{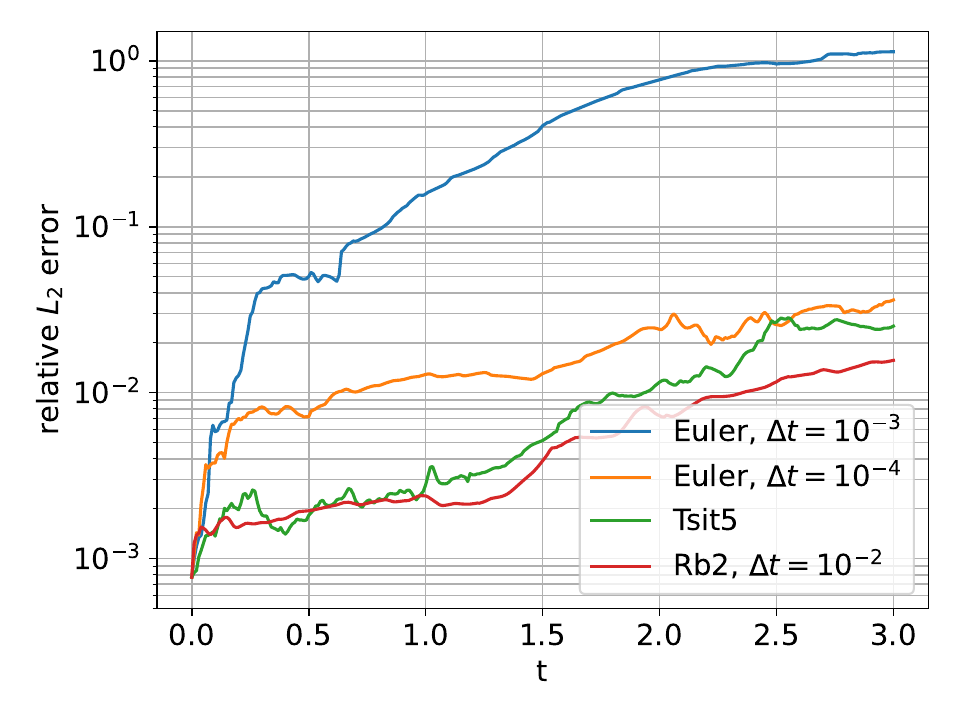}}
\caption{Evolution of error for the Korteweg-de Vries equation.}
\label{fig:kdVerror}
\end{figure}

The results in Figure \ref{fig:kdVsol} show that we visually approximate the solution well. At the end of the simulation time, we observe small oscillation errors for the explicit time stepping methods. This is also reflected by the error plots in Figure \ref{fig:kdVerror}. The explicit Euler scheme with a time step of t=0.001 fails to capture the solution trajectory accurately, while all other time stepping schemes produce results with 1-3 \% error. Even more strikingly, the Rosenbrock method produces the lowest error with the largest time step. Increasing the depth of the neural network as compared to the width seems to produce better results, but we observe no clear trends, as long as the neural network is large enough.
In conclusion, we observe that the neural network indeed can reproduce transport behavior when using static positional embeddings and the modified Rosenbrock method alleviates the stiffness problem of the explicit solvers on this example.

\subsection{A parametrized nonlinear heat equation}
In this example, we consider a heat equation with a cubic nonlinearity and a parametrized initial condition with Dirichlet boundary conditions.
\begin{equation}
    \begin{split}
        & \dfrac{\partial u}{\partial t}   = \dfrac{\partial^2 u}{\partial x^2}  -16 u^3, \\
        & u(0,t) = u(1,t) = 1, \\
       & u(x,0;\alpha_1, \alpha_2) = 1 + \alpha_1 \sin(\pi x) + \alpha_2 \sin(3 \pi x), \\
      & \alpha_1, \alpha_2 \in [-0.5, 0.5].
    \end{split}
\end{equation}
We consider two different positional embeddings 
\begin{align*}
    \Phi_2 &=[ \sin(\pi x), \sin( 2\pi x)] \\
    \Phi_4 &= [\sin(\pi x), \sin( 2\pi x), \sin( 3\pi x), \sin(4\pi x)],
\end{align*}
corresponding to the first eigenfunctions of the Laplace problem with Dirichlet boundary conditions.
Choosing a constant lifting function $u_\text{l}(x,t) =1$, we model the PDE solution as
\begin{equation}
    u(x,t) = u_\text{NN} (\Phi_j(x),\vec{\alpha};\theta(t))  - u_\text{NN} (\vec{0},\vec{\alpha};\theta(t)) + u_\text{l}(x,t)
\end{equation}
We point out that the initial condition does not lie in the span of $\Phi_2$ and that the richness of solutions across the parameter space decreases over time as all initial conditions converge to the same steady state solution.
We use neural networks with 4 layers and 10 or 20 hidden units each and compare training to the initial condition (10000 points for 40000 iterations with ADAM) and the trainingfree scheme.
The reference solutions are computed with a FE method on an equidistant grid of 1000 elements with Lagrange P2 elements, using implicit Euler for the time stepping with a step size of $\Delta t= 10^{-4}$.

\begin{figure}
 \centering
 \subcaptionbox{t=0}
         {\includegraphics[width=0.35\textwidth]{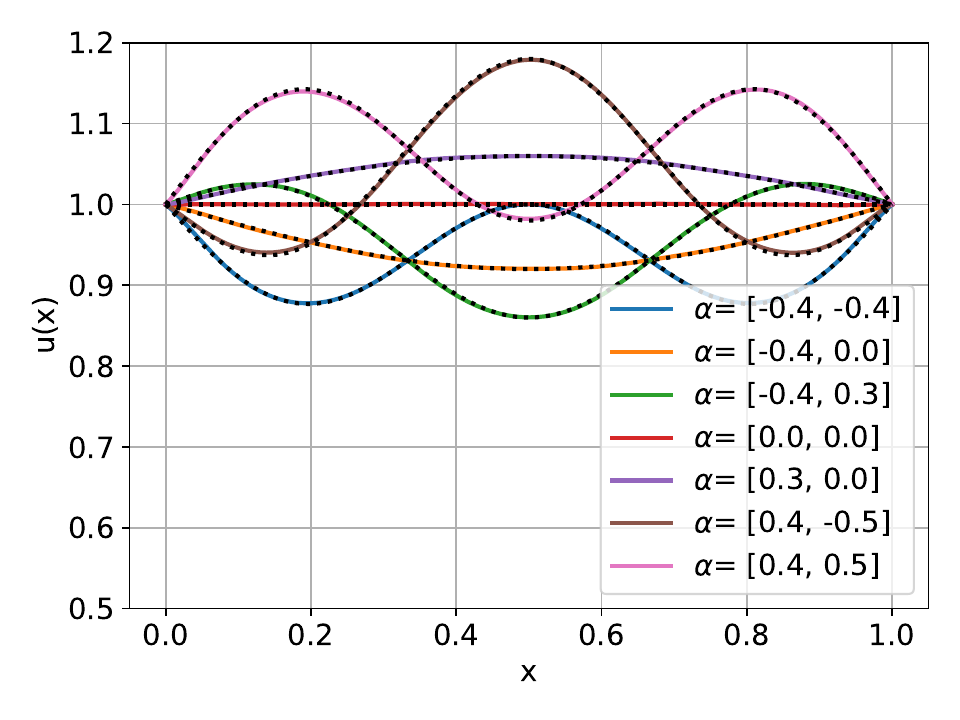}}
  \subcaptionbox{t= 0.002}
         {\includegraphics[width=0.35\textwidth]{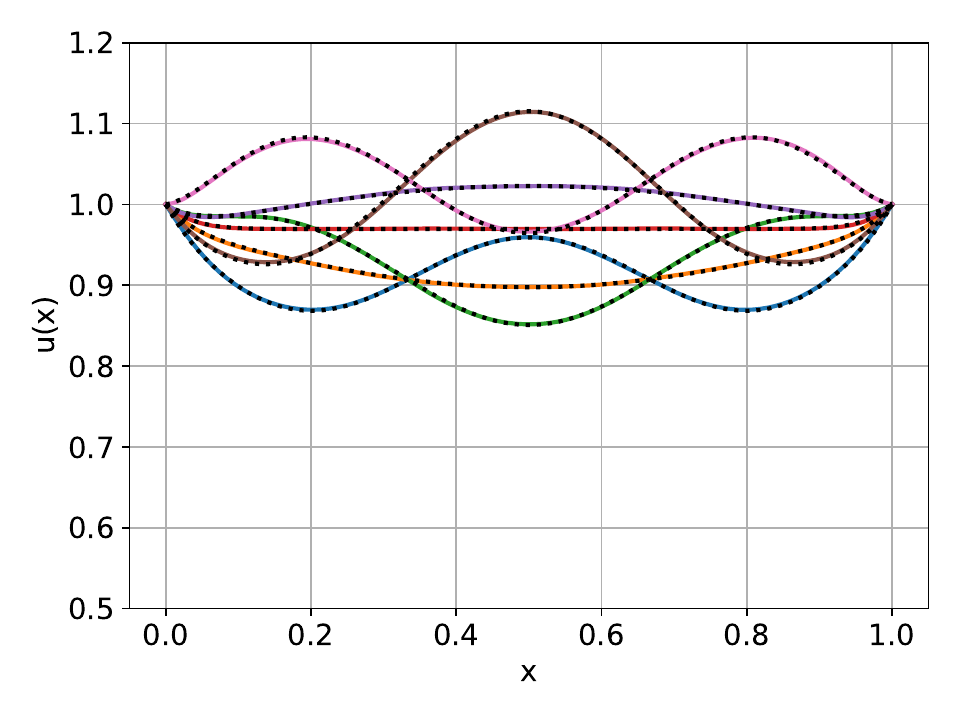}}
  \subcaptionbox{t=0.005}
       {  \includegraphics[width=0.35\textwidth]{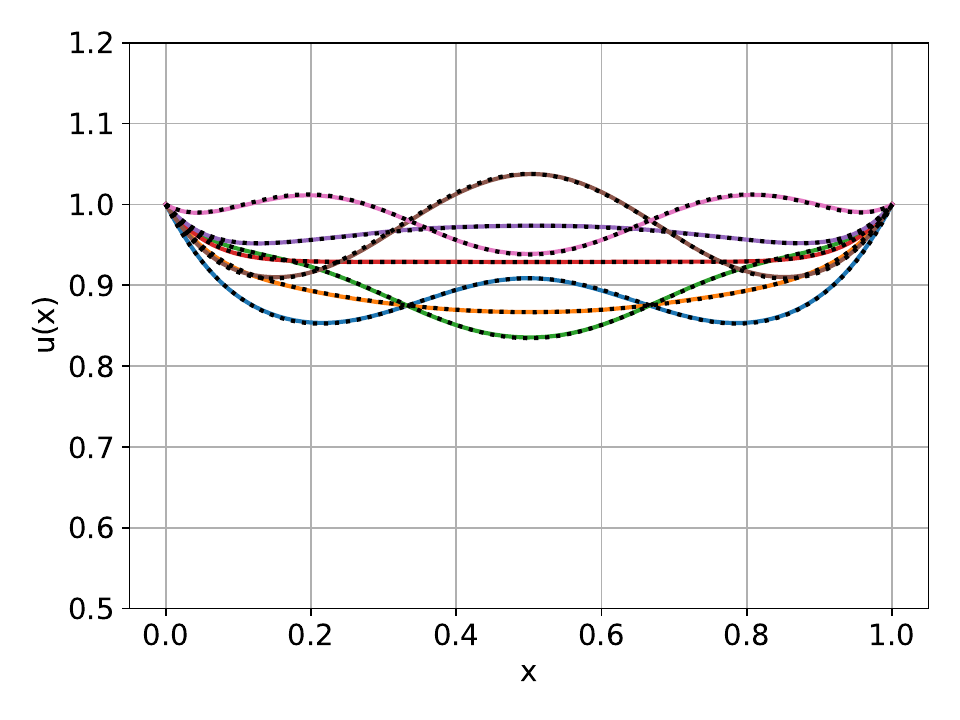}}
    \subcaptionbox{t=0.02}
       {  \includegraphics[width=0.35\textwidth]{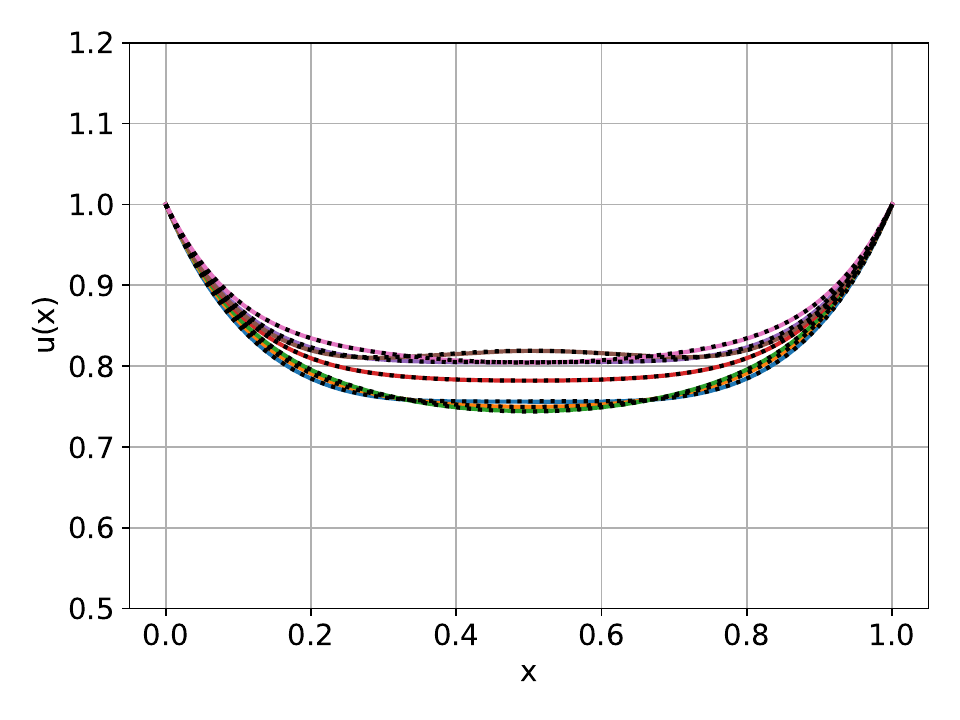}}
    \subcaptionbox{t=0.04}
       {  \includegraphics[width=0.35\textwidth]{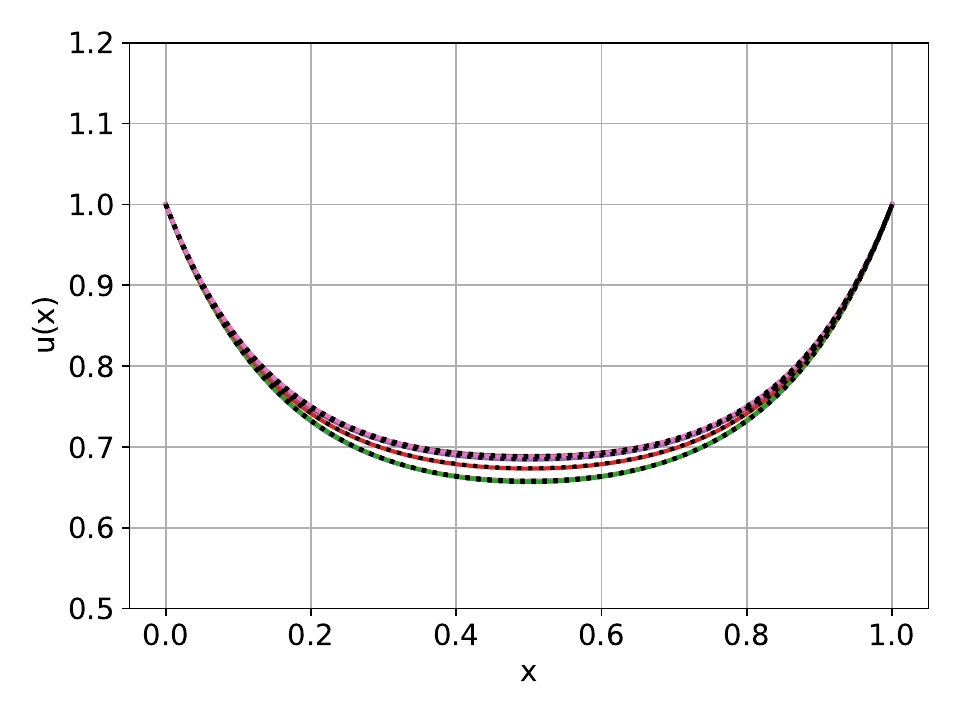}}
    \subcaptionbox{t= 0.1}
       {  \includegraphics[width=0.35\textwidth]{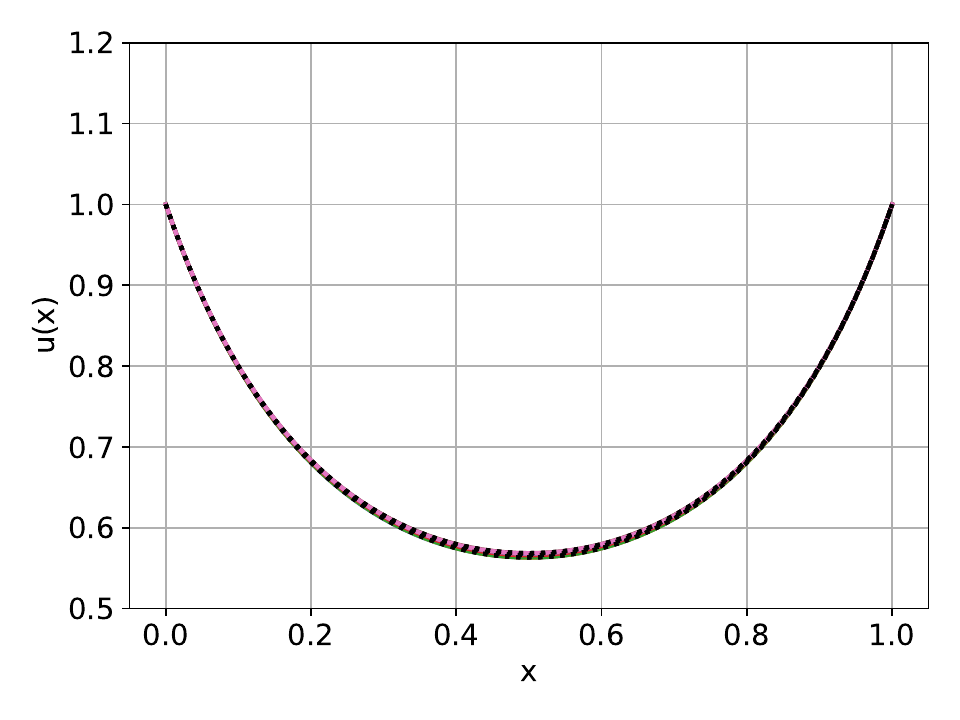}}
\caption{Time evolution of the parametrized heat equation for different parameter values. The black-dotted points indicate the reference FE solution. }
\label{fig:heat1Dsol}
\end{figure}

\begin{figure}

 \centering
 \subcaptionbox{NN: $\Phi_2$, $n_\text{hid}=10$.}
         {\includegraphics[width=0.45\textwidth]{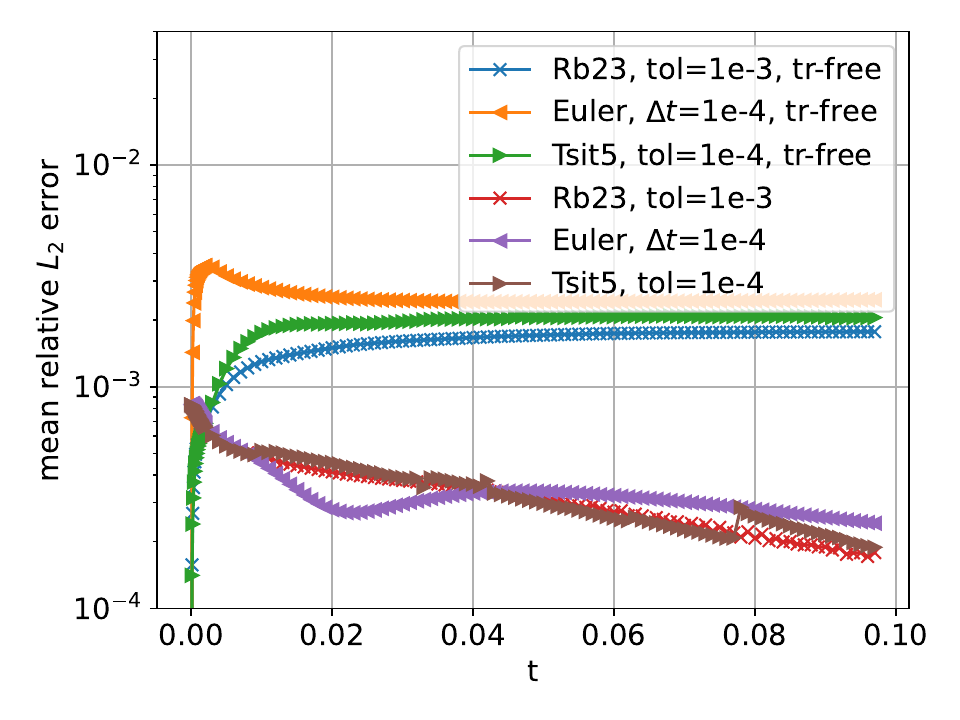}}
  \subcaptionbox{NN: $\Phi_4$, $n_\text{hid}=20$.}
         {\includegraphics[width=0.45\textwidth]{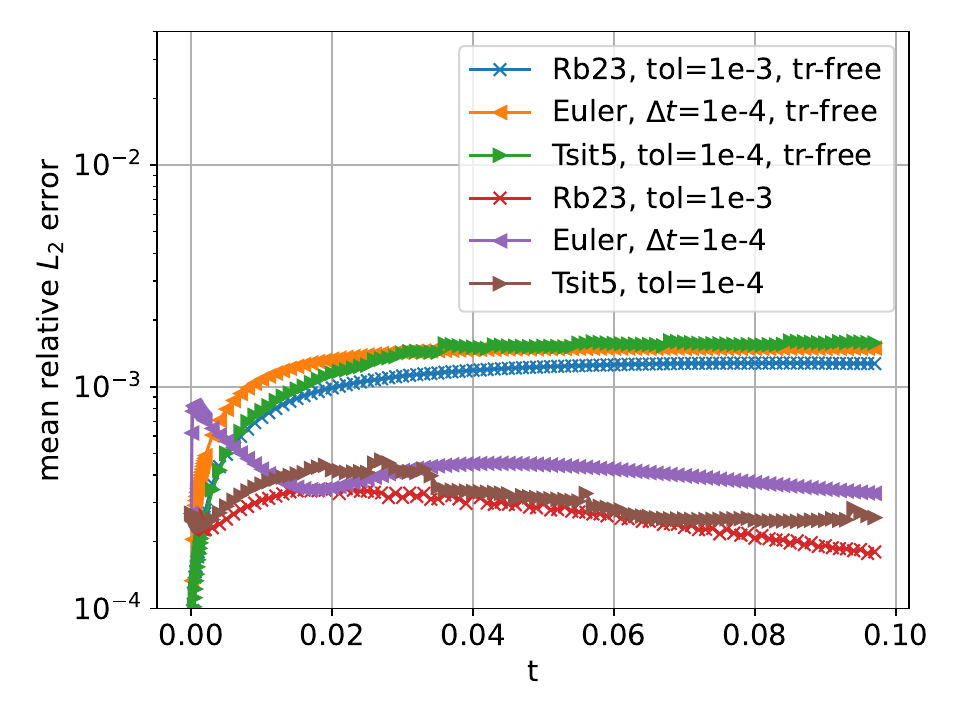}}
\caption{Evolution of error for the parametrized heat equation in 1D. "tr-free" refers to the trainingfree approach \eqref{eq:trainfree}.}
\label{fig:heat1Derror}
\end{figure}

The plots in Figure \ref{fig:heat1Dsol} show that we capture the solution well across a range of different parameter values and there is no visual difference to the FE reference solution.
We compute the mean relative error across a grid of parameters with 121 points, shown in Figure \ref{fig:heat1Derror}.  Comparing the two embeddings in (a) and (b) of Figure \ref{fig:heat1Derror}, we observe negligible differences which demonstrates that  we can leverage the nonlinearity of the neural network to fit functions that are not contained in the linear span of the embedding. We also confirm that for this simple heat equation example, the "trainingfree" approach produces acceptable errors, which are however higher compared to the neural networks which were trained on the initial condition.

We also note that the Rosenbrock method with adaptive time-stepping produces the lowest error. The tolerance of the time step adaptation wase chosen one order of magnitude  lower for Tsit5 to produce comparable results. We conjecture that this is due to the additional error sources in the neural network update equation: while error control is still possible qualitatively, not all quantitative rules carry over directly from traditional Runge-Kutta methods and the error estimate in the parameter update may not always be directly indicative of the error in the physical solution.

\subsection{ Parametrized nonlinear advection-diffusion problem}
We consider  a nonlinear advection-diffusion equation
\begin{equation}
    \begin{split}
     &\dfrac{\partial u}{\partial t} = 0.1 \nabla \cdot \left(\left(1+ \alpha_1 \sin(2\pi x_1)\right) \nabla u\right)  + 4 \begin{bmatrix} \cos(\pi\alpha_2) u \\ \sin(\pi \alpha_2) u \end{bmatrix}  \cdot  \nabla u \\
     &            \nabla u(\vec{x},t) \cdot \vec{n} = 0 \quad \text{for } x \in \partial \Omega \\
     &             u(\vec{x},0) = \sin^2(\pi x_1) \sin^2(\pi x_2) \\
     &             \alpha_1, \alpha_2  \in [-0.5, 0.5],
    \end{split}
\end{equation}
with the physical domain being the unit square.
As this problem uses Neumann boundary conditions, we consider the first three eigenfunctions of the Neumann case for the Laplace operator for the positional embedding:
\begin{equation*}
    \Phi =[ \cos(\pi x_1), \cos(\pi x_2), \cos(\pi x_1) \cos(\pi x_2)],
\end{equation*}
where we have dropped the constant eigenfunction $\phi_0=1$, as it does not enrich the function space - the bias term of the first layer should be able to model any necessary constant terms. 

In this example, all solutions start from the same initial condition and then evolve according to the parametric dependence, i.e. we expect the solution space to become richer/more diverse as time increases. Here $\alpha_1$ controls the variation in the diffusive term, while $\alpha_2$ controls the direction of the nonlinear advective term. 
We train for the initial condition on 10000 randomly sampled points in the domain for 40000 iterations and compare again with the trainingfree approach, which appears especially appealing here: as the initial condition does not depend on the parameters, we train the neural network to produce identical values across the whole parameter space. \\

\begin{figure}

 \centering
 \subcaptionbox{Trained to initial condition.}
         {\includegraphics[width=0.45\textwidth]{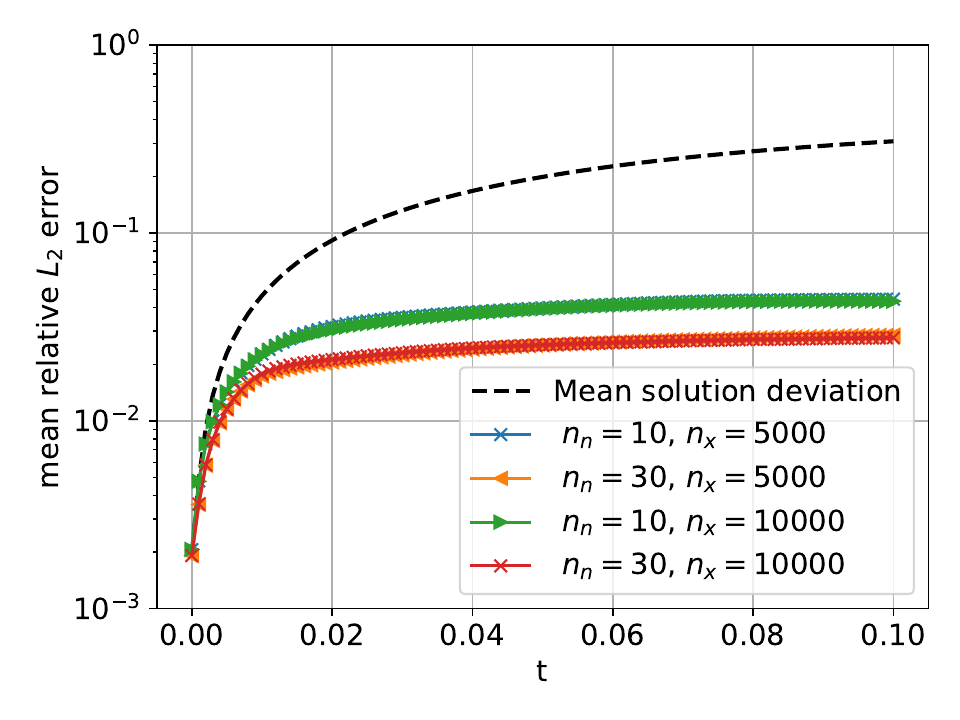}}
  \subcaptionbox{Training free.}
         {\includegraphics[width=0.45\textwidth]{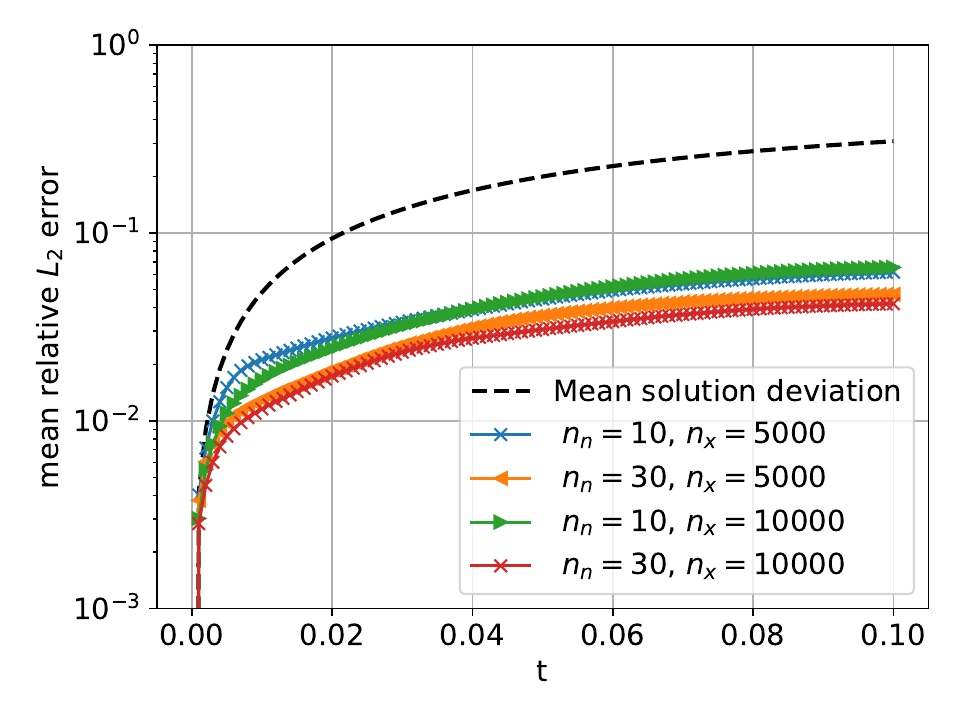}}
\caption{Evolution of error for the parametrized nonlinear advection diffusion equation.}
\label{fig:adv_diff_error}
\end{figure}

In Figure \ref{fig:adv_diff_error}, we compare the error of the two approaches for  neural networks with 4 layers with different number of hidden units $n_n$ and sampling points $n_x$. We plot the mean solution deviation introduced in \eqref{eq:meandev}
to demonstrate how the solution space becomes richer over time. We note that for both the trained and training-free approach the error increases  until t= 0.02. From that point on, the error for the trained neural network increases substantially less than for the trainingfree approach. We conclude that for the best possible performance, training to the initial condition is preferable, while the trainingfree approach may be an acceptable time saving measure for short time spans. We also note that there is only minimal decrease in error when using $n_x= 10000$ vs. $n_x=5000$ integration points, which indicates that we have sufficiently sampled the solution space. 

\begin{figure}
 \centering
 \subcaptionbox{$\vec{\alpha}= [-0.4, -0.4]$ }
         {\includegraphics[width=0.2\textwidth]{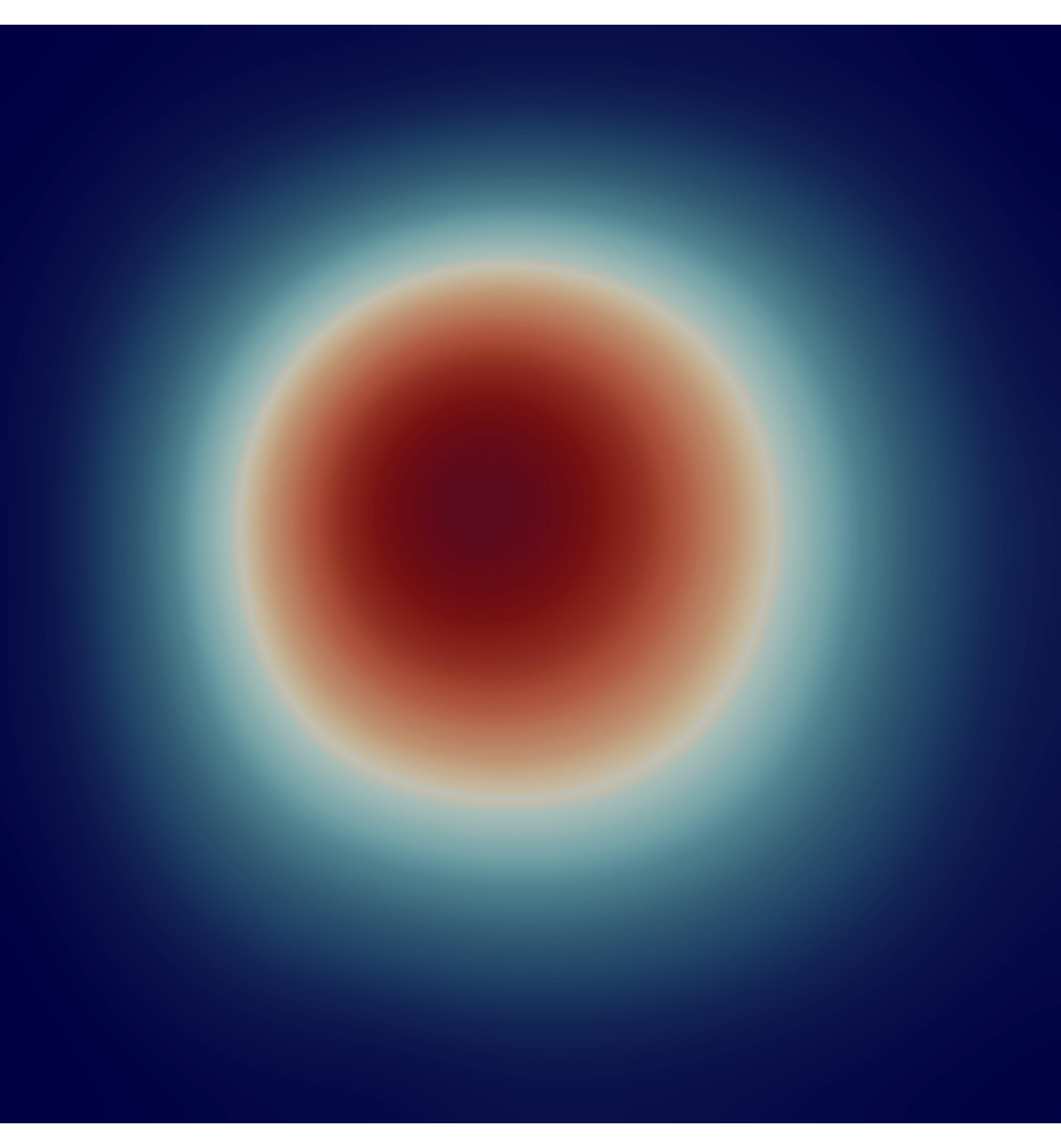}}
 \subcaptionbox{$\vec{\alpha}= [-0.4, 0]$ }
         {\includegraphics[width=0.2\textwidth]{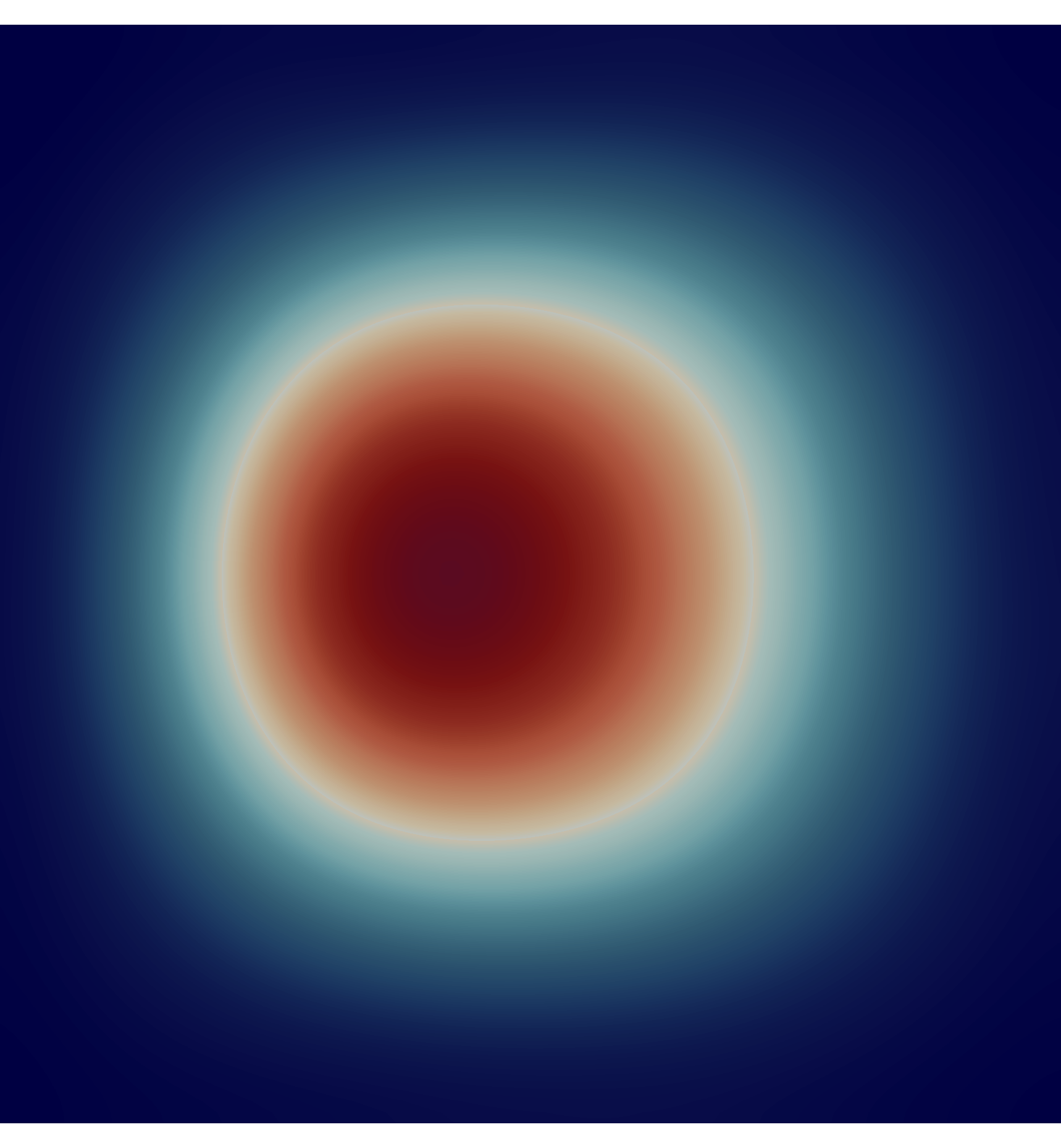}}
 \subcaptionbox{$\vec{\alpha}= [-0.4, 0.4]$ }
         {\includegraphics[width=0.2\textwidth]{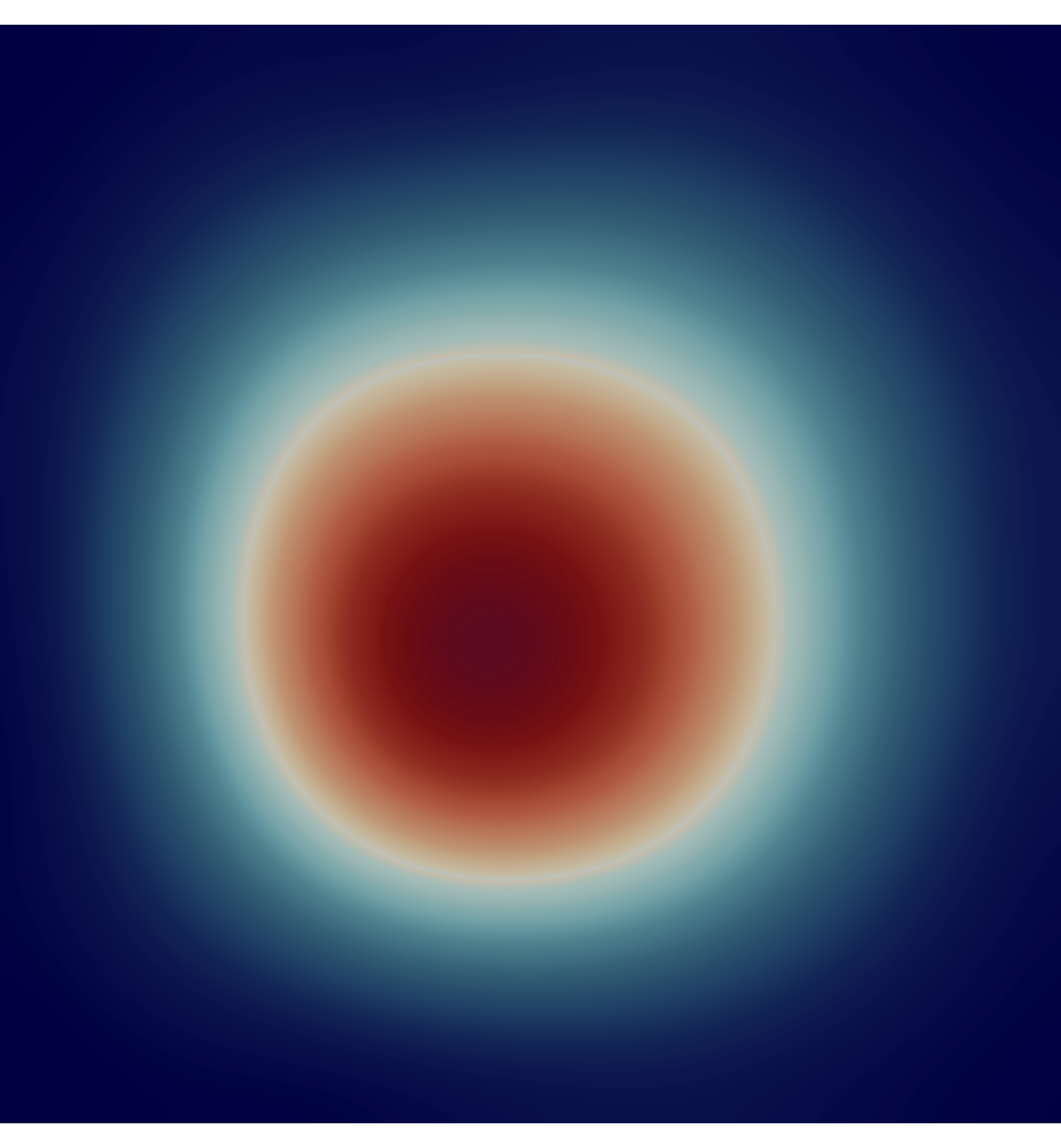}}\\
 \subcaptionbox{$\vec{\alpha}= [0, -0.4]$ }
         {\includegraphics[width=0.2\textwidth]{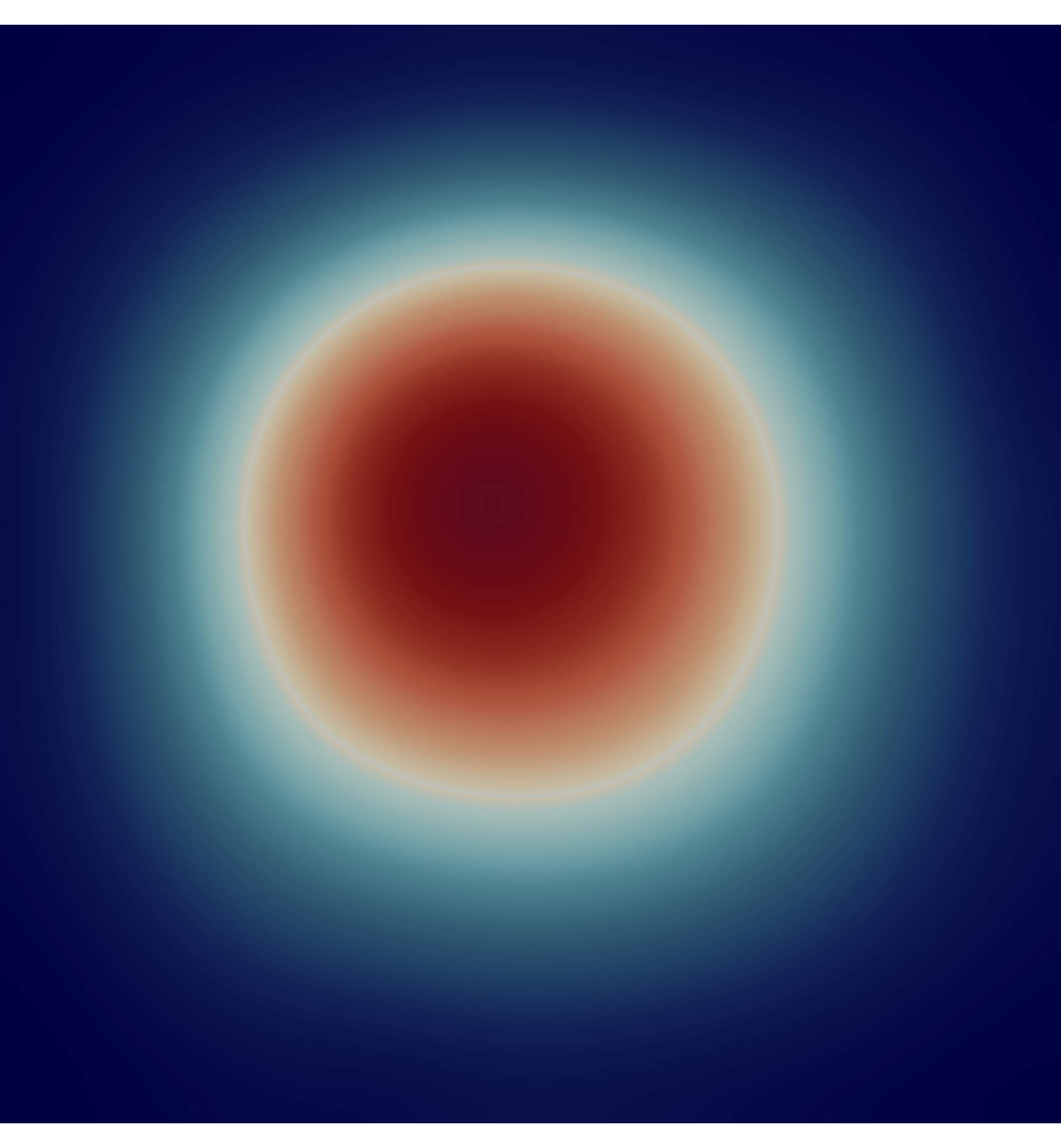}}
 \subcaptionbox{$\vec{\alpha}= [0, 0]$ }
         {\includegraphics[width=0.2\textwidth]{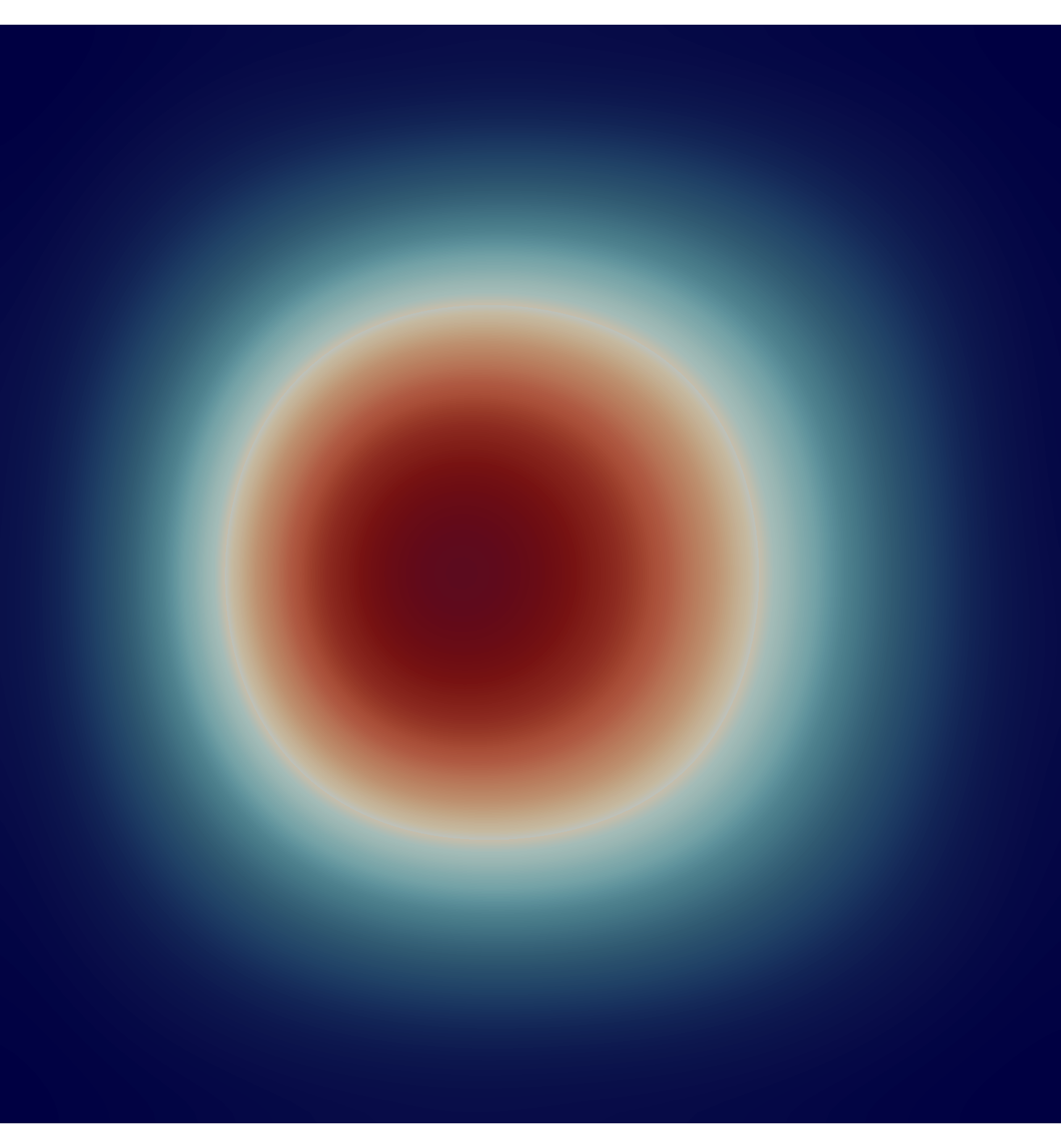}}
 \subcaptionbox{$\vec{\alpha}= [0, 0.4]$ }
         {\includegraphics[width=0.2\textwidth]{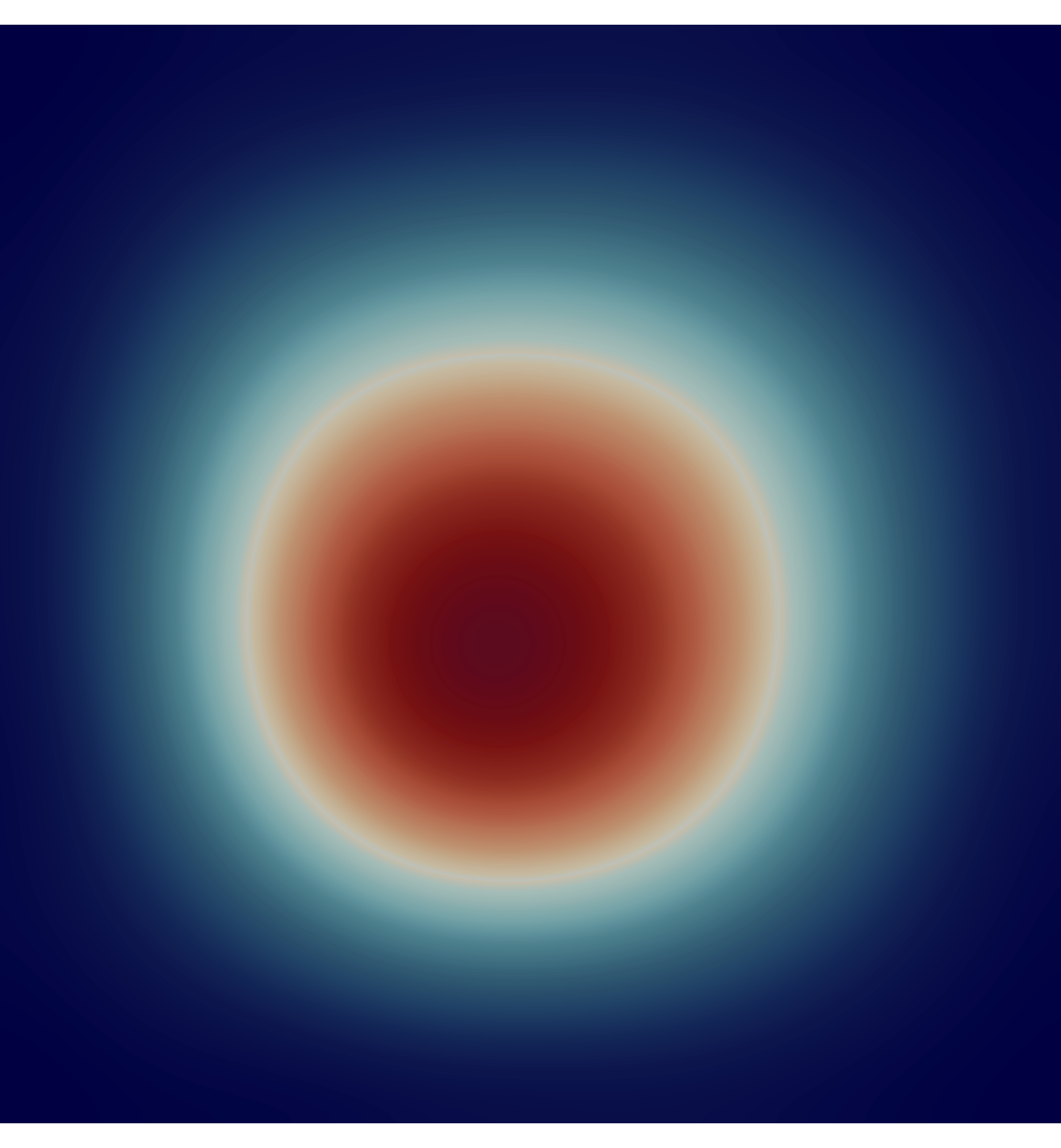}}\\
 \subcaptionbox{$\vec{\alpha}= [0.4, -0.4]$ }
         {\includegraphics[width=0.2\textwidth]{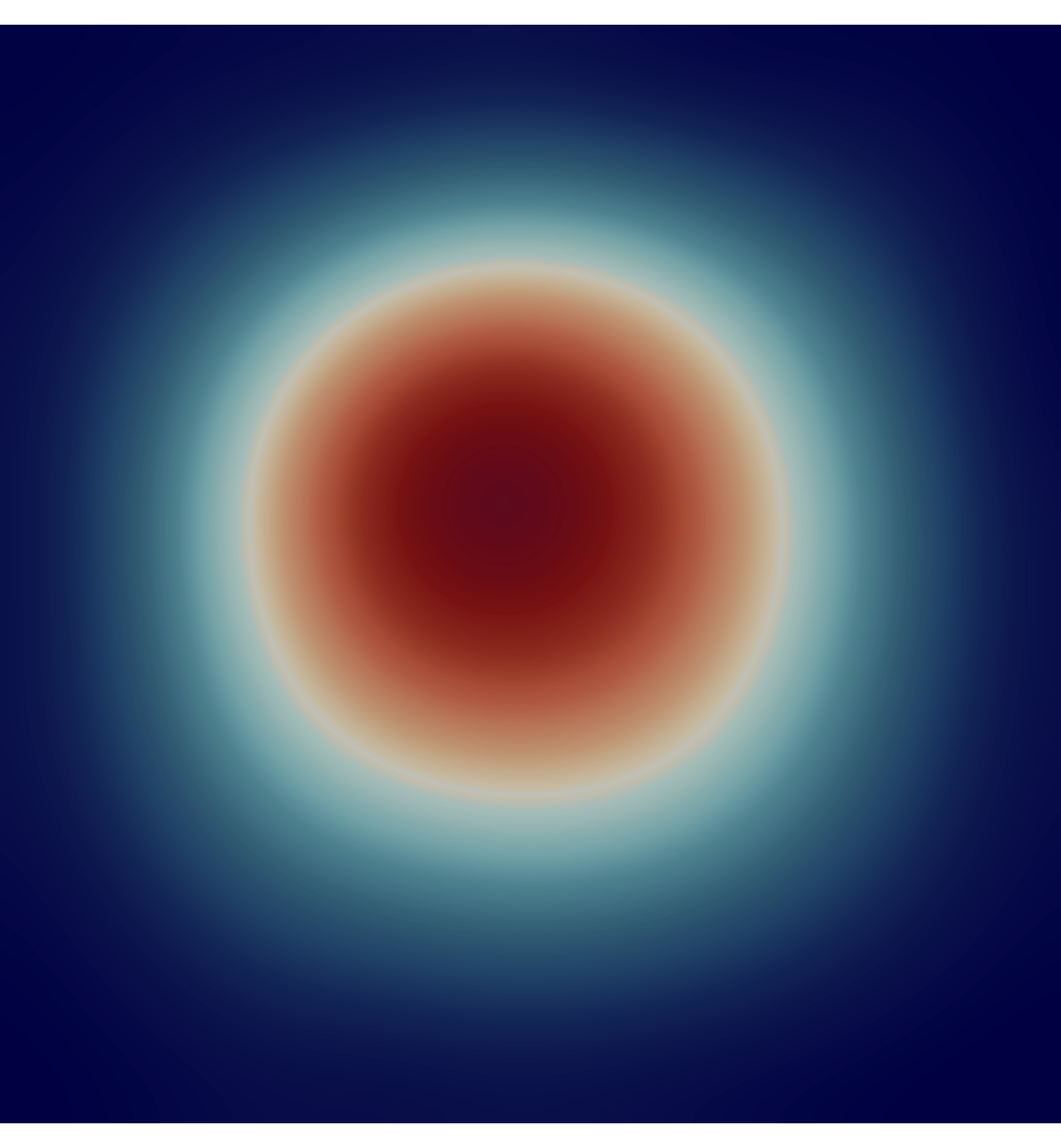}}
 \subcaptionbox{$\vec{\alpha}= [0.4, 0]$ }
         {\includegraphics[width=0.2\textwidth]{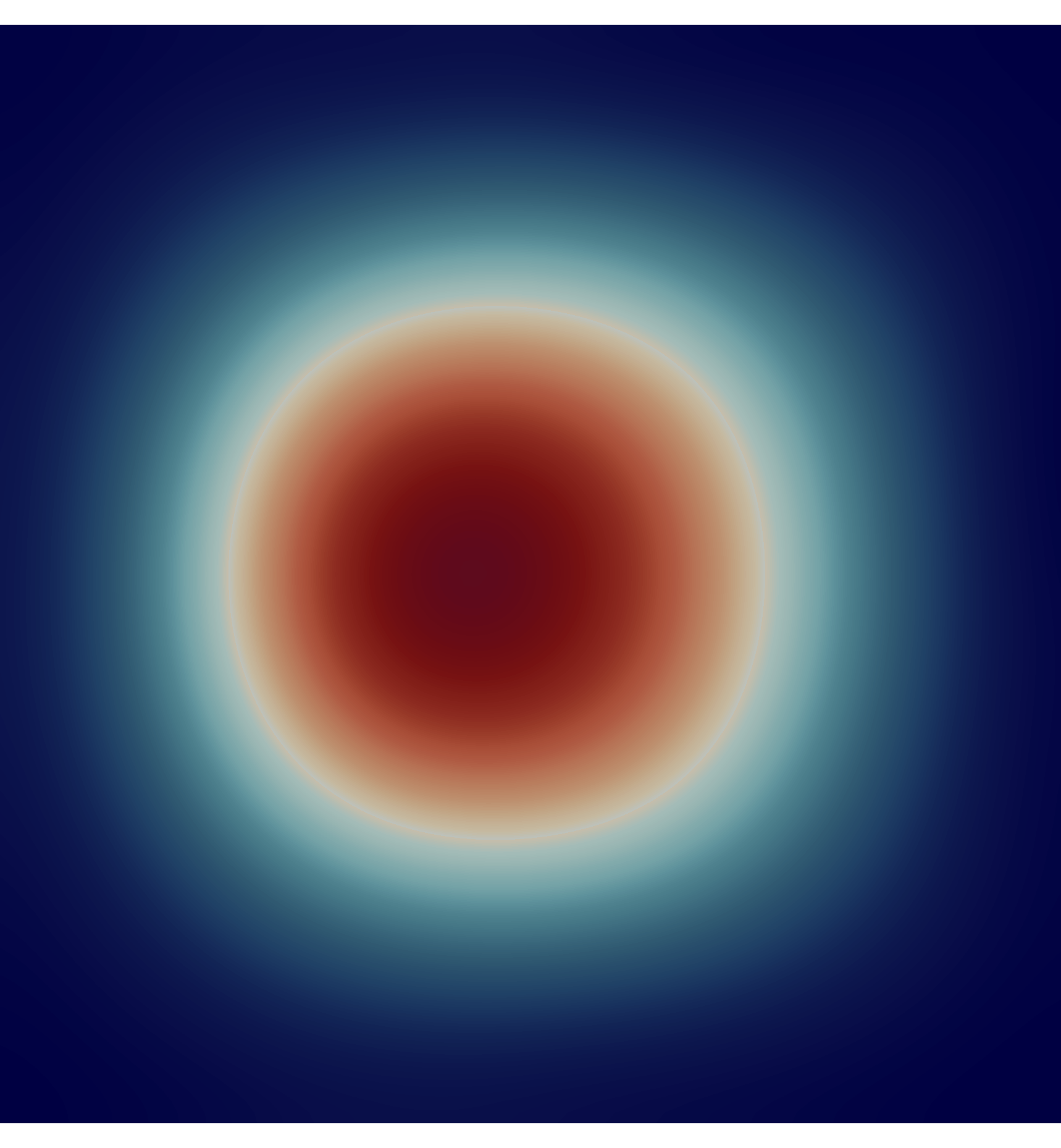}}
 \subcaptionbox{$\vec{\alpha}= [0.4, 0.4]$ }
         {\includegraphics[width=0.2\textwidth]{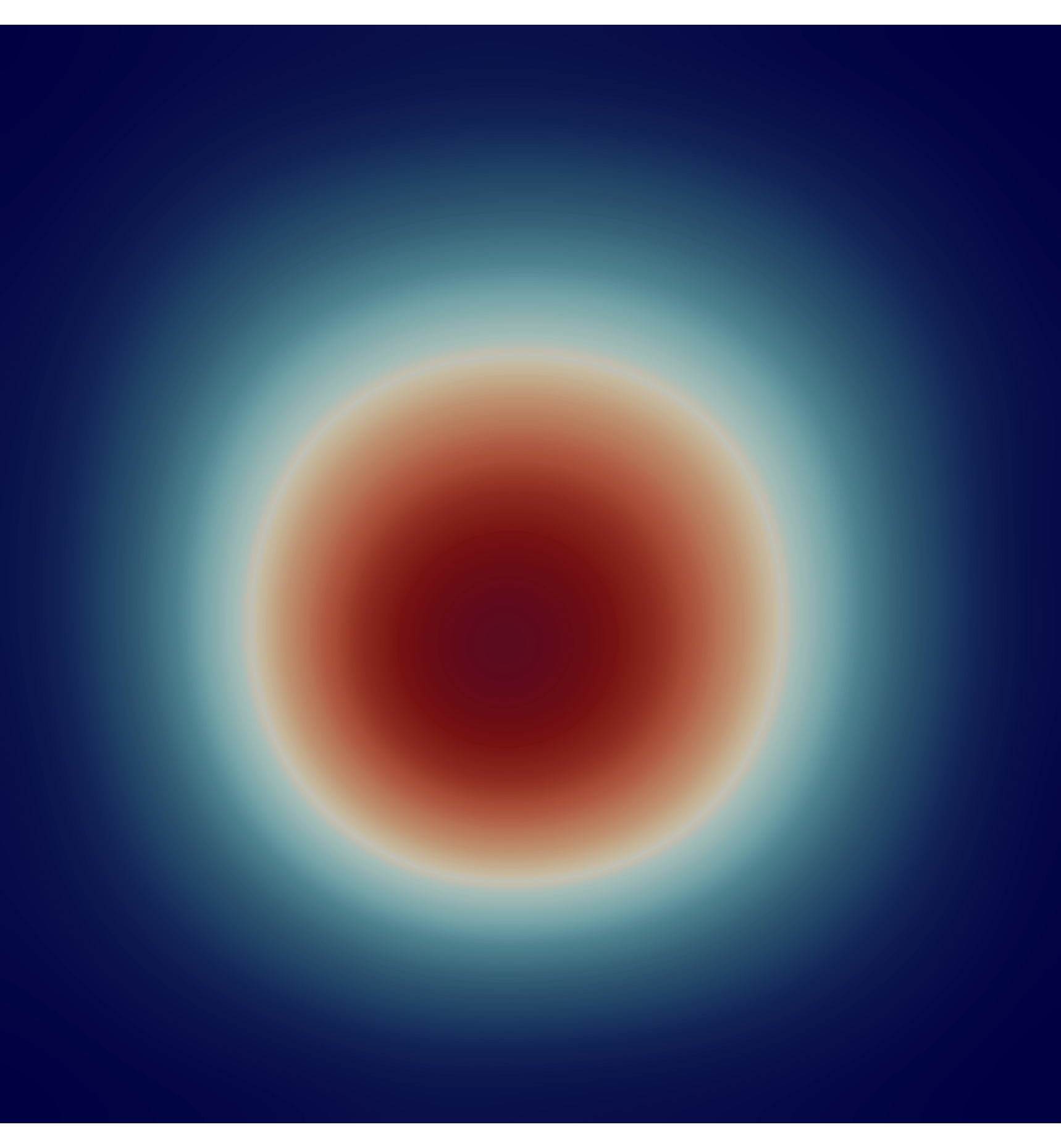}}
         \centering
        \includegraphics[width=0.6\textwidth]{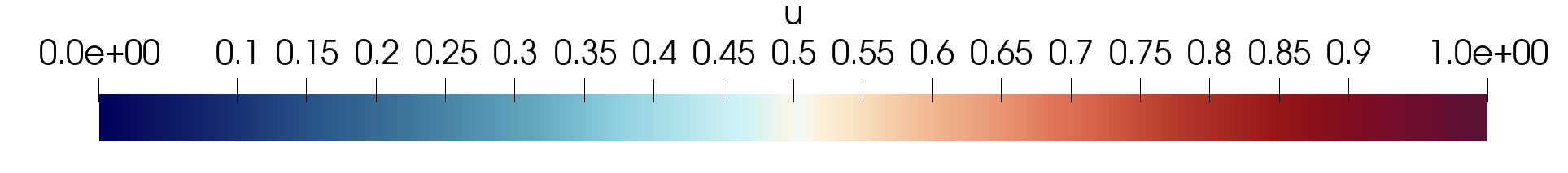}

\caption{Solution snapshots of the nonlinear advection diffusion problem at  time t= 0.02  for different parameter values.}
\label{fig:adv_diff_sol1}
\end{figure}

\begin{figure}
 \centering
 \subcaptionbox{$\vec{\alpha}= [-0.4, -0.4]$ }
         {\includegraphics[width=0.2\textwidth]{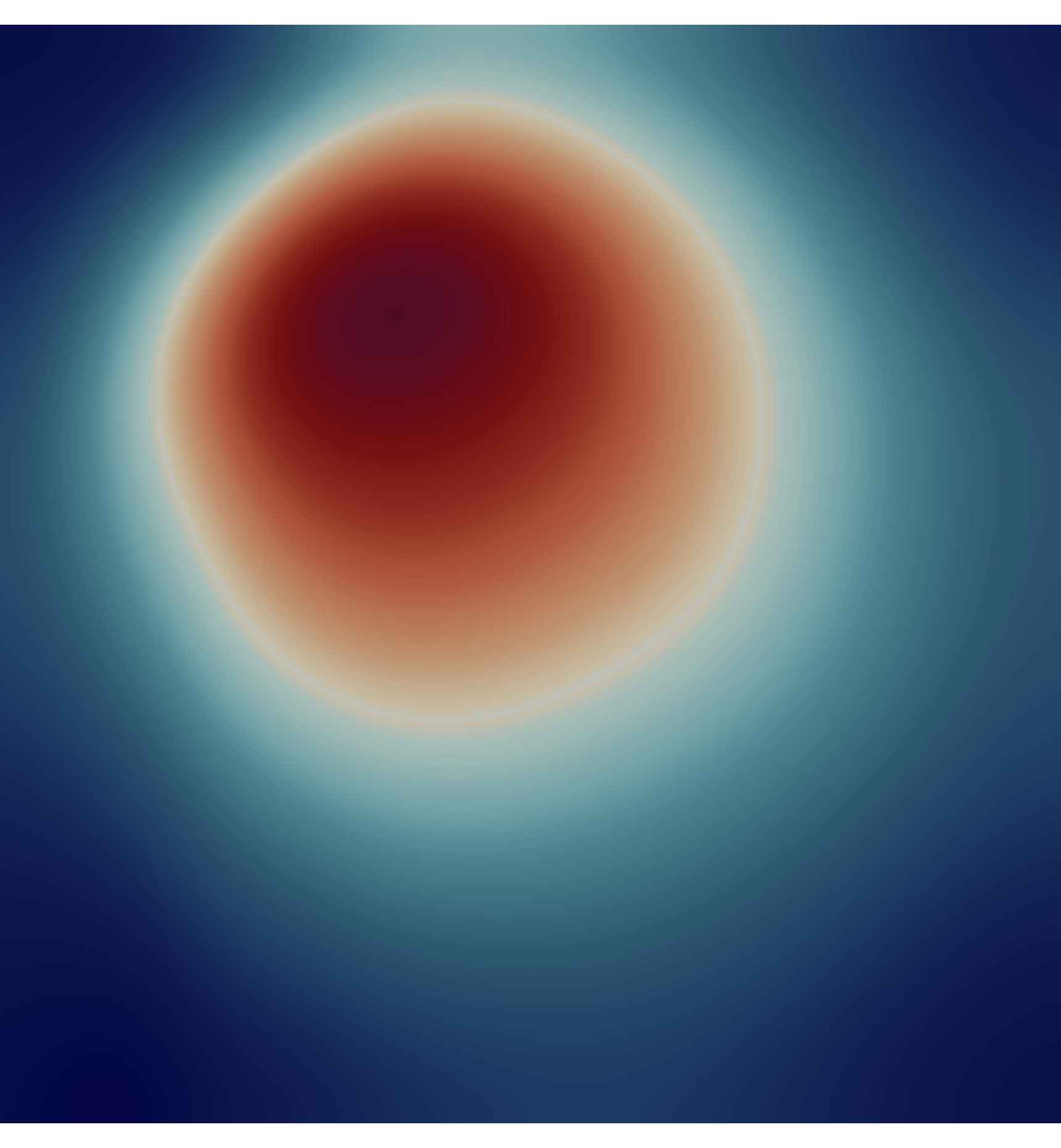}}
 \subcaptionbox{$\vec{\alpha}= [-0.4, 0]$ }
         {\includegraphics[width=0.2\textwidth]{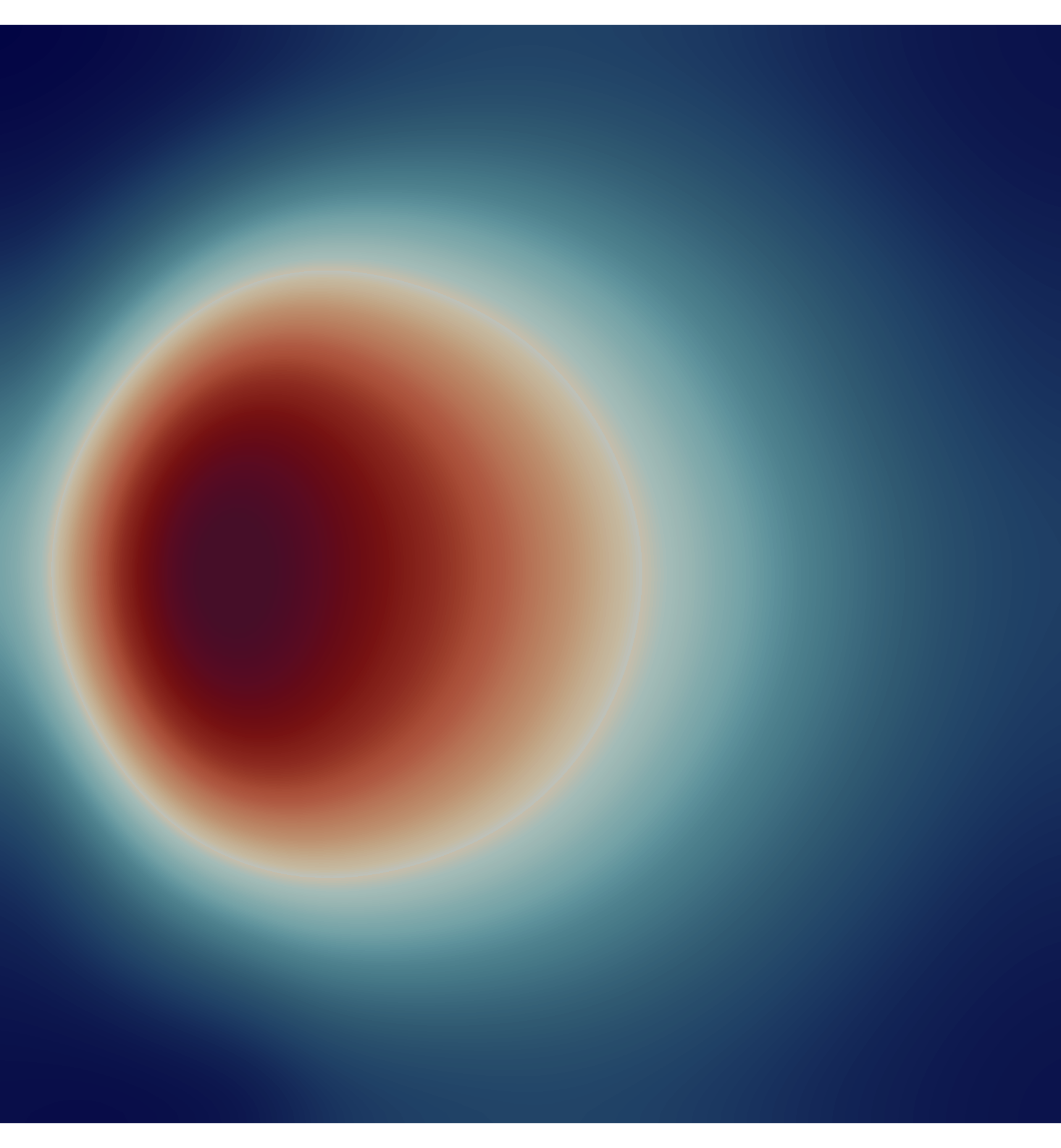}}
 \subcaptionbox{$\vec{\alpha}= [-0.4, 0.4]$ }
         {\includegraphics[width=0.2\textwidth]{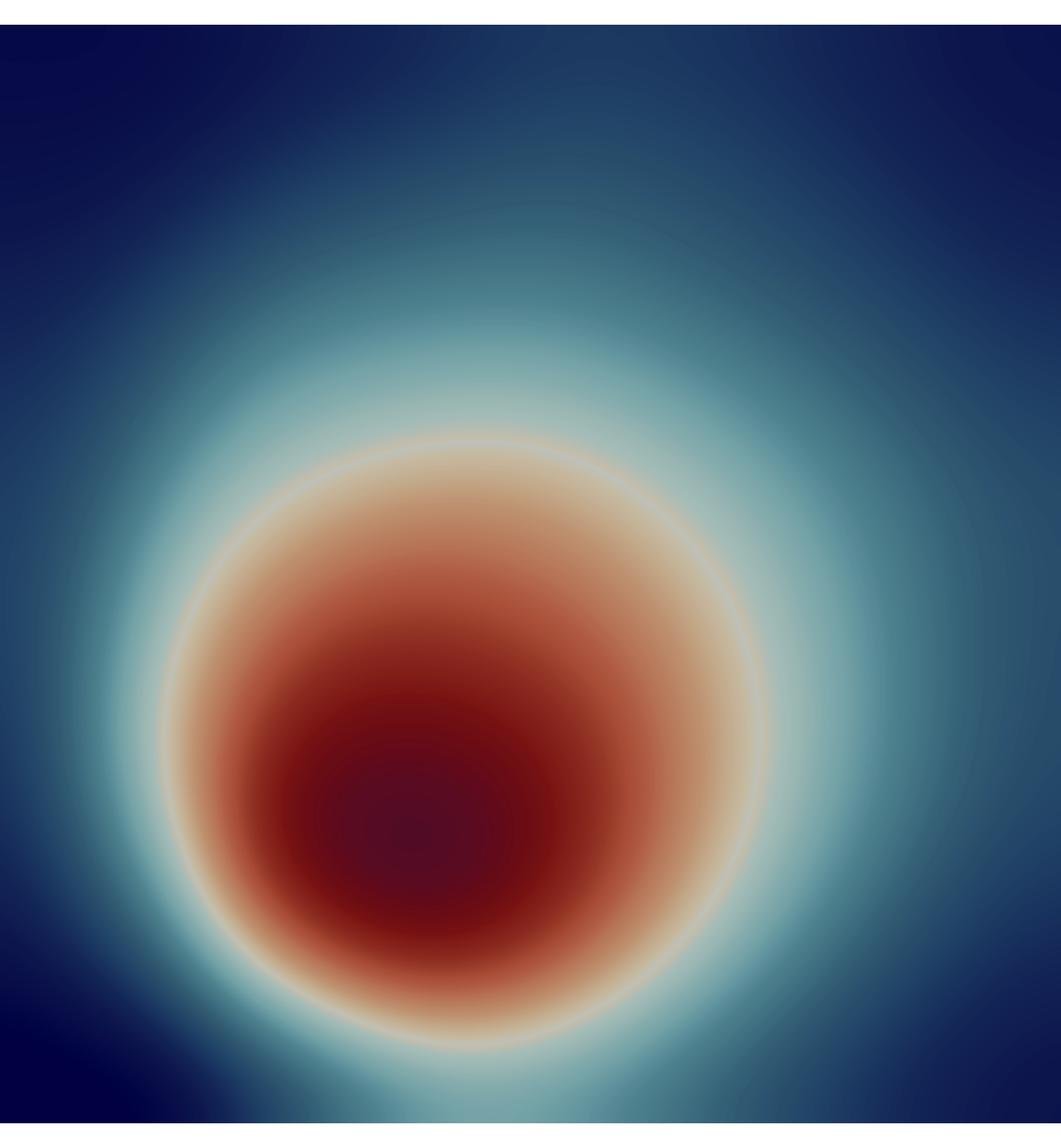}}\\
 \subcaptionbox{$\vec{\alpha}= [0, -0.4]$ }
         {\includegraphics[width=0.2\textwidth]{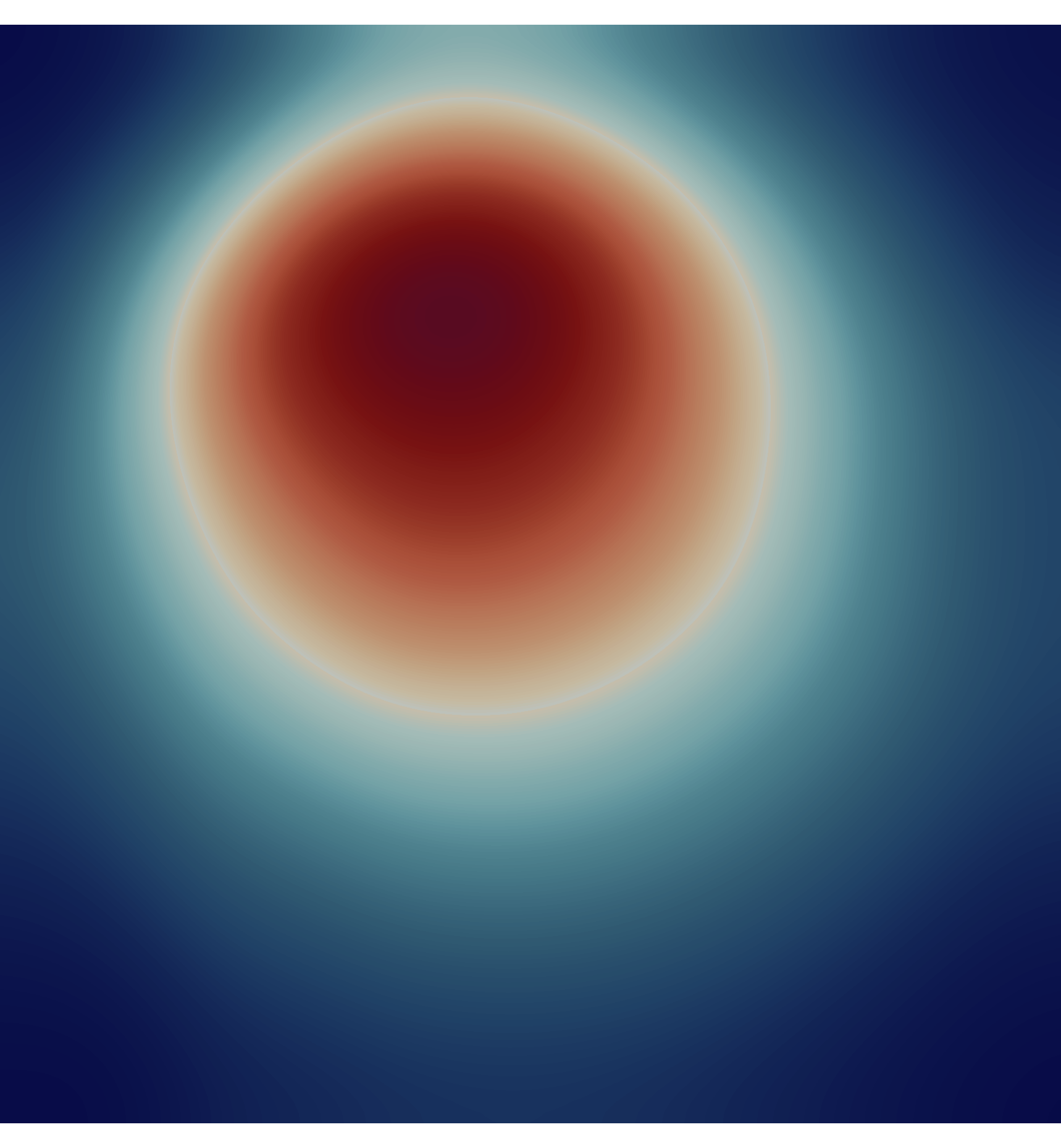}}
 \subcaptionbox{$\vec{\alpha}= [0, 0]$ }
         {\includegraphics[width=0.2\textwidth]{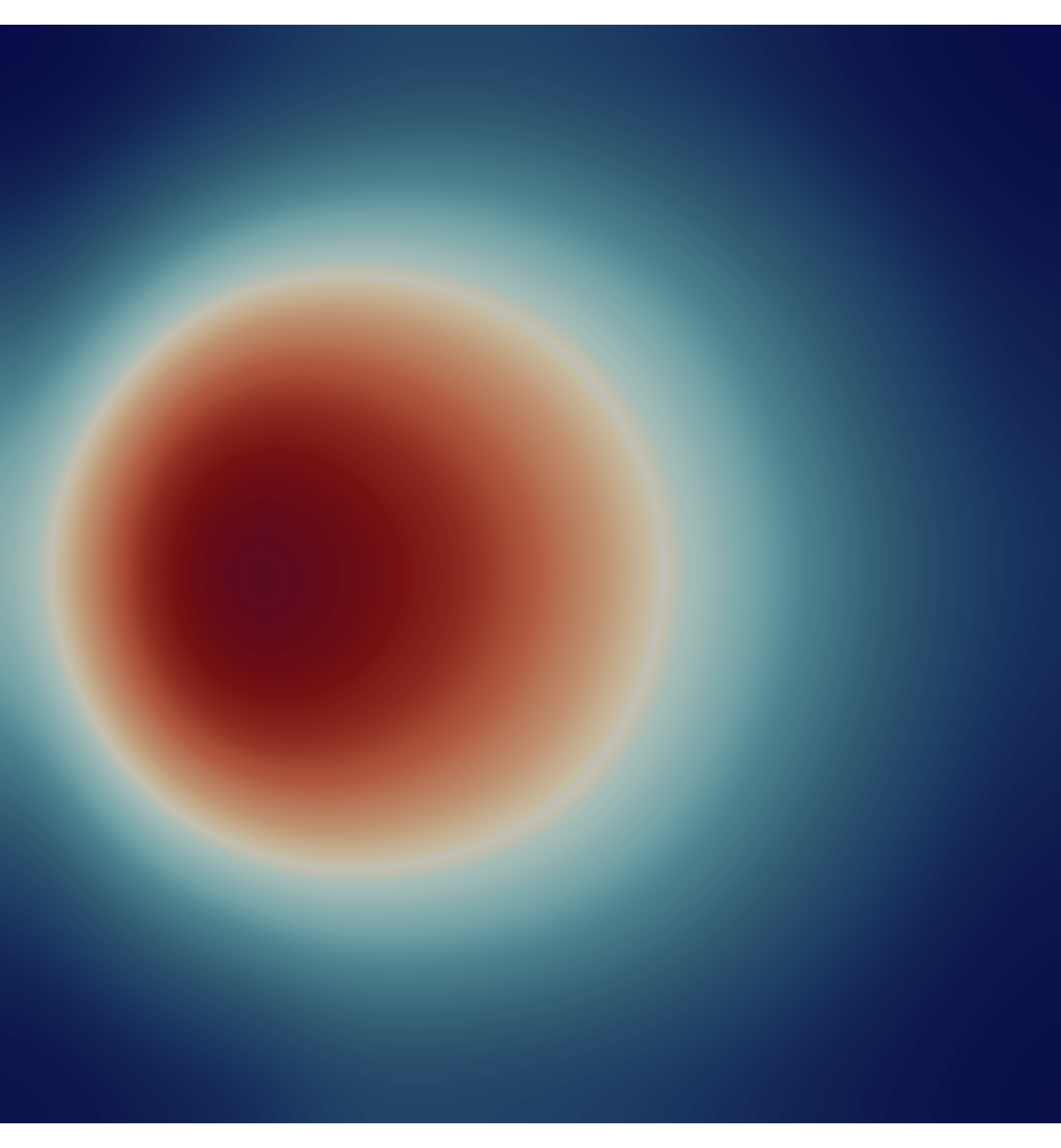}}
 \subcaptionbox{$\vec{\alpha}= [0, 0.4]$ }
         {\includegraphics[width=0.2\textwidth]{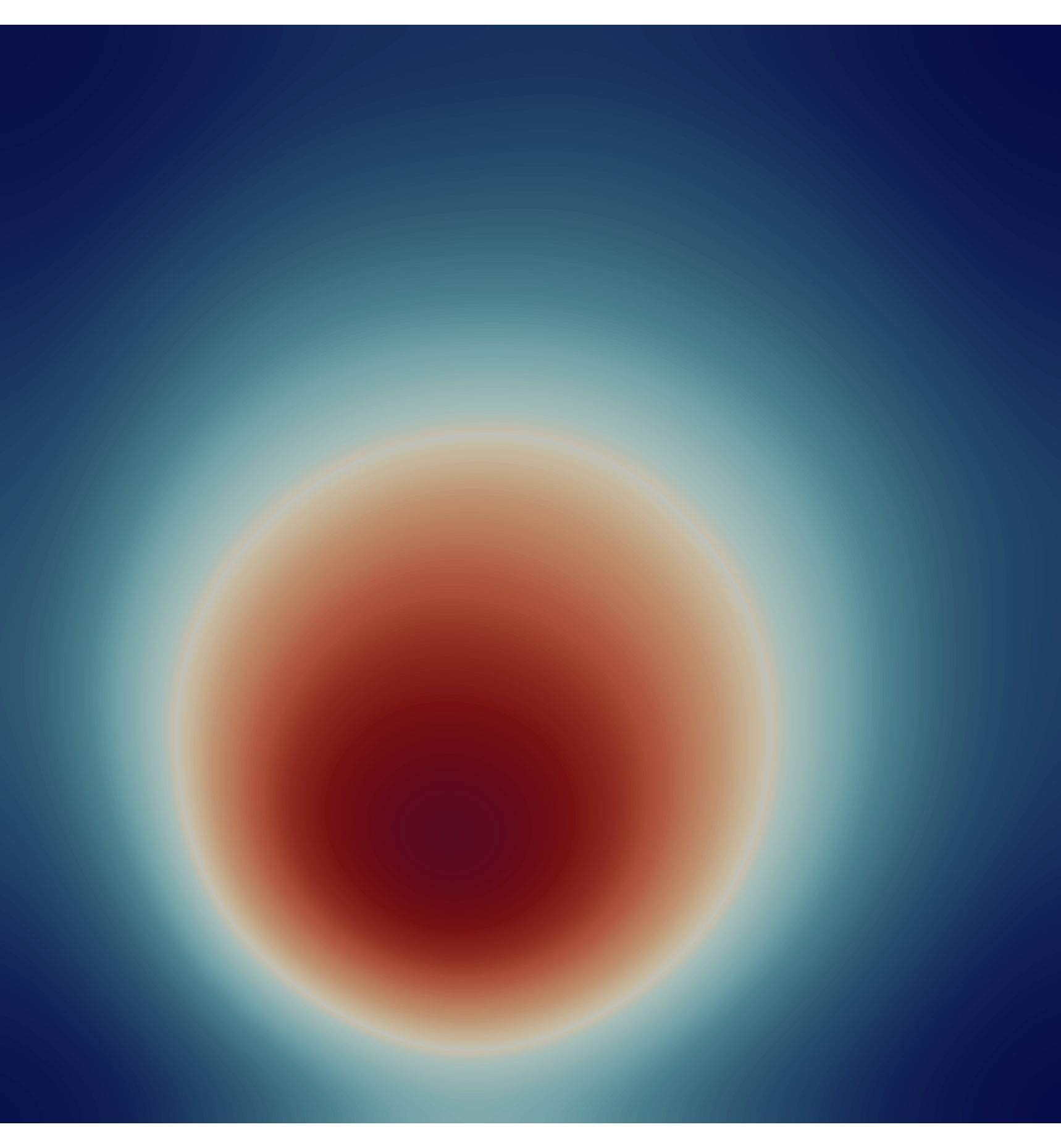}}\\
 \subcaptionbox{$\vec{\alpha}= [0.4, -0.4]$ }
         {\includegraphics[width=0.2\textwidth]{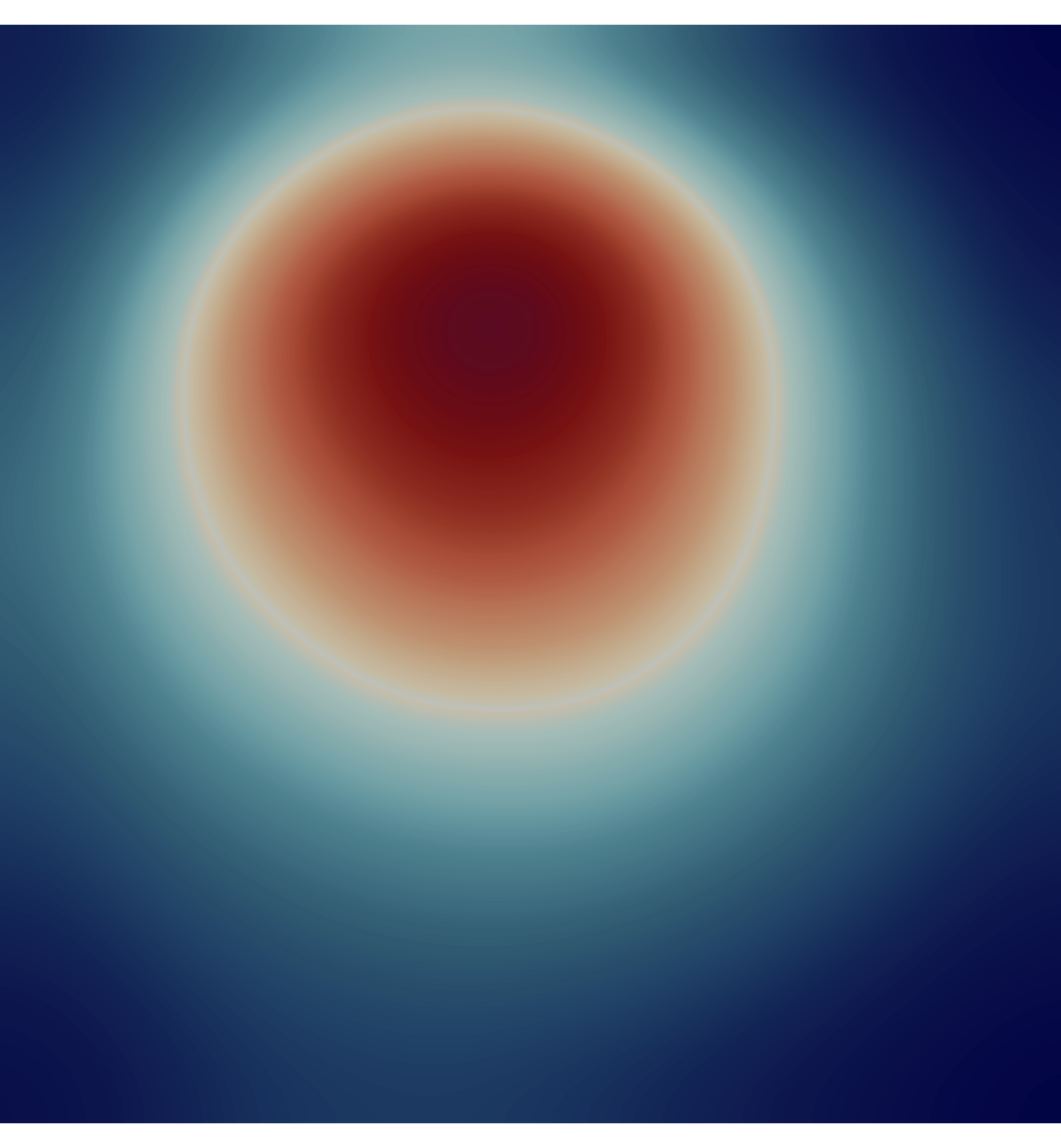}}
 \subcaptionbox{$\vec{\alpha}= [0.4, 0]$ }
         {\includegraphics[width=0.2\textwidth]{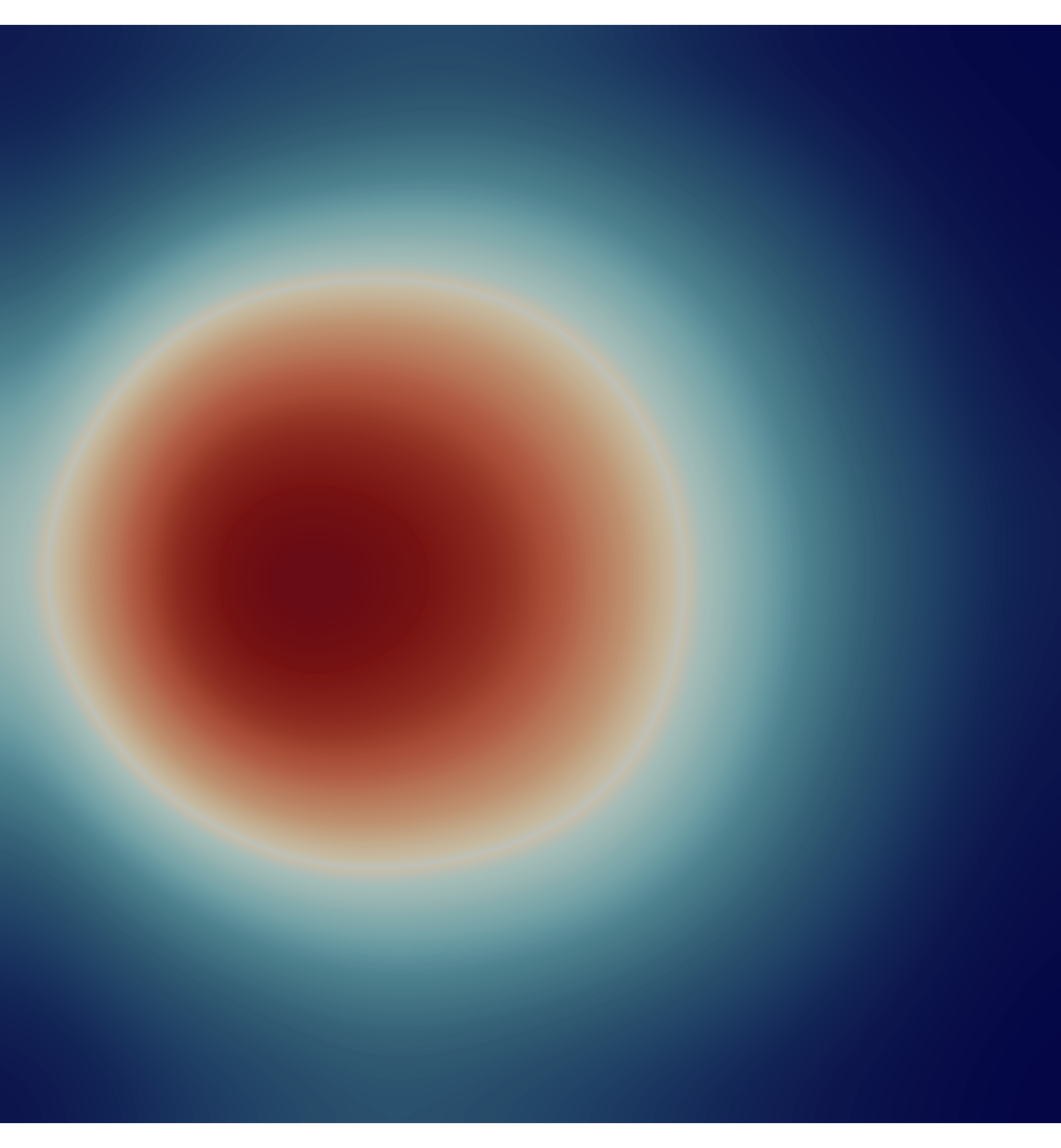}}
 \subcaptionbox{$\vec{\alpha}= [0.4, 0.4]$ }
         {\includegraphics[width=0.2\textwidth]{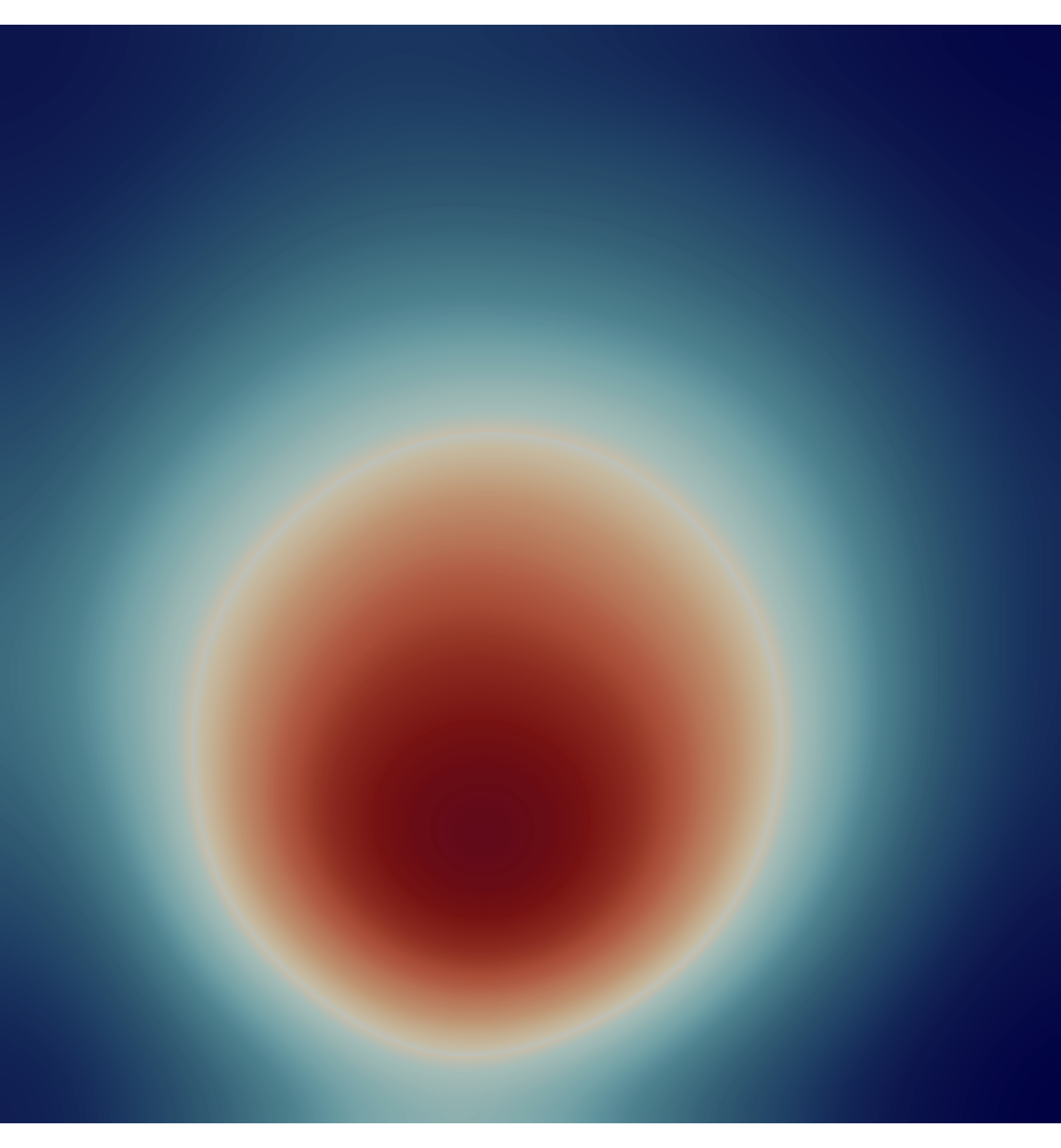}}
         \centering
        \includegraphics[width=0.6\textwidth]{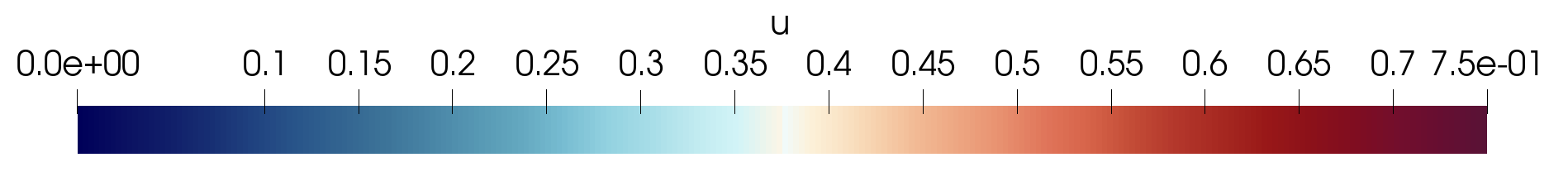}

\caption{Solution snapshots of the nonlinear advection diffusion problem at the final simulation time t=0.1 for different parameter values.}
\label{fig:adv_diff_sol2}
\end{figure}

It is noteworthy that when the error stabilizes around t= 0.02, there is hardly any visual difference between the different parameter values (Figure \ref{fig:adv_diff_sol1}), whereas the difference is clear at the final integration time (Figure \ref{fig:adv_diff_sol2}). This shows that the neural network and its tangent space are rich enough to handle increasingly different parametric solutions without a blow-up of the error.

\subsection{Advection-Diffusion on a domain with holes}
\subsubsection{General setup}
\begin{figure}
 \centering
 \subcaptionbox{Initial condition for $\vec{\alpha}$=[0.0,0.0].}
         {\includegraphics[width=0.32\textwidth]{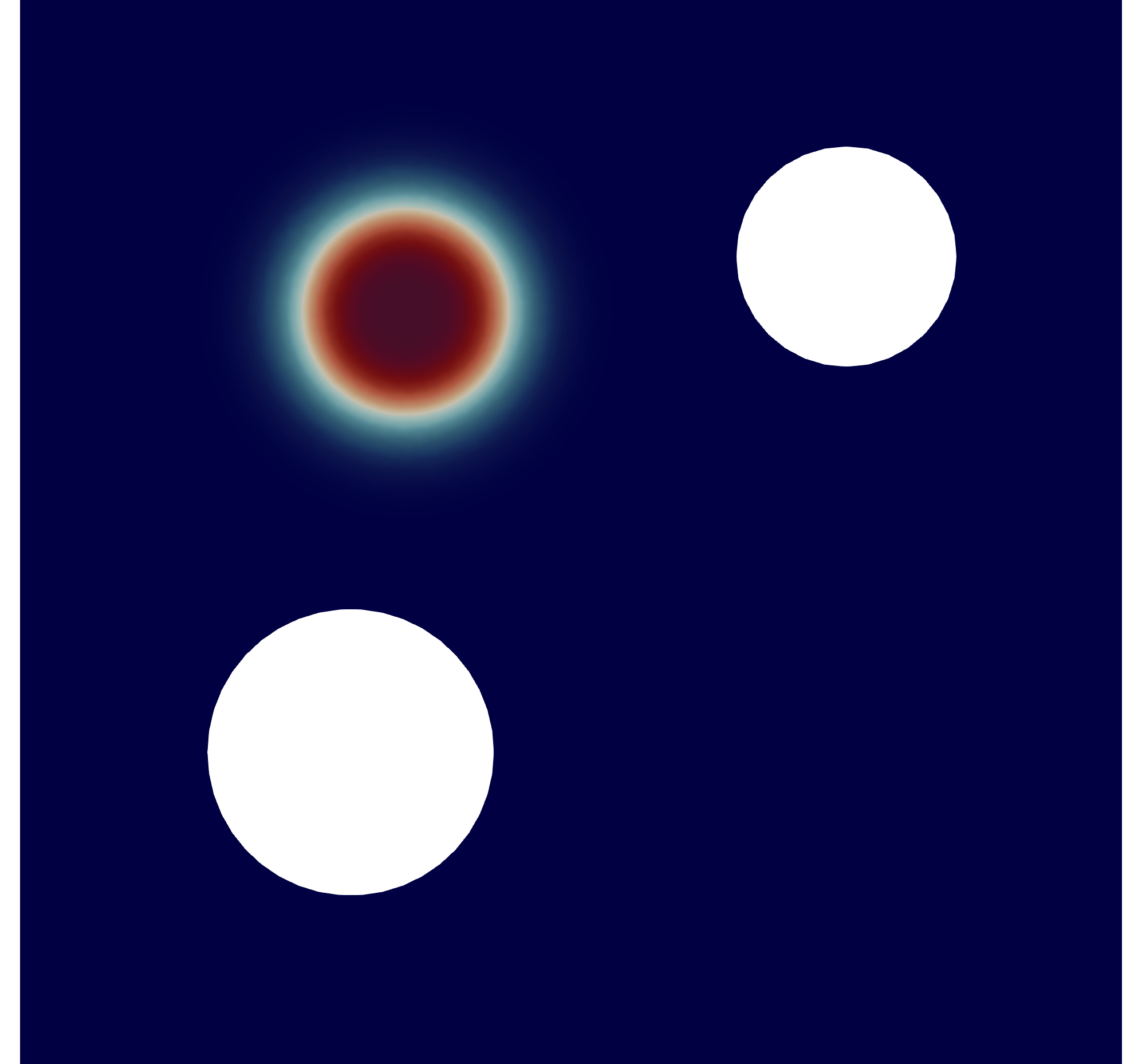}}
  \subcaptionbox{Velocity field $\vec{w}(x)$.}
         {\includegraphics[width=0.32\textwidth]{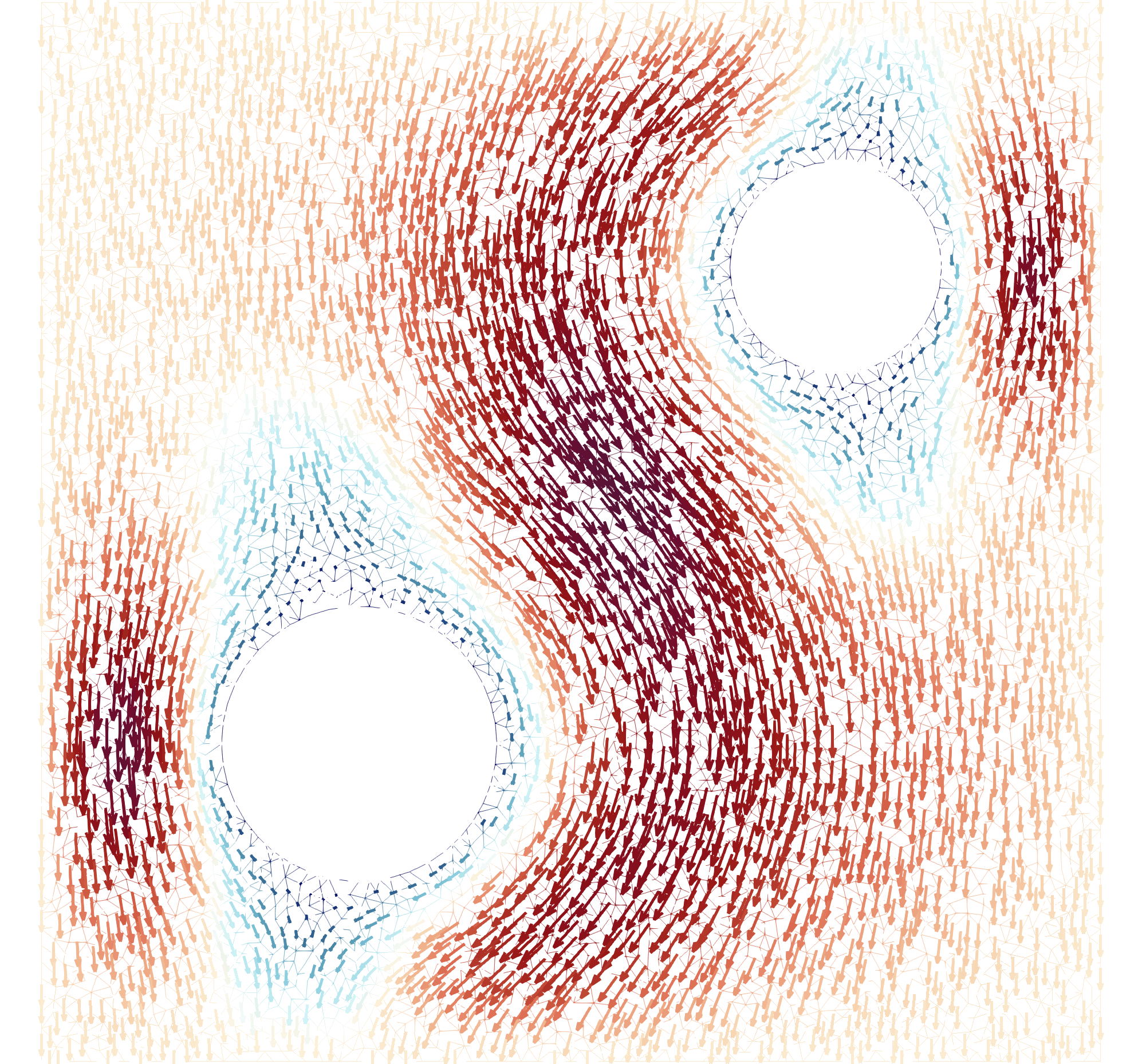}}
\caption{Setting for advection-diffusion problem on a domain with holes.}
\label{fig:setting_adv_diff_holes}
\end{figure}

We consider an advection-dominated equation on the unit squares with two circular holes:
\begin{equation}
\label{eq:adv_diff_holes}
    \begin{split}
             \dfrac{\partial u}{\partial t} &= 0.001 \nabla^2 u  +   \vecf{w}(\vec{x}) \cdot \nabla u, \\
            u(\vec{x},0) &= \text{sech}(-100((x_1-0.35.-0.1 \alpha_1)^2 +  (x_2-0.7-0.1\alpha_2)^2))^2,\\
                  \alpha_1, \alpha_2 & \in [-0.5, 0.5].
    \end{split}
\end{equation}
The initial condition represents a blob of concentration that  spreads according to the advective field $\vecf{w}$ and the small diffusive term.
We impose homogeneous Dirichlet boundary conditions on top and Neumann boundary conditions for the remaining domain boundary. The velocity field $\vecf{w}(\vec{x})$ is shown in Figure \ref{fig:setting_adv_diff_holes} and is precomputed numerically by solving the steady-state Navier Stokes equations, see \ref{sec:appvel}  for details. We compute the reference solutions to \eqref{eq:adv_diff_holes} on the same mesh with a FE method using P2 elements and an implicit Euler method with time step $\Delta t=$ 1e-4.

We first consider a single forward run and then compute the case for the parametrized initial condition. 
We consider positional embeddings of size $n_\Phi=4$, $n_\Phi=10$ and $n_\Phi= 15$, corresponding to the first $n_\Phi$ eigenfunctions of the Laplace operator which are computed numerically on the same FE mesh using P3 elements, see Figure \ref{fig:eigenfunc} in \ref{sec:eigenfuncs} for a visualization. As the solution converges to 0 towards the end of the simulation, we normalize the $L_2$ -error with the norm of the initial condition for this experiment.

\subsubsection{Single parameter case}
\begin{figure}
 \centering
 \subcaptionbox{Full mesh - no sampling}
         {\includegraphics[width=0.45\textwidth]{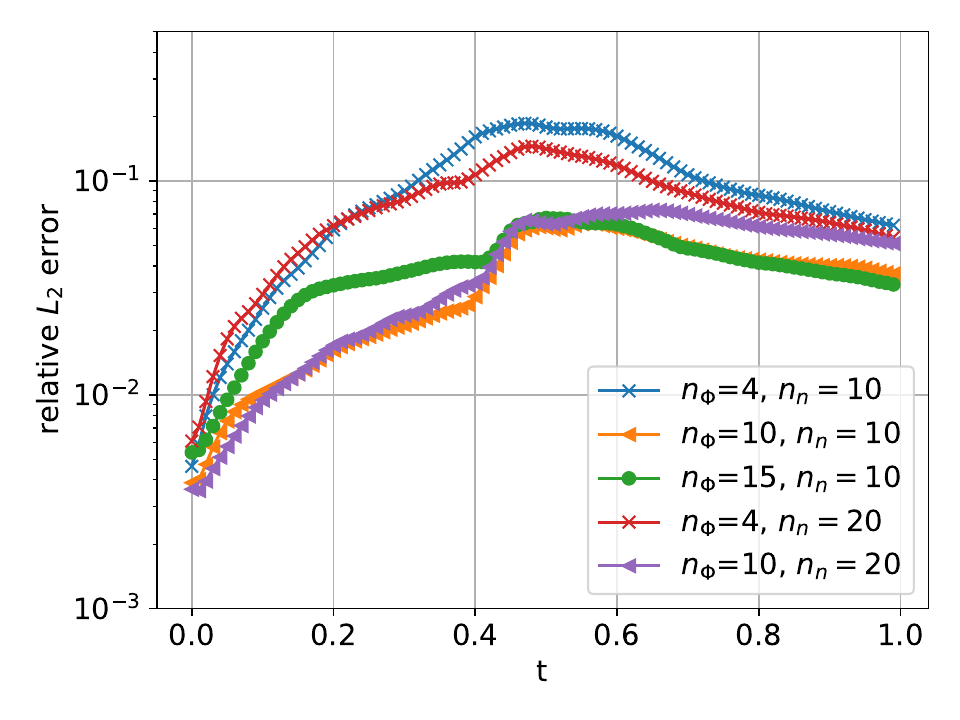}}
  \subcaptionbox{$n_s =$ 500}
         {\includegraphics[width=0.45\textwidth]{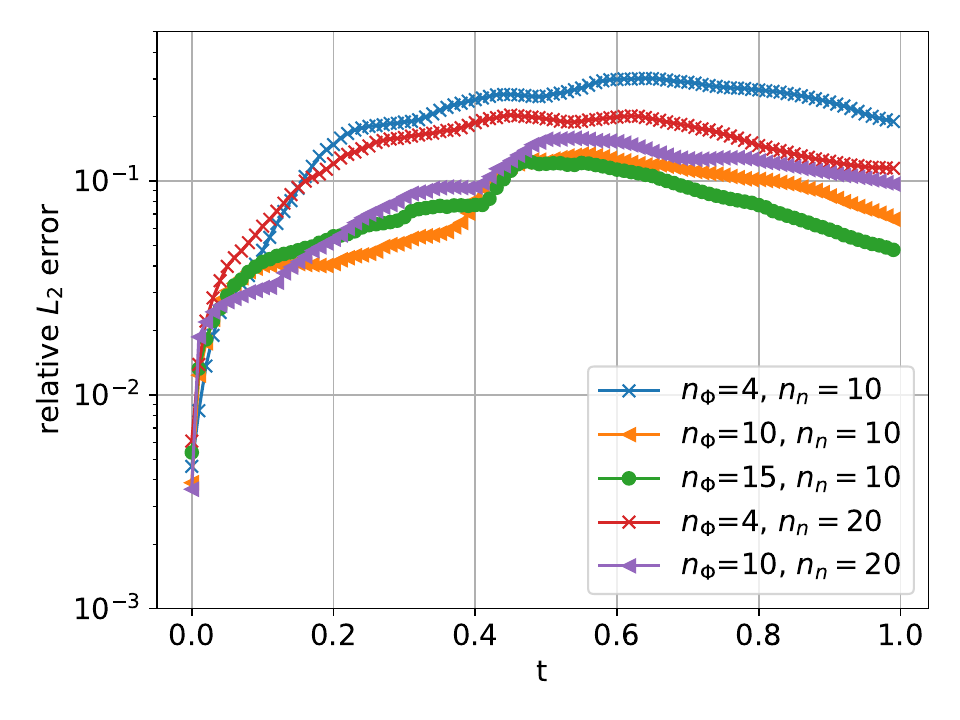}}
   \subcaptionbox{$n_s=$ 1000}
         {\includegraphics[width=0.45\textwidth]{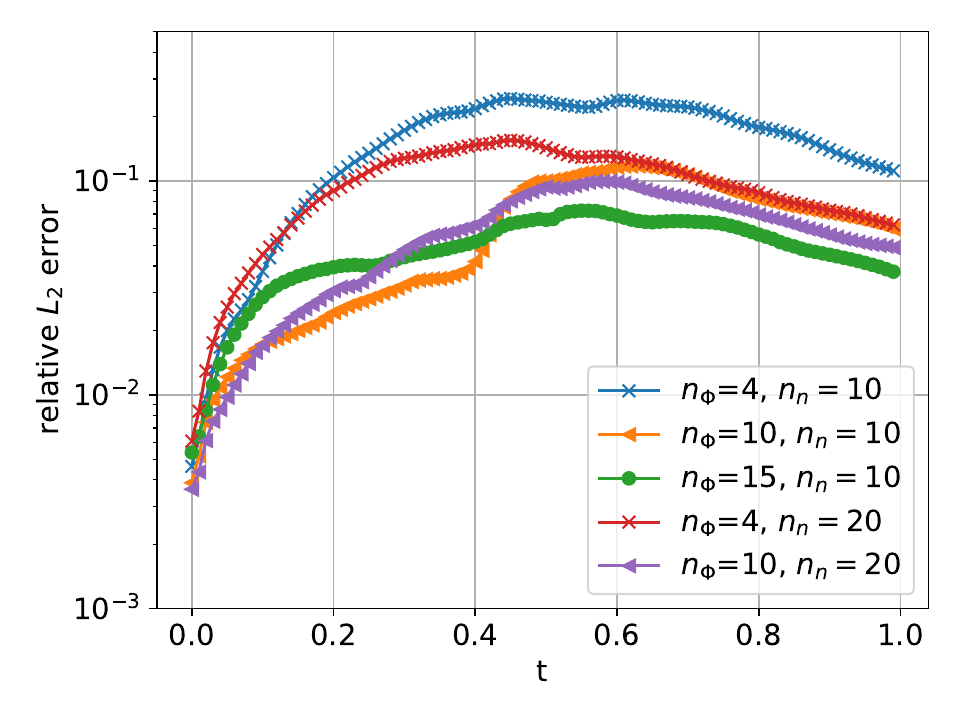}}
   \subcaptionbox{ $n_s=$ 2000}
         {\includegraphics[width=0.45\textwidth]{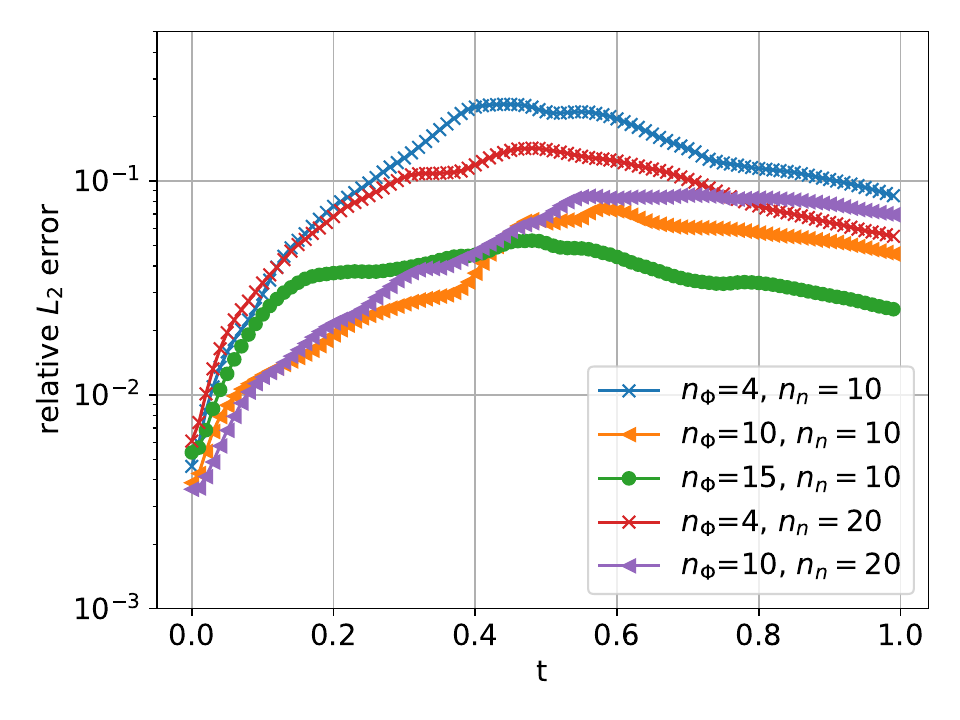}}
\caption{Error plots for the advection-diffusion problem on the domain with holes. We compare the evaluation on the full-mesh, versus subsampling $n_s$ points in every integration step. Each error plot shows varying architectures with $n_\Phi$ harmonic features and $n_n$ number of hidden neurons per layer.}
\label{fig:error_adv_diff_holes}
\end{figure}

\begin{figure}
 \centering
         {\includegraphics[width=0.22\textwidth]{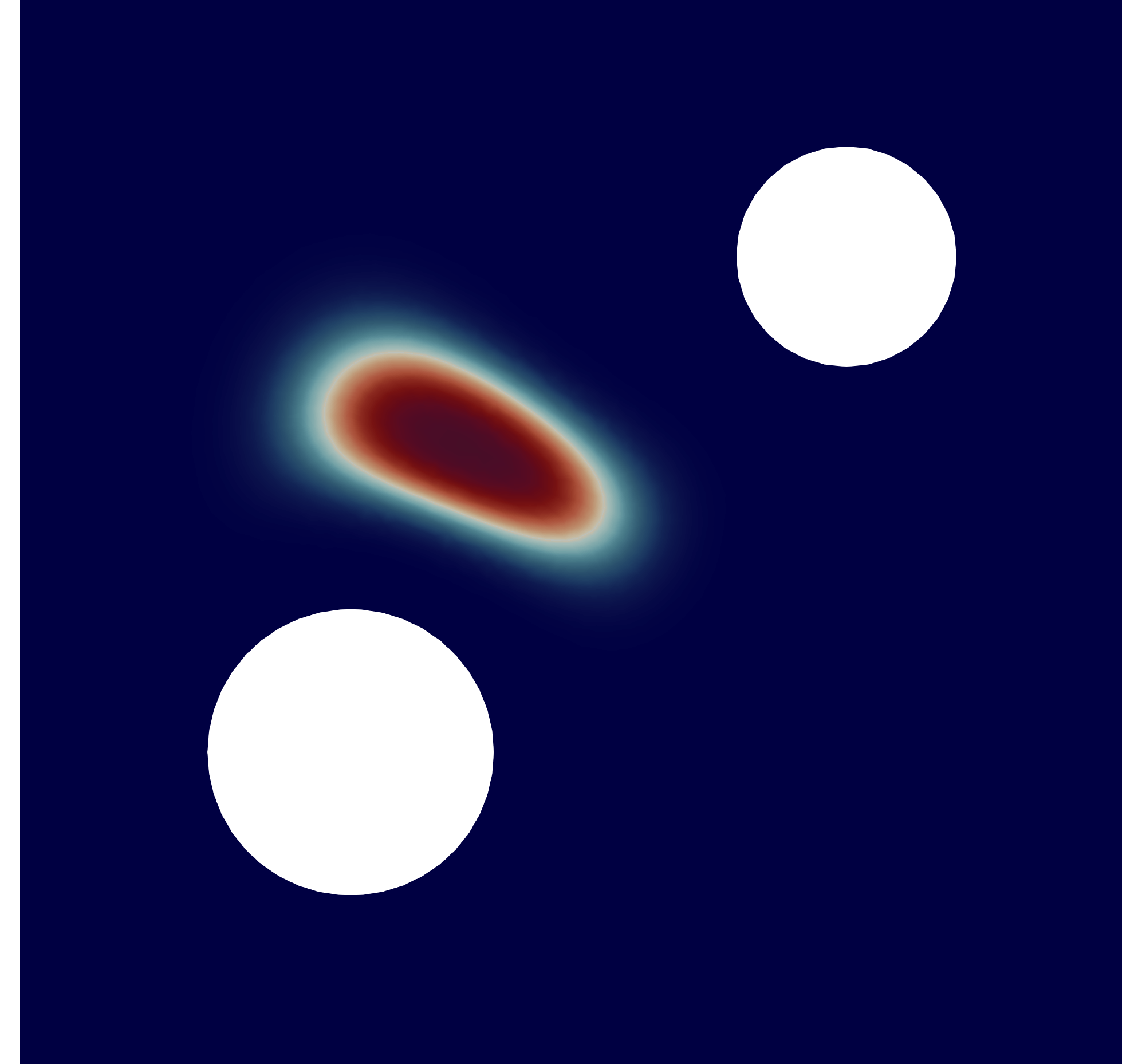}}
         {\includegraphics[width=0.22\textwidth]{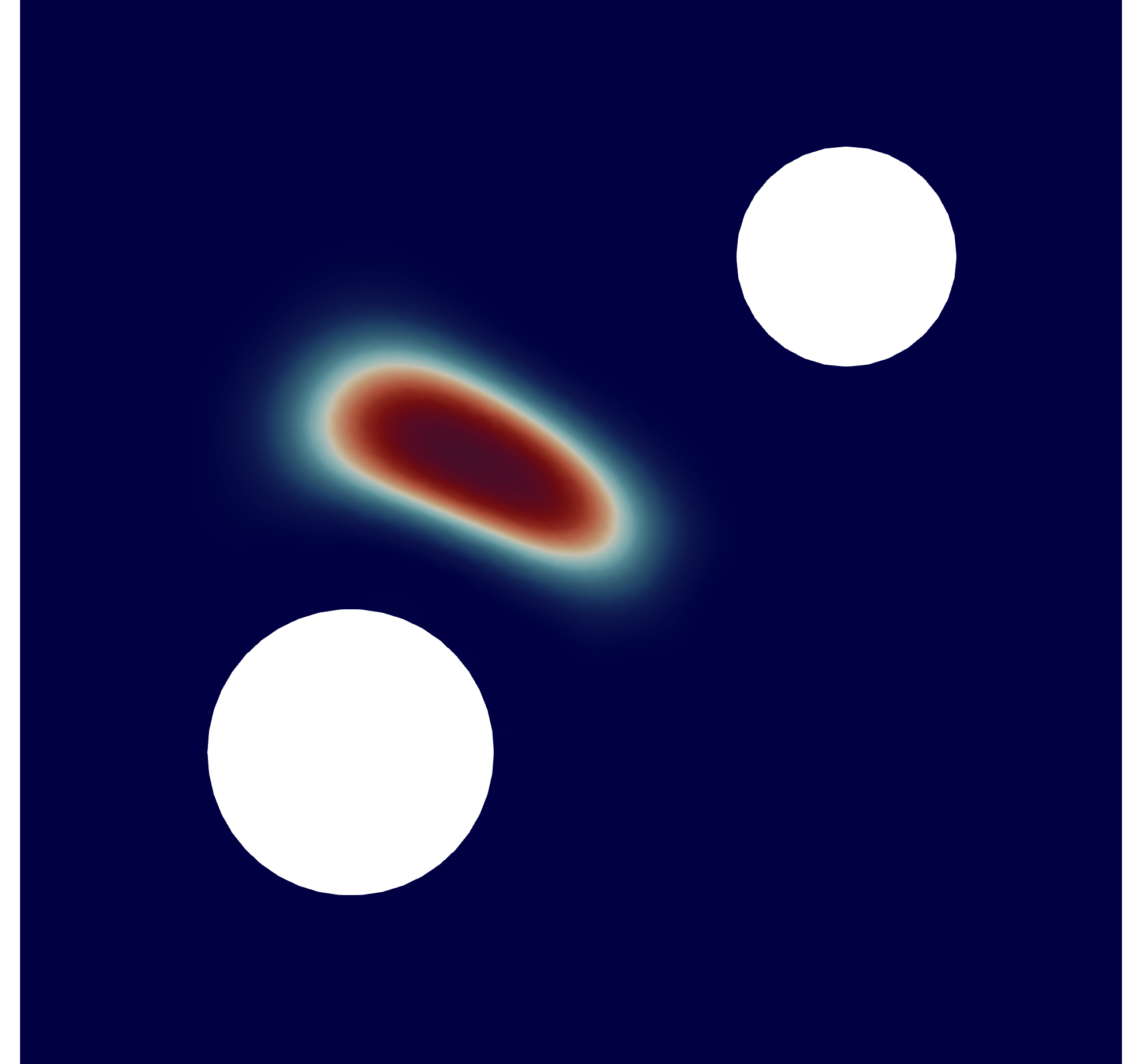}}
         {\includegraphics[width=0.22\textwidth]{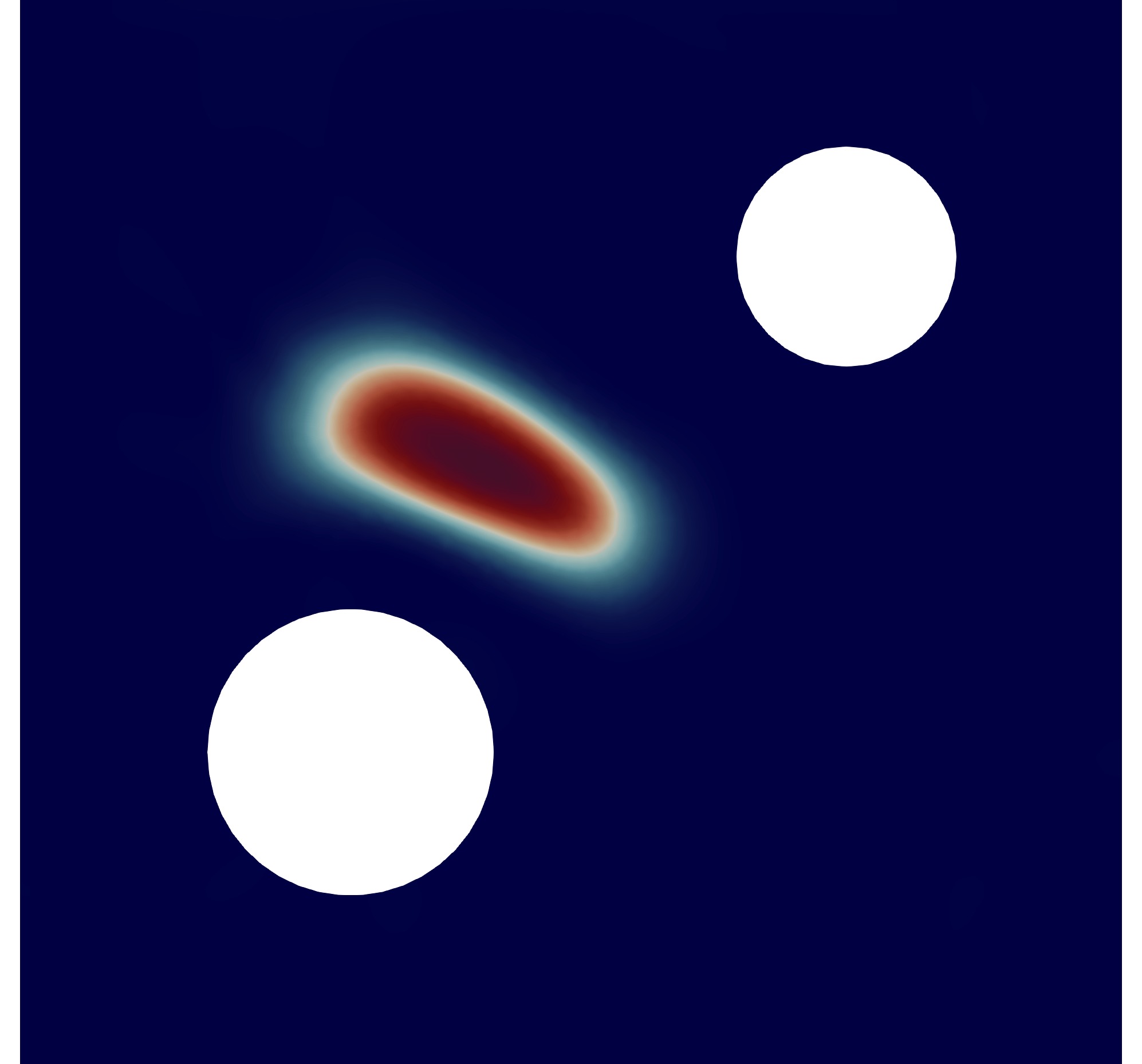}}
         {\includegraphics[width=0.22\textwidth]{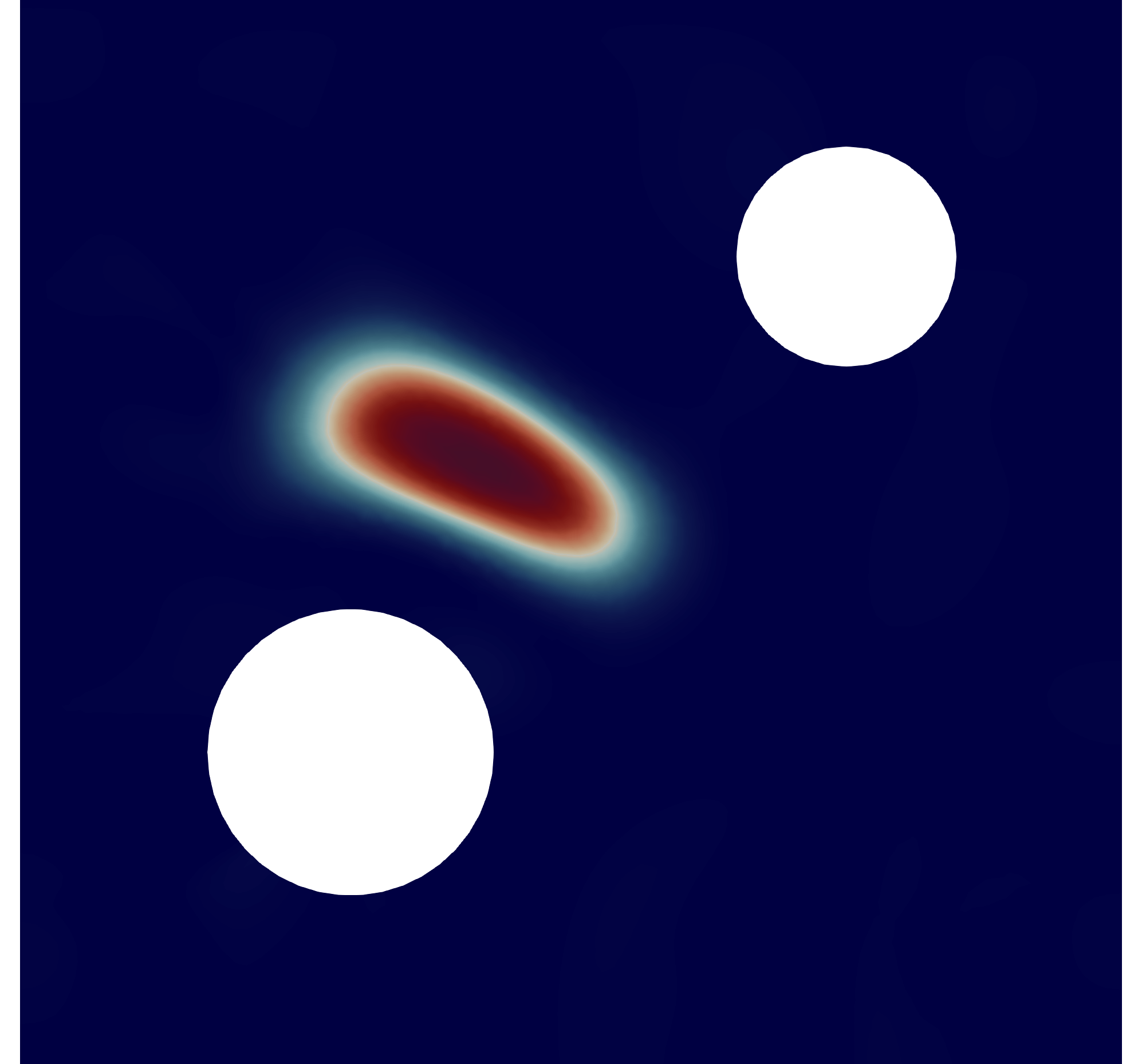}}
         {\includegraphics[width=0.22\textwidth]{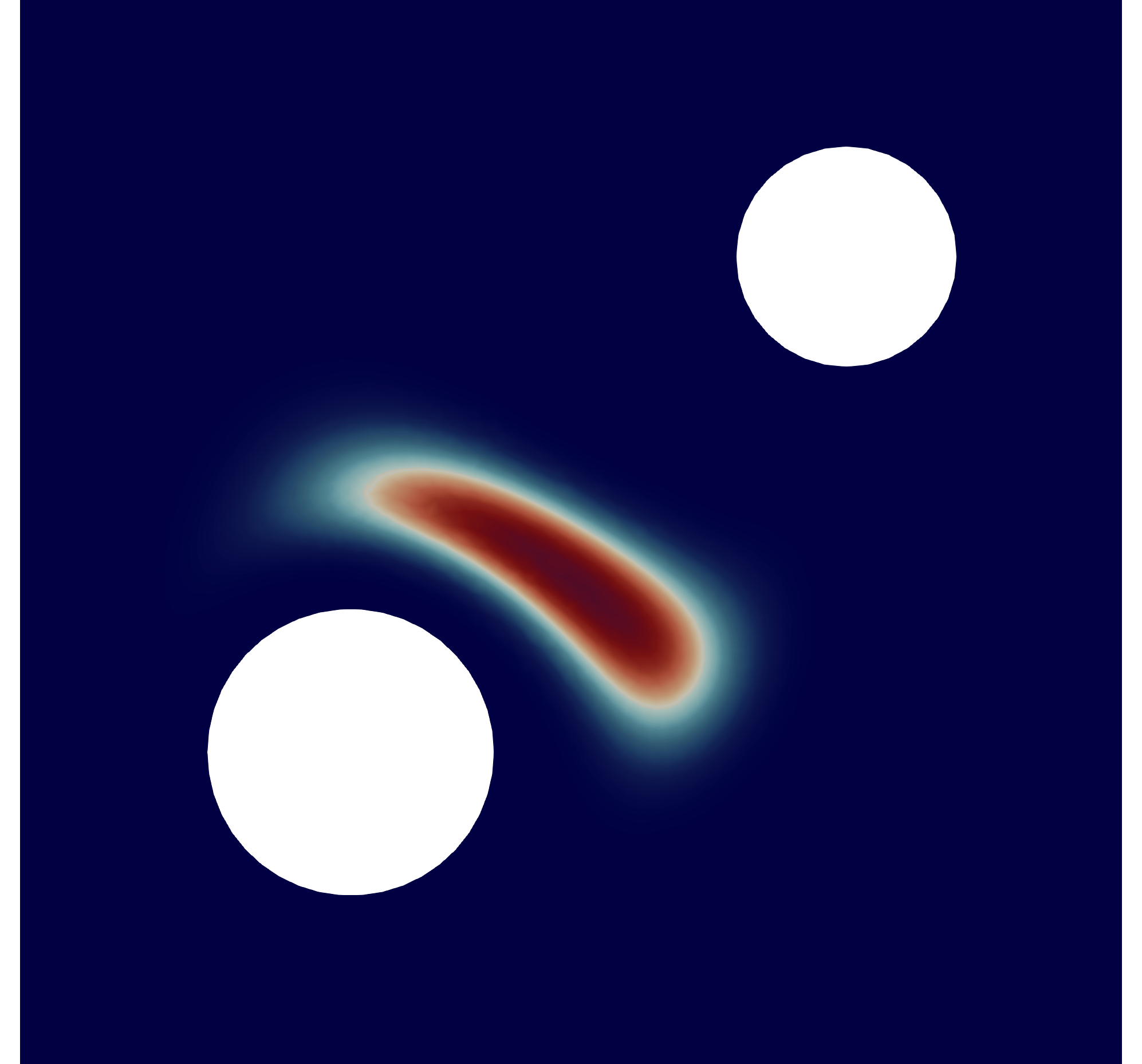}}
         {\includegraphics[width=0.22\textwidth]{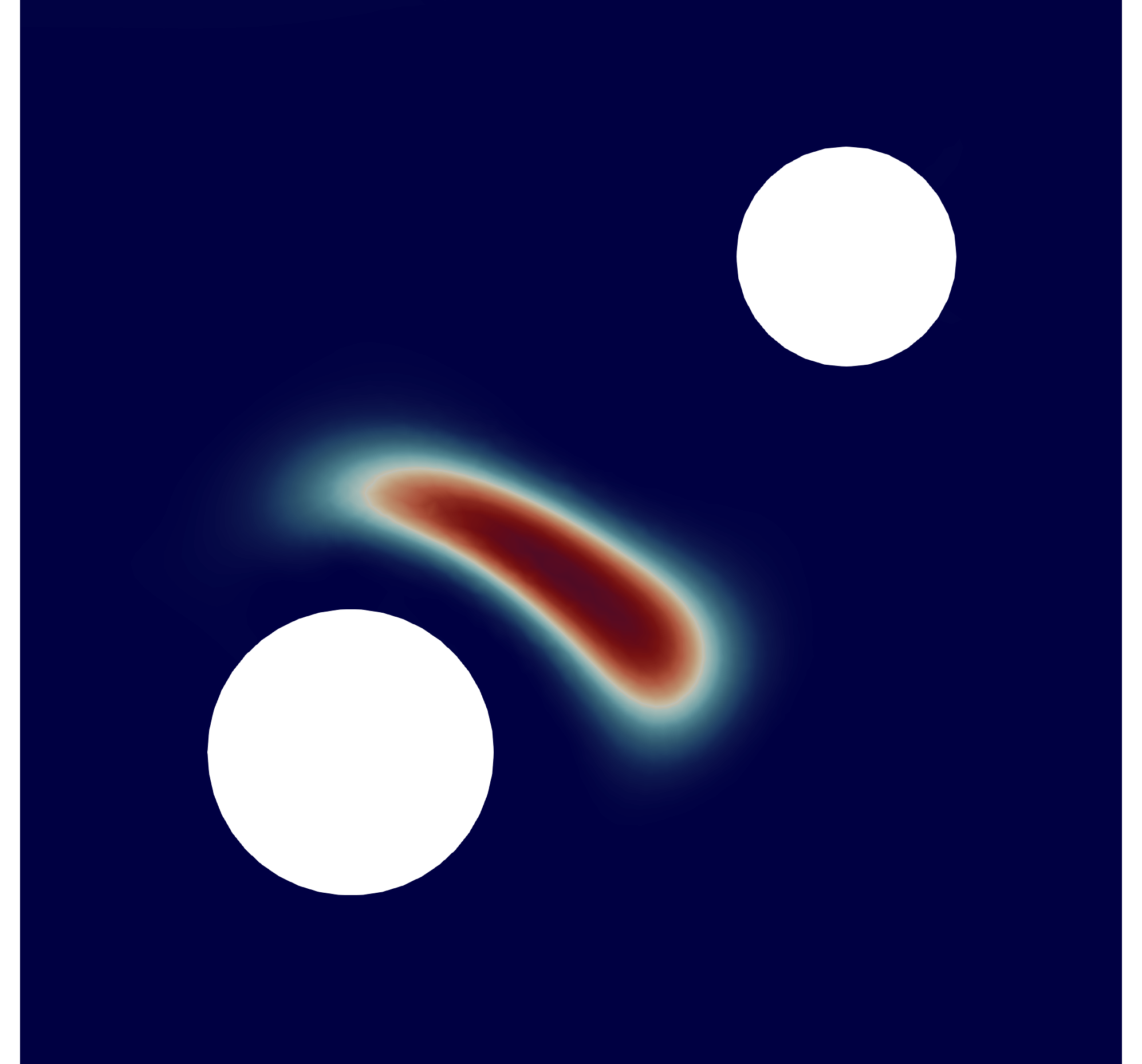}}
         {\includegraphics[width=0.22\textwidth]{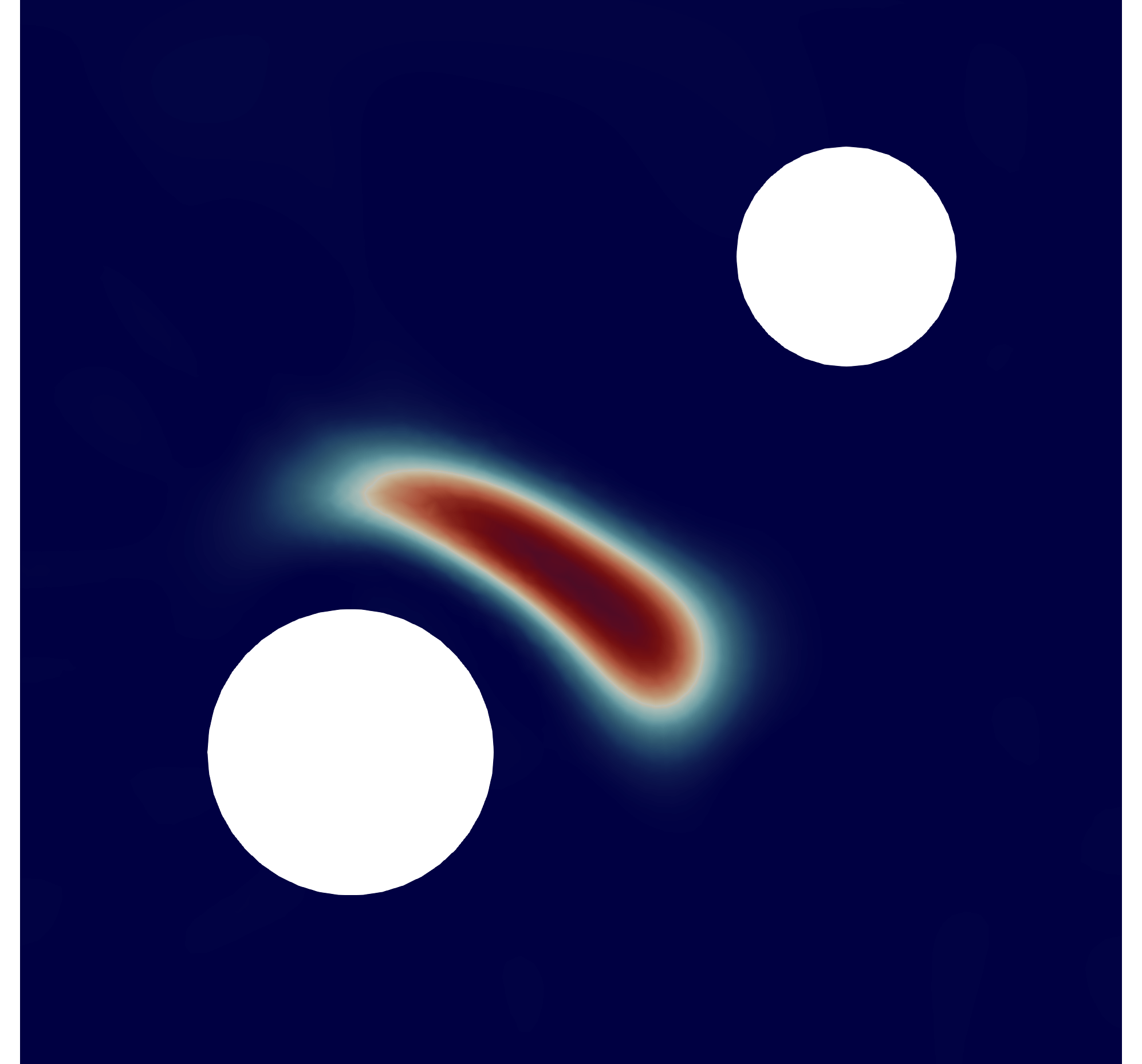}}
         {\includegraphics[width=0.22\textwidth]{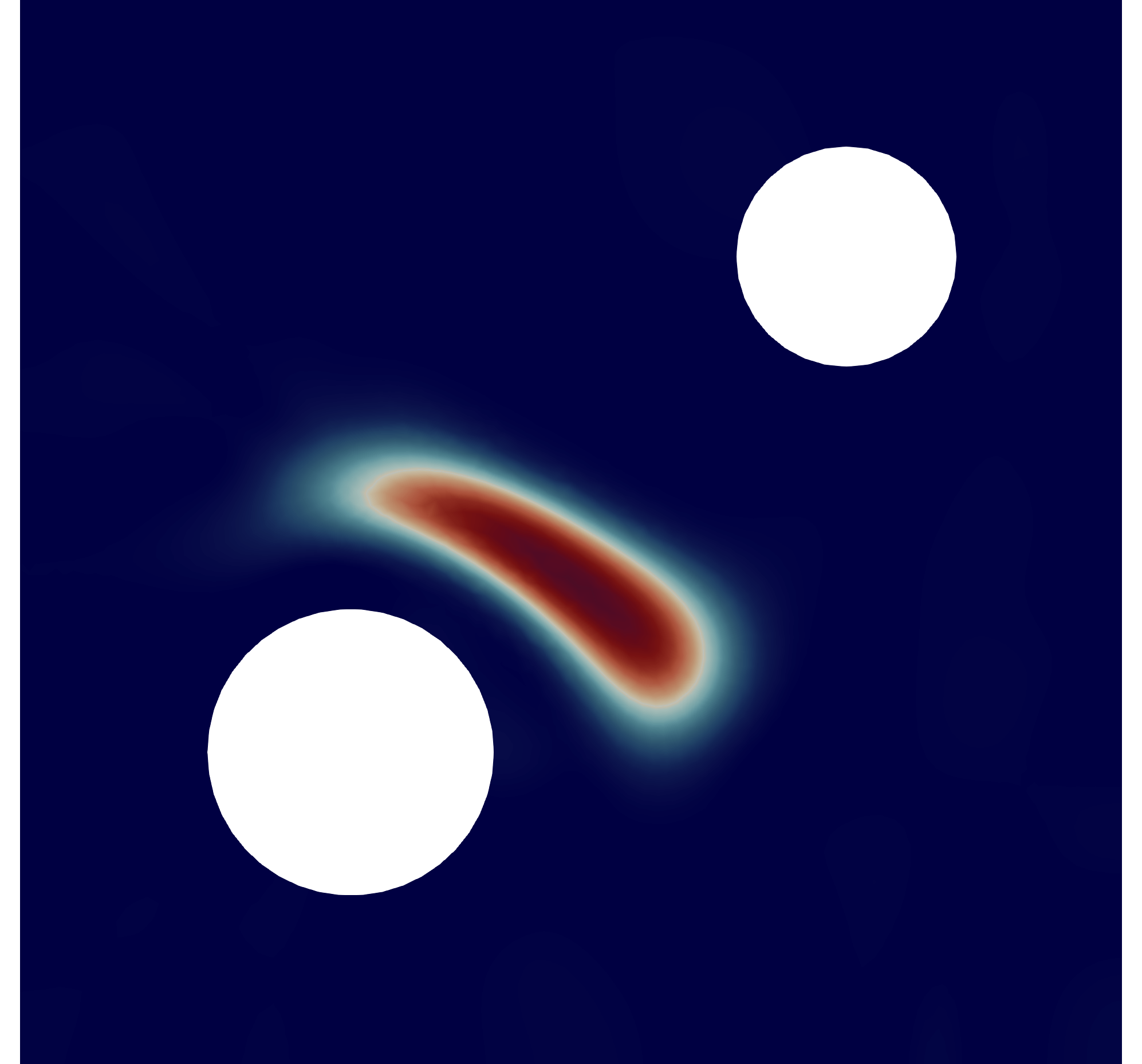}}
         {\includegraphics[width=0.22\textwidth]{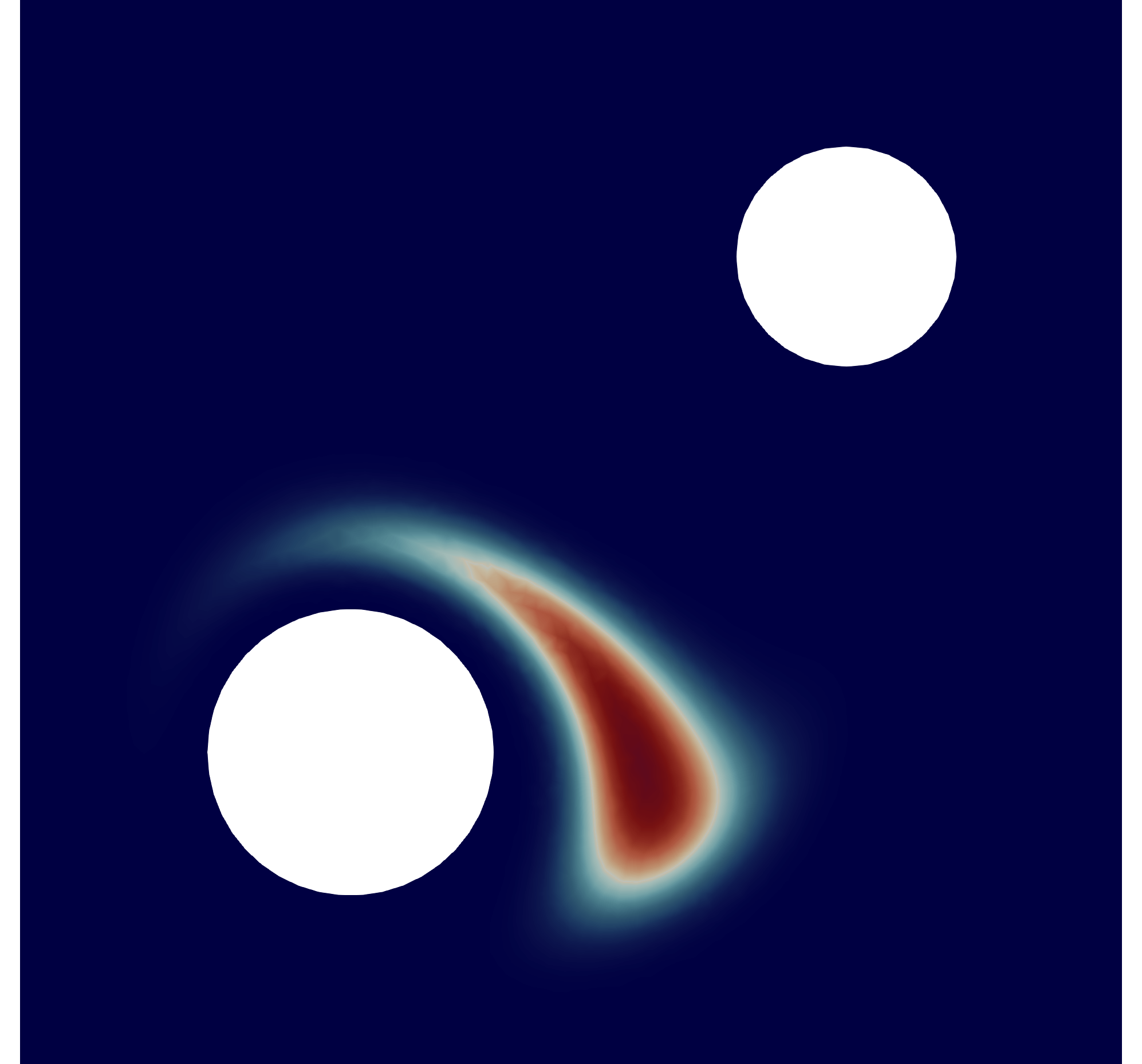}}
         {\includegraphics[width=0.22\textwidth]{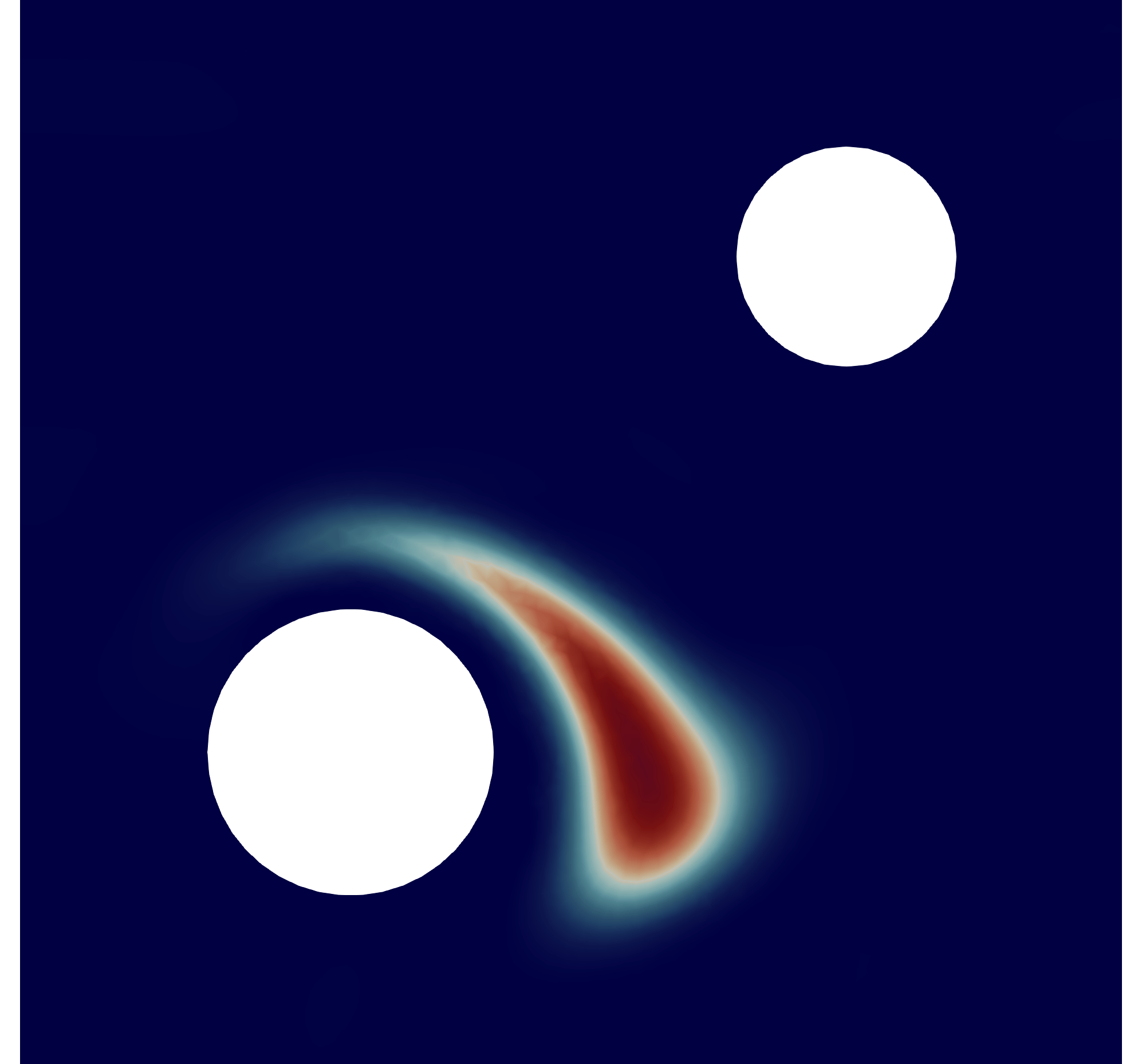}}
         {\includegraphics[width=0.22\textwidth]{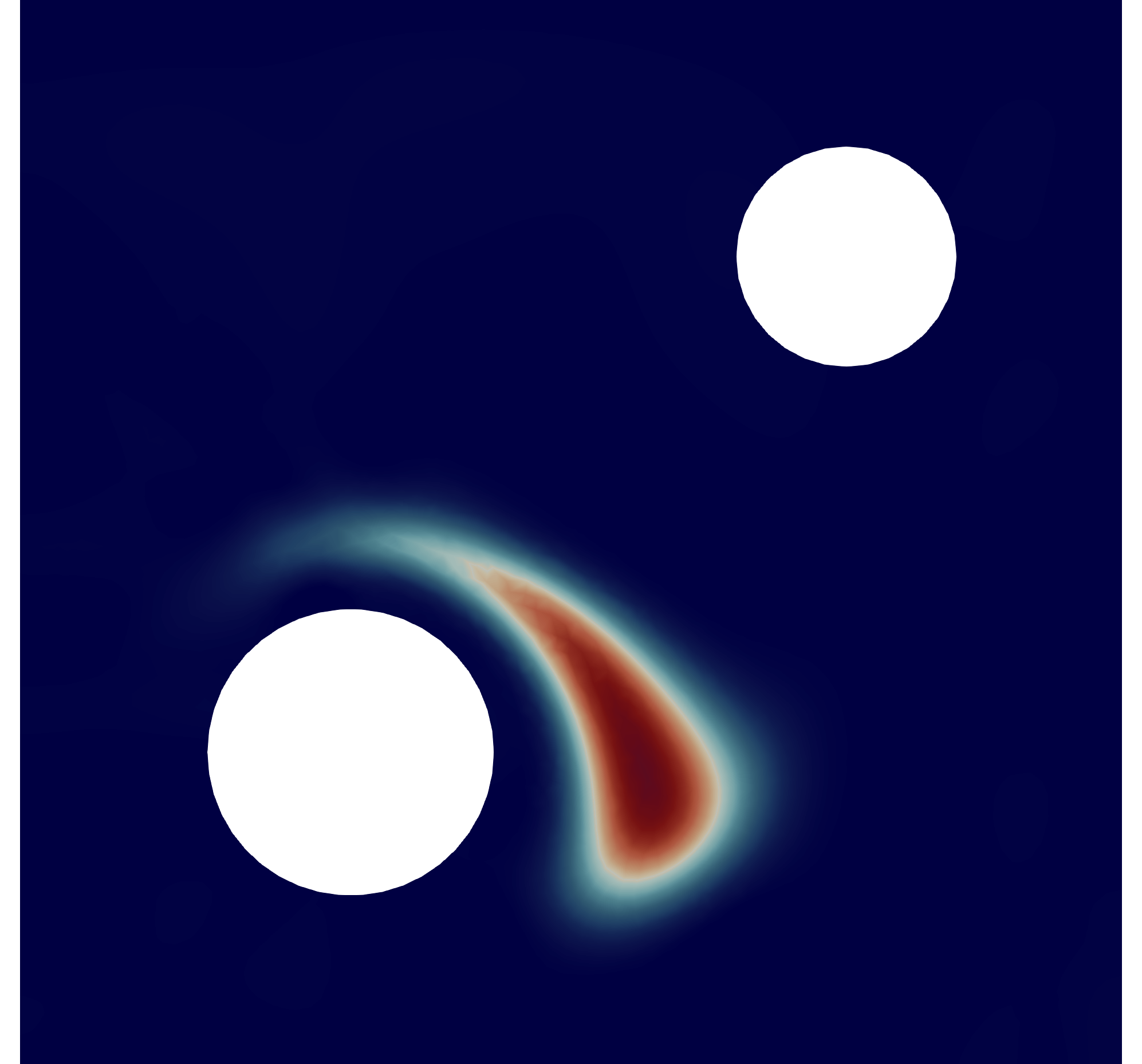}}
         {\includegraphics[width=0.22\textwidth]{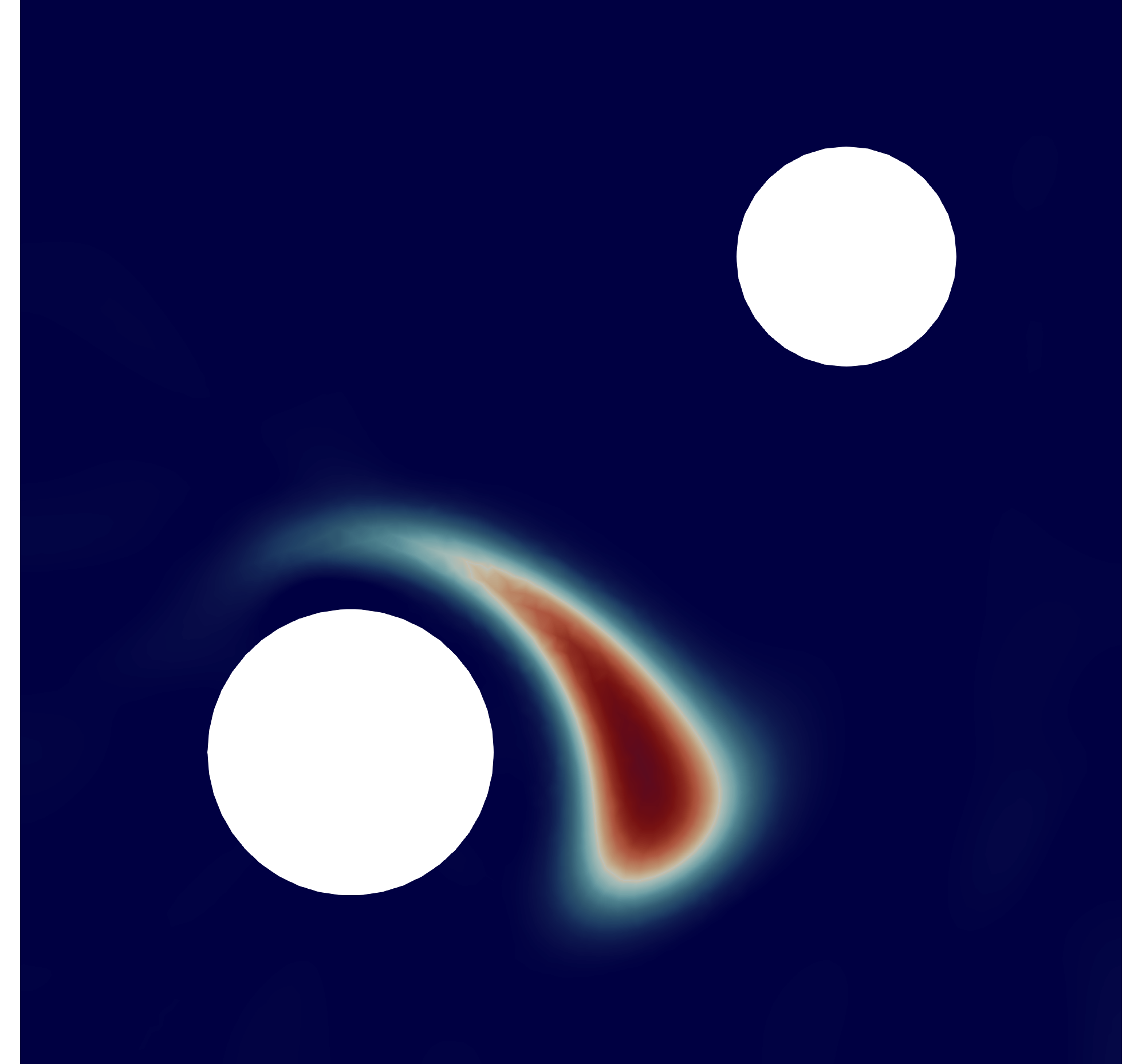}}
         {\includegraphics[width=0.22\textwidth]{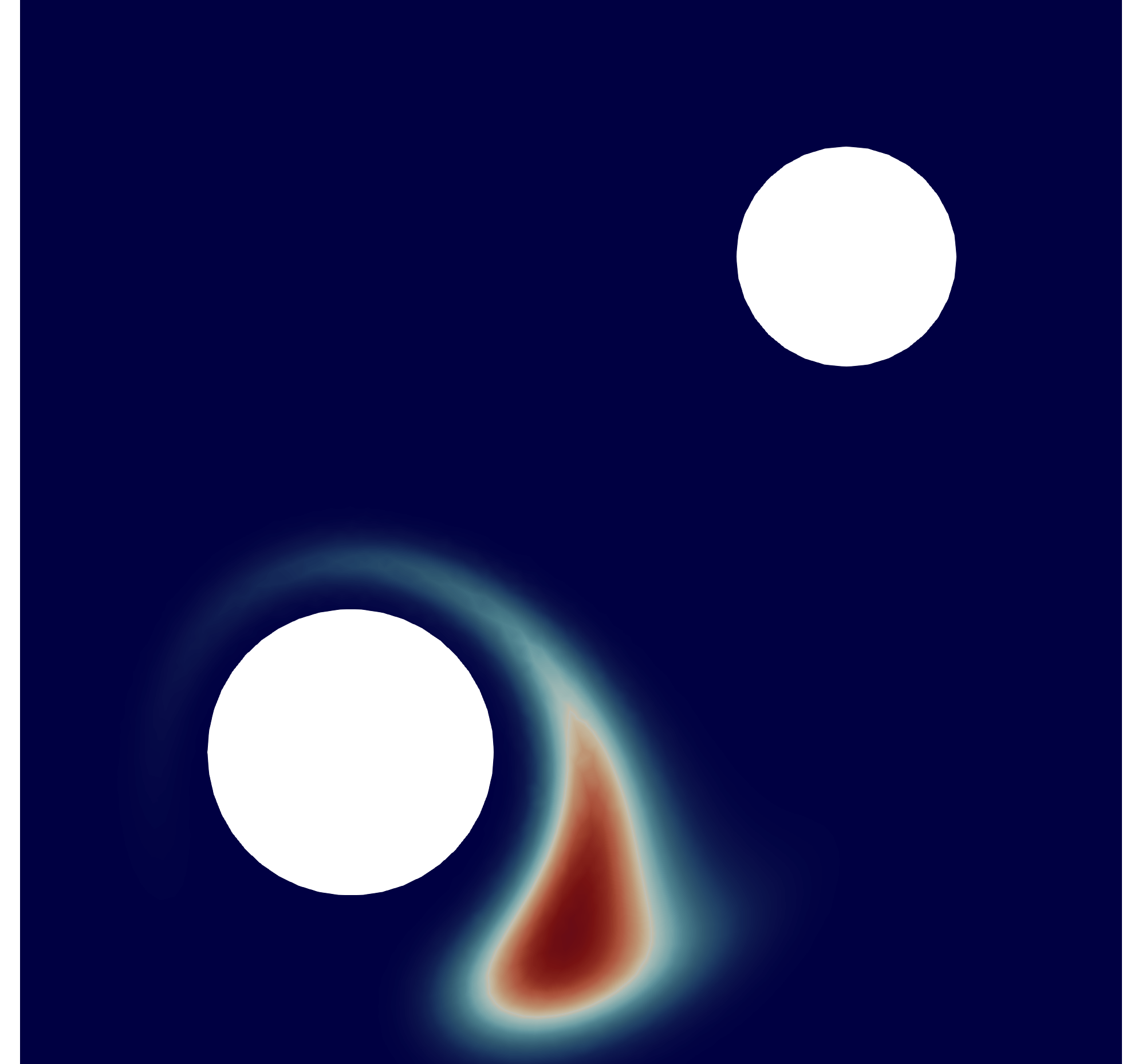}}
         {\includegraphics[width=0.22\textwidth]{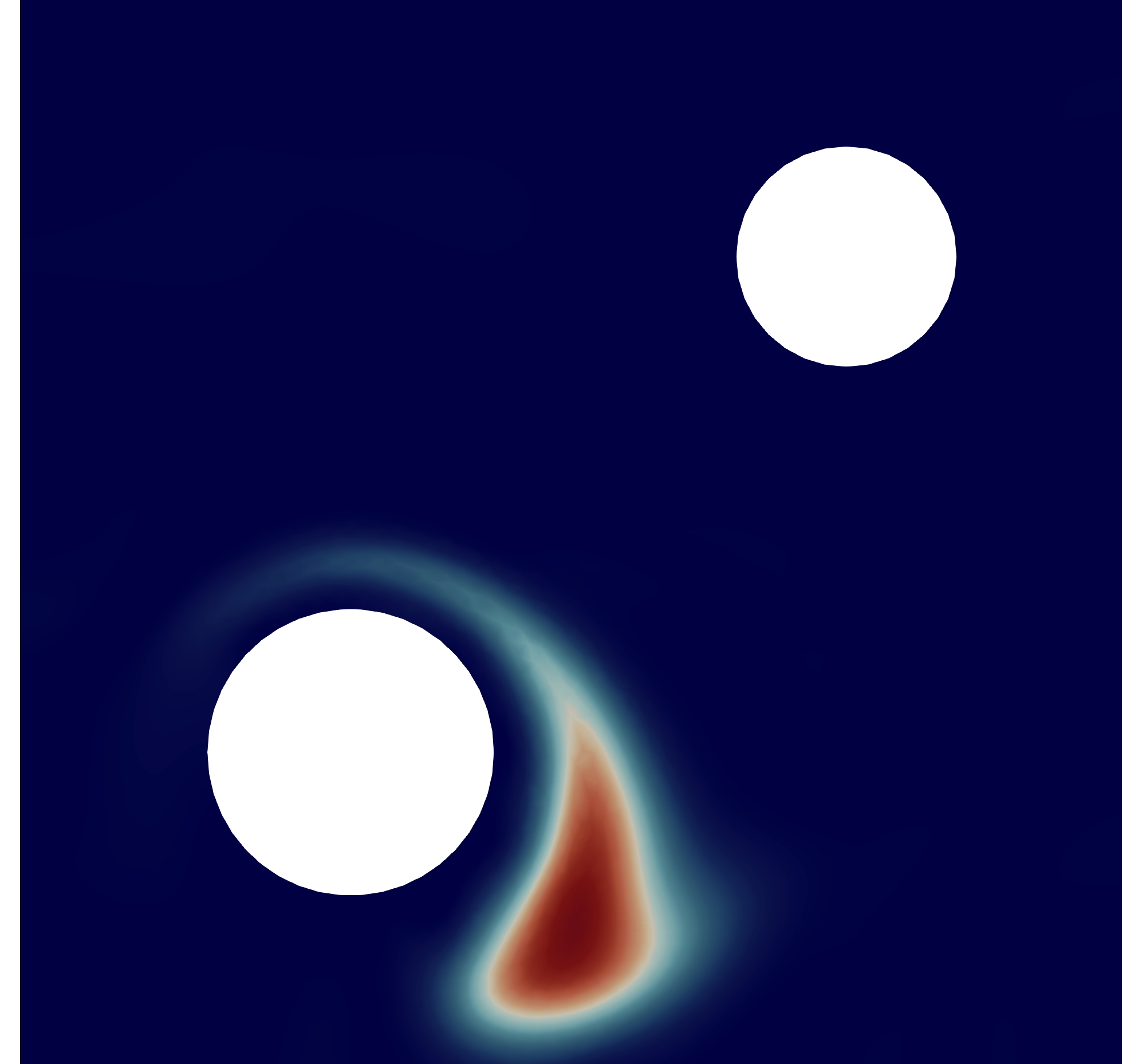}}
         {\includegraphics[width=0.22\textwidth]{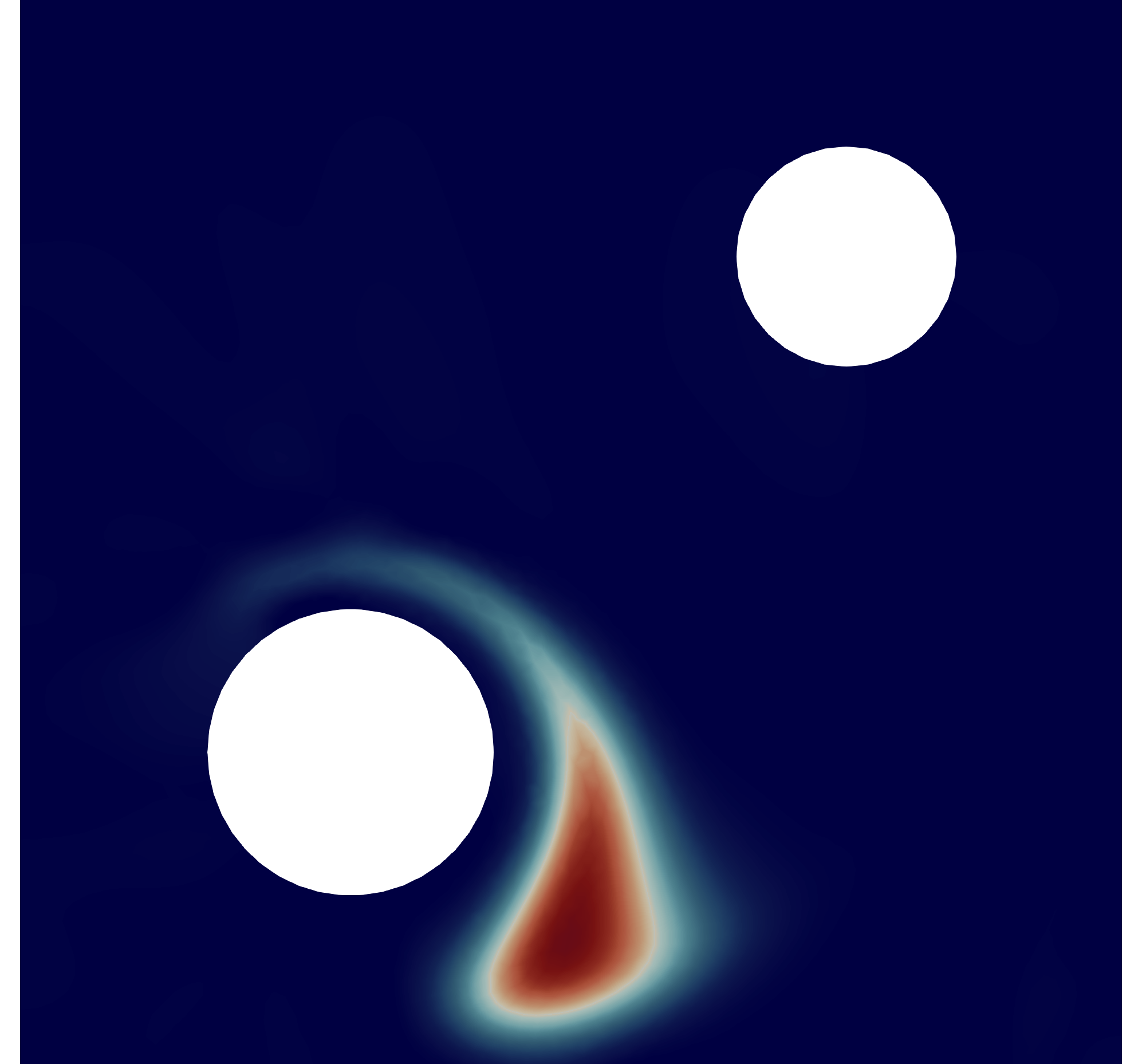}}
         {\includegraphics[width=0.22\textwidth]{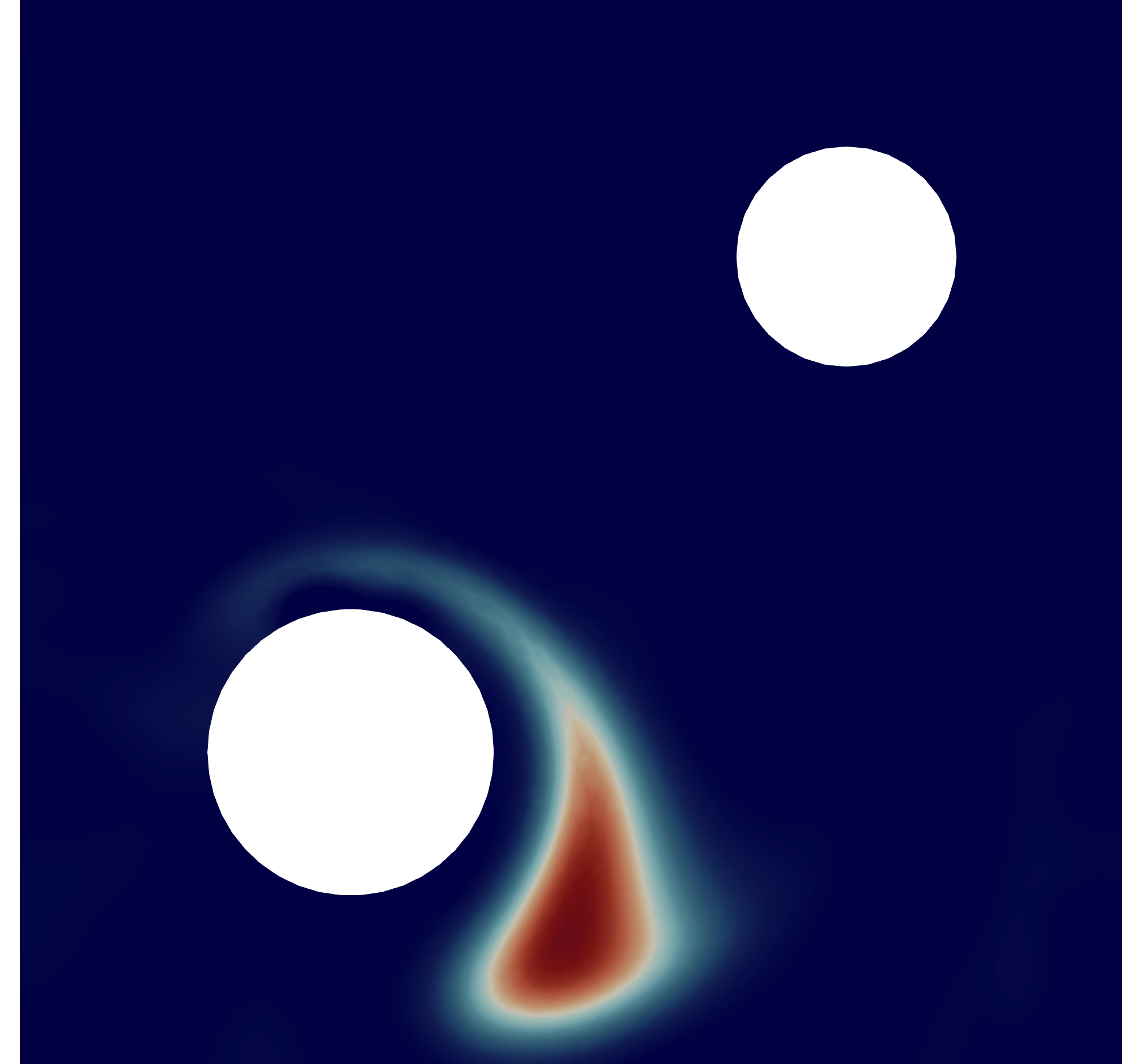}}
        {\includegraphics[width=0.22\textwidth]{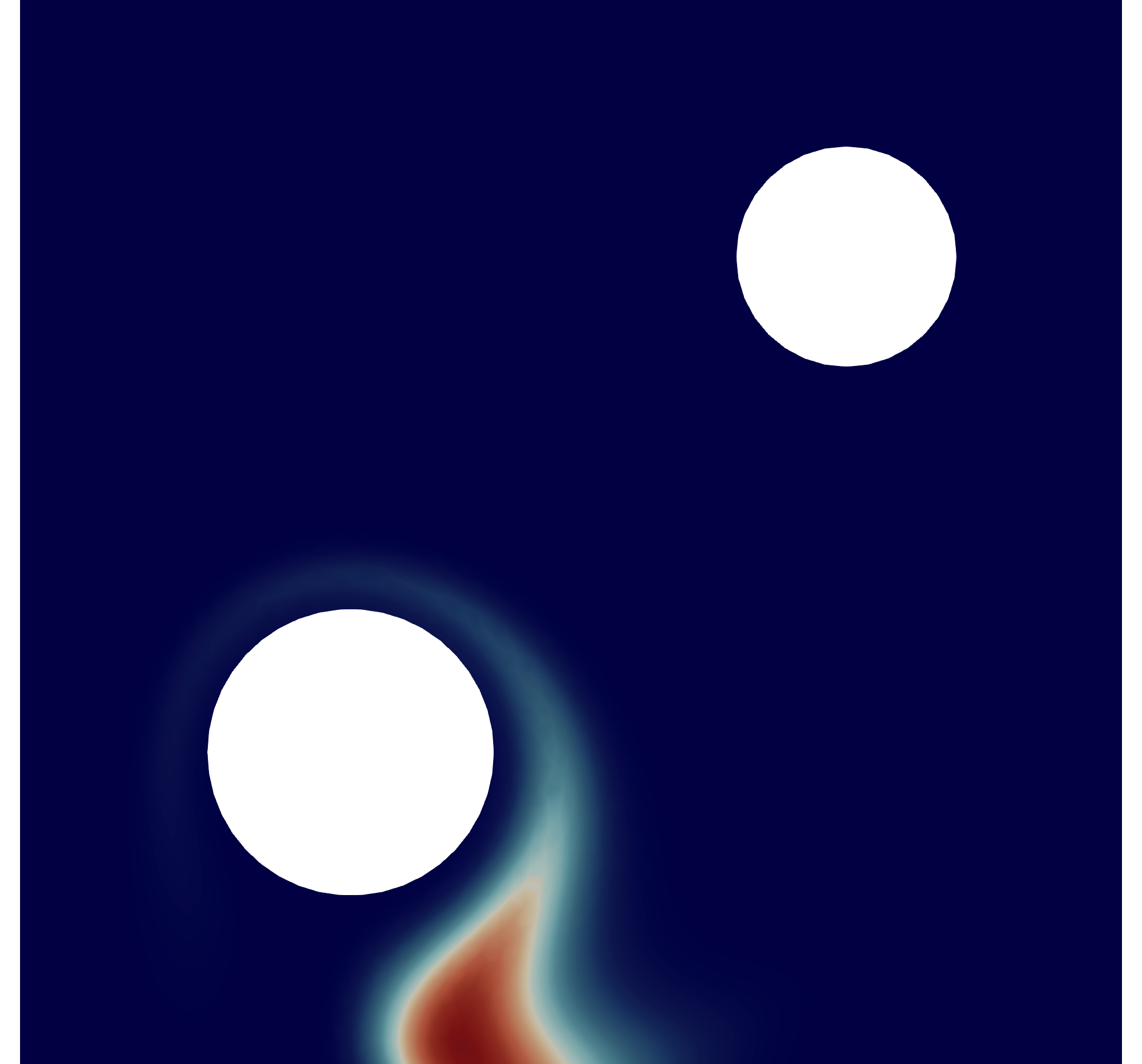}}
         {\includegraphics[width=0.22\textwidth]{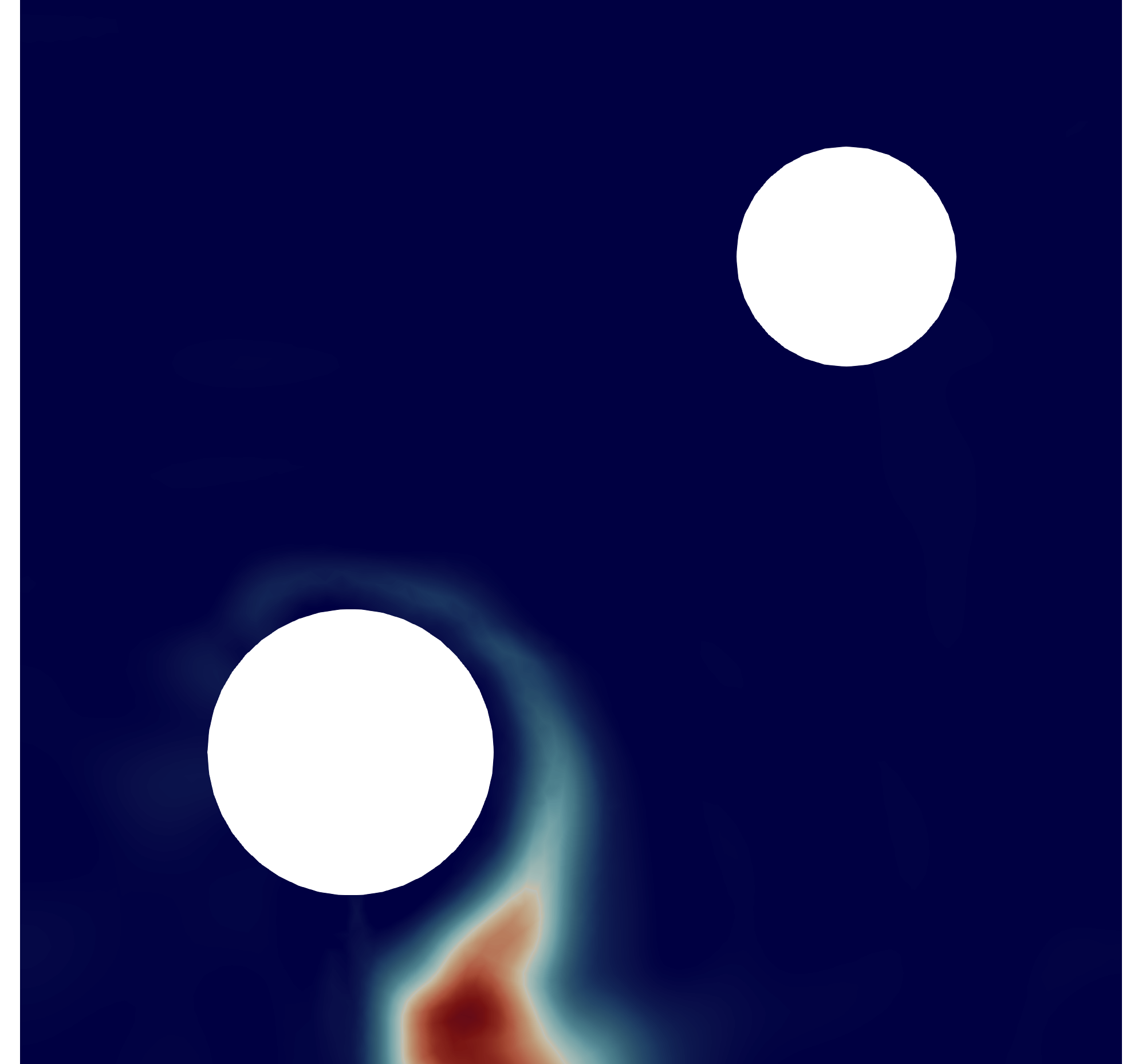}}
         {\includegraphics[width=0.22\textwidth]{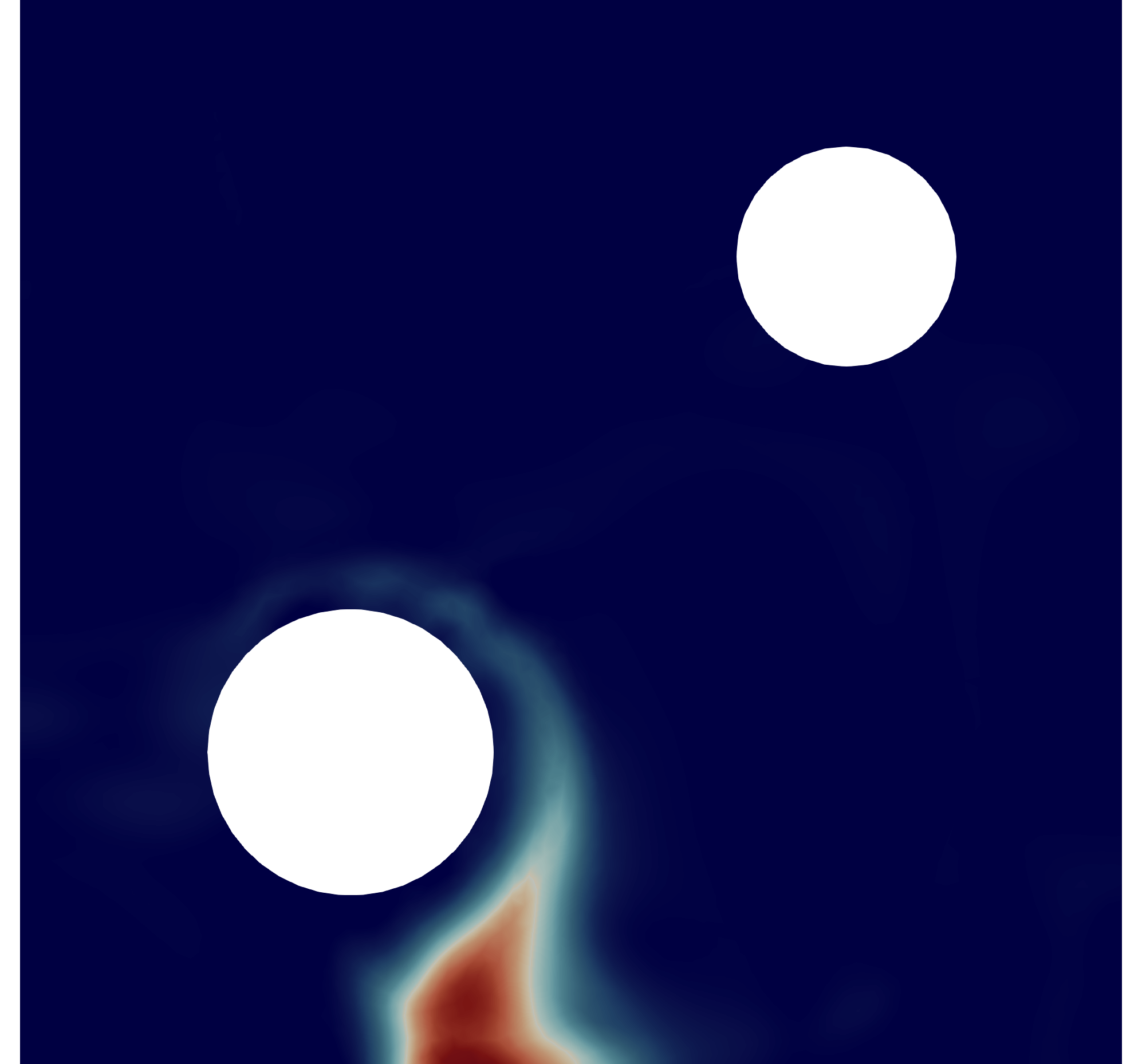}}
         {\includegraphics[width=0.22\textwidth]{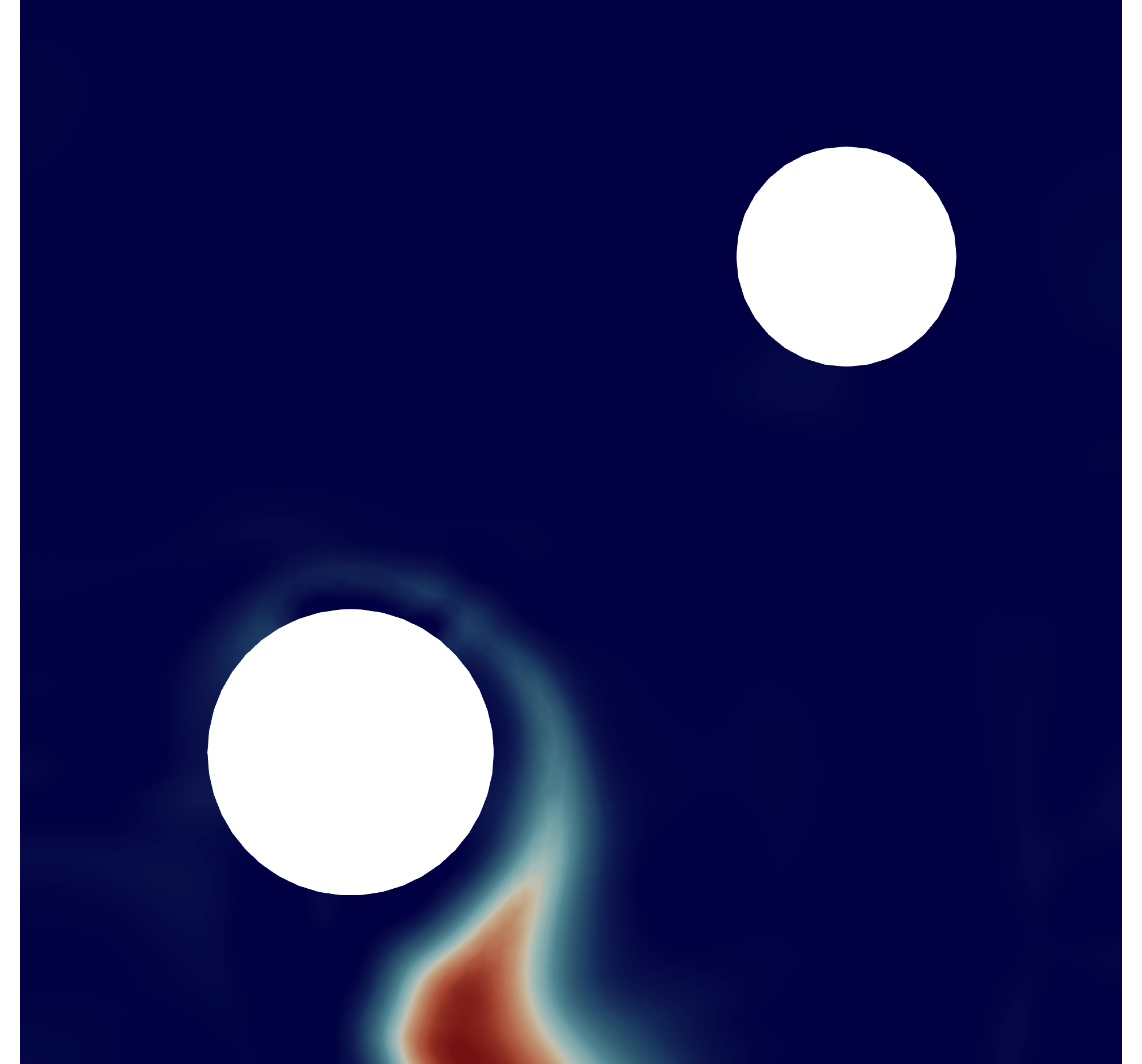}}
 \subcaptionbox{FE reference}
                 {\includegraphics[width=0.22\textwidth]{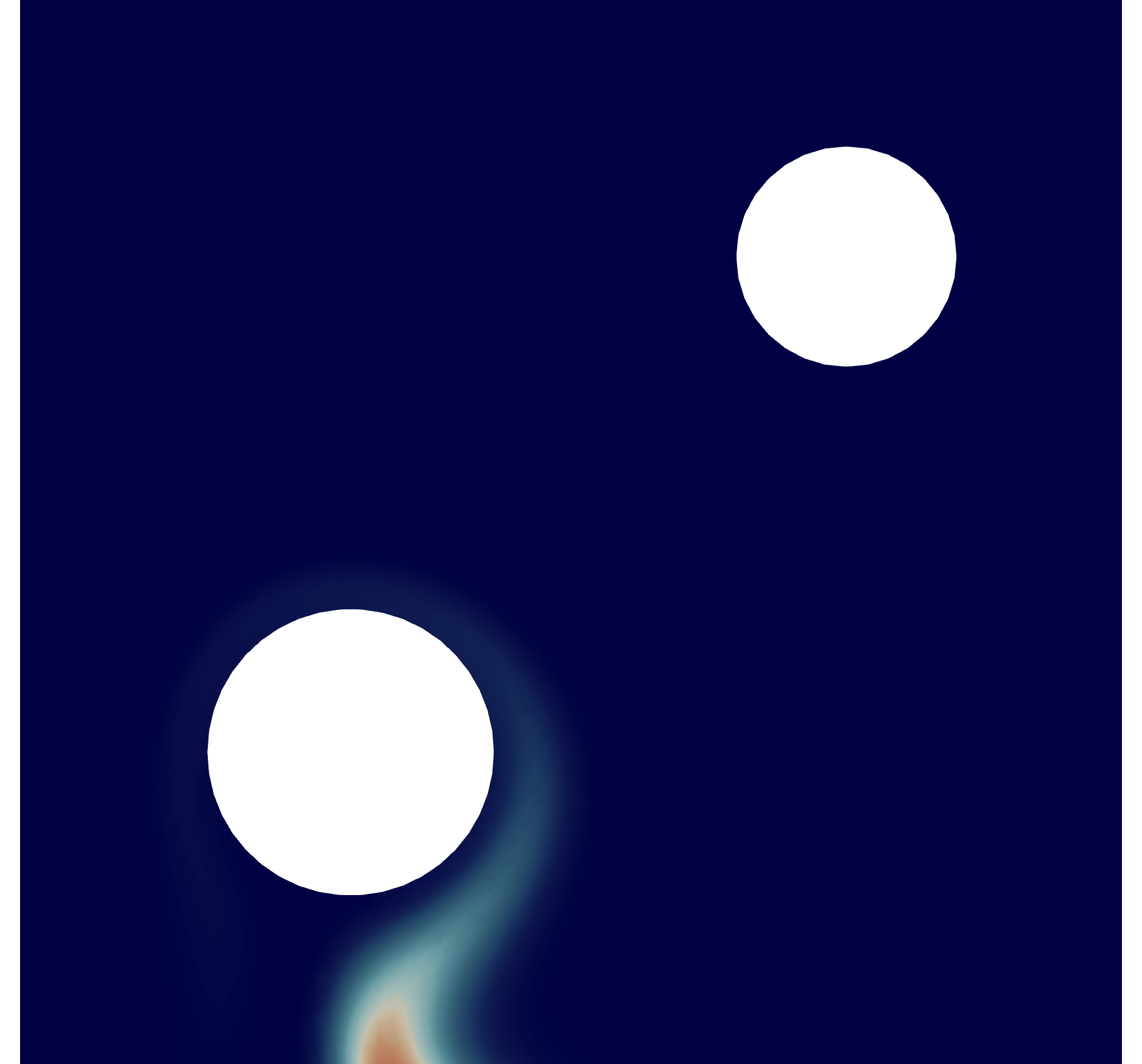}}
  \subcaptionbox{$n_\Phi=10$, full mesh}
         {\includegraphics[width=0.22\textwidth]{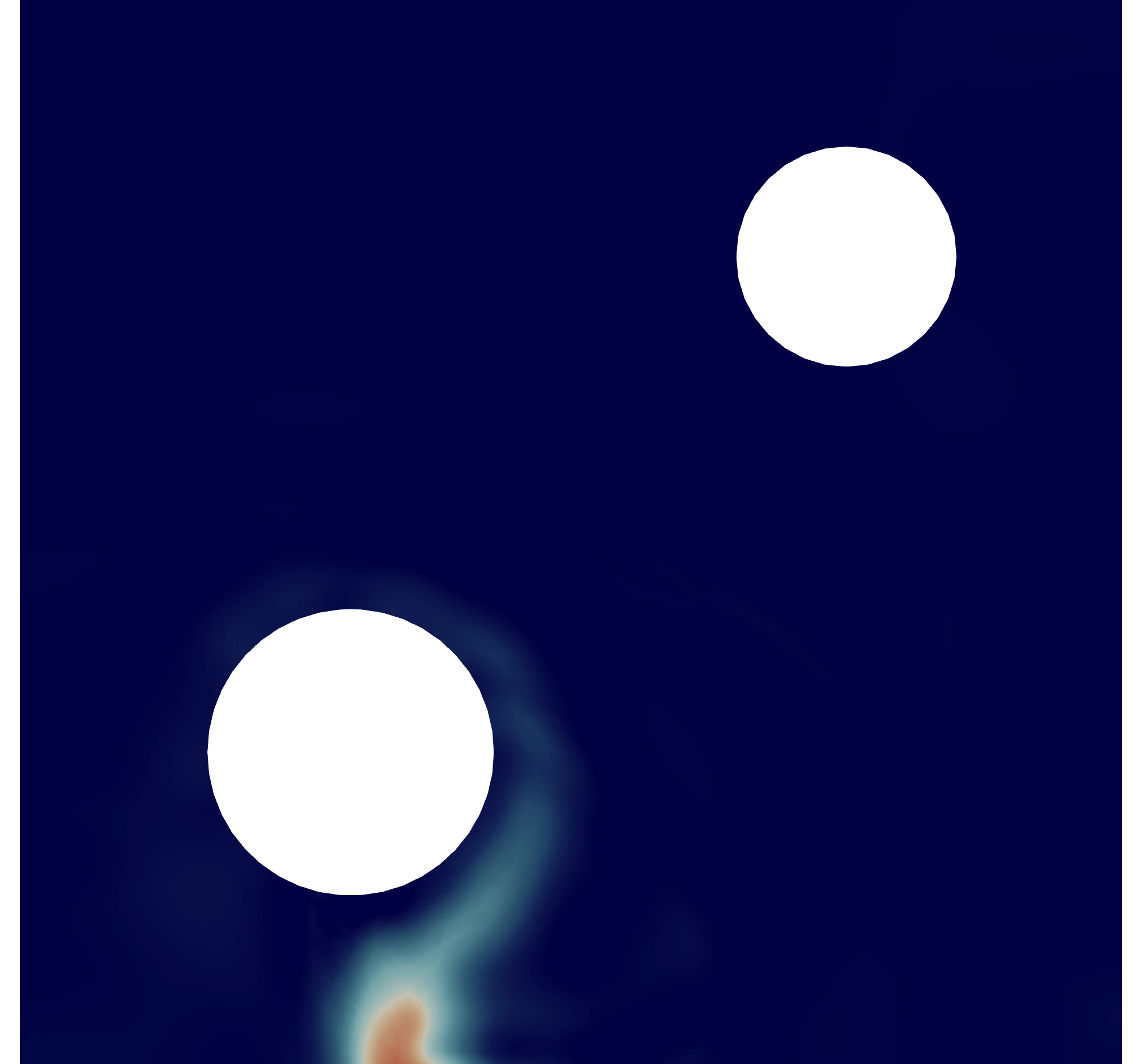}}
   \subcaptionbox{$n_\Phi=15$, full mesh }
         {\includegraphics[width=0.22\textwidth]{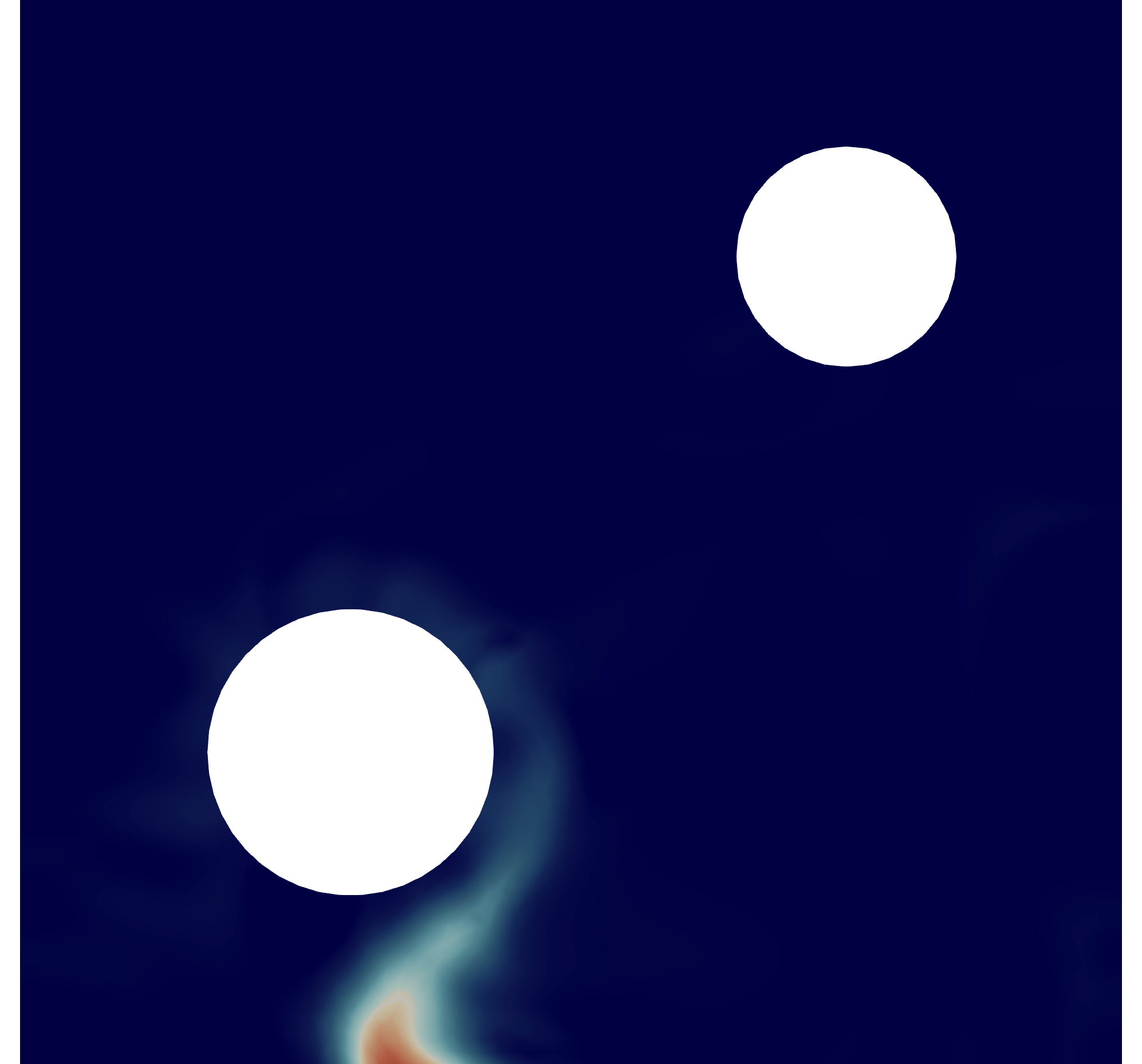}}
   \subcaptionbox{$n_\Phi=15$, $n_s$=1000 }
         {\includegraphics[width=0.22\textwidth]{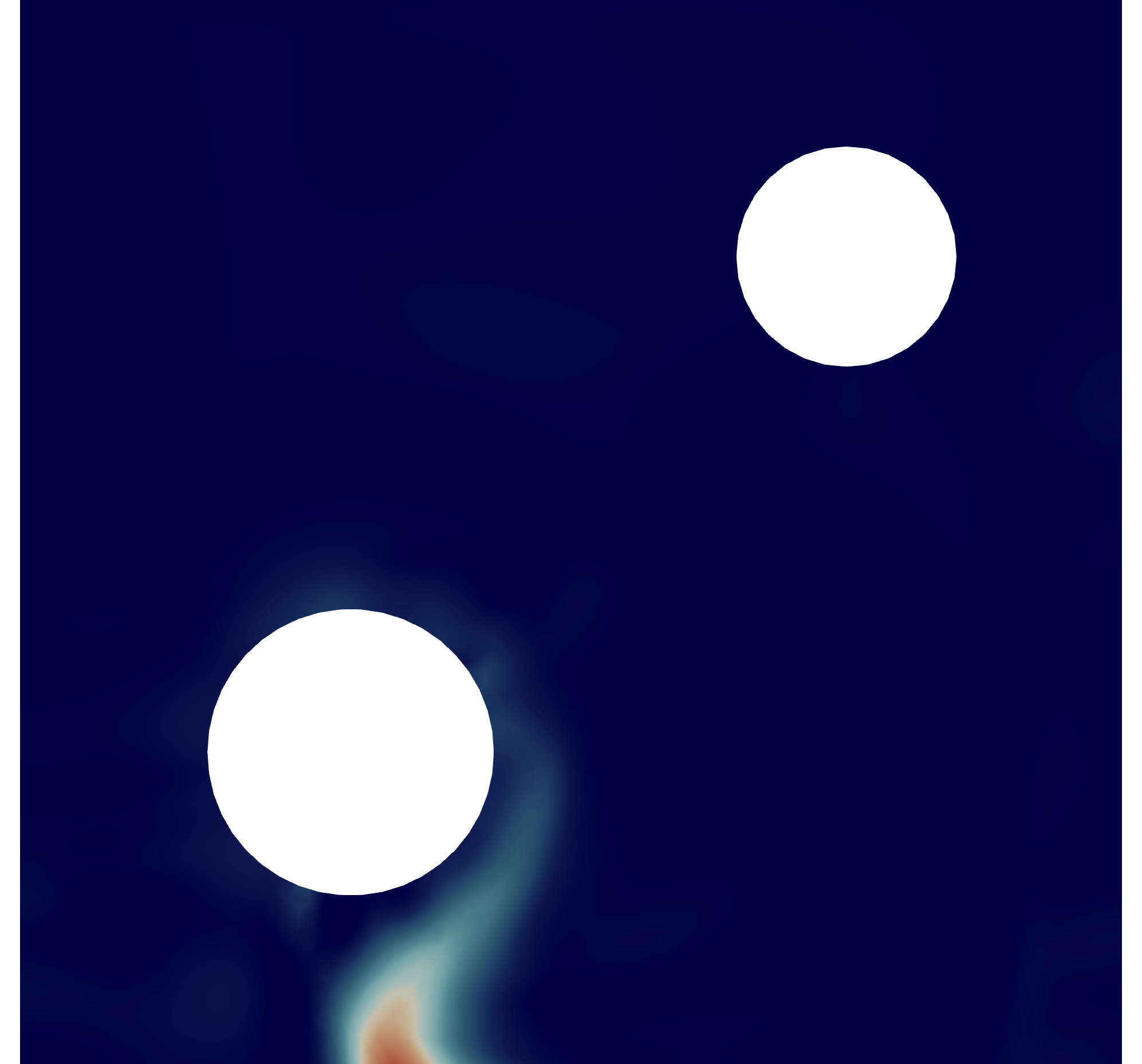}}
\caption{Solution at t= [0.1, 0.2, 0.3, 0.4, 0.5, 0.6] (from top to bottom) for the reference solution and two different archtiectures, as well as the subsampled case, for the advection-diffusion problem on a domain with holes. The neural networks have four layers with 10 hidden units.}
\label{fig:snapshots_adv_diff_holes}
\end{figure}

All neural networks used in this experiment have 4 layers with 10 or 20 hidden units.
We train for the initial condition for $\vec{\alpha}= [0.0,0.0]$ on the $n_x= 6205$ nodes of the mesh for 40000 iterations.
We use an Euler scheme with a time step size of $\Delta t= 10^{-3}$ to run the simulation until time t=1,  but we only visualize the solution until t=0.6, see Figure \ref{fig:snapshots_adv_diff_holes}, when most of the concentration has left the domain. \\
In the error plots in Figure \ref{fig:error_adv_diff_holes} (a), we observe that using only $n_\Phi=4$ features leads to errors larger than 10\%, while using  $n_\Phi=10$ or $n_\Phi=15$ features leads to much better approximations. Increasing the number of hidden layers  only slightly reduces the error compared to increasing the number of features in the positional embedding. This is a striking difference to the cases on simpler domains, where the eigenfunctions were simple trigonometric functions. \\
Taking a closer look at the solution snapshots in Figure \ref{fig:snapshots_adv_diff_holes}, we observe that there is hardly any visual difference to the reference solution from t=0.1 to t=0.4, while we lose both the shape of the thin tail and the exact shape of the solution on the Neumann boundary in the last two snapshots. 

 The figures also show that large parts of the computational domain contain no information about the problem. We therefore  aim to increase computational efficiency by sub-sampling the points on the mesh according to the size of the advective term as described in Algorithm \ref{alg:pseudoalgactive}. As can be seen from both the error plot in Figure \ref{fig:error_adv_diff_holes} and the solution (Figure \ref{fig:snapshots_adv_diff_holes}), the quality of the solution does not significantly degrade when only 1000 points are sampled. On the other hand, 500 points are not enough to recover the same quality of the solution, as errors are considerably higher. We remark that for the larger neural networks, the update equation is underdetermined as the number of network parameters exceeds the number of sampling points. The Krylov solver implementation can deal with this scenario without further modification. For 2000 and 1000 samples, errors are in the same range as for the full mesh. With these encouraging results we turn to the parametrized problem, where we cannot afford to consider the full mesh for every parameter value.

\subsubsection{Parametrized initial condition}
We use 4 layer neural networks with 10 or 20 hidden units and  10, 15 and 20 features in the positional embedding to account for the increased complexity of the solution manifold in the parametrized case. 
We consider the parametrized initial condition and train the neural network on 40000 points across the combined physical and parameter domain with a batch size of 1024 for 8000 epochs. We again use forward Euler with a time step size of $\Delta t= 10^{-3}$ and compute the parametrized solutions until t=0.6.

\begin{figure}
 \centering

         {\includegraphics[width=0.19\textwidth]{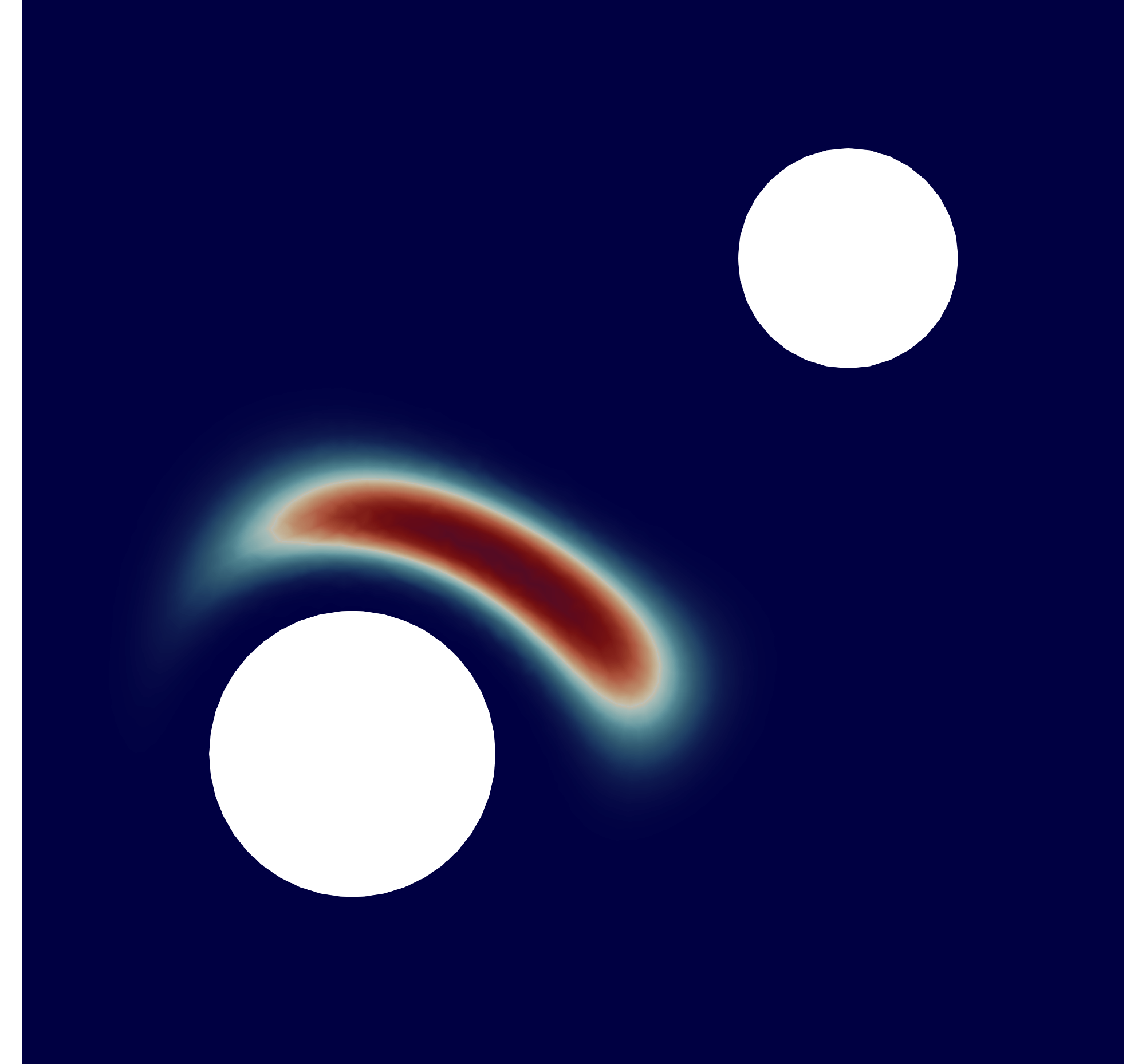}}
  \subcaptionbox{t=0.2, FE reference}
         {\includegraphics[width=0.19\textwidth]{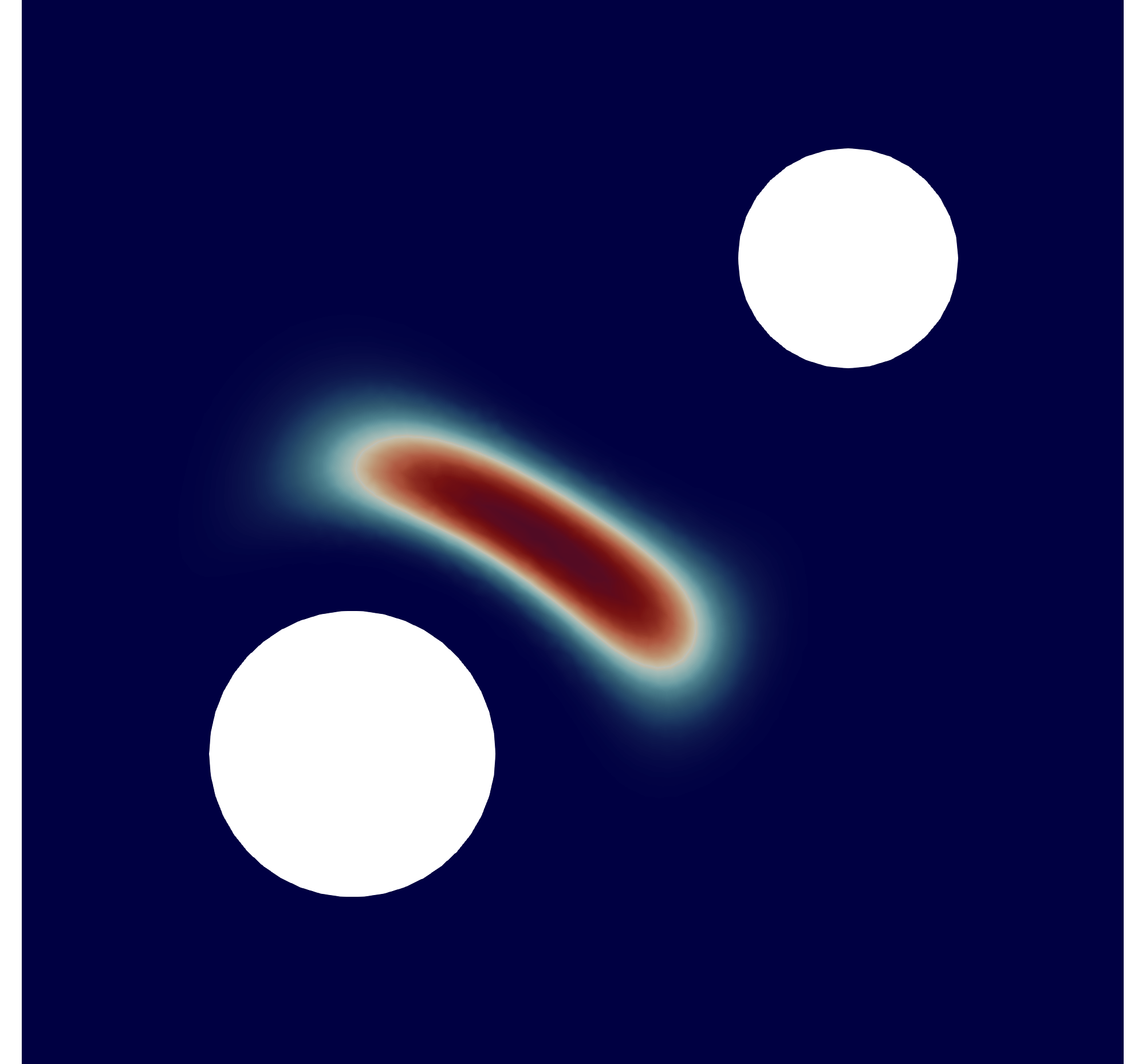}}
         {\includegraphics[width=0.19\textwidth]{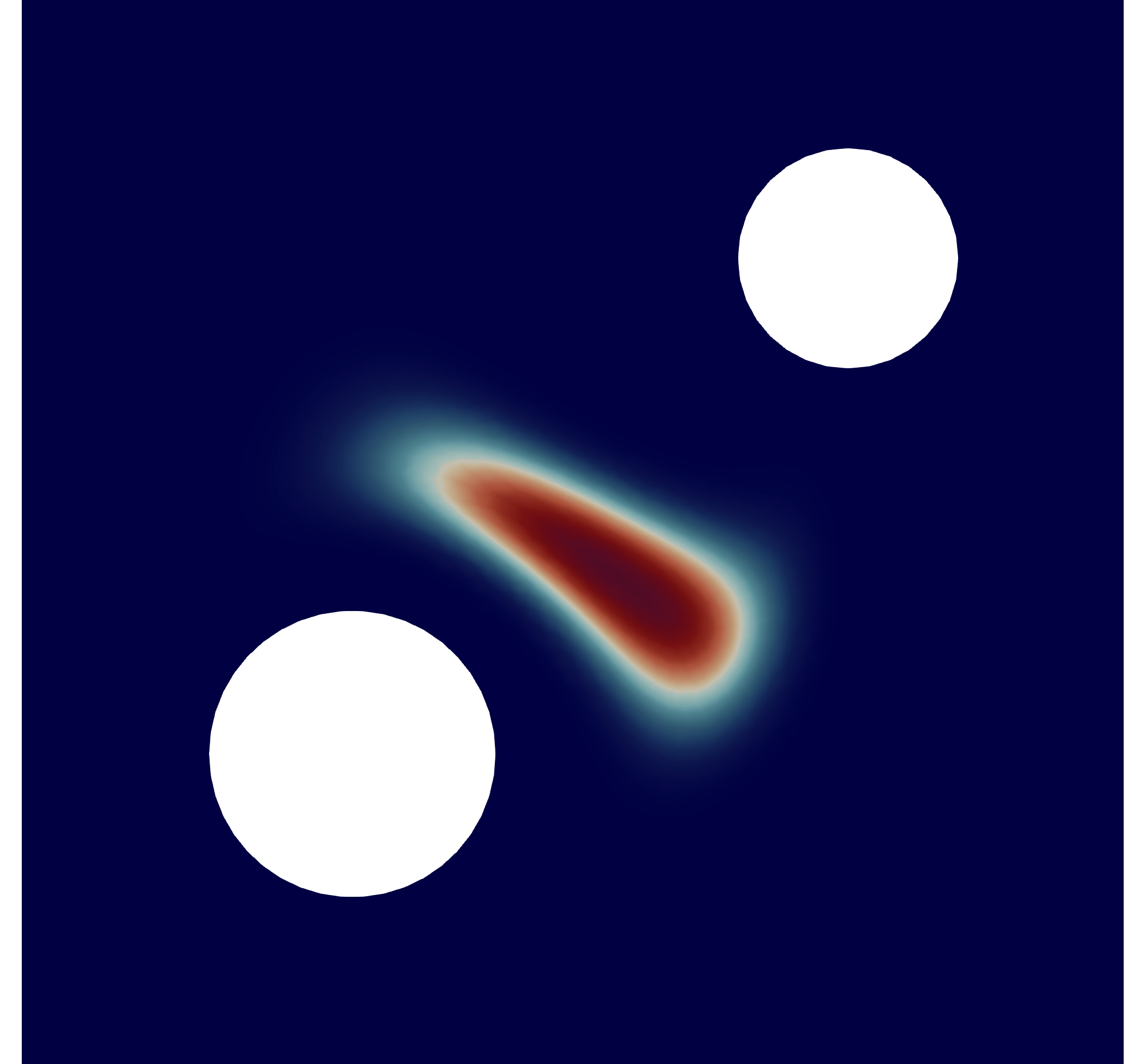}}\\
         {\includegraphics[width=0.19\textwidth]{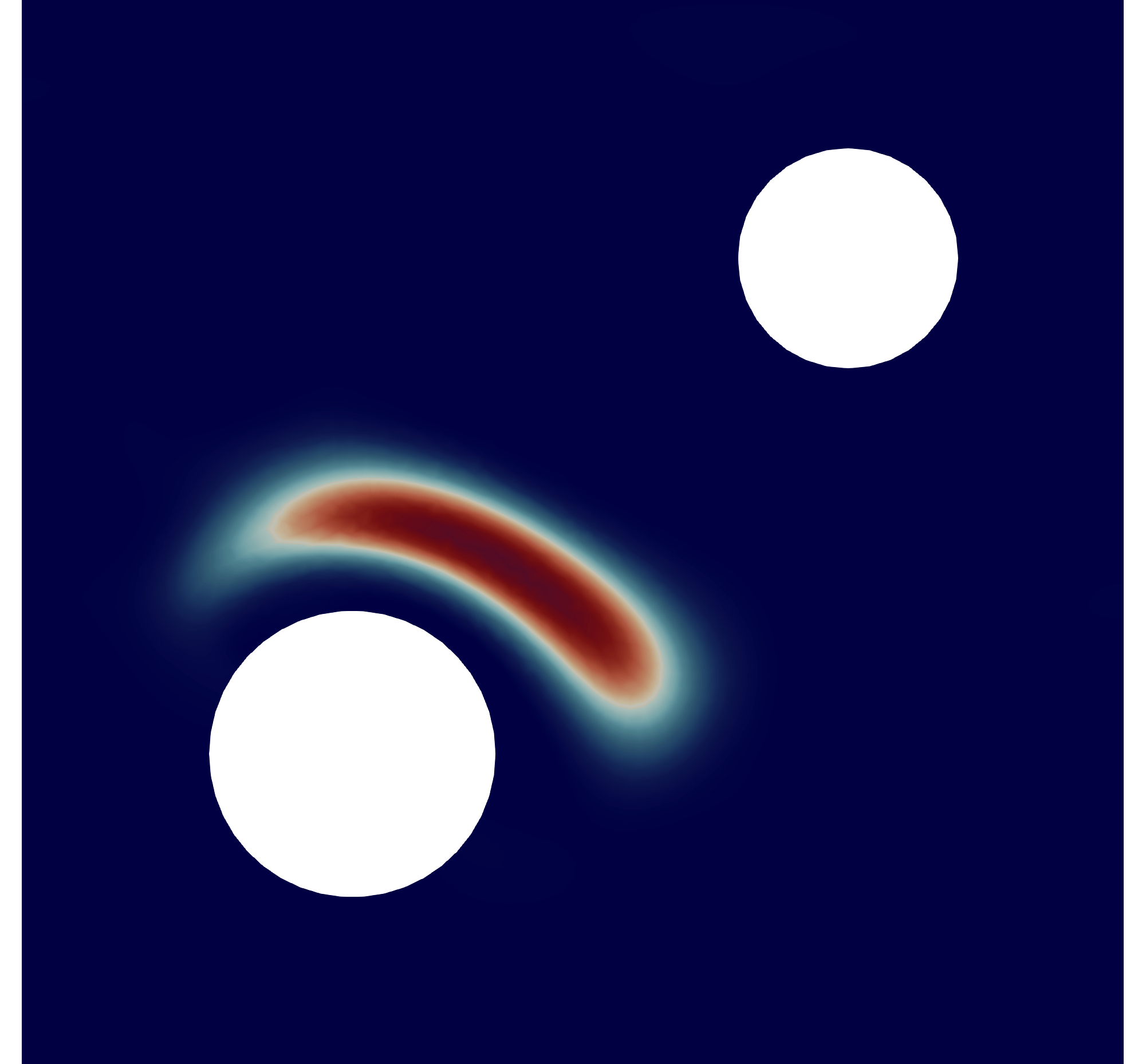}}
 \subcaptionbox{t=0.2, NN solution}
         {\includegraphics[width=0.19\textwidth]{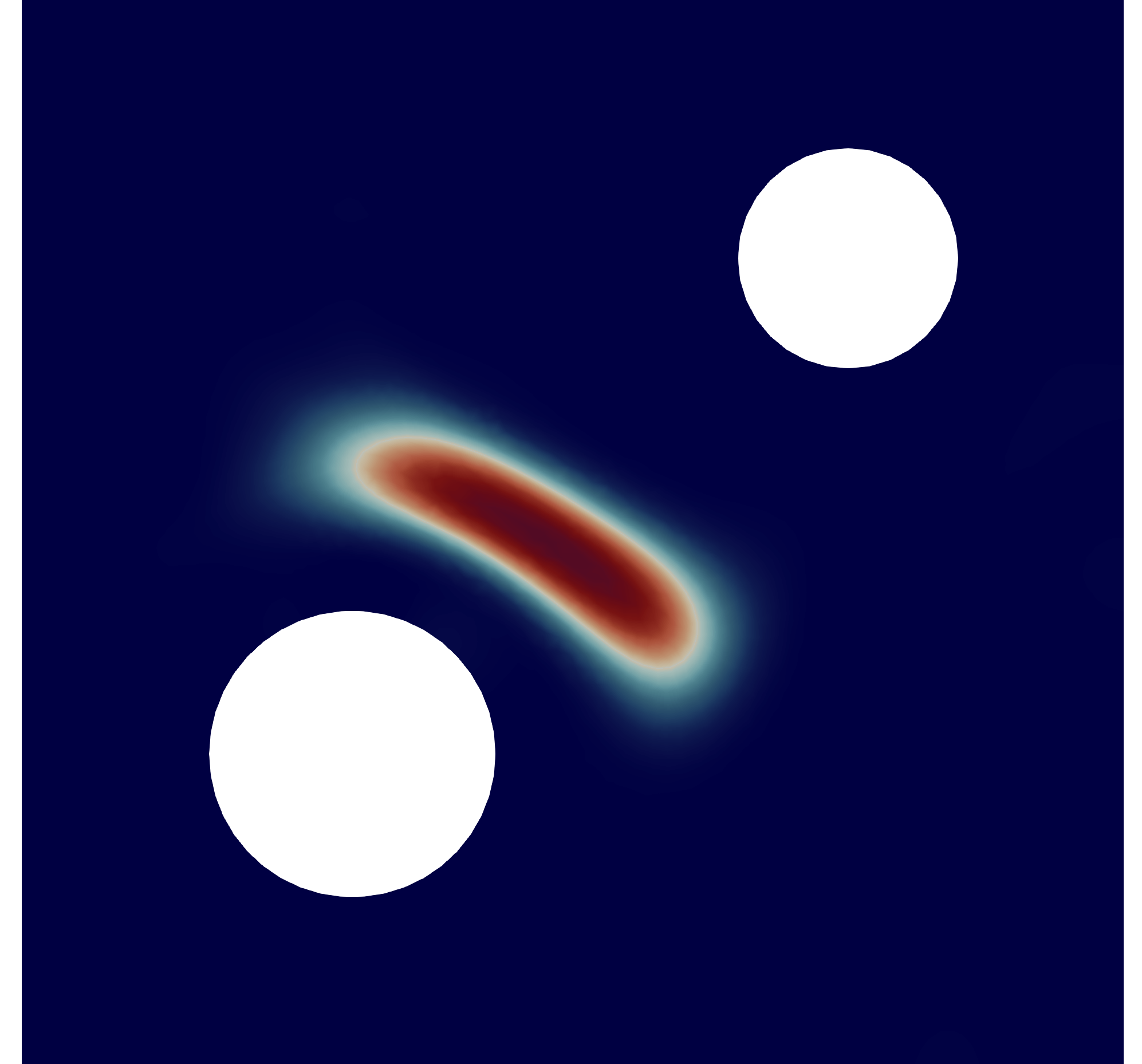}}
         {\includegraphics[width=0.19\textwidth]{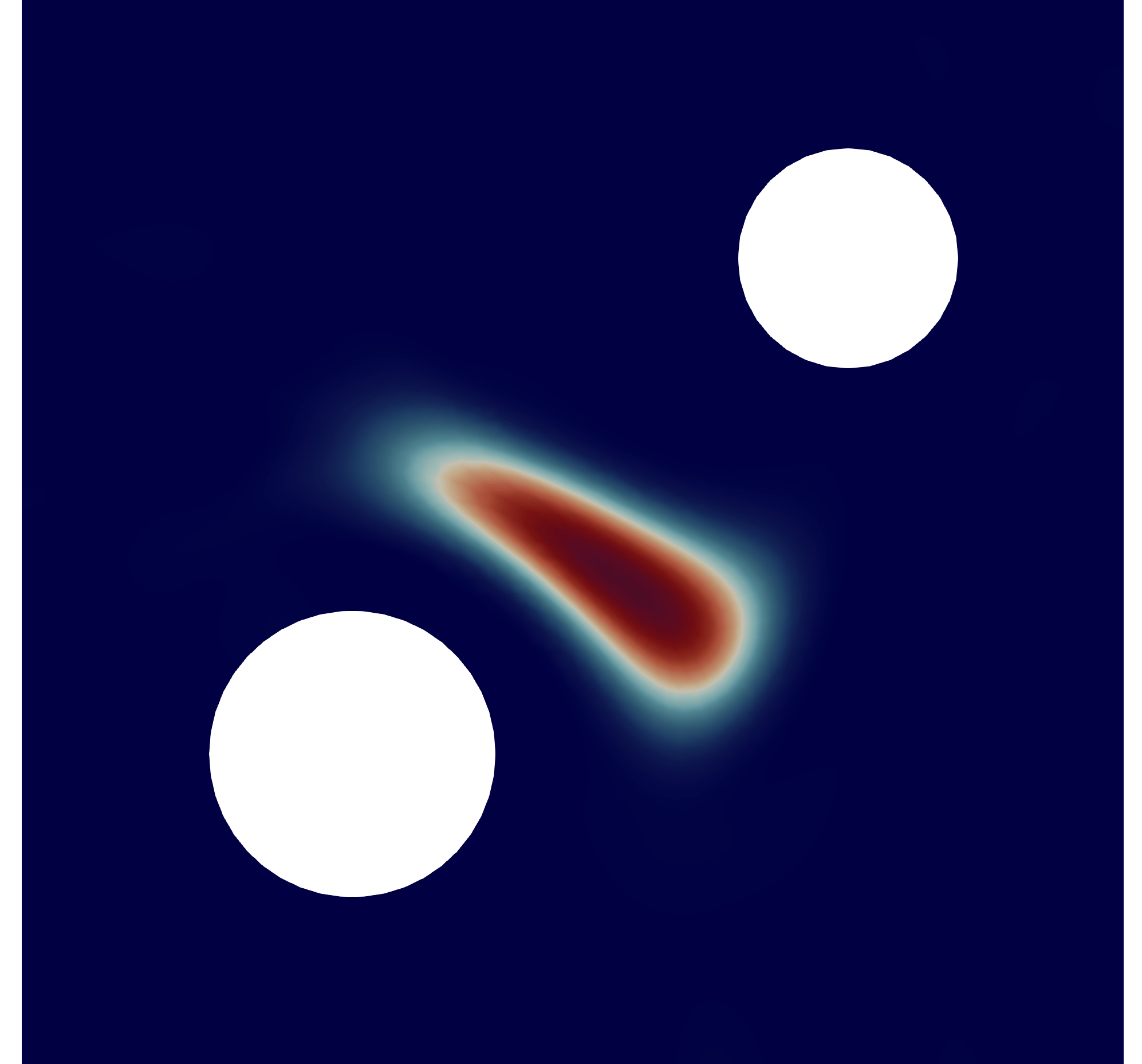}}\\

         {\includegraphics[width=0.19\textwidth]{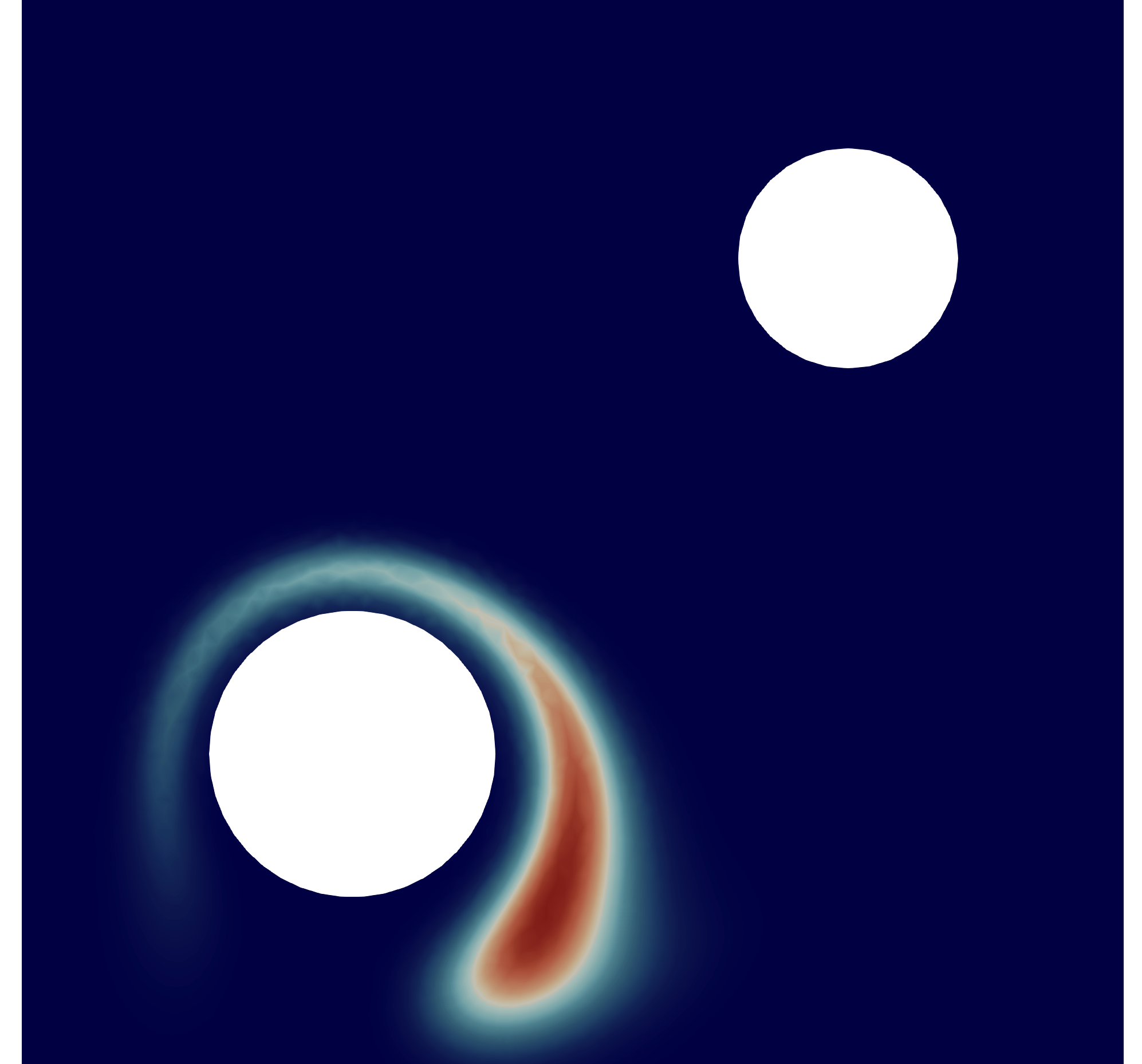}}
  \subcaptionbox{t=0.4, FE reference}
         {\includegraphics[width=0.19\textwidth]{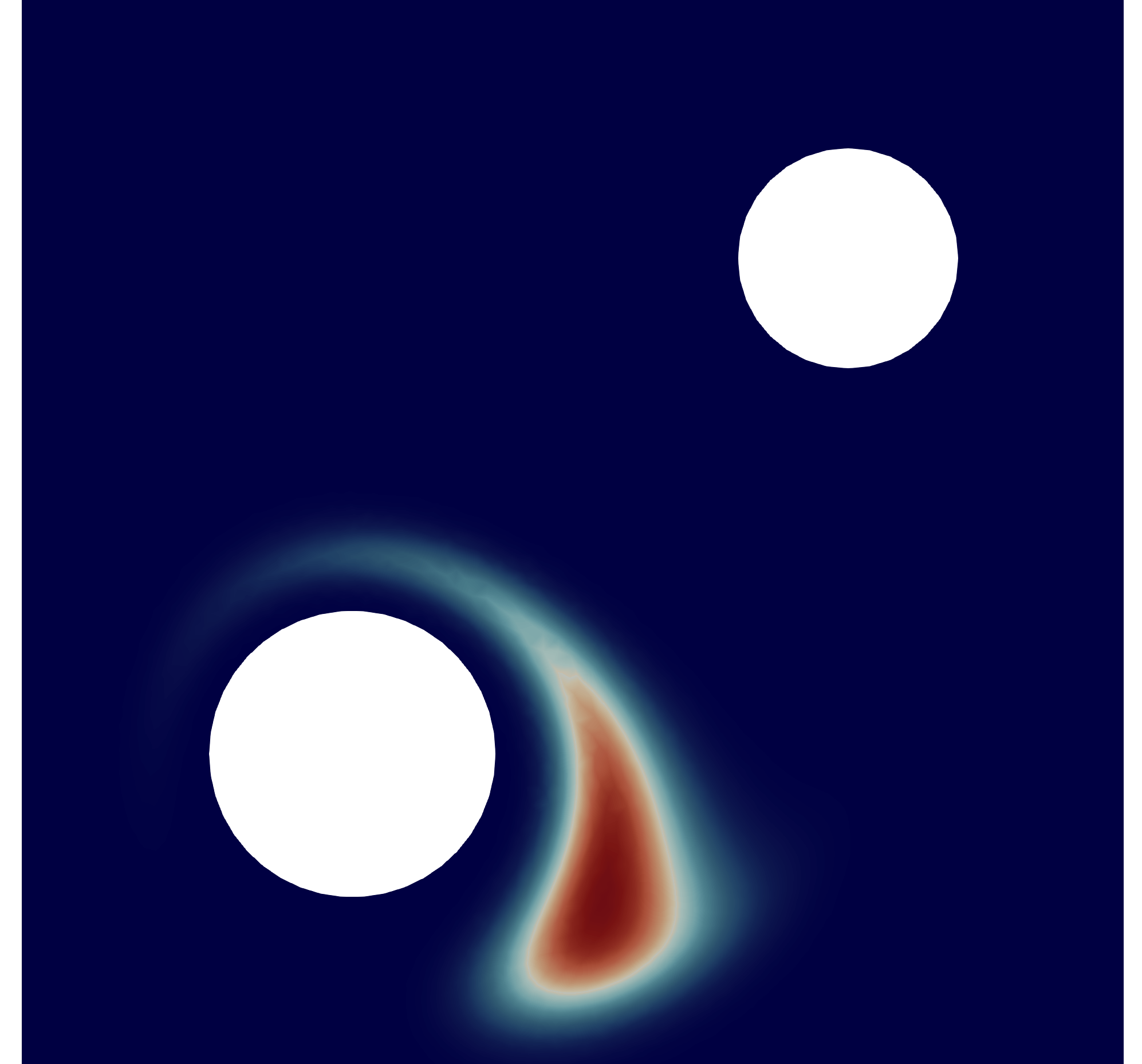}}
         {\includegraphics[width=0.19\textwidth]{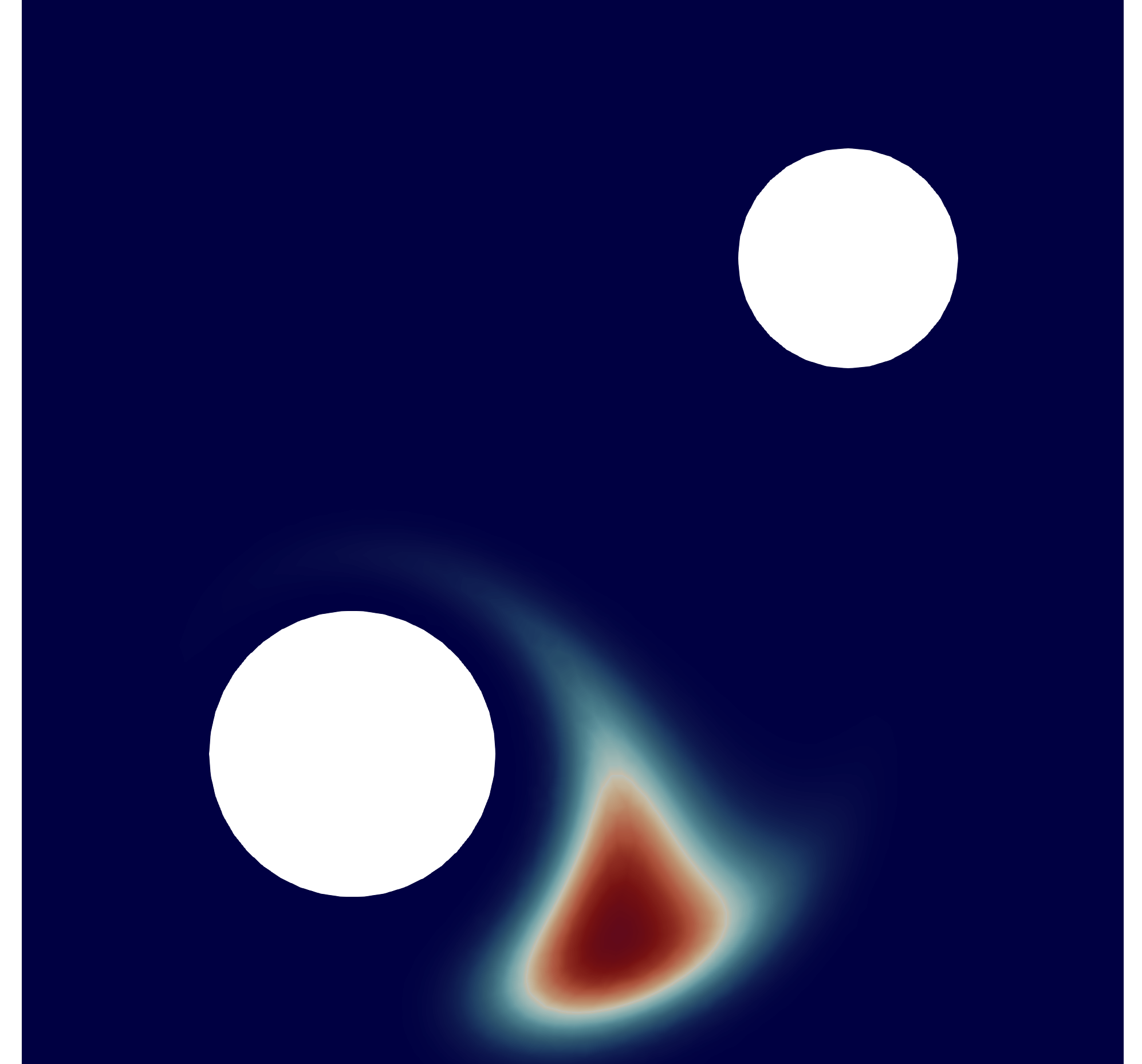}}\\
         {\includegraphics[width=0.19\textwidth]{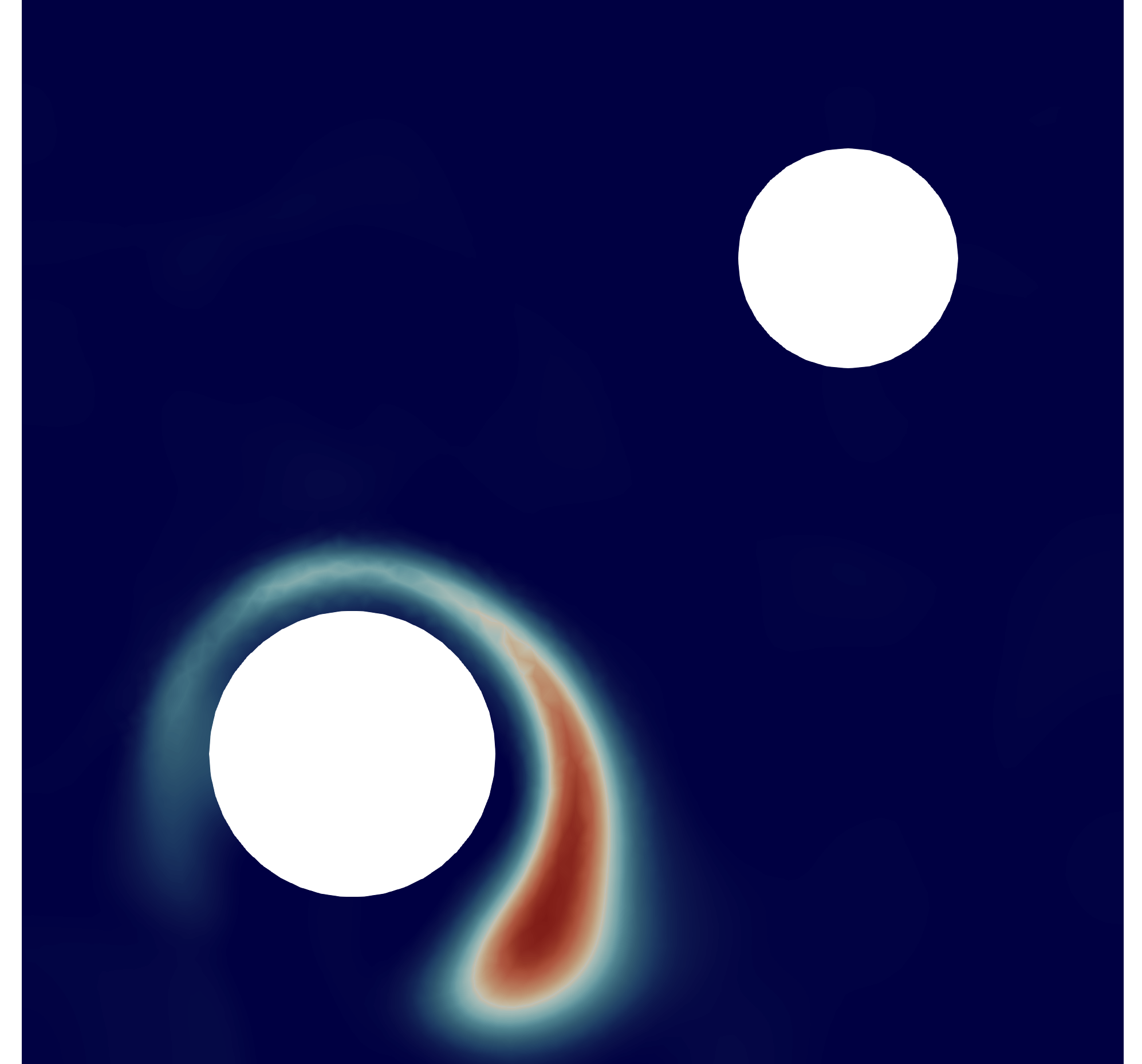}}
  \subcaptionbox{t=0.4, NN solution}
         {\includegraphics[width=0.19\textwidth]{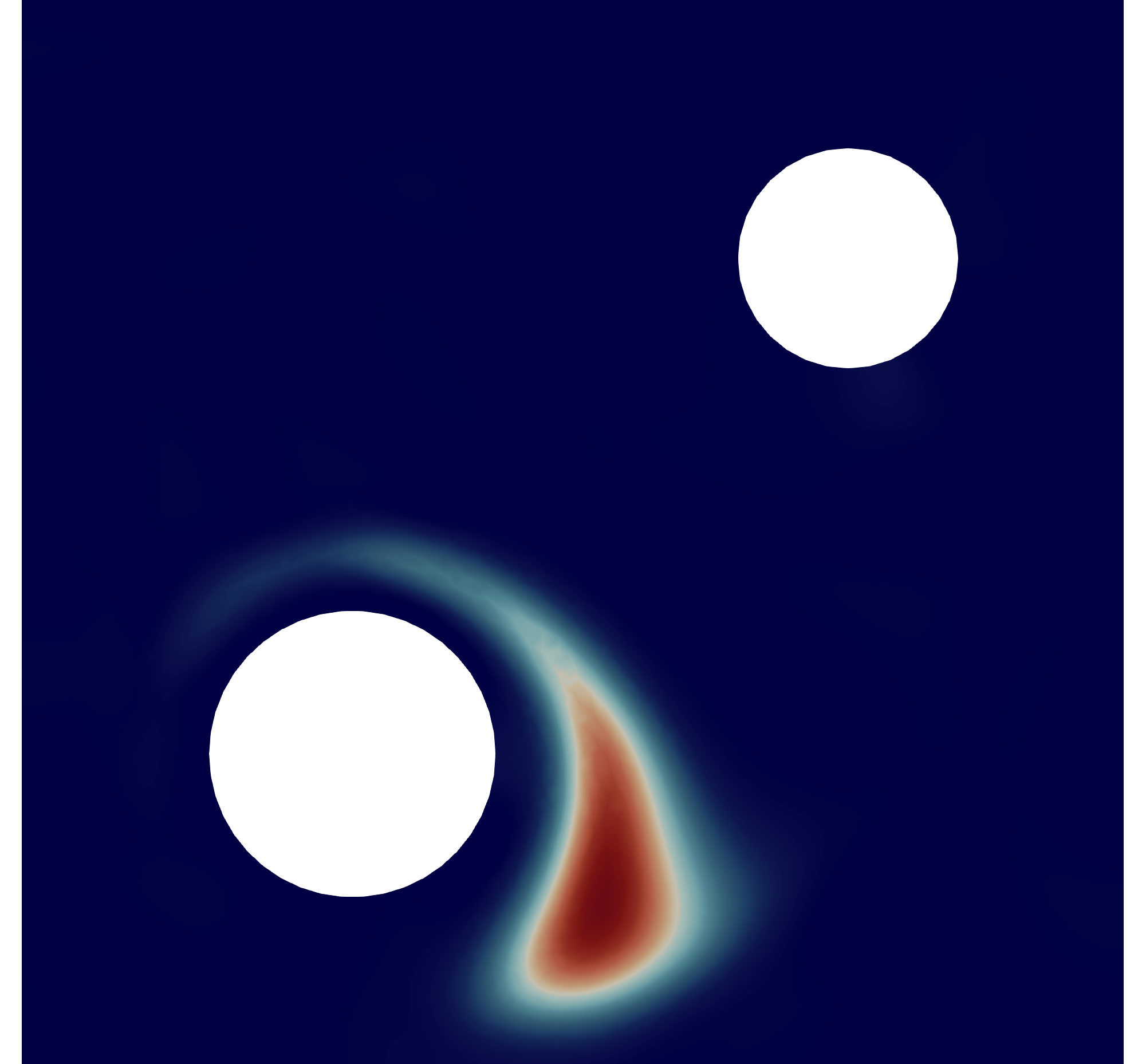}}
         {\includegraphics[width=0.19\textwidth]{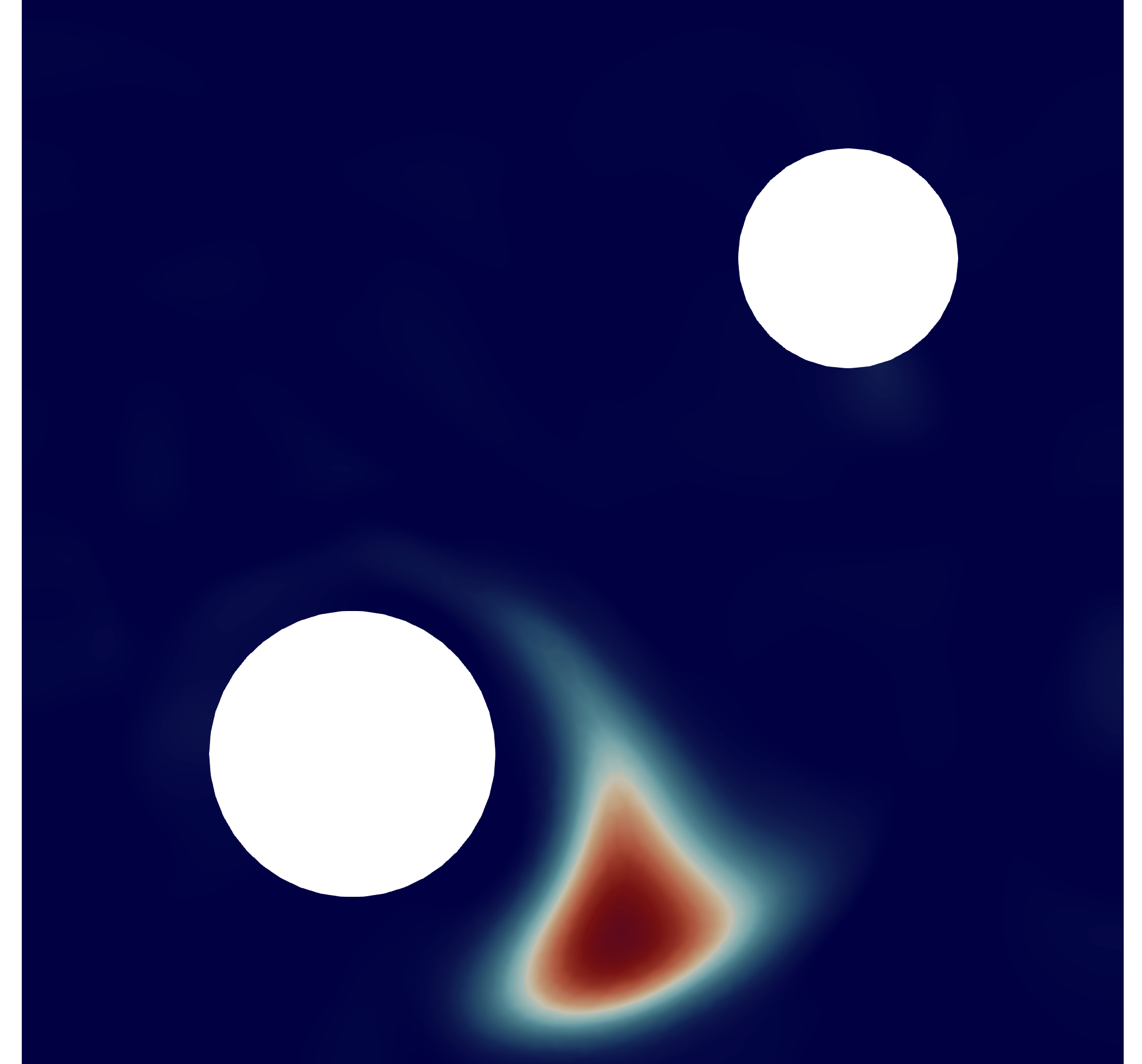}} \\
         {\includegraphics[width=0.19\textwidth]{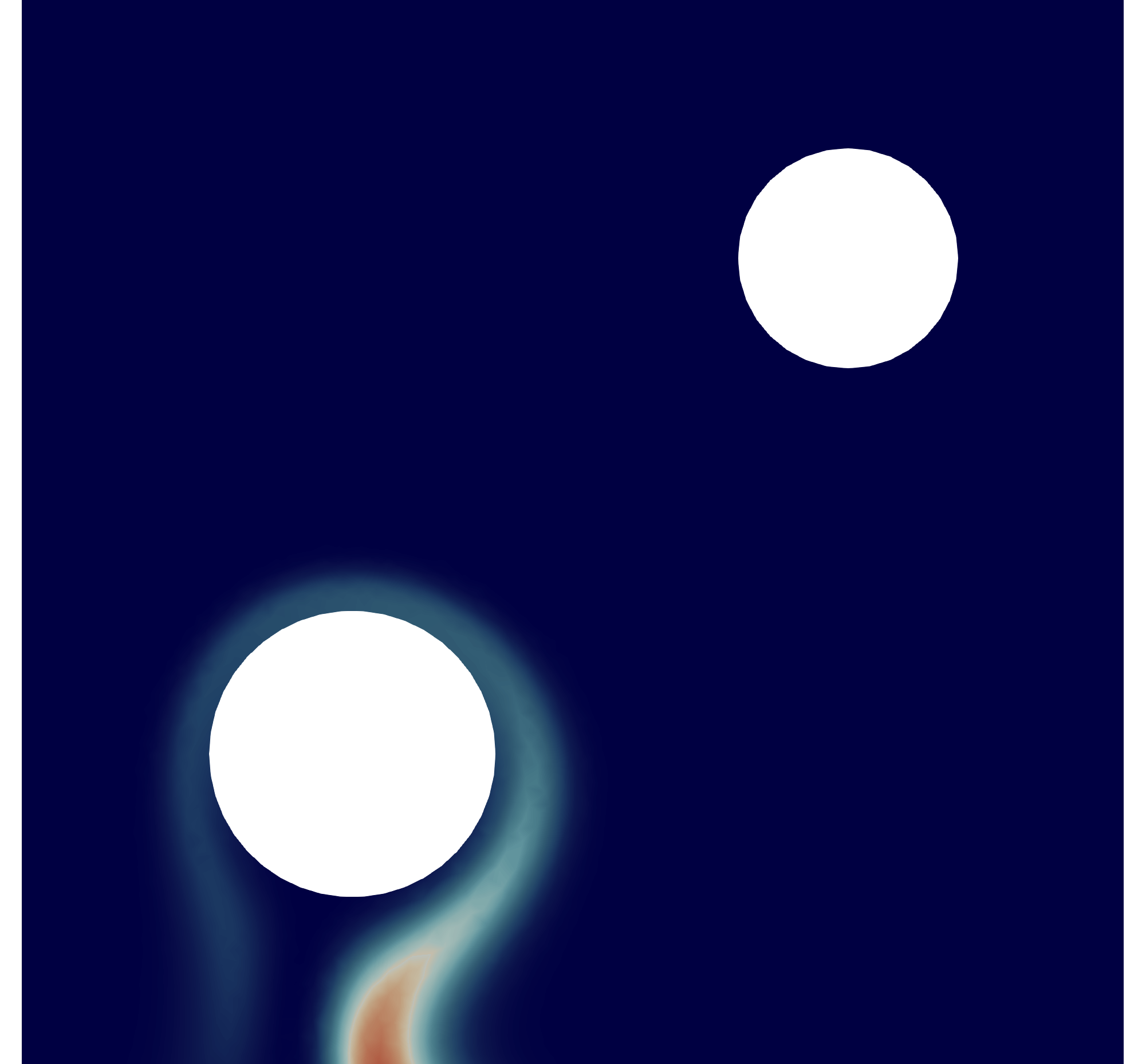}}
  \subcaptionbox{t=0.6, FE reference}
         {\includegraphics[width=0.19\textwidth]{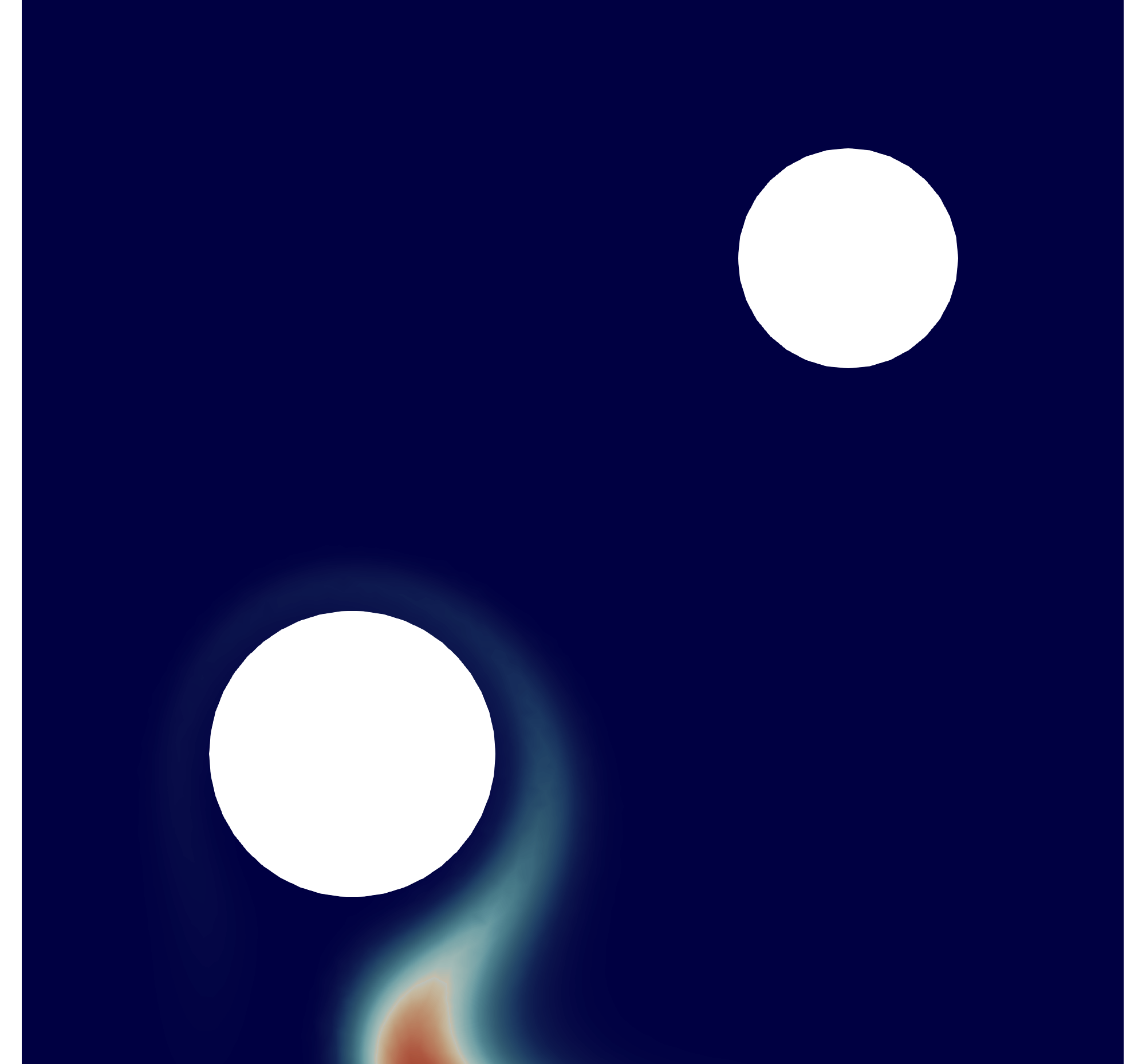}}
         {\includegraphics[width=0.19\textwidth]{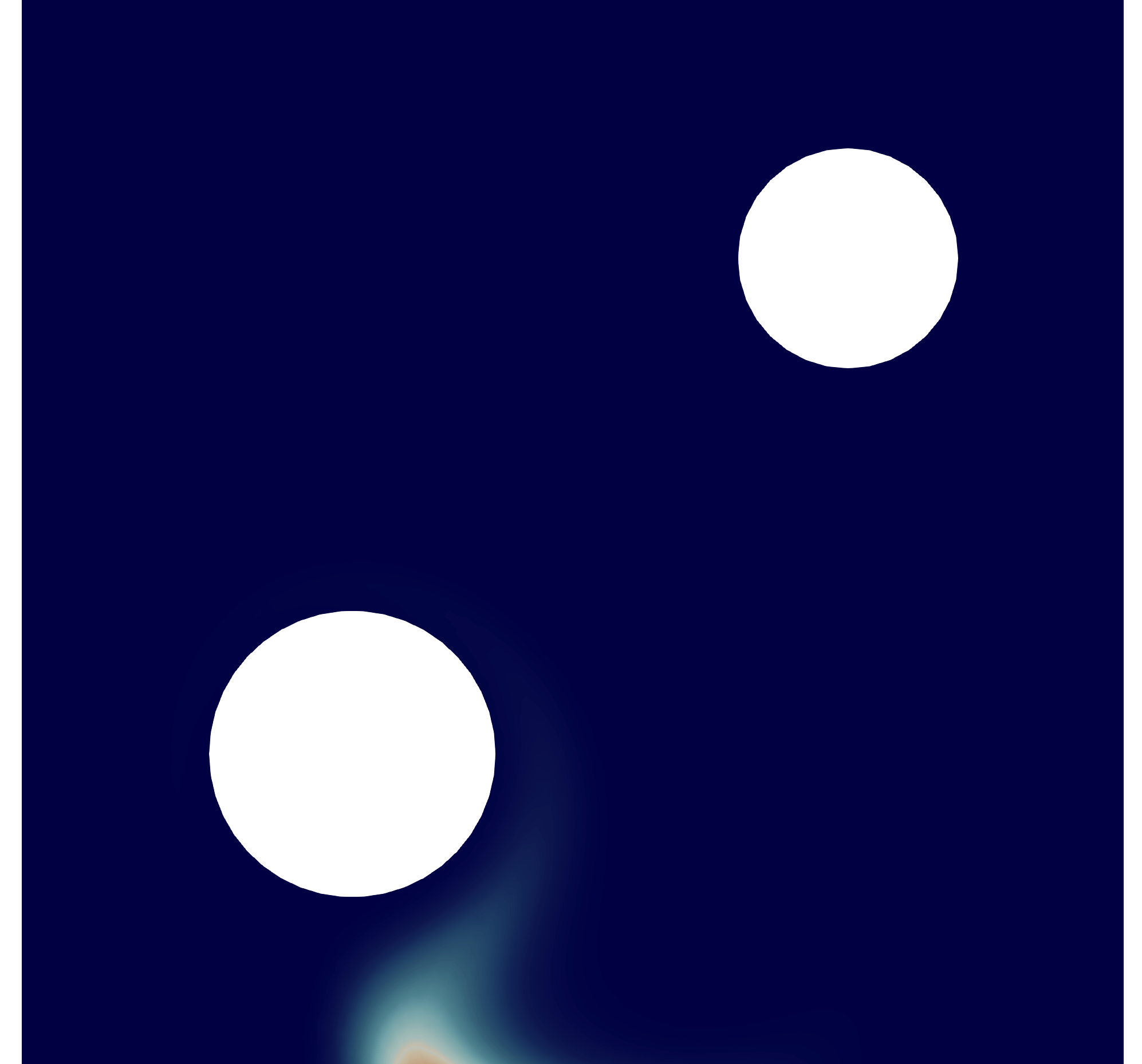}}\\
         {\includegraphics[width=0.19\textwidth]{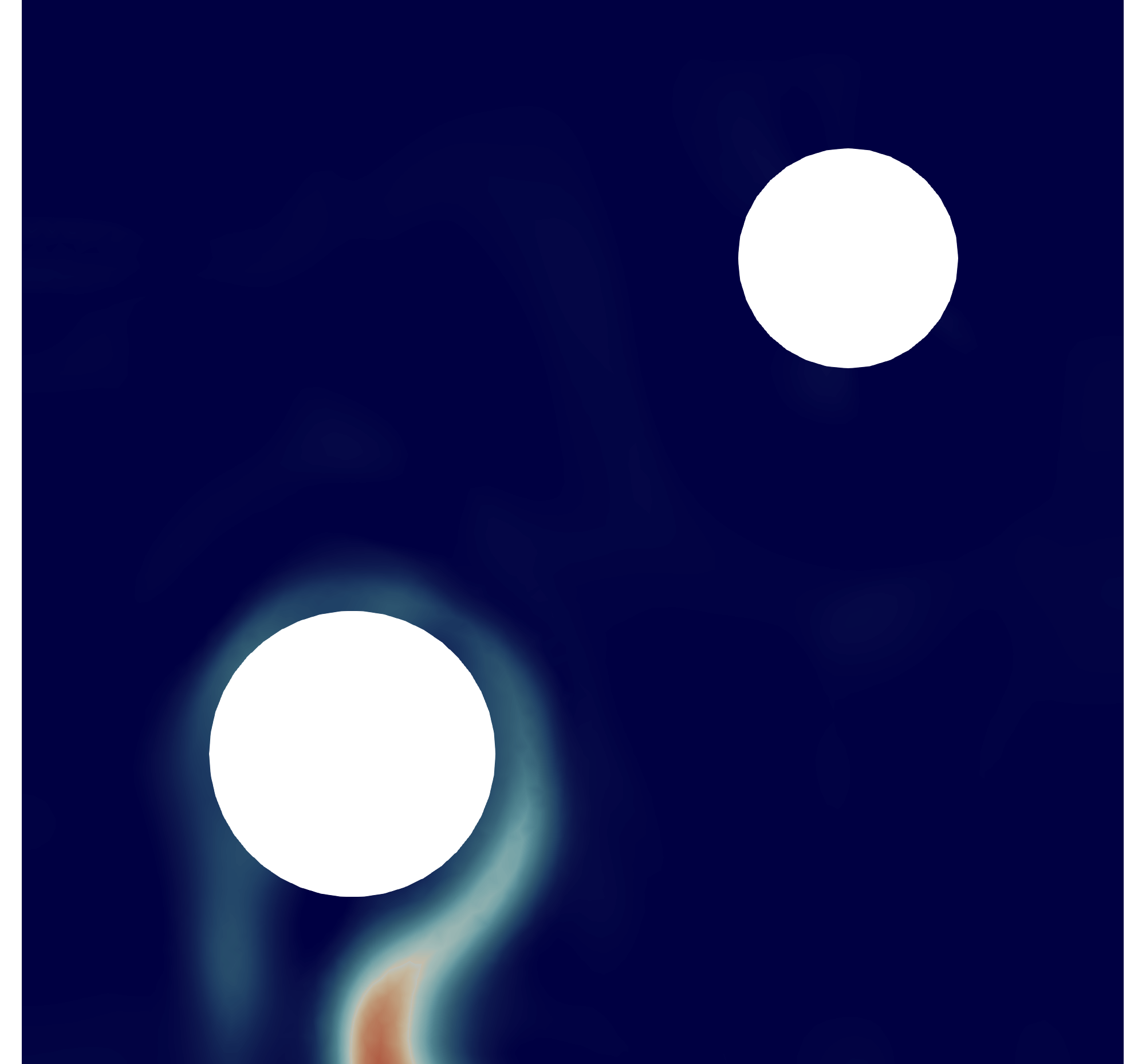}}
  \subcaptionbox{t=0.6, NN solution}
         {\includegraphics[width=0.19\textwidth]{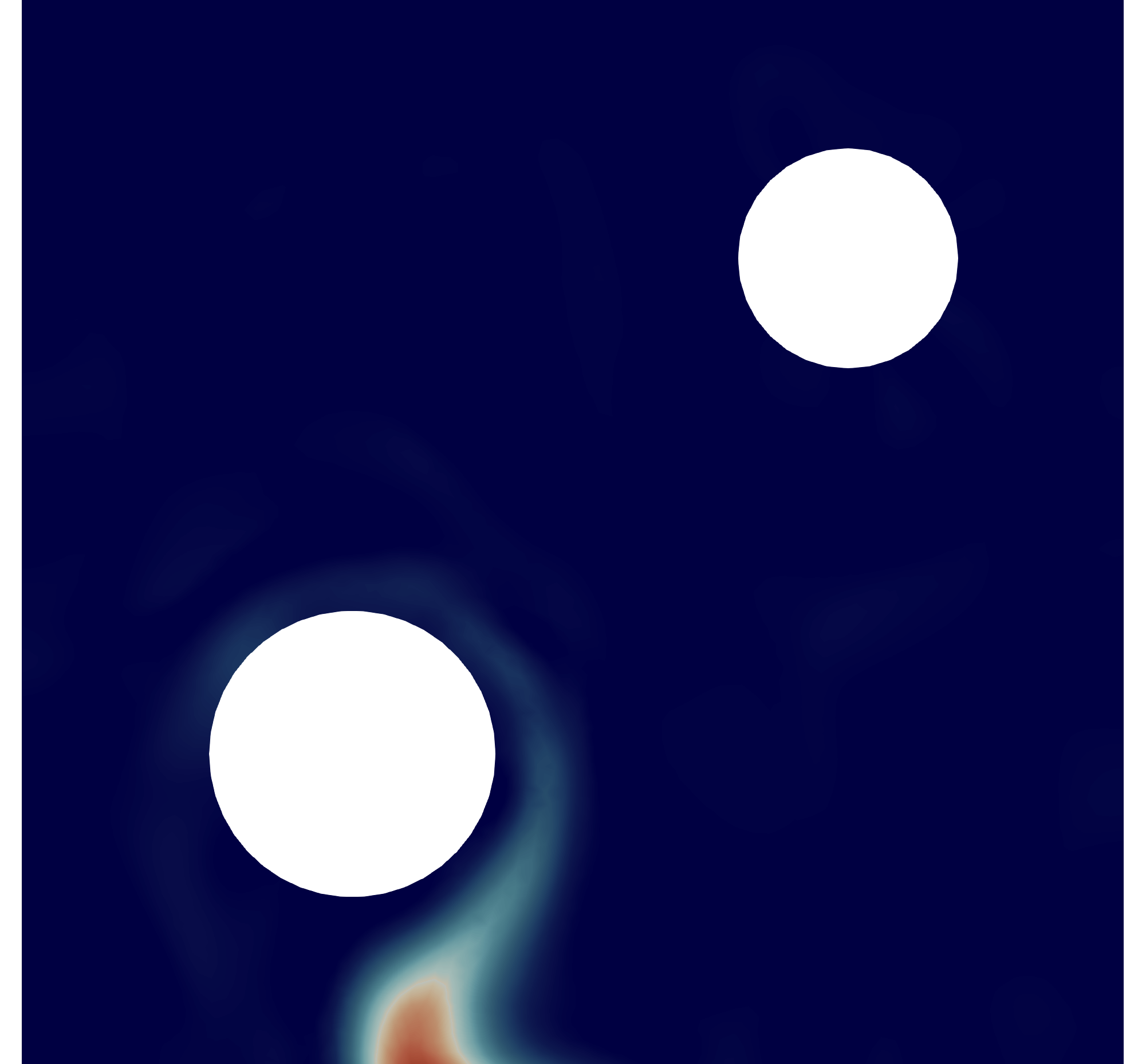}}
         {\includegraphics[width=0.19\textwidth]{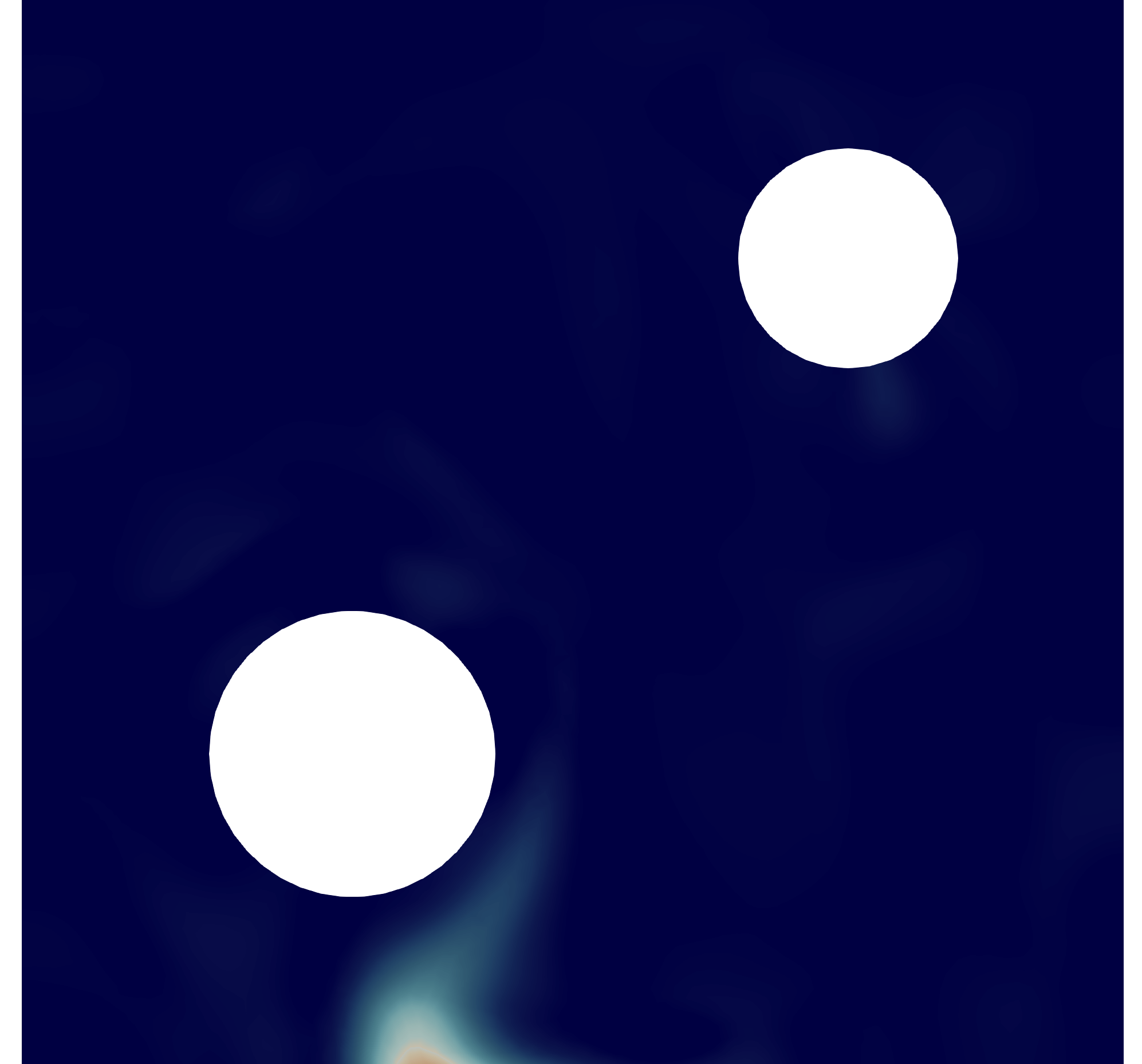}}
 
\caption{Solution of the advection diffusion problem on a domain with holes for three different parameter values (left to right) $\vec{\alpha_1} =[-0.4, -0.4]$, $\vec{\alpha_2} =[0,0.4]$ and $\vec{\alpha_3} =[0.4,0.4]$ of the initial condition. The neural network uses $n_\Phi=20$ input features and has four layers with 20 hidden units. 8000 points are sampled for the update equation in each step.}
\label{fig:snapshots_adv_diff_holes_p}
\end{figure}

We use the sampling algorithm for a discrete points of 40000 candidate points, for which the sampling criterion (absolute value of advective term) is evaluated. We then subsample 5000 and 8000 points respectively which are used to evaluate the update equation. In Figure \ref{fig:snapshots_adv_diff_holes_p}, one  sees the solution for different parameter values and how they evolve differently according to their starting position.
Considering the error plot, we observe a larger error in this challenging case, however, the neural network solutions still provide valuable information.

\begin{figure}
 \centering
   \subcaptionbox{Linear scale}
         {\includegraphics[width=0.45\textwidth]{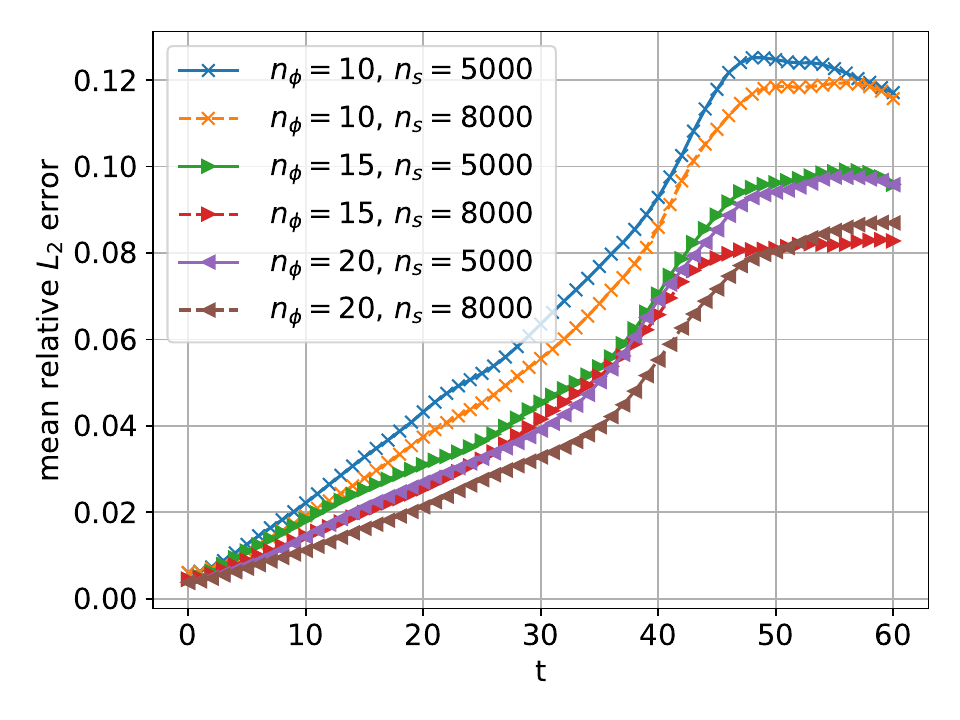}}
   \subcaptionbox{ Log scale}
         {\includegraphics[width=0.45\textwidth]{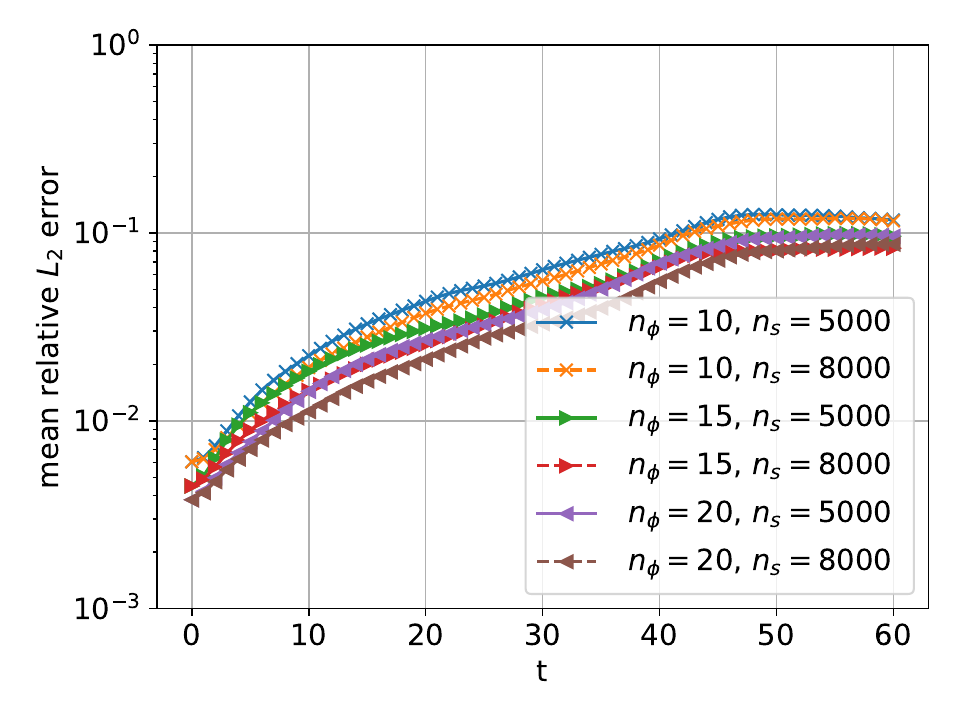}}
\caption{Error plots for the advection-diffusion problem on the domain with holes and the parametrized initial condition.}
\label{fig:error_adv_diff_holes_p}
\end{figure}

\section{Conclusion}
We have successfully demonstrated the effectiveness of evolutional deep neural networks (EDNNs) for solving parametric  time-dependent PDEs on domains with geometric structure.
Our work proposes significant steps towards employing EDNNs  to solve PDEs in real-world scenarios. 

By introducing positional embeddings based on eigenfunctions of the Laplace-Beltrami operator, we achieve intrinsic encoding of geometric properties and automatic enforcement of Dirichlet, Neumann and periodic boundary conditions, resulting in a simplified and better-conditioned learning problem for both static and dynamic PDEs.
We have shown that even a small number of feature functions and neural networks of moderate size can effectively handle transport-dominated problems. 

Utilizing Krylov solvers instead of direct assembly has allowed us to scale the EDNN approach to larger neural networks, while our active sampling scheme has further improved computational efficiency. Moreover, the use of linearly implicit Rosenbrock methods allows us to effectively handle stiff PDE problems. We have further highlighted that training can be completely eliminated from the time-stepping scheme, albeit at the cost of larger approximation errors.

For single-query scenarios, traditional numerical methods will likely continue to outperform EDNNs, but in the many-query case, our approach demonstrates competitiveness as the computation of the solution across the parameter domain can happen in one single time integration.  Additionally, the solution is available as a functional model with respect to the problem parameters, so that sampling or computing derivatives, for e.g. sensitivity studies, becomes inexpensive once time integration has been performed. 

To fully understand the theoretical properties of EDNNs, additional investigations are necessary. Particularly, we need a further understanding of  error control and how the approximation capabilities of neural networks affect  their tangent spaces along the solution trajectory. Preconditioning, although unexplored in this study, holds potential as a beneficial technique for the Krylov solver.
While our experiments focus on moderate-sized NNs, in future work the same methodology could be applied to more complex neural network structures. Furthermore, the idea of harmonic features can potentially be  integrated into other neural network structures that currently rely on Fourier embeddings, e.g. the Fourier Neural Operator.
In conclusion, our work contributes towards  utilizing EDNNs for real-world PDE applications and we believe that some of the ideas introduced in this work can be applied to improve training in other neural-network based methods in PDE surrogate modelling.
 
\clearpage
\appendix

\section{Initial condition for the Korteweg-de-Vries equation}
\label{app:kdV}
The exact solution for the two soliton collision problem can be described by the following construction:
\begin{equation}
u(x,t) = 2 \dfrac{\partial^2 \log(f(x,t))}{\partial x^2},
\end{equation}
where 
\begin{equation}
\begin{split}
    f  &= 1 + e^{\eta_1} + e^{\eta_2} + A e^{\eta_1 +\eta_2}, \\
    \eta_i &= k_i x - k_i^3 t + \eta_i ^{(0)}, \\
    A &= \left (\dfrac{k_1 - k_2}{k_1+ k_2} \right)^2,
\end{split}
\end{equation}
with parameters $k_1 = 1$, $k_2 = \sqrt{5}$, $\eta_1^{(0)} =0$, $\eta_2^{(0)} =10.73$.
The initial condition is given by evaluating the exact solution at $t=0$:  $u_0 = u(x,0)$.
We point out that this solution is defined on the infinite domain, while the computational domain used here is the interval (-20,20). More details can be found in \cite{taha_analytical_1984}.

\section{Velocity field computation for the domain with holes}
\label{sec:appvel}
To obtain a physically valid velocity field, we solve the steady states Navier-Stokes equation:

\begin{equation}
    \begin{split}
      &  \frac{-1}{\text{Re}} \nabla^2 \vecf{w} + \nabla q + \vecf{w} \cdot \nabla \vecf{w} =0, \\
       & \nabla \cdot \vecf{w} =0,\\
       & \vecf{w} = \vec{g_1}   \quad x \in \partial\Omega_\text{top},\\ 
       &\vecf{w} = \vec{g_2}  \quad x \in  \partial \Omega \backslash \partial \Omega_\text{top}
    \end{split}
    \label{eq:NS_pre}
\end{equation}
By setting the Reynolds number $\text{Re}= 100$, the solution will be in the laminar flow regime.
The flow is driven by an inflow Dirichlet boundary condition $\vec{g_1}= (0,-1)$  on the top. A no-slip boundary condition $\vec{g_2}=(0,0)$ is imposed on all remaining parts of the domain boundary. The pressure is set to 1 in one corner of the domain, to obtain a well-posed problem.  \\
Equation \eqref{eq:NS_pre} above are solved numerically in FEnics. The equations are discretized by a mixed element formulation, using P2-Lagrange elements for the velocity field and P1-Lagrange elements for the pressure. The nonlinear system is then solved with a Newton solver. 

\section{Eigenfunctions for the domain with holes}
\label{sec:eigenfuncs}
\begin{figure}[h]
 \centering
    \subcaptionbox{$\phi_1$}
{\includegraphics[width=0.19\textwidth]{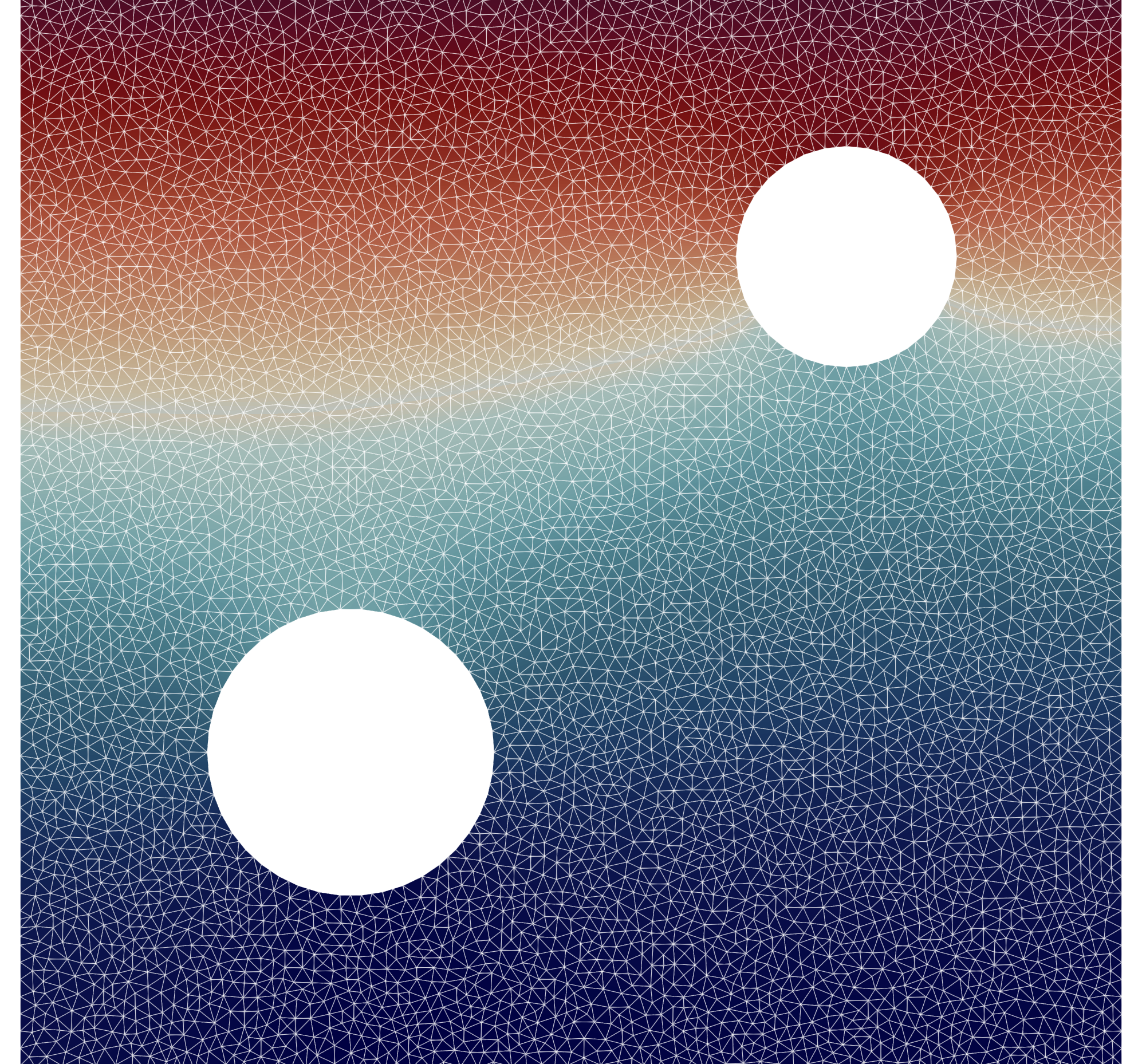}}
    \subcaptionbox{$\phi_2$}
{\includegraphics[width=0.19\textwidth]{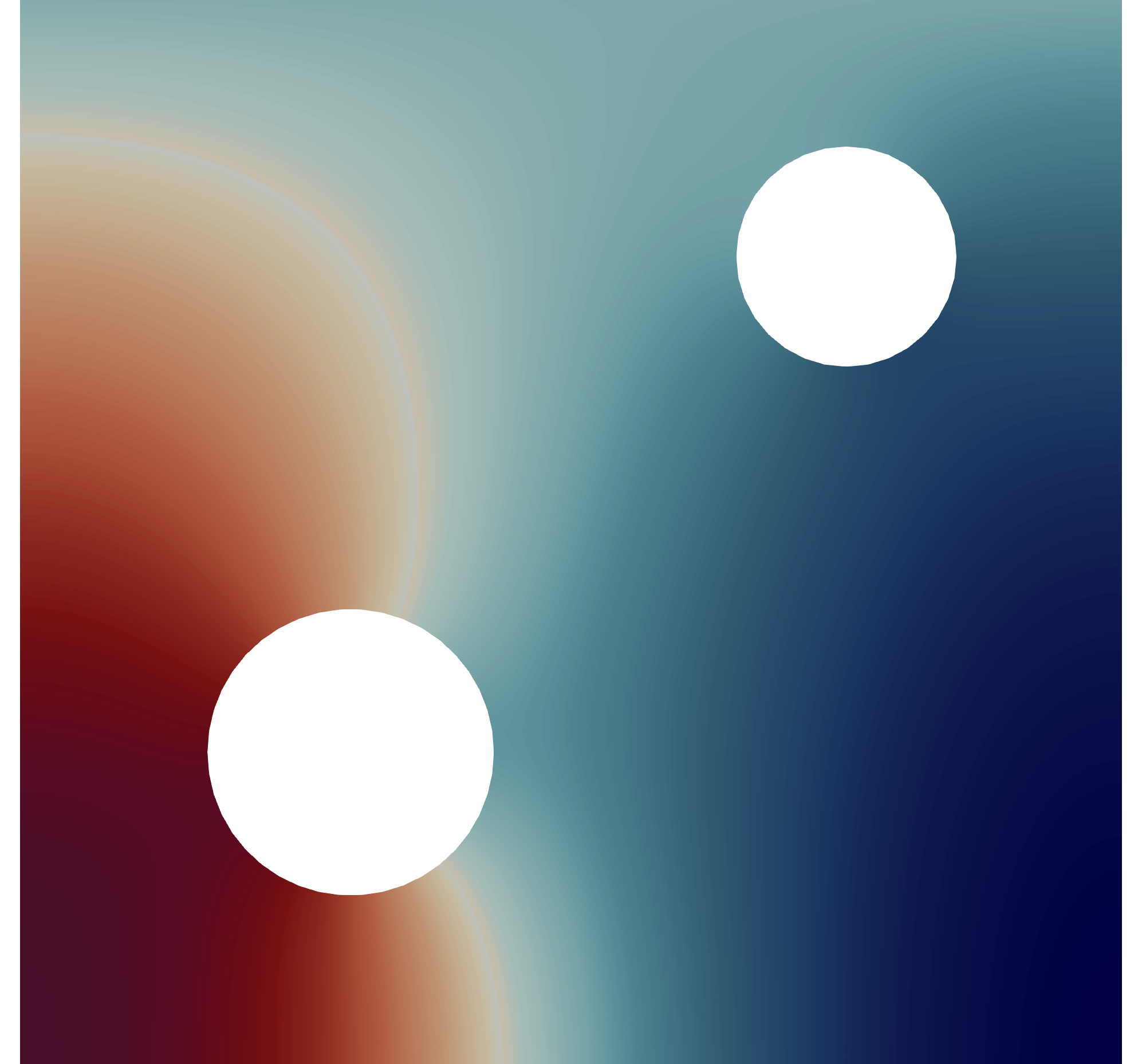}}
    \subcaptionbox{$\phi_3$}
{\includegraphics[width=0.19\textwidth]{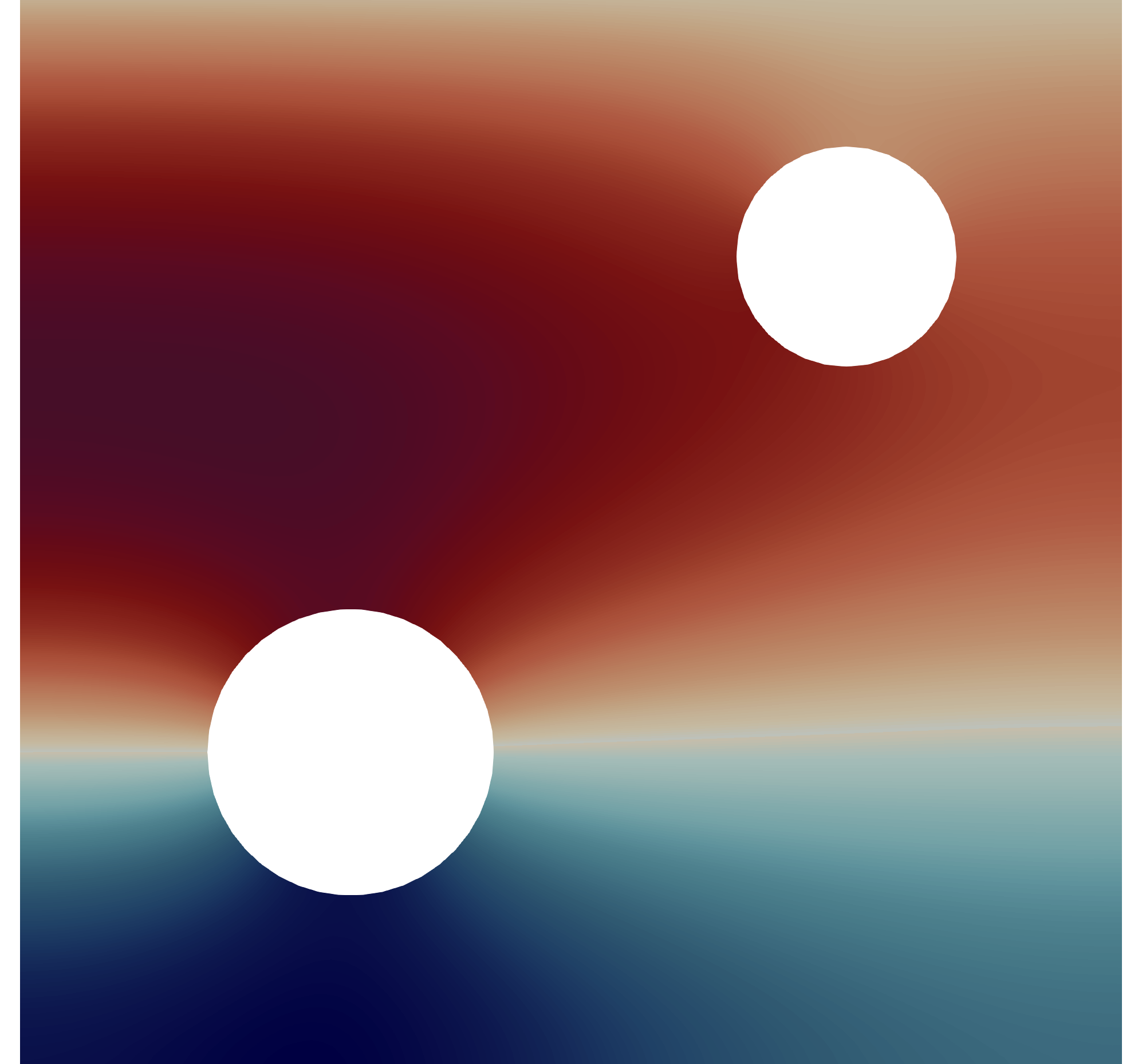}}
    \subcaptionbox{$\phi_4$}
{\includegraphics[width=0.19\textwidth]{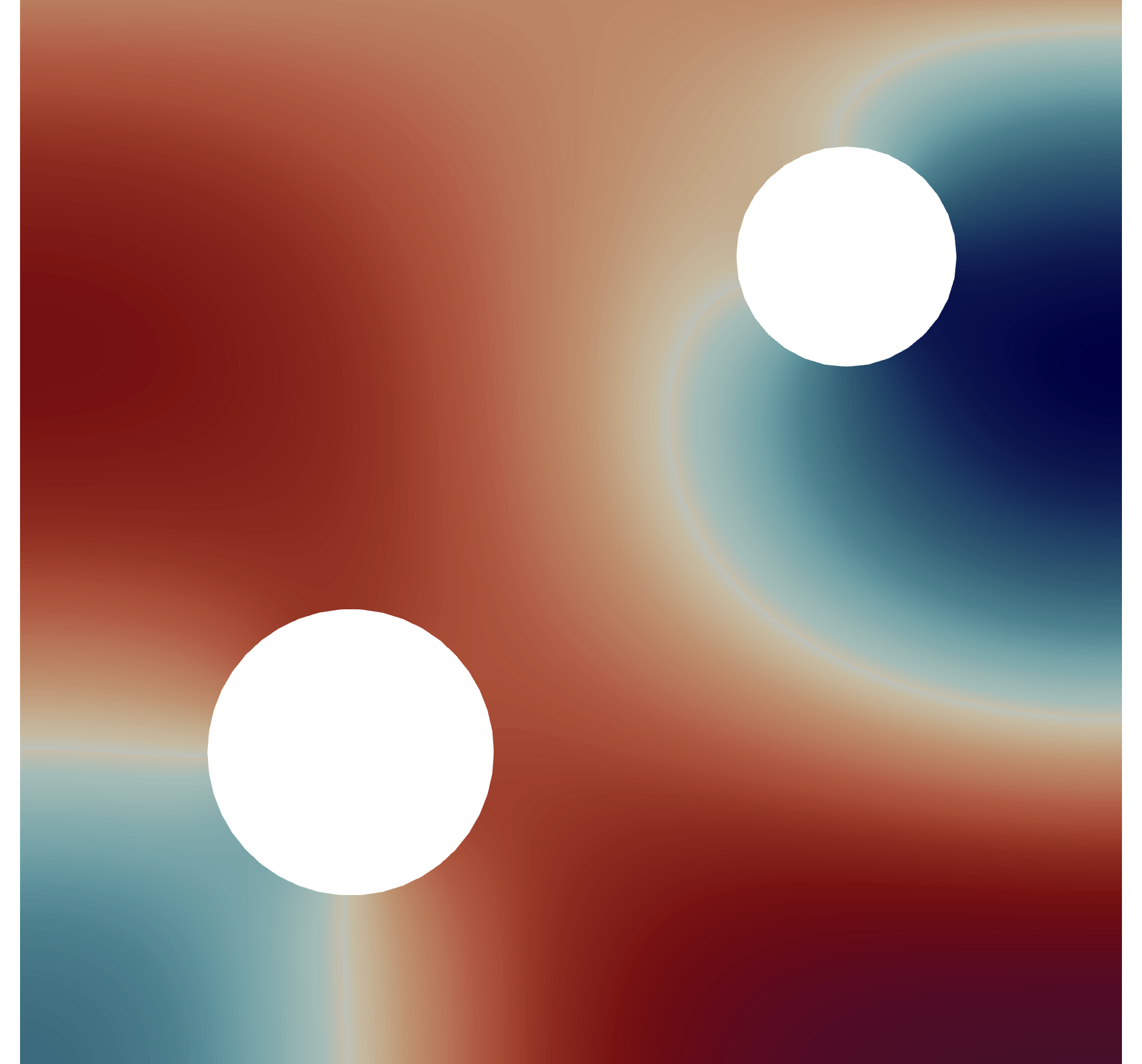}}
    \subcaptionbox{$\phi_5$}
{\includegraphics[width=0.19\textwidth]{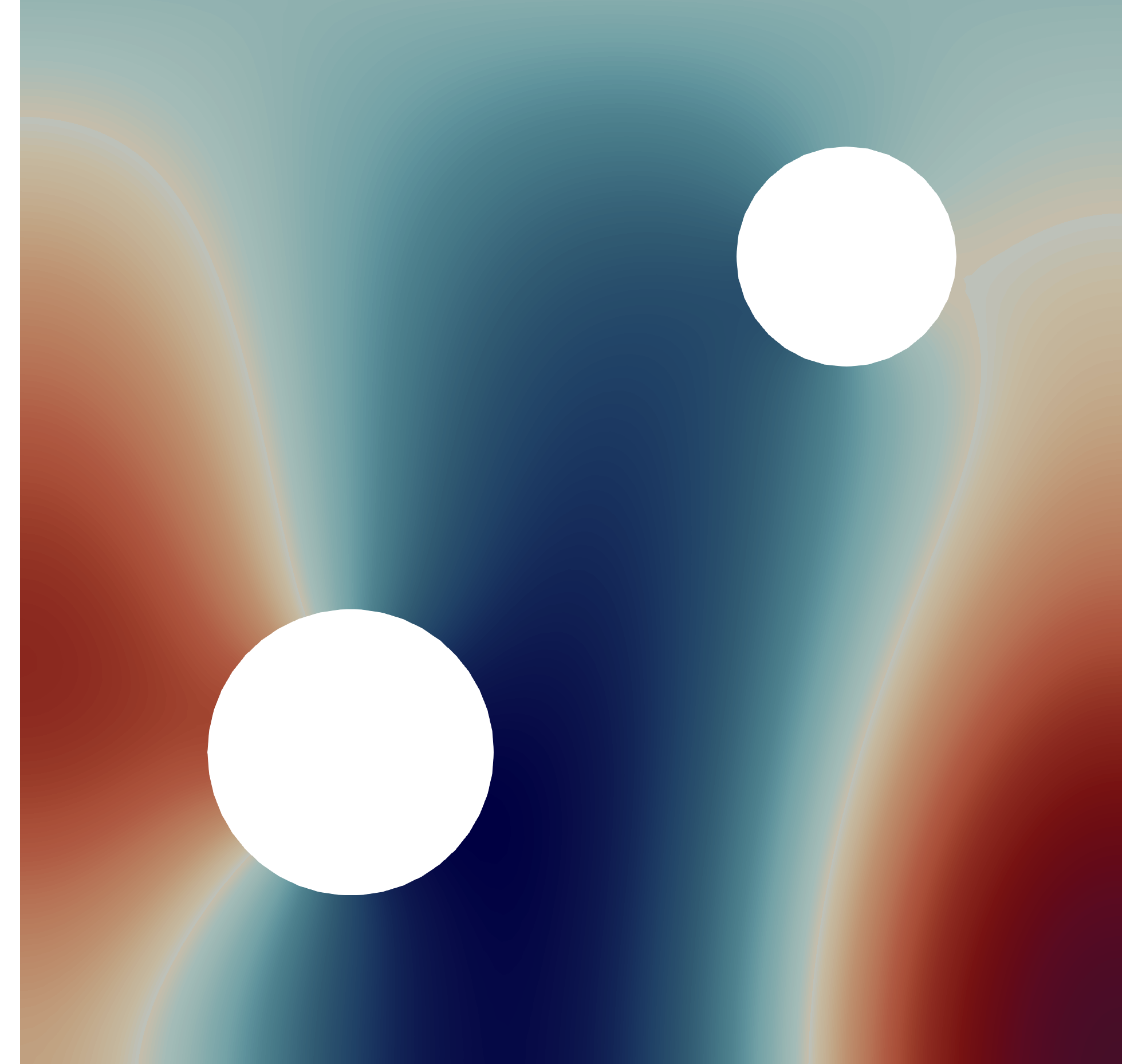}}
    \subcaptionbox{$\phi_6$}
{\includegraphics[width=0.19\textwidth]{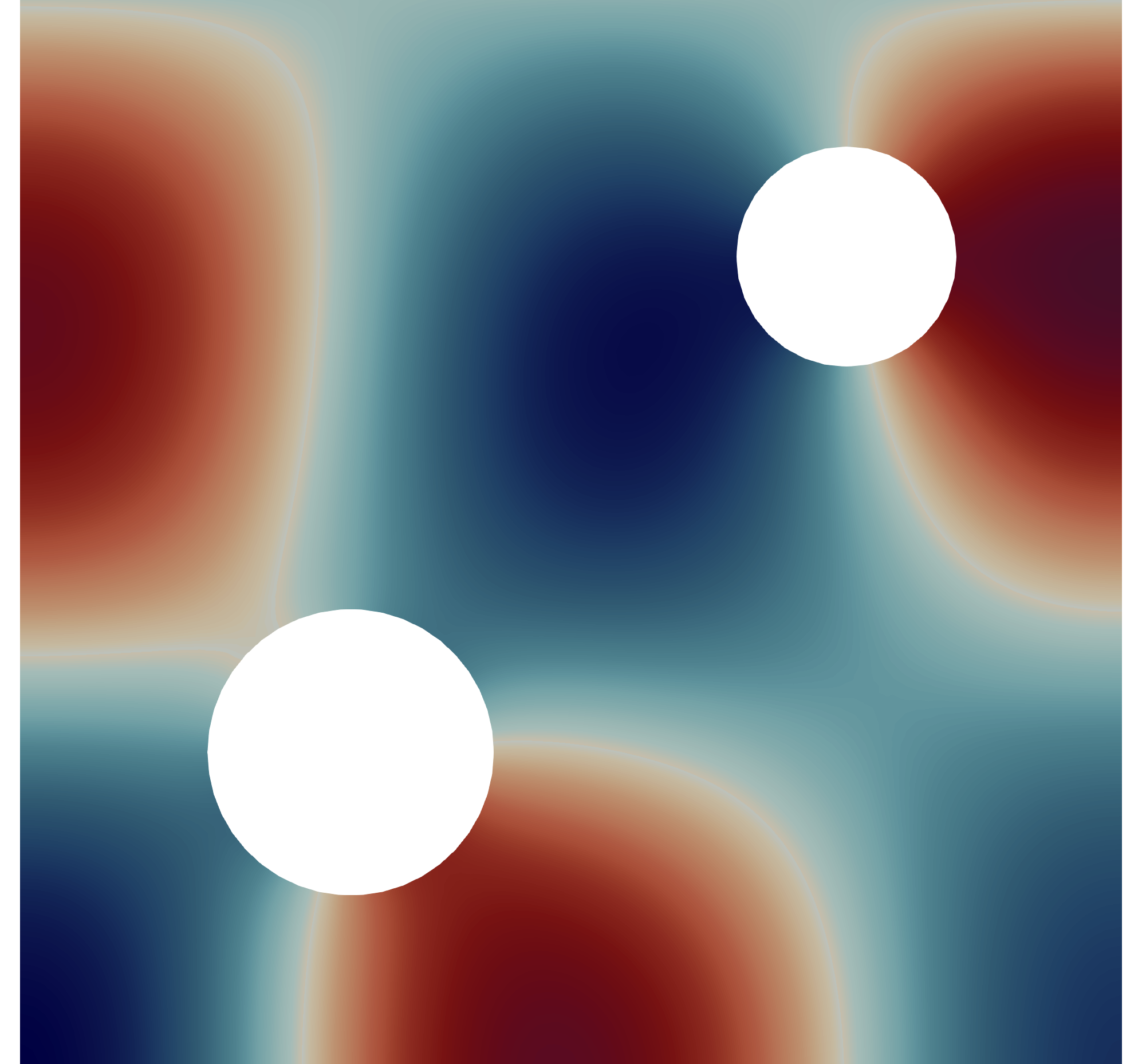}}
    \subcaptionbox{$\phi_7$}
{\includegraphics[width=0.19\textwidth]{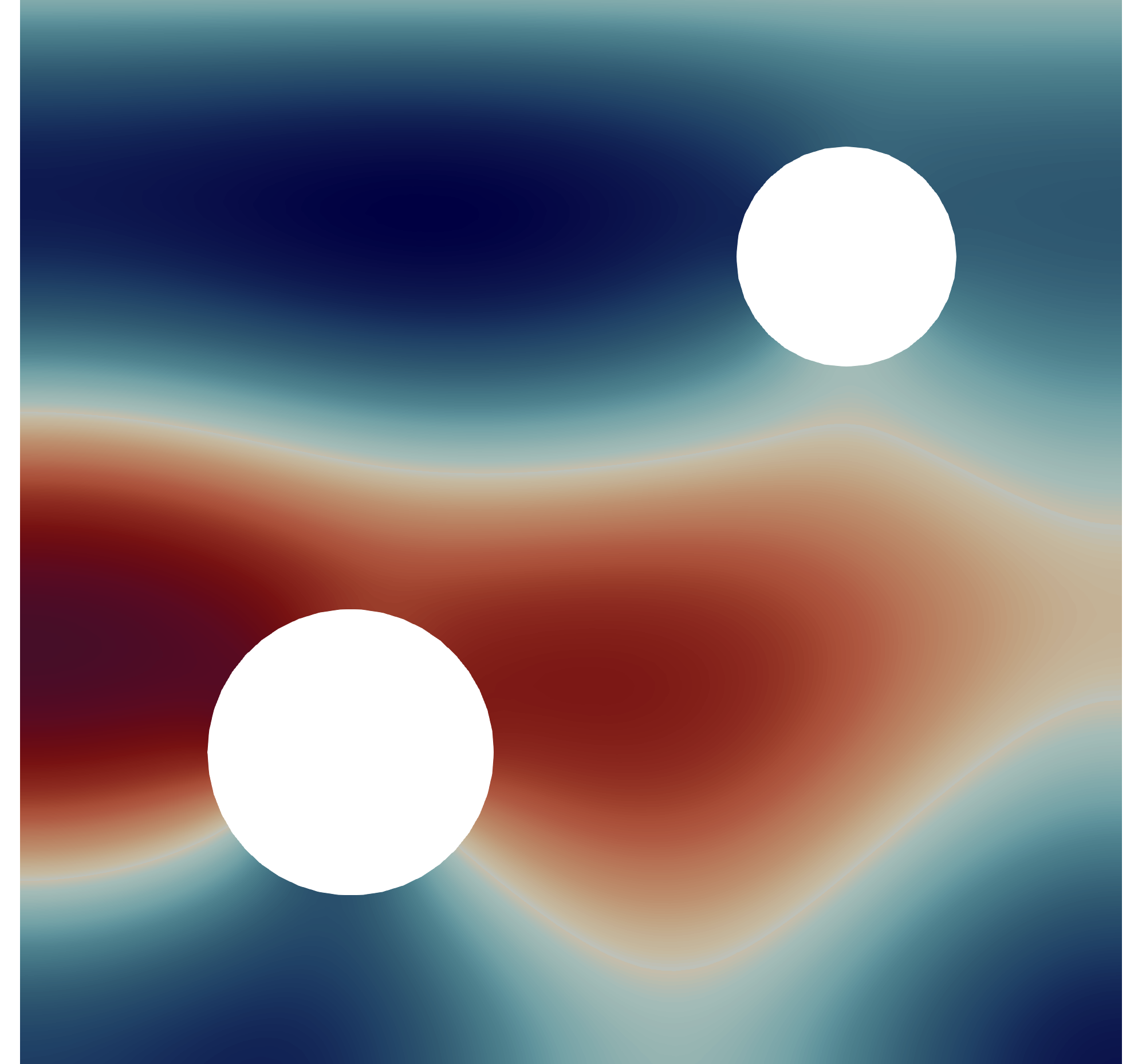}}
    \subcaptionbox{$\phi_8$}
{\includegraphics[width=0.19\textwidth]{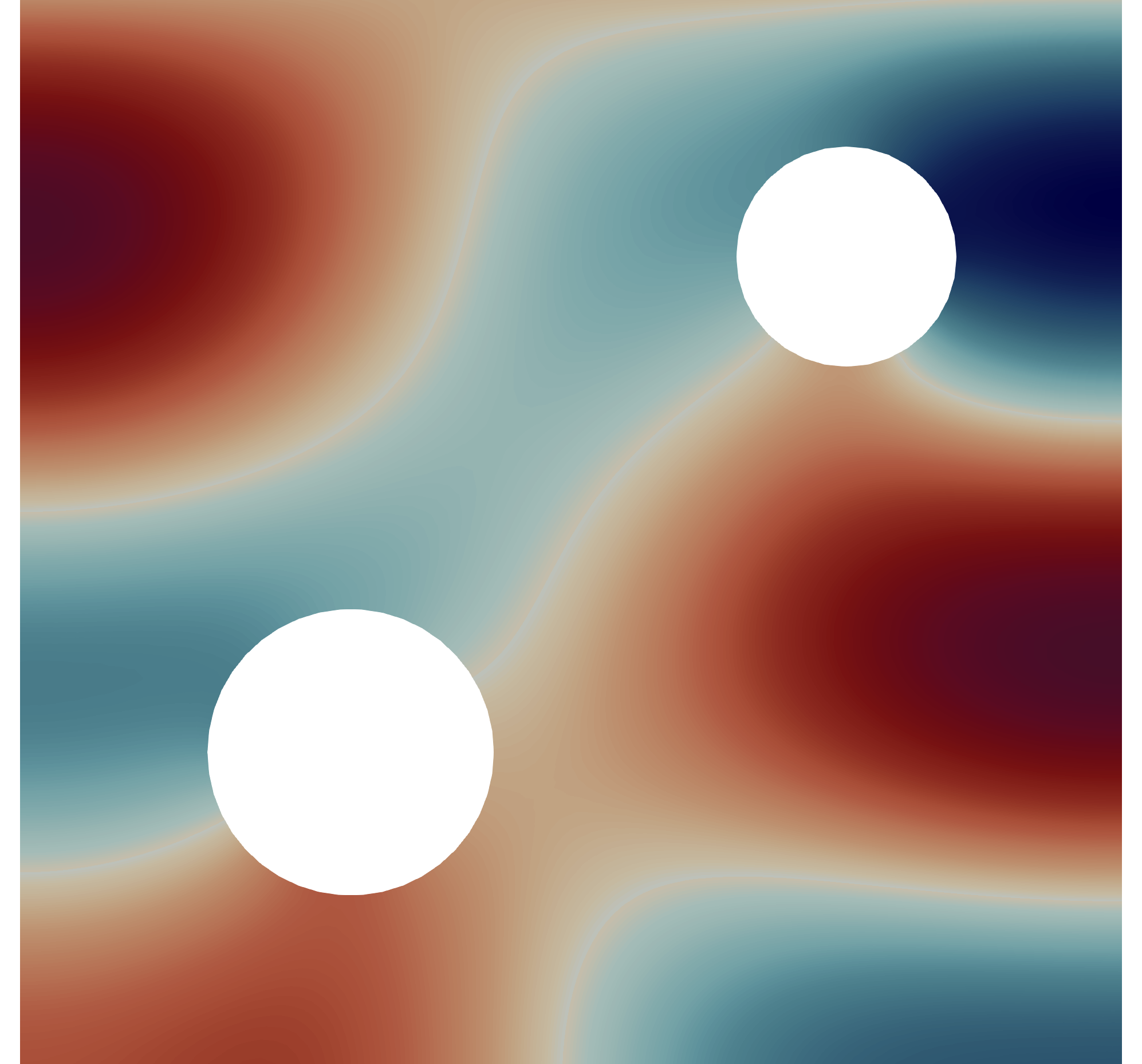}}
    \subcaptionbox{$\phi_9$}
{\includegraphics[width=0.19\textwidth]{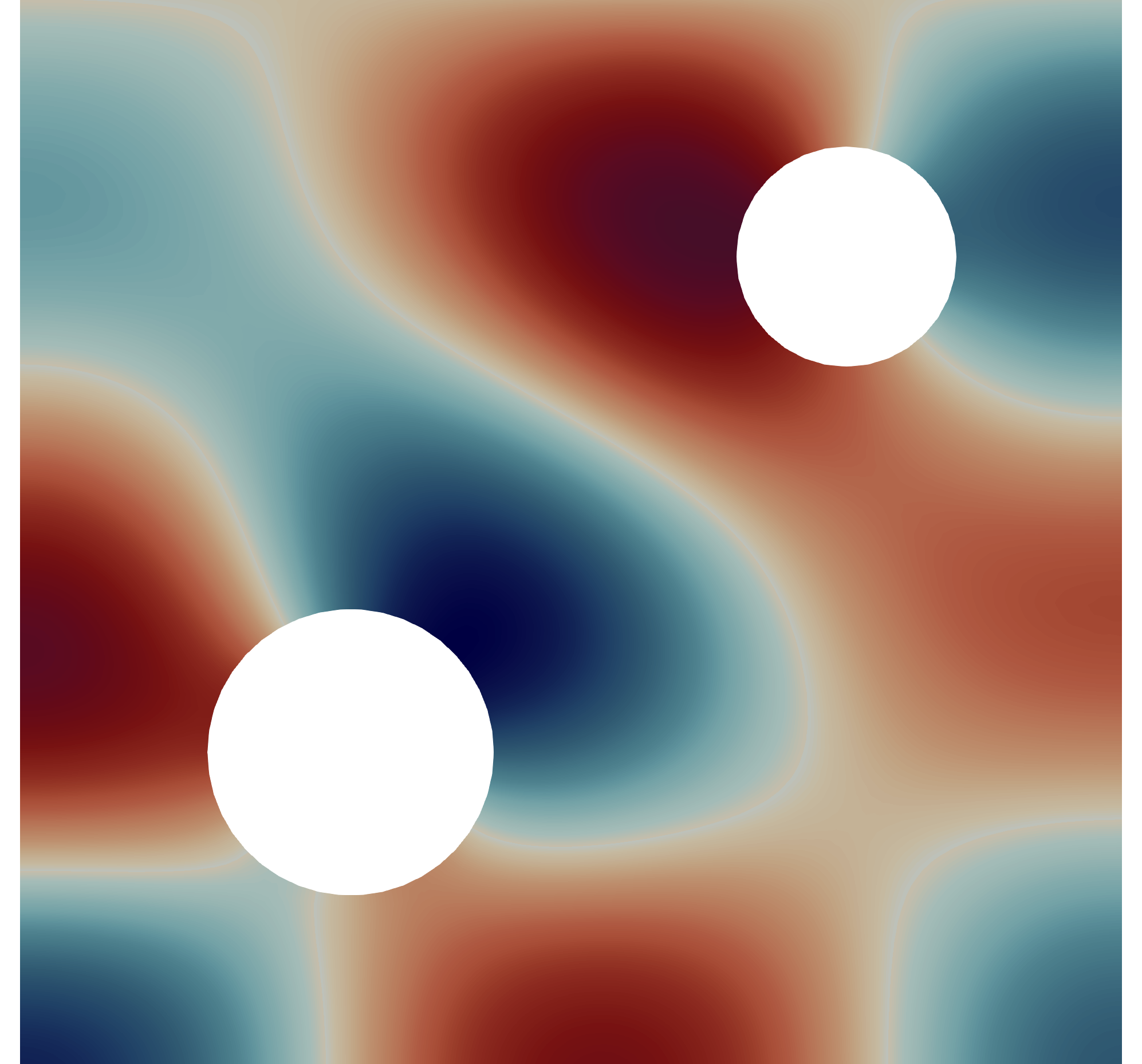}}
    \subcaptionbox{$\phi_{10}$}
{\includegraphics[width=0.19\textwidth]{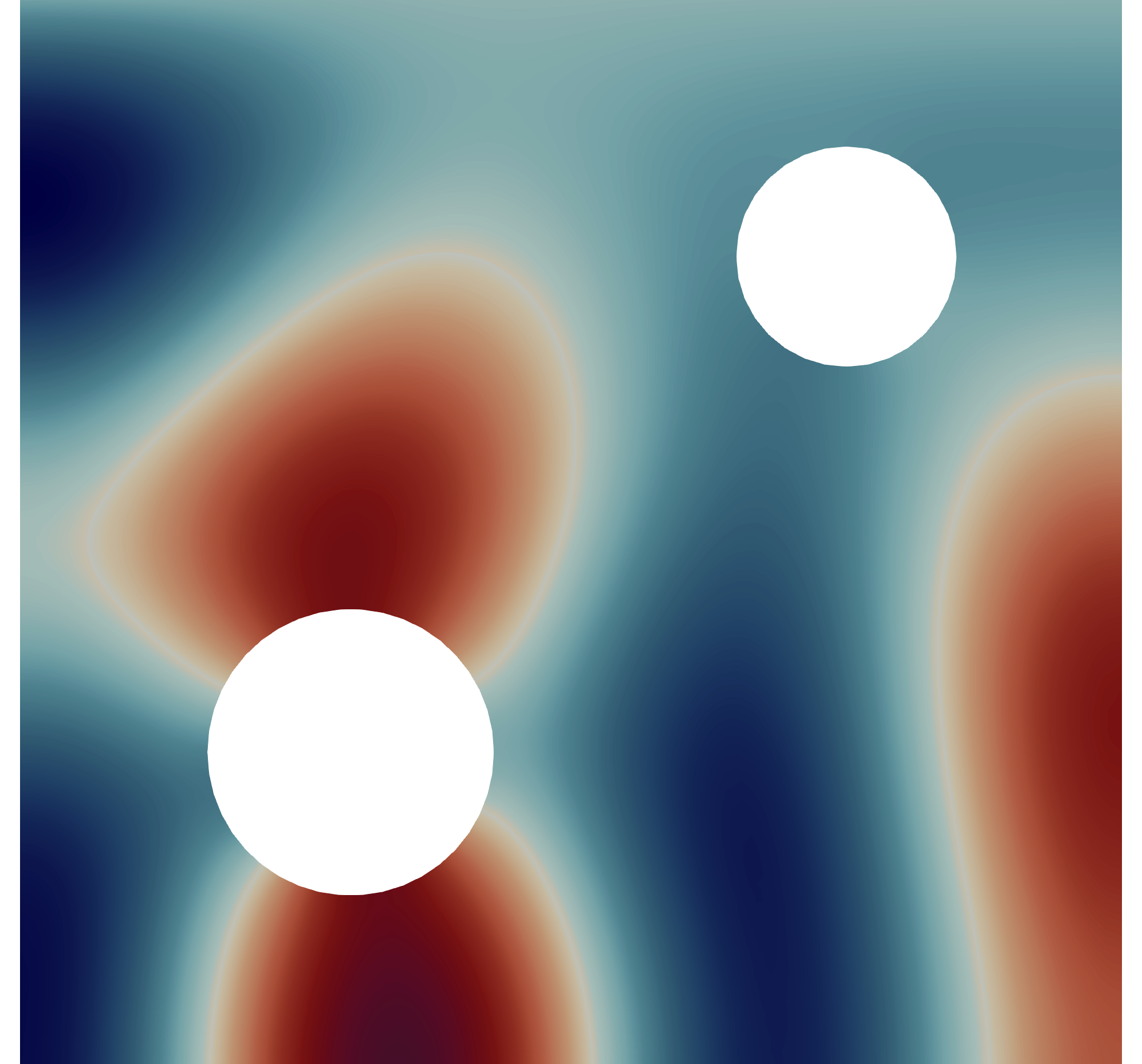}}

\caption{The first ten eigenfunctions that are used for the positional embedding on the domain with a hole. We plot the FE mesh on the first mode to visualize the resolution.}
\label{fig:eigenfunc}
\end{figure}
\newpage
 \bibliographystyle{elsarticle-num} 
 \bibliography{literature_pos_enc}

\end{document}